\documentclass[12pt]{amsart}
\usepackage{latexsym,amssymb,amsrefs,graphicx}
\usepackage{mathrsfs}
\usepackage[all]{xy}
\newtheorem{lemma}{Lemma}
\numberwithin{lemma}{section}

\newtheorem{claim}[lemma]{Claim}
\newtheorem{proposition}[lemma]{Proposition}
\newtheorem{theorem}[lemma]{Theorem}
\newtheorem{corollary}[lemma]{Corollary}
\newtheorem{definition}[lemma]{Definition}

\theoremstyle{definition}

\newtheorem{example}[lemma]{Example}
\newtheorem{remark}[lemma]{Remark}

\newcommand{\proofstart}{\smallskip\noindent Proof:\ }

\newcommand{\proofend}{\begin{flushright}
                       $\Box$
                       \end{flushright}}

\begin{document}

\title[Twisted Morse complexes]{Twisted Morse complexes}

\author{Augustin Banyaga}
\address{Department of Mathematics \\
         Penn State University \\
         University Park, PA 16802}
\email{banyaga@math.psu.edu}

\author{David Hurtubise}
\address{Department of Mathematics and Statistics\\
         Penn State Altoona\\
         Altoona, PA 16601-3760}
\email{Hurtubise@psu.edu}

\author{Peter Spaeth}
\address{Department of Mathematics and Statistics\\
         Penn State Altoona\\
         Altoona, PA 16601-3760}
\email{pwspaeth@gmail.com}



\begin{abstract}
In this paper we study Morse homology and cohomology with local coefficients, i.e. ``twisted'' Morse homology and cohomology, on closed finite dimensional smooth manifolds.  We prove a Morse theoretic version of Eilenberg's Theorem, and we prove isomorphisms between twisted Morse homology, Steenrod's CW-homology with local coefficients for regular CW-complexes, and singular homology with local coefficients. By proving Morse theoretic versions of the Poincare Lemma and of the de Rham Theorem, we show that twisted Morse cohomology with coefficients in a local system determined by a closed 1-form is isomorphic to the Lichnerowicz cohomology obtained by deforming the de Rham differential by the 1-form. We demonstrate the effectiveness of twisted Morse complexes by using them to compute Lichnerowicz cohomology, to compute obstructions to spaces being associative H-spaces, and to compute Novikov numbers. 
\end{abstract}

\maketitle

\eject


\vspace*{\fill}

\begin{center}
{\bf\large Preface}
\end{center}

\bigskip

This paper is a discourse on Morse homology and cohomology with local coefficients, i.e. ``twisted''
Morse homology and cohomology. The Morse-Smale-Witten chain complex of a smooth function 
$f:M \rightarrow \mathbb{R}$ with coefficients in a bundle of abelian groups $G$ over a closed
finite dimensional smooth Riemannian manifold $(M,\mathsf{g})$ is constructed 
(Definition \ref{twistedboundary}), and an Invariance Theorem (Theorem \ref{homologyindependence}),
which shows that the homology of the twisted Morse-Smale-Witten chain complex is independent of the
Morse-Smale pair $(f,\mathsf{g})$ and depends only on the isomorphism class of $G$, is proved. 

An isomorphism between the twisted Morse homology of $(M,G)$ and the singular homology of $M$ with 
coefficients in $G$ (Theorem \ref{twistedsame}) is proved by comparing with Steenrod's CW-chain
complex with coefficients in $G$, defined for regular CW-complexes \cite{SteHom}. The existence 
of a Morse-Smale pair $(f,\mathsf{g})$ whose unstable manifolds determine a regular CW-structure
on the manifold is demonstrated (Theorem \ref{regularMorse}), and a Morse theoretic version of 
Eilenberg's Theorem (Theorem \ref{EilenbergMorse}), which relates homology with local coefficients
to equivariant homology, is proved. 

The Morse-Smale-Witten cochain complex with coefficients in a bundle of groups $G$ is
defined (Definition \ref{twistedcochain}) and studied in the case where the local system $G=e^\eta$
has fiber $\mathbb{R}$ and the homomorphisms are defined by integrating a closed 1-form $\eta$
along paths (Example \ref{etasystem}). The twisted Morse cohomology with coefficients in 
$e^\eta$ is shown to be isomorphic to the Lichnerowicz cohomology $H^\ast_{-\eta}(M)$ 
(an invariant used to study locally conformal symplectic manifolds) by proving $\eta$-twisted 
Morse theoretic versions of the Poincar\'e Lemma (Lemma \ref{MSWPoincare}) and of the 
de Rham Theorem (Theorem \ref{twistedDeRham}). 

Twisted Morse complexes are effective tools for computing homology and cohomology with
local coefficients, because they often have fewer generators than a Steenrod CW-chain complex
with local coefficients. To demonstrate this, the Lichnerowicz cohomology of a surface of 
genus two is computed explicitly with respect to all closed 1-forms on the surface 
(Example \ref{genus2Lich}), and following a result due to Albers, Frauenfelder, and Oancea
\cite{AlbLoc}, twisted Morse complexes are used to show certain spaces are not
associative H-spaces (Corollary \ref{Hspaceobstruction}, Examples \ref{genus2notH} 
and \ref{RP2notH}). The last section of this paper contains a discussion of the Novikov 
homology of a closed 1-form and how Novikov homology can be viewed as a special case 
of twisted Morse homology. Twisted Morse complexes are used to explicitly compute the 
Novikov numbers of various surfaces with respect to all closed 1-forms on the surfaces 
(Examples \ref{Novcircle}\,--\,\ref{Novsurface2}).

\bigskip

\begin{flushright}
\begin{tabular}{l}
Augustin Banyaga\\
David Hurtubise\\
Peter Spaeth
\end{tabular}
\end{flushright}

\vspace*{\fill}

\eject


\vspace*{\fill}

\tableofcontents

\vspace*{\fill}

\eject


\section{Introduction}

Morse homology is a homology theory defined using a smooth Morse-Smale function 
$f:M \rightarrow \mathbb{R}$ on a smooth Riemannian manifold $(M,\mathsf{g})$. The $k$-th 
chain group $C_k(f)$ is defined using the critical points of index $k$, and the boundary 
operator is defined by counting with sign the number of gradient flow lines between 
critical points of relative index one. The Morse Homology Theorem says that on a closed finite 
dimensional smooth manifold the Morse homology with integer coefficients is 
isomorphic to the singular homology of the manifold with integer coefficients, cf. \cite{BanLec},
\cite{MilLec} \cite{SchMor}.

\smallskip
In more detail, a function $f:M\rightarrow \mathbb{R}$ is called a \textbf{Morse} function 
if all of its critical points are \textbf{non-degenerate}, i.e. the \textbf{Hessian} of $f$ at every 
critical point is non-degenerate, cf. Definition 3.1 of \cite{BanLec}. The assumption that 
$p \in Cr(f)$ is a non-degenerate critical point implies that the \textbf{stable manifold} 
$W^s(p)$ and the \textbf{unstable manifold} $W^u(p)$, defined by
\begin{eqnarray*}
W^s(p) & = & \{ x\in M | \lim_{t \rightarrow \infty} \varphi_t(x) = p \}\\
W^u(p) & = & \{ x\in M | \lim_{t \rightarrow -\infty} \varphi_t(x) = p \},
\end{eqnarray*}
are embedded open disks, where $\varphi_t$ is the 1-parameter group of diffeomorphisms generated by
the negative of the gradient vector field. The \textbf{index} of $p$, defined as the dimension of
the subspace of $T_pM$ on which the Hessian is negative definite, coincides with the dimension of 
the disk $W^u(p)$, cf. Theorem 4.2 of \cite{BanLec}. 

If all the stable and unstable manifolds of $f$ intersect transversally, i.e.
$$
W^u(q) \pitchfork W^s(p)
$$
for all $p,q \in Cr(f)$, then the pair $(f,\mathsf{g})$ is called a \textbf{Morse-Smale pair}, 
and $f$ is said to be \textbf{Morse-Smale} or to satisfy the \textbf{Morse-Smale transversality}
condition on $(M,\mathsf{g})$. When $(f,\mathsf{g})$ is a Morse-Smale pair, 
$$
W(q,p) = W^u(q) \cap W^s(p) \subset M
$$ 
is either empty or an embedded submanifold of dimension $\lambda_q - \lambda_p$, where
$\lambda_q$ and $\lambda_p$ are the indexes of $q$ and $p$ respectively, 
cf. Proposition 6.2 of \cite{BanLec}. In particular, $W(q,p)$ is a one dimensional manifold
when $\lambda_q - \lambda_p = 1$, and if $M$ is compact, then $W(q,p)$ consists of the images
of finitely many gradient flow lines from $q$ to $p$, cf. Corollary 6.29 of \cite{BanLec}.
Thus, when $M$ is compact and the pair $(f,\mathsf{g})$ is Morse-Smale we can count the number of 
gradient flow lines between any two critical points of relative index one. This is sufficient to define
the Morse-Smale-Witten chain complex with coefficients in $\mathbb{Z}_2$, but to define the complex
with integer coefficients we need to consider orientations.

Arbitrarily choosing orientations on the unstable manifolds of a Morse-Smale pair $(f,\mathsf{g})$
determines orientations on the normal bundles of the stable manifolds, which then determines signs
associated to the gradient flow lines between critical points of relative index one 
(see Section \ref{pathcomponents}). With these signs we can define the Morse-Smale-Witten
chain complex over $\mathbb{Z}$.

\begin{definition}[Morse-Smale-Witten chain complex]
Let $f:M \rightarrow \mathbb{R}$ be a smooth Morse-Smale function on a closed smooth Riemannian
manifold $(M,\mathsf{g})$ of dimension $m < \infty$. Fix orientations on the unstable manifolds
of $(f,\mathsf{g})$. The \textbf{Morse-Smale-Witten chain complex with coefficients in $\mathbb{Z}$}
is defined to be the chain complex $(C_\ast(f),\partial_\ast)$, where $C_k(f)$ is the free abelian
group generated by the critical points of $f$ of index $k$ for all $k=0,\ldots, m$, and 
$\partial_k:C_k(f) \rightarrow C_{k-1}(f)$ is defined on a generator $q \in C_k(f)$ to be
$$
\partial_k(q)\ = \sum_{p \in Cr_{k-1}(f)} \left(\sum_{\nu \in \mathcal{M}(q,p)} \epsilon(\nu)
\right) p,
$$ 
where $\epsilon(\nu) = \pm 1$ denotes the sign associated to the unparameterized gradient flow
line $\nu \in \mathcal{M}(q,p) = W(q,p)/\mathbb{R}$ by the orientations.
\end{definition}

\noindent
Using the structure of the compactified moduli spaces $\overline{\mathcal{M}}(r,p)$, where
$\lambda_r - \lambda_p = 2$ one can show that $(\partial_\ast)^2 = 0$, 
cf. Lemma \ref{boundarysquared}. Thus, $(C_\ast(f),\partial_\ast)$ is a chain complex.

\smallskip
The following theorem is well-known and can be proved using many different approaches
\cite{BanLec} \cite{HelPui} \cite{MilLec} \cite{QinOnm} \cite{SalMor}
\cite{SchMor} \cite{WitSup}.

\begin{theorem}[Morse Homology Theorem]\label{MorseHomology}
Let $f:M \rightarrow \mathbb{R}$ be a smooth Morse-Smale function on a closed finite dimensional
smooth Riemannian manifold $(M,\mathsf{g})$. Then the homology of the Morse-Smale-Witten chain
complex with coefficients in $\mathbb{Z}$ is isomorphic to the singular homology of $M$ with coefficients 
in $\mathbb{Z}$.
\end{theorem}

Morse homology can be extended in several ways. For instance, it can be extended
to infinite dimensional manifolds, it can be extended to include critical submanifolds, 
it can be extended to manifolds with boundaries or corners, and it can be extended to 
include local coefficient systems. Moreover, these extensions may overlap, 
i.e. one can consider infinite dimensional versions of Morse homology with local coefficients 
or Morse homology on manifolds with boundaries or corners with local coefficients.

\smallskip
In previous work, the first two authors explored various ways to extend Morse homology to 
include critical submanifolds \cite{BanMor} \cite{BanCas} \cite{HurThr}. In this paper
we discuss how to extend Morse homology to include local coefficient systems, which we
refer to as \textbf{twisted Morse homology}. The discussion is limited to Morse functions 
on closed finite dimensional smooth manifolds.  Morse-Bott functions and  manifolds with boundaries
or corners are not discussed outside of this introduction. Even so, there are several interesting
applications of Morse homology with local coefficients that fit within the scope of this paper.

For instance, a relationship between twisted Morse cohomology and \textbf{Lichnerowicz cohomology}
is proved in Section \ref{CoLich}, and we show how to use twisted Morse cohomology to compute 
the Lichnerowicz cohomology of a surface with respect to any closed 1-form on the surface
in Section \ref{ParLich}. In Section \ref{Hspace} we discuss, following a result due to Albers, 
Frauenfelder, and Oancea \cite{AlbLoc}, how Morse homology with local coefficients in certain 
rank one $R$-modules serves as an obstruction to a space having an \textbf{associative H-space
structure}, and in Section \ref{NovHomology} we discuss how the Novikov Principle implies that 
the \textbf{Novikov homology} of a closed 1-form can be viewed as a special case of Morse homology 
with coefficients in a local system whose fiber is a Novikov ring. Applying this result, we use 
twisted Morse complexes to explicitly compute the \textbf{Novikov numbers} of closed 1-forms
on certain manifolds.

As for other applications, there are at least two areas of mathematics that use Morse homology
with local coefficients on manifolds with boundaries or corners which should be mentioned. 
Although this paper does not discuss twisted Morse homology on manifolds with boundaries or 
corners outside of this introduction, the results, techniques, and proofs used in this paper
should extend to that context with the right definitions and assumptions. 

The first area that deserves mention is \textbf{Floer homology} for 3 and 4-manifolds, and in particular,
the approach based on the Seiberg-Witten monopole equations \cite{KroMon}. In fact, Section 2.4 
of \cite{KroMon} outlines an approach for extending Morse homology to manifolds with boundary where
the gradient of the function is tangent to the boundary, and Section 2.7 of \cite{KroMon} gives an
outline of how to extend Morse homology to include local coefficients. Kronheimer and Mrowka do not 
include proofs in Section 2.7 of their book because, ``In the main part of this book, the Floer 
homology of a 3-manifold will be constructed by taking these constructions of Morse theory and 
repeating them in an infinite dimensional setting." \cite[p. 1]{KroMon} We note that Theorems
\ref{homologyindependence} and \ref{twistedsame} of this paper prove the claims contained 
in Section 2.7 of \cite{KroMon} within the context of closed finite dimensional smooth manifolds.

The second area that deserves mention is \textbf{symplectic cohomology}, and more specifically, the proof
of Viterbo's Theorem given by Abouzaid \cite{AboSym}. Viterbo's Theorem asserts that there is an isomorphism
between the twisted homology of the free loop space $\mathcal{L}Q$ of a closed differentiable manifold
$Q$ and the symplectic cohomology of its cotangent bundle $T^\ast Q$. The proof given by Abouzaid models
the free loop space $\mathcal{L}Q$ as a direct limit of spaces of piecewise geodesics 
$\mathcal{L}^r_{\delta^r} Q$, which are finite dimensional smooth manifolds with corners. The twisted 
Morse homology of Morse functions on the finite dimensional approximations $\mathcal{L}^r_{\delta^r} Q$
whose gradient flows point outward at the boundary is used in Chapter 11  of \cite{AboSym}. In fact,
Proposition 3.13 in Chapter 11 of \cite{AboSym} would be a corollary of Theorems \ref{homologyindependence}
and \ref{twistedsame} of this paper extended to Morse-Smale functions on finite dimensional smooth 
manifolds with corners whose gradient flows point outward at the boundary.


\smallskip
We now outline the results and approaches contained in this paper.

\smallskip
A \textbf{local coefficient system} on a topological space $X$ assigns a fiber $G$, which can be a
group, ring, module, or field, to every point in the space, and to every homotopy class of paths
rel endpoints it assigns a homomorphism in the appropriate category between the fibers of the 
endpoints. We take as our starting point a \textbf{bundle of abelian groups}, which 
is a local coefficient system in the category of abelian groups (Definition \ref{bundleofgroups}).
Most of the results in Sections \ref{MorComSec}\,--\,\ref{singularCW} are proved for bundles of 
abelian groups, but the proofs carry over without change to bundles of rings, modules, and fields, 
cf. Remark \ref{modulebundle}. A bundle of abelian groups $G$ over a pointed space $(X,x_0)$ 
induces a representation
$$
\pi_1(X,x_0) \times G_{x_0} \rightarrow G_{x_0}
$$
on the fiber $G_{x_0}$ over the basepoint $x_0$, which can be viewed as an isomorphism class 
of local coefficient systems (Theorem \ref{actionsystem}). 

The definition of the boundary operator $\partial^G_\ast$ in the twisted Morse-Smale-Witten 
chain complex depends on the homomorphisms associated to the gradient flow lines by the
local coefficient system (Definition \ref{twistedboundary}), but the homology of the twisted
Morse-Smale-Witten chain complex only depends on the isomorphism class of the local coefficient
system (Theorem \ref{homologyindependence}). The distinction between a bundle of abelian groups and 
its isomorphism class is examined in detail in Section \ref{ExamplesSec} with Examples
\ref{circle} and \ref{projectivespace}, where we compute the twisted Morse homology of $S^1$ 
and $\mathbb{R}P^2$ with respect to several different local coefficient systems. 

In Section \ref{EMTheorem} we prove the (untwisted) Morse Homology Theorem for a Covering Space
(Theorem \ref{MorHomCover}) and the \textbf{Morse Eilenberg Theorem} (Theorem \ref{EilenbergMorse}). 
The Morse Eilenberg Theorem relates the twisted Morse homology of a Morse-Smale pair
$(f,\mathsf{g})$ with coefficients in a bundle of abelian groups $G$ to the homology 
of the Morse-Smale-Witten chain complex of $(f,\mathsf{g})$ pulled back to the universal 
cover $\widetilde{M}$ of the manifold and twisted by the representation determined by the 
local coefficient system
$$
\pi_1(M,x_0) \times G_{x_0} \rightarrow G_{x_0}
$$
and the action of $\pi_1(M,x_0)$ on the universal cover $\widetilde{M}$ by deck transformations.
One consequence of the Morse Eilenberg Theorem is that the twisted Morse homology of a
Morse-Smale pair $(f,\mathsf{g})$ can be computed using the representation determined by the
local coefficient system, instead of the homomorphisms associated to the gradient
flow lines by the local coefficient system. This shows that the twisted Morse homology of a 
Morse-Smale pair $(f,\mathsf{g})$ only depends on the isomorphism class of $G$, a fact that is proved 
separately using an entirely different approach in the next section.

Section \ref{continuation} is devoted to proving an \textbf{Invariance Theorem} 
(Theorem \ref{homologyindependence}), which shows that the homology of the Morse-Smale-Witten
chain complex of a Morse-Smale pair $(f,\mathsf{g})$ with coefficients in a bundle of abelian
groups $G$ depends only on the isomorphism class of $G$ and not at all on the Morse-Smale pair
$(f,\mathsf{g})$. The proof of Theorem \ref{homologyindependence} is fairly standard within 
Morse homology and Floer homology, except for the necessary additions required to take into account
local coefficient systems. Theorem \ref{homologyindependence} allows us to pick any 
Morse-Smale pair $(f,\mathsf{g})$ when computing twisted Morse homology, a fact 
that we use throughout the rest of the paper.

Section \ref{singularCW} is devoted to proving the \textbf{Twisted Morse Homology Theorem}
(Theorem \ref{twistedsame}), which says that the homology of the twisted
Morse-Smale-Witten chain complex of $(M,G)$ is isomorphic to the singular homology
of $M$ with coefficients in $G$. In \cite{BanLec} the first two authors presented a proof
of the (untwisted) Morse Homology Theorem using \textbf{Conley index theory} following \cite{SalMor}. 
The Conley index approach uses a filtration of index pairs associated to the gradient flow
which has properties similar to a CW filtration. So, with the Conley index approach it is
unnecessary to address the question of whether or not the unstable manifolds of a Morse-Smale
function determine a CW-structure -- a deep result that is usually proved with additional
assumptions on the Riemannian metric (see the references in Remark \ref{signdegree}). 

Rather than extending the Conley index approach to the Morse Homology Theorem to include local
coefficients, we choose instead in this paper to work directly with the CW-structures determined 
by certain Morse-Smale pairs $(f,\mathsf{g})$. The Invariance Theorem (Theorem \ref{homologyindependence})
shows that there is no loss of generality to the Twisted Morse Homology Theorem 
(Theorem \ref{twistedsame}) if it is proved for a restricted class of Morse-Smale pairs 
$(f,\mathsf{g})$, as long as at least one pair $(f,\mathsf{g})$ in the restricted class exists. 
The main advantage to working with CW-structures, instead of Conley index filtrations, is that we 
can directly compare the twisted Morse-Smale-Witten boundary operator to 
\textbf{Steenrod's CW-boundary operator} with local coefficients
(Lemma \ref{boundarysame}). This approach also illuminates certain technical issues that arise when
extending CW-homology to include local coefficient systems. Specifically, Steenrod's CW-chain
complex with local coefficients is not defined for all CW-complexes. Instead, one must restrict
to a subcategory of CW-complexes where Steenrod's CW-boundary operator is well-defined \cite{SteHom}.
We choose to restrict to the category of \textbf{regular CW-complexes}, which are CW-complexes where the 
attaching maps are all homeomorphisms (Definition \ref{regularCW}).

Regular CW-complexes are much more restrictive than general CW-complexes. For instance, the
CW-structure determined by the usual height function on $S^1$ is not regular (Example \ref{circle}),
although it is possible to deform $S^1$ so that the unstable manifolds of a height function do 
determine a regular CW-structure (Example \ref{deformedcircle}). Similarly, the unstable manifolds 
of a tilted height function on a torus do not determine a regular CW-structure, cf. Example 7.11 
of \cite{BanLec}. However, on every closed finite dimensional smooth manifold $M$ it is always 
possible to find a Morse-Smale pair $(f,\mathsf{g})$ whose unstable manifolds do determine a
regular CW-structure on $M$ (Theorem \ref{regularMorse}). This result is proved by
starting with a triangulation of $M$, fine enough so that each simplex fits inside a coordinate chart,
and then constructing a Morse-Smale pair $(f,\mathsf{g})$ whose unstable manifolds mimic the
triangulation. We expect that this result should be of independent interest in \textbf{combinatorial Morse
theory}, cf. Remark \ref{CombMorse}.

The isomorphism between the twisted Morse homology of a Morse-Smale pair $(f,\mathsf{g})$ whose 
unstable manifolds determine a regular CW-structure and the homology of Steenrod's CW-chain 
complex with local coefficients is proved in Lemma \ref{boundarysame}. We also include for 
completeness a detailed proof of the fact that the homology of Steenrod's CW-chain complex with local 
coefficients is isomorphic to singular homology with local coefficients (Lemma \ref{CWsingular}).
These two lemmas, together with the Invariance Theorem (Theorem \ref{homologyindependence}) 
and Theorem \ref{regularMorse}, constitute our proof of the Twisted Morse Homology Theorem
(Theorem \ref{twistedsame}).

A key step in the proof of Lemma \ref{boundarysame} involves relating the signed count of 
the number of gradient flow lines of a Morse-Smale pair $(f,\mathsf{g})$ to the degree of 
the attaching map in the associated CW-complex, i.e. $\#\mathcal{M}(q,p) = [e^k_q:e^{k-1}_p]$.
This result follows for Morse-Smale pairs $(f,\mathsf{g})$ (with possibly some additional 
assumptions on $\mathsf{g}$) using the manifolds with corners structure on compactified 
moduli spaces of gradient flow lines, cf. Remark \ref{signdegree}, but it is also easy to see 
directly for the specific Morse-Smale pair $(f,\mathsf{g})$ we construct in Theorem \ref{regularMorse}, 
cf. Remark \ref{signdegreespecial}. Thus, the proof of the Twisted Morse Homology Theorem
(Theorem \ref{twistedsame}) given in Section \ref{singularCW} is mostly independent of the 
manifolds with corners structure on the compactified moduli spaces; although the proof of 
the Invariance Theorem (Theorem \ref{homologyindependence}) in Section \ref{continuation} 
does rely heavily on the manifolds with corners structures on various compactified moduli spaces
of gradient flow lines.

In Section \ref{CoLich} we examine some aspects of \textbf{twisted Morse cohomology}.
After defining a general Morse-Smale-Witten cochain complex with coefficients
in a bundle of abelian groups $G$ (Definition \ref{twistedcochain}), we quickly limit our 
discussion to the local coefficient system $G = e^\eta$ defined by a closed 1-form
(Example \ref{etasystem}), i.e. to an $\eta$-twisted Morse-Smale-Witten cochain complex
 (Definition  \ref{etacochain}).  A closed 1-form $\eta$ on a smooth manifold $M$ can be 
used to twist both the Morse-Smale-Witten coboundary operator and also the differential 
in the de Rham complex via
$$
d_\eta \xi = d\xi + \eta \wedge \xi,
$$
where $\xi \in \Omega^\ast(M,\mathbb{R})$. The resulting cohomology groups are known as
the \textbf{Lichnerowicz cohomology} groups $H^\ast_\eta(M)$ (Section \ref{LichLCS}).

Lichnerowicz cohomology is an invariant used to study \textbf{locally conformal symplectic (LCS)}
manifolds.  An LCS manifold is a smooth manifold with a smooth nondegenerate 2-form $\Omega$
which becomes closed locally when multiplied by a smooth positive function (Definition \ref{LCSDef}).
Associated to every LCS form $\Omega$ there is a closed 1-form $\eta$, known as the
\textbf{Lee form}, which satisfies the equation
$$
d\Omega = -\eta \wedge \Omega
$$
(Proposition \ref{LCSLee}). The de Rham cohomology class of the Lee form $\eta$ only depends on 
the conformal class of the LCS form $\Omega$ (Proposition \ref{LCSConformal}), and 
the $\eta$-twisted Morse homology and the $\eta$-twisted Morse cohomology only depend
on the de Rham cohomology class of $\eta$ (Corollary \ref{etadeRhamclass} and
Corollary \ref{cohomologyetainvariant}). Thus, the $\eta$-twisted Morse homology and cohomology
groups are invariants of the conformal class of a locally conformal symplectic form 
$\Omega$ with associated Lee form $\eta$ (Corollary \ref{etaconformalinv}).

Our main result in Section \ref{CoLich} is the \textbf{$\eta$-Twisted Morse de Rham Theorem} 
(Theorem \ref{twistedDeRham}), which says that the $\eta$-twisted Morse cohomology groups 
are isomorphic to the Lichnerowicz cohomology groups defined by $-\eta$. Our proof of
Theorem \ref{twistedDeRham} relies on the \textbf{$\eta$-twisted Morse Poincar\'e Lemma} 
(Lemma \ref{MSWPoincare}) and is modeled on the proof of the de Rham Theorem found in
Section V.9 of \cite{BreTop}. We conclude our discussion of twisted Morse cohomology in 
Section \ref{sheafSec} where we discuss its relationship with \textbf{sheaf cohomology}.
In particular, we discuss a result due to Vaisman \cite{VaiRem} that says that the Lichnerowicz 
cohomology groups $H^\ast_{-\eta}(M)$ are isomorphic to the cohomology of the manifold 
$M$ with coefficients in a sheaf $\mathcal{F}_\eta(M)$ (Theorem \ref{sheafiso}). 

\smallskip
The isomorphisms proved in Section \ref{MorComSec}\,--\,\ref{CoLich} show that homology
with local coefficients and Lichnerowicz cohomology on closed finite dimensional smooth manifolds
can be computed using finitely generated chain and cochain complexes.  This is enough to establish 
the \textbf{invariance of the twisted Euler number} for homology with local coefficients in a bundle of 
finitely generated free $R$-modules  (Theorem \ref{Eulerinvariant}) and Lichnerowicz cohomology  
(Corollary \ref{coEulerindex}) using the \textbf{Euler-Poincar\'e Theorem}
(Theorem \ref{EulerPoincare}) for modules over a principle ideal domain $R$.
Moreover, twisted Morse complexes and cochain complexes often have fewer generators than 
Steenrod's CW-chain complex for regular CW-complexes (compare Examples \ref{circle} 
and \ref{deformedcircle}). Thus, twisted Morse chain and cochain complexes are superb tools
for computing homology and cohomology with local coefficients, a topic we turn to in
Section \ref{App}.

We begin Section \ref{App}  by reviewing a result due to Le\'on, L\'opez, Marrero, and Padr\'on \cite{LeoOnt}
that says that if a nonzero closed $1$-form $\eta$ on a closed smooth manifold $M$ is parallel with 
respect to some Riemannian metric on $M$, then the Lichnerowicz cohomology $H^\ast_\eta(M)$ 
vanishes (Theorem \ref{Leo4.5}). Thus, the $\eta$-twisted Morse cohomology groups serve as obstructions 
to the existence of \textbf{nonzero closed parallel $1$-forms} on closed smooth manifolds 
(Corollary \ref{notparallel}), as well as invariants for the conformal class of a LCS 
form (Corollary \ref{etaconformalinv}). This shows why it might be useful to compute $\eta$-twisted Morse 
cohomology groups, and in Example \ref{genus2Lich} we use twisted Morse cochain complexes to explicitly 
compute the Lichnerowicz cohomology of a surface $S$ of genus two with respect to all closed 1-forms 
$\eta$ on the surface. The invariance of the $\eta$-twisted Euler number and the small number of generators
required for the twisted Morse cochain complex provide strong constraints on the solution. In fact, 
these constraints imply that there are only 4 possibilities for the Lichnerowicz cohomology groups 
$H^\ast_\eta(S)$, no matter which closed 1-form $\eta\in \Omega^1_{\text{cl}}(S,\mathbb{R})$ is 
considered, and our computation shows that 2 of these 4 possibilities actually occur. The approach 
used in Example \ref{genus2Lich} should extend to many other surfaces and manifolds that 
have ``conducive'' CW-structures.

In Section \ref{Hspace} we present a proof of a result due to Albers, Frauenfelder, and Oancea \cite{AlbLoc}
that says that if $X$ is a path connected \textbf{associative H-space} with a local coefficient system 
$\mathcal{L}$ of rank one $R$-modules satisfying certain conditions, then the singular homology of $X$ 
with coefficients in $\mathcal{L}$ vanishes (Proposition \ref{singularvanish}). Thus, the Twisted Morse
Homology Theorem (Theorem \ref{twistedsame}) implies that Morse homology with coefficients in certain local
coefficient systems $\mathcal{L}$ serve as obstructions to a closed smooth manifold $M$ having an associative
H-space structure (Corollary \ref{Hspaceobstruction}). Moreover, the invariance of the twisted Euler number
(Corollary \ref{Eulerclassical}) implies that the Euler number serves as an obstruction to the
existence of an associative H-space structure on  a closed smooth manifold $M$ with 
$H^1_{\text{dR}}(M;\mathbb{R}) \neq 0$ (Corollary \ref{Eulerobstruction}). In Examples 
\ref{circleHspace}\,--\,\ref{RP2notH} we consider the twisted Morse homology of various manifolds 
with coefficients in local systems $\mathcal{L}$ meeting the conditions listed in 
Proposition \ref{singularvanish}. In particular, we use twisted Morse complexes to shows that a surface
of genus two is not an associative H-space, and $\mathbb{R}P^n$ is not an associative H-space when $n$ 
is even (Example 1 of \cite{AlbLoc}).

Section \ref{NovHomology} contains an exposition of Novikov homology \cite{NovMul} \cite{NovThe},
which extends Morse homology from exact 1-forms to closed 1-forms, and its relationship with
twisted Morse homology. Novikov homology is more general than (untwisted) Morse homology; however,
the Novikov Principle (Theorem \ref{NovPrin}) implies that Novikov homology is a special case of
twisted Morse homology (Corollary \ref{NovTwist}). The Novikov chain complex is constructed by pulling
back a closed one form $\zeta \in \Omega^1_{\text{cl}}(M,\mathbb{R})$ to a covering space of the manifold
$M$ where $\zeta$ becomes exact. The flow of the pullback of a generic closed 1-form $\zeta$ to a
covering space where it is exact will be the gradient flow of some Morse-Smale pair. However, if the
\textbf{rank} of $\zeta$ is greater than zero, then the pullback of $\zeta$ will have an infinite number
of zeros (Corollary \ref{rank1cyclic}).

To account for the infinite number of zeros introduced when pulling back the form to a covering
space the \textbf{Novikov ring} (Definition \ref{Novikovring}), which consists of countable 
``half infinite'' Laurent series with integer coefficients and exponents in $\mathbb{R}$, is used.
In Definition \ref{zetasystem} we introduce a local coefficient system $\mathcal{L}_\zeta$
of rank one Nov-modules associated to a closed 1-form $\zeta$. The local coefficient system 
$\mathcal{L}_\zeta$ is suitable for use with twisted Morse homology, and it is in the isomorphism
class of local coefficient systems commonly used for Novikov homology (Claim \ref{zetarep}). 

Using twisted Morse homology with coefficients in the local system $\mathcal{L}_\zeta$ we compute
the \textbf{Novikov numbers} (Definition \ref{NovNum}) of all closed 1-forms on various manifolds
(Examples \ref{Novcircle}\,--\,\ref{Novsurface2}). The Novikov numbers $b_k([\zeta])$ and 
$q_k([\zeta])$ for $k=0,\ldots, m$, which generalize the Betti numbers and the torsion numbers 
of a manifold (Proposition \ref{Novexact}), satisfy the \textbf{Novikov inequalities} 
(Theorem \ref{Novinequalities}), which generalize the Morse inequalities, cf. Theorem 3.33 
of \cite{BanLec}. Hence, the Novikov inequalities give lower bounds on the number of zeros of 
a closed Morse 1-form $\zeta$. In Example \ref{Novtorus} we note that the local system of rank one
Nov-modules $\mathcal{L}_\zeta$ satisfies the conditions listed in Proposition \ref{singularvanish}
whenever $\zeta$ is a non-exact closed 1-form, and since a torus is an associative H-space this 
implies that the Novikov numbers of a non-exact closed 1-form on a torus vanish. Our last
example is a concrete example where the Novikov inequalities give a non-trivial lower bound
on the number of zeros of a non-exact closed 1-form. In Example \ref{Novsurface2} we use a
twisted Morse complex with coefficients in the local system $\mathcal{L}_\zeta$ to compute the 
Novikov numbers of a surface $S$ of genus two with respect to any closed 1-form 
$\zeta\in \Omega^1_{\text{cl}}(S,\mathbb{R})$. Our computation shows that every non-exact closed
Morse 1-form on $S$ must have at least 2 zeros with Morse index 1.


\section{The Morse complex with twisted coefficients}\label{MorComSec}

In this section we construct the Morse-Smale-Witten chain complex with coefficients in 
a bundle of abelian groups, and we present some examples where the homology of
this twisted chain complex can be computed explicitly.

\subsection{Local coefficients}
A bundle of abelian groups $G$ over a topological space $X$ is a functor from the
fundamental groupoid of $X$ to the category of abelian groups. More explicitly, we have
the following.

\begin{definition}\label{bundleofgroups}
A {\bf bundle of abelian groups $G$} over a topological space $X$ associates to every point
$x\in X$  an abelian group $G_x$ and to every continuous path $\gamma:[0,1] \rightarrow X$
a homomorphism $\gamma_\ast:G_{\gamma(1)} \rightarrow G_{\gamma(0)}$ such that the following
conditions are satisfied.
\begin{enumerate}
\item If two paths $\gamma_1,\gamma_2:[0,1] \rightarrow X$ from $x\in X$ to $y \in X$
      are homotopic rel endpoints, then the homomorphisms from $G_y$ to $G_x$ associated
      to $\gamma_1$ and $\gamma_2$ are the same, i.e. $(\gamma_1)_\ast = (\gamma_2)_\ast$.
\item If $\gamma:[0,1] \rightarrow X$ is constant, then $\gamma_\ast$ is the identity.
\item If $\gamma_1,\gamma_2:[0,1] \rightarrow X$ are paths with $\gamma_1(1) = \gamma_2(0)$,
      then $(\gamma_1\gamma_2)_\ast = (\gamma_1)_\ast \circ (\gamma_2)_\ast$, where
      $\gamma_1\gamma_2$ denotes the concatenation of $\gamma_1$ and $\gamma_2$.
\end{enumerate}
\end{definition}

\noindent
Letting $\gamma_2(t) = \gamma_1(1-t)$ in (3), we see that the above conditions imply that
the homomorphism $\gamma_\ast$ associated to a path $\gamma$ is in fact an isomorphism.

\smallskip\noindent
\textbf{Note:} If $G$ is any abelian group and $\gamma_\ast$ is the identity map
for all paths $\gamma:[0,1] \rightarrow X$, then associating $G = G_x$ to every point $x\in X$
determines a \textbf{constant} bundle of abelian groups.

\begin{definition}\label{isobundles}
Suppose that $G_1$ and $G_2$ are both bundles of abelian groups over a topological
space $X$. If there exists a family of isomorphisms $\Phi:G_1 \rightarrow G_2$
such that for every continuous path $\gamma:[0,1] \rightarrow X$ the diagram
$$
\xymatrix{
(G_1)_{\gamma(1)} \ar[r]^{\gamma_\ast^{G_1}} \ar[d]_{\Phi_{\gamma(1)}} & (G_1)_{\gamma(0)} 
  \ar[d]^{\Phi_{\gamma(0)}} \\
(G_2)_{\gamma(1)} \ar[r]^{\gamma_\ast^{G_2}}                           & (G_2)_{\gamma(0)}
}
$$
commutes, then $G_1$ and $G_2$ are said to be \textbf{isomorphic}.
\end{definition}

\smallskip\noindent
\textbf{Note:} A bundle of abelian groups that is isomorphic to a constant bundle 
is called \textbf{simple}. A bundle of abelian groups $G$ is simple if and only if 
for any $x,y \in X$ the homomorphism $\gamma_\ast$ is independent of the path
$\gamma$ from $x$ to $y$.

\smallskip
A bundle of abelian groups $G$ over a pointed space $(X,x_0)$ induces a representation
$$
\pi_1(X,x_0) \times G_{x_0} \rightarrow G_{x_0},
$$ 
i.e. $[\gamma] \in \pi_1(X,x_0)$ determines an isomorphism  $\gamma_\ast:G_{x_0}
\rightarrow G_{x_0}$ and $(\gamma_1\gamma_2)_\ast = (\gamma_1)_\ast \circ (\gamma_2)_\ast$
for any $[\gamma_1],[\gamma_2] \in \pi_1(X,x_0)$. The following converse is well known, 
cf. Theorem VI.1.11 and VI.1.12 of \cite{WhiEle}.
\begin{theorem}\label{actionsystem}
Let $X$ be a connected topological space with a basepoint $x_0$, and let $G_0$ be an abelian 
group on which $\pi_1(X,x_0)$ operates. Then there exists a bundle of abelian groups $G$ 
over $X$ with $G_{x_0} = G_0$ which induces the operation of $\pi_1(X,x_0)$ on $G_0$, and
the bundle $G$ is unique up to isomorphism.
\end{theorem}

\begin{remark}
The boundary operator in the twisted Morse-Smale-Witten chain complex
(Definition \ref{twistedboundary}) depends on the choice of the bundle of abelian groups $G$.
However, we will show that the homology of the twisted Morse-Smale-Witten complex only
depends on the isomorphism class of $G$, i.e. only on the representation
$
\pi_1(X,x_0) \times G_{x_0} \rightarrow G_{x_0}
$ 
(Theorem \ref{EilenbergMorse} and Theorem \ref{homologyindependence}).
\end{remark}


\begin{example}[The local coefficient system $e^\eta$ determined by a closed $1$-form]
\label{etasystem}
Let $\eta \in \Omega^1_{\text{cl}}(M,\mathbb{R})$ be a closed smooth real valued $1$-form on 
a  finite dimensional smooth manifold $M$. To each point $x \in M$ associate the
additive abelian group $\mathbb{R}$, and to each smooth path $\gamma:[0,1] \rightarrow M$
associate the homomorphism $\gamma_\ast:\mathbb{R}_{\gamma(1)} \rightarrow \mathbb{R}_{\gamma(0)}$
defined by
$$
\gamma_\ast(s) = e^{\int_1^0 \gamma^\ast(\eta)} \cdot s \quad \text{for all } s \in \mathbb{R}.
$$
Since every continuous path in $M$ is homotopic rel endpoints to a smooth path,
Stokes' Theorem shows that this defines a bundle of (additive) $\mathbb{R}$ groups over
$M$. The above definition of $\gamma_\ast$ extends to paths $\gamma:\overline{\mathbb{R}}
\rightarrow M$ using any diffeomorphism $\overline{\mathbb{R}} \approx [0,1]$.  We will denote
this flat line bundle by $e^\eta$.

\begin{claim}\label{etaisomorphic}
If $\eta_1,\eta_2 \in \Omega^1_{\text{cl}}(M,\mathbb{R})$ are in the same de Rham cohomology
class, then $e^{\eta_1}$ is isomorphic to $e^{\eta_2}$.
\end{claim}

\proofstart
By assumption there exists a smooth function $h:M \rightarrow \mathbb{R}$ with
$\eta_1 - \eta_2 = dh$.  Define a family of isomorphisms $\Phi:e^{\eta_1} \rightarrow
e^{\eta_2}$ by $\Phi_x(s) = e^{-h(x)} \cdot s$ for all $x \in M$ and $s \in \mathbb{R}$.
Then the following diagram commutes for any path $\gamma:[0,1] \rightarrow \mathbb{R}$
$$
\xymatrix{
\mathbb{R} \ar[rr]^{\times e^{\int_1^0 \gamma^\ast(\eta_1)}} 
 \ar[d]_-{\times e^{-h(\gamma(1))}} & & \mathbb{R} 
 \ar[d]^-{\times e^{-h(\gamma(0))}} \\
\mathbb{R} \ar[rr]^{\times e^{\int_1^0 \gamma^\ast(\eta_2)}} & & \mathbb{R}
}
$$
because $e^{\int_1^0 \gamma^\ast(\eta_1)} = e^{\int_1^0 \gamma^\ast(\eta_2+ dh)}
= e^{\int_1^0 \gamma^\ast(\eta_2)} e^{h(\gamma(0)) - h(\gamma(1))}$.
\proofend

\end{example}

\begin{remark}\label{modulebundle}
Note that $e^\eta$ is not only a bundle of abelian groups with respect to addition on the fibers,
but it is also a bundle of commutative rings, a rank one bundle of $\mathbb{R}$-modules, 
and a flat $\mathbb{R}$-vector bundle of dimension one, also known as a {flat line bundle}.
That is, the homomorphisms $\gamma_\ast(s) = e^{\int_1^0 \gamma^\ast(\eta)} \cdot s$ are
isomorphisms in each of those categories.  Moreover, the notion of isomorphic bundles from Definition
\ref{isobundles} carries over to those other categories by requiring $\Phi$ to be an
isomorphism in the category under consideration, and the proof of Claim \ref{etaisomorphic}
carries over to those other categories without modification.
\end{remark}

\begin{remark}\label{aeta}
We can also define a local coefficient system $a^\eta$ for any $a > 0$ by replacing 
$e^{\int_1^0 \gamma^\ast(\eta)}$ with $a^{\int_1^0 \gamma^\ast(\eta)}$ in the preceding
example.  The proof of Claim \ref{etaisomorphic} is easily modified for the system $a^\eta$.
However, $a_1^\eta$ may not be isomorphic to $a_2^\eta$ if $a_1 \not= a_2$.
\end{remark}


\subsection{Path components of compactified moduli spaces}\label{pathcomponents}
Let $f:M \rightarrow \mathbb{R}$ be a smooth Morse-Smale function on a closed finite 
dimensional smooth Riemannian manifold $M$. 
Let $Cr(f) = \{p \in M |\, df_p = 0 \}$ denote the set of critical points of $f$, let
$Cr_k(f) \subset Cr(f)$ denote the set of critical points of index $k$, and for
any $p \in Cr(f)$ let $\lambda_p$ denote the index of $p$. For $p,q \in Cr(f)$ let 
$W^u(q)\subset M$ be the unstable manifold of $q$, $W^s(p) \subset M$ the stable manifold
of $p$, and define
$$
W(q,p) = W^u(q) \cap W^s(p) \subset M.
$$
If this space is nonempty, then one says that $q$ is succeeded by $p$, i.e. $q \succeq p$. 
In this case, $W(q,p)$ is a noncompact smooth manifold of dimension $\lambda_q - \lambda_p$.
Choosing orientations for the unstable manifolds $W^u(q)$ for all $q \in Cr(f)$ determines
an orientation on $W(q,p)$ for all $p,q \in Cr(f)$ via the short exact sequence
$$
\xymatrix{
0 \ar[r] & T_\ast W(q,p) \ar@{^{(}->}[r] & T_\ast W^u(q) |_{W(q,p)} \ar[r] & 
\nu_\ast(W(q,p), W^u(q))|_{W(q,p)} \ar[r] & 0,
}
$$
where the fibers of the normal bundle are canonically isomorphic to $T_p W^u(p)$ via the
gradient flow. Taking a quotient by the action of $\mathbb{R}$ given by the flow of 
$-\nabla f$ then gives a smooth manifold
$$
\mathcal{M}(q,p) = W(q,p)/\mathbb{R}
$$
of dimension $\lambda_q - \lambda_p - 1$, which we orient as follows. For any regular
value $y$ between $f(p)$ and $f(q)$ we identify $\mathcal{M}(q,p) = W(q,p) \cap f^{-1}(y)$,
and for any $x \in W(q,p) \cap f^{-1}(y)$ we declare $B_x$ to be a positive basis for
$T_x \mathcal{M}(q,p)$ if and only if $(-(\nabla f)(x),B_x)$ is a positive basis for
$T_x W(q,p)$.  (See Section 6.1 of \cite{QinOnm} or Proposition 3.10 of \cite{WebThe}
for more details.)


The moduli space $\mathcal{M}(q,p)$ has a compactification $\overline{\mathcal{M}}(q,p)$
consisting of the piecewise gradient flow lines from $q$ to $p$, which can be given the
structure of a smooth manifold with corners \cite{BFKOnt} \cite{BurOnt}
\cite{LatExi} \cite{QinOnm} \cite{QinAna}. As in Section 6.1 of \cite{QinOnm}, we orient
the (codimension) $1$-stratum using the convention that an outward pointing normal vector
field followed by a positive basis for a tangent space of 
$\partial^1\overline{\mathcal{M}}(q,p)$ should be a positive basis for a tangent space
of $\overline{\mathcal{M}}(q,p)$.

A piecewise gradient flow line from $q$ to $p$ can be identified with its image in $M$, which
is an element of $\mathcal{P}^c(M)$, the space of all nonempty closed subsets of $M$ with the
Hausdorff topology. This identification is compatible with the topology of the smooth manifold
with corners $\overline{\mathcal{M}}(q,p)$ in the sense that the map that sends an element of
$\nu \in \overline{\mathcal{M}}(q,p)$ to its image $Im(\nu)$ is a homeomorphism onto its image
$Im(\overline{\mathcal{M}}(q,p))$ in $\mathcal{P}^c(M)$ \cite{BanCas} \cite{HurFlo}.

Denote the path component of $\nu \in \overline{\mathcal{M}}(q,p)$ by 
$\overline{\mathcal{M}}(q,p;[\nu])$. Let $(\gamma_1, \ldots , \gamma_l)$ be a sequence
of gradient flow lines with 
\begin{eqnarray*}
\lim_{t \rightarrow -\infty}\gamma_1(t) & = & q,\\
\lim_{t \rightarrow \infty}\gamma_j(t) & = & \lim_{t \rightarrow -\infty}\gamma_{j+1}(t)
\text{ for all }j = 1, \ldots, l-1,\\
\lim_{t \rightarrow \infty} \gamma_l(t) & = & p,
\end{eqnarray*}
and denote the corresponding sequence of unparameterized gradient flow lines by
$(\nu_1,\ldots, \nu_l)$.  We will write $[(\nu_1,\ldots, \nu_l)] = [\nu]$ to indicate
that the image of the piecewise gradient flow line $(\nu_1,\ldots, \nu_l)$
$$
Im(\nu_1,\ldots, \nu_l) = 
Im(\gamma_1,\ldots, \gamma_l) = 
\bigcup_{j=1}^l \gamma_j(\overline{\mathbb{R}}) \in \mathcal{P}^c(M)
$$
is in the same path component as $Im(\nu)$ in $\mathcal{P}^c(M)$.  

\begin{lemma}\label{componentboundary}
Let $r,p \in Cr(f)$. If $\nu \in \mathcal{M}(r,p)$, then the closure of
$\mathcal{M}(r,p;[\nu])$ in $\overline{\mathcal{M}}(r,p)$ consists of the piecewise
gradient flow lines from $r$ to $p$ that are in the same path component as $\nu$.
Moreover, when $\lambda_r - \lambda_p = 2$ we have
$$
\partial^1\overline{\mathcal{M}}(r,p;[\nu]) \ = \partial\overline{\mathcal{M}}(r,p;[\nu]) \ =
(-1) \bigcup_{\stackrel{r \succeq q \succeq p}
{[\nu]=[(\nu_1,\nu_2)]}} \mathcal{M}(r,q;[\nu_1]) \times \mathcal{M}(q,p;[\nu_2])
$$
as oriented manifolds. Thus when $\lambda_r - \lambda_p = 2$, 
$$
\sum_{r \succeq q \succeq p} \sum_{\stackrel{[\nu]=[(\nu_1,\nu_2)]}{(\nu_1,\nu_2) \in
\mathcal{M}(r,q) \times \mathcal{M}(q,p)}} \epsilon(\nu_1) \epsilon(\nu_2) = 0
$$
where $\epsilon(\nu_j) = \pm 1$ is the sign of the zero dimensional oriented manifold
$\nu_j$ for $j=1,2$.
\end{lemma}

\proofstart
The boundary of the smooth manifold with corners $\overline{\mathcal{M}}(r,p)$ consists of
the broken gradient flow lines from $r$ to $p$, and the interior $\mathcal{M}(r,p)$ consists
of the (nonbroken) gradient flow lines from $r$ to $p$. Thus, the closure of the path component
$\mathcal{M}(r,p;[\nu])$ consists of the piecewise gradient flow lines in 
$\overline{\mathcal{M}}(r,p)$ that are in same path component as $\nu$. 

The second statement follows immediately from Theorem 3.6 of \cite{QinOnm} and 
Theorem 8.1 of \cite{QinAna}. The last statement follows because 
$\overline{\mathcal{M}}(r,p;[\nu])$ is a compact connected oriented $1$-dimensional
smooth manifold with boundary when $\lambda_r - \lambda_p = 2$ (thus diffeomorphic to
$S^1$ or $[0,1]$), and the oriented sum of the signs associated to the boundary points of a
$1$-dimensional compact smooth manifold is zero.
\proofend

\begin{remark}
Formulas similar to those in the preceding lemma can be found in Lemma 3.4 of \cite{AusMor},
Proposition 5.2 of \cite{BFKOnt}, Theorem 2 of \cite{BurOnt}, Sections 2.14 and 2.15
of \cite{LatExi}, Lemma 4.3 of \cite{SchMor}, and Proposition 4.4 of \cite{WebThe}. 
However, we will following the sign conventions in \cite{QinOnm} and \cite{QinAna}. 
\end{remark}


\subsection{Twisting the Morse-Smale-Witten boundary operator}

Let $G$ be a bundle of abelian groups over $M$ and let $\gamma^{\nu}:[0,1] \rightarrow M$
be a continuous path from $p$ to $q$ whose image coincides with the image of some element 
$\nu \in \overline{\mathcal{M}}(q,p;[\tilde{\nu}])$, where $q,p \in Cr(f)$ and $\tilde{\nu}
\in \overline{\mathcal{M}}(q,p)$. Lemma \ref{componentboundary} implies that any two
paths representing elements of a path component $\overline{\mathcal{M}}(q,p;[\tilde{\nu}])$
are homotopic rel endpoints, and hence condition (1) of Definition \ref{bundleofgroups}
implies that there is a well-defined homomorphism $\gamma^{\nu}_\ast: G_q \rightarrow G_p$
which is independent of the element $\nu \in \overline{\mathcal{M}}(q,p;[\tilde{\nu}])$.

\begin{definition}[Twisted Morse-Smale-Witten chain complex]\label{twistedboundary}
Let $f:M \rightarrow \mathbb{R}$ be a smooth Morse-Smale function on a closed 
smooth Riemannian manifold $(M,\mathsf{g})$ of dimension $m < \infty$. Fix orientations on the
unstable manifolds of $(f,\mathsf{g})$, and let $G$ be a bundle of abelian groups over $M$. The
\textbf{Morse-Smale-Witten chain complex with coefficients in $G$}
is defined to be the chain complex $(C_\ast(f;G),\partial_\ast^G)$ where
$$
C_k(f;G) \stackrel{def}{=} \left.\left\{\sum_{q \in Cr_k(f)} gq \right|\ g \in G_q \right\} 
\approx \bigoplus_{q \in Cr_k(f)} G_q
$$
for all $k=0,\ldots ,m$, and the homomorphism $\partial_k^G:C_k(f;G) \rightarrow 
C_{k-1}(f;G)$ is defined on an elementary chain $gq \in C_k(f;G)$ to be
$$
\partial_k^G(gq)\ \ = \sum_{p \in Cr_{k-1}(f)} \sum_{\nu \in \mathcal{M}(q,p)} \epsilon(\nu)
\gamma^{\nu}_\ast(g)p, 
$$
where $\gamma^{\nu}:[0,1] \rightarrow M$ is any continuous path from $p$ to $q$ whose
image coincides with the image of $\nu \in \mathcal{M}(q,p)$ and $\epsilon(\nu) = \pm 1$ 
is the sign determined by the orientation on $\mathcal{M}(q,p)$. 
\end{definition}

\smallskip\noindent
Note: Since $G_p$ is abelian, there is an action $\mathbb{Z} \times G_p \rightarrow G_p$ that
sends $-1 \cdot g$ to the inverse of $g$.  Also, if $G=\mathbb{Z}$ is a constant bundle, then
$\gamma^{\nu}_\ast = id$ for all $p,q\in Cr(f)$ and the above reduces to the 
Morse-Smale-Witten chain complex with coefficients in $\mathbb{Z}$, cf. Chapter 7 
of \cite{BanLec}.  We will sometimes drop the $G$ in the notation
$\partial^G_\ast$ when the choice of the bundle of coefficients is clear.

\begin{definition}[$\eta$-Twisted Morse-Smale-Witten chain complex]\label{etatwisted}
Let $f:M \rightarrow \mathbb{R}$ be a smooth Morse-Smale function on a closed finite
dimensional smooth Riemannian manifold $(M,\mathsf{g})$. Fix orientations on the unstable
manifolds of $(f,\mathsf{g})$, and let $\eta \in \Omega^1_{\text{cl}}(M,\mathbb{R})$. 
The Morse-Smale-Witten chain complex with coefficients in the local system $e^\eta$ is called
the  \textbf{$\eta$-twisted Morse-Smale-Witten chain complex}. In this case, $C_k(f;e^\eta)
\approx C_k(f) \otimes \mathbb{R}$, where $C_k(f)$ is the free abelian group generated by
the critical points $q$ of index $k$, and the homomorphism $\partial_k^\eta:C_k(f) \otimes
\mathbb{R} \rightarrow C_{k-1}(f) \otimes \mathbb{R}$ is given on a critical point
$q \in Cr_k(f)$ by
$$
\partial_k^\eta(q)\ \ = \sum_{p \in Cr_{k-1}(f)} \sum_{\nu \in \mathcal{M}(q,p)}
\epsilon(\nu) \exp{\left(\int_{\overline{\mathbb{R}}} \gamma^\ast_\nu(\eta) \right)} p,
$$
where $\gamma_\nu:\overline{\mathbb{R}} \rightarrow M$ is any gradient flow line
from $q$ to $p$ parameterizing $\nu \in \mathcal{M}(q,p)$ and $\epsilon(\nu) = 
\pm 1$ is the sign determined by the orientation on $\mathcal{M}(q,p)$. 
\end{definition}

\smallskip\noindent
\textbf{Note:} Example \ref{etasystem} shows that the $\eta$-twisted Morse-Smale-Witten
chain complex is a special case of the general twisted Morse-Smale-Witten chain
complex in Definition \ref{twistedboundary}.

\begin{lemma}\label{boundarysquared}
The pair $(C_\ast(f;G),\partial_\ast^G)$ is a chain complex, i.e. 
$(\partial_\ast^G)^2 = 0$.
\end{lemma}

\proofstart
Let $r \in Cr(f)$ with $\lambda_r = k+1$ for some $k=1,\ldots,m-1$, where $m = \text{dim }M$.
For any $g \in G_r$ we have
$$
\begin{array}{lcl}
\partial_k^G(\partial_{k+1}^G(gr)) & = & 
      \displaystyle \partial^G_k\left(\sum_{q \in Cr_k(f)} \sum_{\nu_1 \in 
      \mathcal{M}(r,q)} \epsilon(\nu_1) \gamma^{\nu_1}_\ast(g) q \right)\\
& = & \displaystyle \sum_{q \in Cr_k(f)} \sum_{\nu_1 \in \mathcal{M}(r,q)} 
      \epsilon(\nu_1)  \partial^G_k \left( \gamma^{\nu_1}_\ast(g)q \right)\\
& = & \displaystyle \sum_{q \in Cr_k(f)} \sum_{\nu_1 \in \mathcal{M}(r,q)} 
      \epsilon(\nu_1) \sum_{p \in Cr_{k-1}(f)} \sum_{\nu_2 \in \mathcal{M}(q,p)} 
      \epsilon(\nu_2) \gamma^{\nu_2}_\ast(\gamma^{\nu_1}_\ast(g)) p\\
& = & \displaystyle \sum_{p \in Cr_{k-1}(f)} \sum_{q \in Cr_k(f)} 
      \sum_{\nu_1 \in \mathcal{M}(r,q)} \sum_{\nu_2 \in \mathcal{M}(q,p)} 
      \epsilon(\nu_1)\epsilon(\nu_2) \gamma^{(\nu_1,\nu_2)}_\ast(g) p.\\
\end{array}
$$

Now consider the coefficient in front of some fixed $p \in Cr_{k-1}(f)$.
\begin{eqnarray*}
\text{coef}(p) 
& = & \sum_{q \in Cr_k(f)} \sum_{(\nu_1,\nu_2) \in \mathcal{M}(r,q) \times \mathcal{M}(q,p)}
      \epsilon(\nu_1)\epsilon(\nu_2)\gamma^{(\nu_1,\nu_2)}_\ast(g).
\end{eqnarray*}
We can group the terms in the above sum according to the various path components
$\overline{\mathcal{M}}(r,p;[\nu])$ and use the fact that the homomorphism
$\gamma^{(\nu_1,\nu_2)}_\ast = \gamma^\nu_\ast$ on each path component to get 
terms of the form
$$
\gamma^\nu_\ast(g) \sum_{q \in Cr_k(f)}
\sum_{\stackrel{[\nu] = [(\nu_1,\nu_2)]}{(\nu_1,\nu_2) \in \mathcal{M}(r,q) \times
\mathcal{M}(q,p)}} \epsilon(\nu_1) \epsilon(\nu_2).
$$
These terms are zero by Lemma \ref{componentboundary}.
\proofend



\subsection{Examples}\label{ExamplesSec}

\begin{example}[A circle]\label{circle}
Consider the height function $f:S^1 \rightarrow \mathbb{R}$ on the unit circle 
$S^1\subset \mathbb{R}^2$ with a critical point $q$ of index $1$ and a critical point 
$p$ of index $0$. Orient the unstable manifold of $q$ clockwise and the unstable manifold 
of $p$ as $+1$.
\begin{figure}[h]
\includegraphics{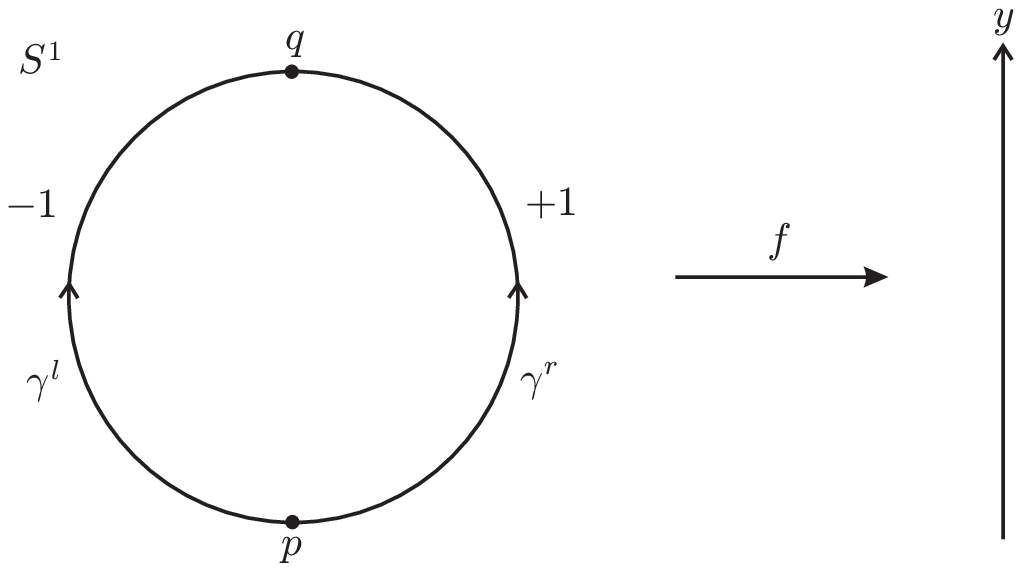}
\end{figure}
The (non-twisted) Morse-Smale-Witten chain complex of $f$ is 
$$
\xymatrix{
0 \ar[r] & C_1(f) \ar[r]^-{\partial_1}\ar@{<->}[d]^{\approx} & C_0(f) \ar@{<->}[d]^{\approx}
\ar[r] & 0\\
0 \ar[r] & <q> \ar[r]^-{\partial_1} & <p> \ar[r] & 0
}
$$
with $\partial_1(q) = 0$ zero since the two gradient flow lines have opposite orientations.
(See Example 7.7 of \cite{BanLec} for more details.)

If $G$ is a bundle of abelian groups over $S^1$ and $g \in G_q$, then
$$
\partial_1^G(gq)  = \left( \gamma^r_\ast(g) - \gamma^l_\ast(g) \right) p
$$
where $\gamma^r$ is a parameterization of the right half of the circle and $\gamma^l$ is
a parameterization of the left half of the circle, both from $p$ to $q$. Thus, for $k=0,1$
$$
H_k((C_\ast(f;G),\partial^G_\ast)) \approx
H_k(S^1;G_q) \text{ if } \gamma^r_\ast(g) = \gamma^l_\ast(g) \text{ for all } g \in G_q.
$$

However, it is possible to have a bundle of abelian groups over $S^1$ where
$\gamma_\ast^r \neq \gamma_\ast^l$. To see this, recall that a representation
$$
\pi_1(S^1,q) \times G_q \rightarrow G_q
$$
determines a bundle of abelian groups over $S^1$ that is unique up to isomorphism.
The bundle of abelian groups $G$ is defined by associating the group $G_q$ to every
point and arbitrarily fixing a homotopy class of paths rel endpoints $[\gamma_{qx}]$ from
$q$ to $x$ for every $x \in S^1$. If $\gamma$ is a path from $x_0 \in S^1$ to $x_1 \in S^1$ and 
$\gamma_{qx_0}$ and $\gamma_{qx_1}$ are paths from $q$ to $x_0$ and from $q$ to $x_1$
in the chosen homotopy classes, then the concatenation $\gamma_{qx_0}\gamma\gamma_{qx_1}^{-1}$
represents an element in $\pi_1(S^1,q)$ and the homomorphism $\gamma_\ast:G_{x_1}
\rightarrow G_{x_0}$ is defined to be the associated homomorphism determined by the representation.
(For more details see the proof of Theorem VI.1.12 in \cite{WhiEle}.)

For instance, suppose $G_q = \mathbb{Z}$ and $\pi_1(S^1,q)$ is identified with
$\mathbb{Z}$ by sending the clockwise generator to $1$.  A representation
$$
\pi_1(S^1,q) \times \mathbb{Z} \rightarrow \mathbb{Z}
$$
is then determined by whether $(1,1) \mapsto 1$ or $(1,1) \mapsto -1$, since the only group
automorphisms of $\mathbb{Z}$ are $\pm$id.
If $(1,1) \rightarrow 1$, then the representation is trivial and the homomorphism associated
to any path is the identity.  So, assume that $(1,1) \mapsto -1$, which implies
that $(n,g) \mapsto (-1)^n g$ for all $(n,g) \in \pi_1(S^1,q) \times \mathbb{Z}$.

For the chosen homotopy classes of paths $[\gamma_{qx}]$ from $q$ to $x \in S^1$ we will
first take the constant path when $x=q$ and choose paths that wrap clockwise $n$-times 
around $S^1$ and then continue clockwise to $x$ whenever $x \neq q$. Taking $x_0=p$ and
$x_1=q$ in the above we see that $\gamma_\ast^r(g) = (-1)^n g$ and $\gamma_\ast^l(g) = 
(-1)^{n+1} g$, and hence $\partial^G_1(gq) = \left((-1)^n - (-1)^{n+1}\right)gp = 2(-1)^n gp$. 
Therefore,
$$
H_1((C_\ast(f;G),\partial^G_\ast)) \approx 0, \quad
H_0((C_\ast(f;G),\partial^G_\ast))  \approx \mathbb{Z}_2.
$$
Alternately, we could choose the constant path when $x=q$ and paths that wrap counterclockwise
$n$-times around $S^1$ and then continue counterclockwise to $x\in S^1$ when $x \neq q$ for the 
homotopy classes of paths. Then $\gamma_\ast^r(g) = (-1)^{-(n+1)}g$ and $\gamma_\ast^l(g) = 
(-1)^{-n} g$, and hence $\partial^G_1(gq) = \left((-1)^{-(n+1)} - (-1)^{-n}\right)gp = 
2(-1)^{-(n+1)}gp$. The reader can verify that alternate choices for the homotopy class of paths
from $q$ to $q$ yield similar boundary operators. Hence, we see that although $\partial_\ast$
depends on the specific bundle of abelian groups $G$ over $S^1$, the homology of $(C_\ast(f;G), 
\partial^G_\ast)$ depends only on the isomorphism class of $G$.

\smallskip
As another example of a bundle of abelian groups over $S^1$, consider a closed $1$-form 
$\eta$ on $S^1$, its associated flat line bundle $e^\eta$, and the associated $\eta$-twisted
Morse-Smale-Witten boundary operator 
$$
\partial_1^\eta(q) = \left(\exp\left(\int_1^0 (\gamma^r)^\ast (\eta)
\right) -\ \exp\left( \int_1^0 (\gamma^l)^\ast (\eta) \right) \right)p.
$$
If $\eta = dh$ is exact, then the integral of $\eta$ along any path from $q$ to $p$ is 
$h(q) - h(p)$. Hence, 
$$
\partial^\eta_1(q) = \left(e^{h(q) - h(p)} -\ e^{h(q) - h(p)} \right) p = 0,
$$
and $H_\ast((C_\ast(f;\mathbb{R}), \partial^\eta_\ast)) = H_\ast(S^1;\mathbb{R})$.
However, if $\eta$ is not exact, then $\int_1^0 (\gamma^r)^\ast (\eta)$ is not equal
to $\int_1^0 (\gamma^l)^\ast (\eta)$. In this case $\partial_1^\eta(q) \neq 0$, and 
$H_k((C_\ast(f;\mathbb{R}), \partial^\eta_\ast)) = 0$ for all $k$. Explicitly, consider 
the form
$$
d\theta = \frac{1}{x^2+y^2}(-y dx + x dy)
$$
and the parameterization of $S^1$ given by $\gamma(t) = (\cos{t},\sin{t})$. Then we have
$$
\int_1^0 (\gamma^r)^\ast(d\theta) = \int_{\pi/2}^{-\pi/2} \sin^2{t} + \cos^2{t}\ dt = -\pi
$$
and
$$
\int_1^0 (\gamma^l)^\ast(d\theta) = \int_{\pi/2}^{3\pi/2} \sin^2{t} + 
\cos^2{t}\ dt = \pi.
$$
Thus, $\partial^\eta_1(q) = (e^{-\pi} - e^{\pi}) p \neq 0$, and $H_k((C_\ast(f) \otimes
\mathbb{R},\partial_\ast^\eta)) \approx 0$ for all $k$.  
\end{example}


\begin{example}[Real projective space]\label{projectivespace}
The real projective space $M = \mathbb{R}P^2$ can be viewed as $S^2 \subset \mathbb{R}^3$
with diametrically opposed points identified or as the closed disk $D^2 \subset \mathbb{R}^2$
with diametrically opposed points on the boundary identified. The Morse-Smale function
$\tilde{f}:S^2 \rightarrow \mathbb{R}$ defined by
$$
\tilde{f}(x_1,x_2,x_3) = x_2^2 + 2x_3^2
$$
satisfies $\tilde{f}(-x_1,-x_2,-x_3) = \tilde{f}(x_1,x_2,x_3)$, and hence it descends to a Morse-Smale
function $f:\mathbb{R}P^2 \rightarrow \mathbb{R}$ with three critical points $p$, $q$, $r$ of 
index $0$, $1$, and $2$ respectively, whose unstable manifolds give a CW-structure with three cells 
$e_0$, $e_1$, $e_2$. (Compare with Example \ref{projectivespacecover} and Example 3.7 of \cite{BanLec}.)
The unstable manifolds, the gradient flow lines, and their orientations are
indicated in the following diagram. The convention for the orientations is given by
the short exact sequence
$$
\xymatrix{
0 \ar[r] & T_\ast W(q,p) \ar@{^{(}->}[r] & T_\ast W^u(q) |_{W(q,p)} \ar[r] & 
\nu_\ast(W(q,p), W^u(q))|_{W(q,p)} \ar[r] & 0
}
$$
where $W(q,p) = W^u(q) \pitchfork W^s(p)$ for any two critical points $q$ and $p$,
and the fibers of the normal bundle are canonically isomorphic to $T_p^u W^u(p)$.

\begin{figure}[h]
\includegraphics{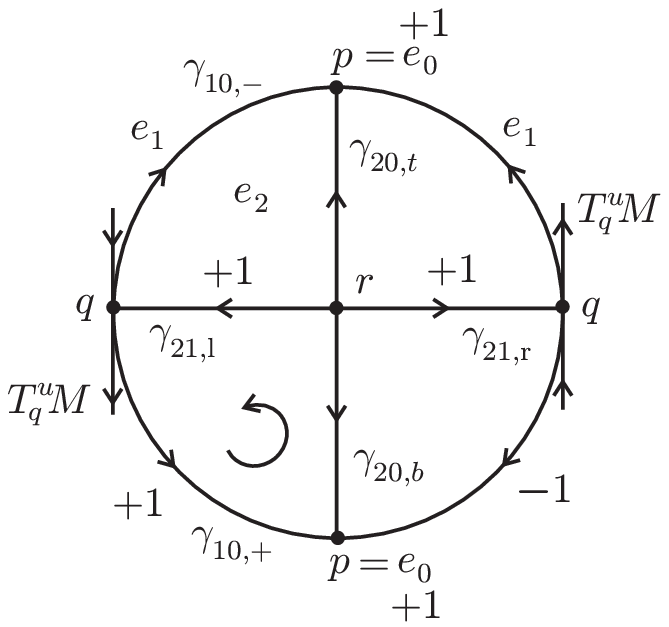}
\end{figure}

The (non-twisted) Morse-Smale-Witten chain complex of $f$ is 
$$
\xymatrix{
0 \ar[r] & C_2(f) \ar[r]^-{\partial_2} \ar@{<->}[d]^{\approx} & C_1(f) \ar[r]^-{\partial_1}
 \ar@{<->}[d]^{\approx} & C_0(f) \ar@{<->}[d]^{\approx} \ar[r] & 0\\
0 \ar[r] & <r> \ar[r]^-{\partial_2} & <q> \ar[r]^-{\partial_1} & <p> \ar[r] & 0
}
$$
with $\partial_2(r) = 2 q$ and $\partial_1(q) = 0$.  Thus, 
$$
H_2((C_\ast(f),\partial_\ast)) \approx 0, \quad H_1((C_\ast(f),\partial_\ast)) 
\approx \mathbb{Z}_2, \quad H_0((C_\ast(f),\partial_\ast)) \approx \mathbb{Z}.
$$ 

Now recall that $\pi_1(\mathbb{R}P^2) \approx \mathbb{Z}_2$, since its universal cover is $S^2$, 
and consider a representation
$$
\pi_1(\mathbb{R}P^2,r) \times \mathbb{Z} \rightarrow \mathbb{Z}.
$$
Denote the gradient flow lines in the above diagram from $r$ to $q$ by $\gamma_{21,r}$
and $\gamma_{21,l}$, from $r$ to $p$ by $\gamma_{20,t}$ and $\gamma_{20,b}$,
and the flow lines from $q$ to $p$ by $\gamma_{10,+}$ and $\gamma_{10,-}$. Denote 
the flow lines in the opposite direction by reversing the subscripts. Choose the constant path at
$r$ and the homotopy classes of paths rel endpoints of $\gamma_{21,r}$ and $\gamma_{20,t}$ for the
bundle of abelian groups $G$ over $\mathbb{R}P^2$ corresponding to the representation. Note that with
this notation the concatenation $\gamma_{21,r}\gamma_{12,l}$ is a generator for 
$\pi_1(\mathbb{R}P^2,r)$.

If $(1,1) \mapsto 1$ under the above action, then the homomorphism associated to
any path is the identity. Hence, the bundle of abelian groups $G$ over $\mathbb{R}P^2$ 
determined by the representation is the constant bundle $\mathbb{Z}$, and the Morse-Smale-Witten
chain complex with coefficients in $G$ reduces to the above (non-twisted) Morse-Smale-Witten
chain complex. So, assuming that the action sends $(1,1) \mapsto -1$, the relevant
homomorphisms in the associated bundle of abelian groups are as follows. 
\begin{eqnarray*}
(\gamma_{12,r})_\ast(1) & = &  +1\\
(\gamma_{12,l})_\ast(1) & = & -1\\
(\gamma_{01,+})_\ast(1) & = & +1\\
(\gamma_{01,-})_\ast(1) & = & -1\\
\end{eqnarray*}
Therefore,
$$
H_2((C_\ast(f;G),\partial^G_\ast)) \approx \mathbb{Z}, \quad
H_1((C_\ast(f;G),\partial^G_\ast)) \approx 0, \quad
H_0((C_\ast(f;G),\partial^G_\ast))  \approx \mathbb{Z}_2.
$$
The reader can verify that different choices for the homotopy classes of paths rel endpoints
from $r$ to $q$ and from $r$ to $p$ yield the same homology.  Hence, the homology of 
$(C_\ast(f;G),\partial^G_\ast)$ only depends on the isomorphism class of the bundle of abelian 
groups over $\mathbb{R}P^2$ when $G_q = \mathbb{Z}$, rather than the specific bundle 
of abelian groups $G$.

\smallskip
Now let $\eta \in \Omega^1_{\text{cl}}(\mathbb{R}P^2,\mathbb{R})$, and consider
the $\eta$-twisted Morse-Smale-Witten chain complex $(C_\ast(f)\otimes \mathbb{R},
\partial_\ast^\eta)$. For the $\eta$-twisted Morse-Smale-Witten boundary operator
we have the following.
\begin{eqnarray*}
\partial_2^\eta(r) & = & \exp{\left(\int_{\overline{\mathbb{R}}} \gamma_{21,r}^\ast(\eta)
\right)} q + \exp{\left(\int_{\overline{\mathbb{R}}} \gamma_{21,l}^\ast(\eta) \right)}q\\
\partial_1^\eta(q) & = & \exp{\left(\int_{\overline{\mathbb{R}}} \gamma_{10,+}^\ast(\eta)
\right)} p - \exp{\left(\int_{\overline{\mathbb{R}}} \gamma_{10,-}^\ast(\eta) \right)} p
\end{eqnarray*}
Hence, $\partial_2^\eta:C_2(f) \otimes \mathbb{R} \rightarrow C_1(f) \otimes \mathbb{R}$
is surjective for any $\eta \in \Omega^1_{\text{cl}}(\mathbb{R}P^2,\mathbb{R})$. If $\eta$
is exact, then it is clear that $\partial_1^\eta = 0$.  If $\eta$ were not exact,
then $\partial_1^\eta:C_1(f) \otimes \mathbb{R} \rightarrow C_0(f) \otimes \mathbb{R}$
would be surjective. However, $H^1(\mathbb{R}P^2;\mathbb{R}) = 0$, and hence 
$\partial_1^\eta$ = 0 for all $\eta \in \Omega^1_{\text{cl}} (\mathbb{R}P^2,\mathbb{R})$.
Thus for any $\eta \in \Omega^1_{\text{cl}} (\mathbb{R}P^2,\mathbb{R})$,
$$
H_2((C_\ast(f)\otimes \mathbb{R},\partial_\ast^\eta)) \approx 0, \quad
H_1((C_\ast(f)\otimes \mathbb{R},\partial_\ast^\eta)) \approx 0, \quad
H_0((C_\ast(f)\otimes \mathbb{R},\partial_\ast^\eta)) \approx \mathbb{R}.
$$
\end{example}


\subsection{Computation of $H_0((C_\ast(f;G),\partial^G_\ast))$ }
Let $G$ be a bundle of abelian groups over a closed smooth Riemannian manifold
$M$ of dimension $m < \infty$, and let $f:M \rightarrow \mathbb{R}$ be
a smooth Morse-Smale function on $M$.  Choose a basepoint $x_0 \in Cr_0(f)$ of $M$.
Recall that if $M$ is connected then the isomorphism class of $G$ is determined by a representation
$$
\pi_1(M,x_0) \times G_{x_0} \rightarrow G_{x_0},
$$
and let $H_{x_0} \subseteq G_{x_0}$ be the subgroup generated by elements of 
the form $g - \gamma_\ast(g)$ where $g \in G_{x_0} \text{ and }[\gamma] \in \pi_1(M,x_0)$.
The following theorem gives the $0$-dimensional twisted Morse homology group of $M$ in
terms of the above action, cf. Theorem VI.3.2 of \cite{WhiEle}.

\begin{theorem}\label{H_0}
If $M$ is connected, then the $0$-dimensional twisted Morse homology group of $M$
is isomorphic to $G_{x_0}/H_{x_0}$, i.e.
$$
H_0((C_\ast(f;G),\partial_\ast^G)) \approx G_{x_0}/H_{x_0}.
$$
\end{theorem}

\proofstart
Let $p_0 \in Cr_0(f)$ and pick any path from $p_0$ to $x_0$. For all critical points
$p_k \in Cr_k(f)$ with $k \geq 2$ the path is homotopic rel endpoints to a path that
does not intersect $W^s(p_k)$, since dim $W^s(p_k) = m-k \leq m-2$, cf. Theorems 5.16
and 5.17 of \cite{BanMor}. The flow of $-\nabla f$ then gives a homotopy rel endpoints
to a path from $p_0$ to $x_0$ that lies within an $\varepsilon$-neighborhood of the
$1$-skeleton of $f$, i.e.
$$
\bigcup_{q \in Cr_1(f)} \overline{W^u(q)}
$$
for any $\varepsilon > 0$. An additional homotopy rel endpoints then produces a path
$\gamma$ from $p_0$ to $x_0$ that lies in the $1$-skeleton of $f$. Thus, every
critical point $p_0 \in Cr_0(f)$ is connected by a path $\gamma$ to $x_0\in Cr_0(f)$
that lies in the $1$-skeleton of $f$.

\smallskip
Now, let $Z_0(f;G) = C_0(f;G) = \text{ker }\partial_0^G$ denote the group of 
$0$-cycles and $B_0(f;G) = \text{Im }\partial_1^G$ the group of $0$-boundaries.
We claim that every element of $Z_0(f;G)/B_0(f;G)$ can be represented by
 an elementary cycle supported at $x_0$, i.e. $g x_0$ for some $g \in G_{x_0}$.  To see this,
let $g_0p_0$ be an elementary cycle and consider a path $\gamma$ contained in the $1$-skeleton
of $f$ connecting $p_0$ to $x_0$.  Then
$$
\text{Im } \gamma = \bigcup_{j=1}^n \overline{W^u(q_j)}
$$
for some critical points $q_j \in Cr_1(f)$.  Number the critical points consecutively along
$\text{Im }\gamma$ starting with $q_1$ as the critical point with $p_0 \in  \overline{W^u(q_1)}$. 
Let $\nu$ be the gradient flow line from $q_1$ to $p_0$, and let $p_1 \in  \overline{W^u(q_1)}$
be the next element of $Cr_0(f)$ along the path $\gamma$ connecting $p_0$ to $x_0$.  Then
$$
g_0 p_0 - \partial^G_1\left( \epsilon(\nu)  (\gamma_\nu)_\ast (g_0)    q_1  \right) 
$$
has support at $p_1$, where $\gamma_\nu$ is any continuous path from $q_1$ to $p_0$ whose image
coincides with the image of $\nu$.  By consecutively subtracting boundary elements of this
type for $q_2, \dots, q_n$ we get a cycle in the same equivalence class of $[g_0p_0] \in Z_0(f;G)/B_0(f;G)$
with support at $x_0$. Applying this algorithm to each term in the sum defining an arbitrary element of 
$Z_0(f;G)$ establishes the claim. 

To finish the proof of the theorem, note that any boundary that is supported at $x_0$ must be the image
under $\partial^G_1$ of a chain $\sum_{j=1}^n g_j q_j \in C_1(f;G)$, where
$$
 \bigcup_{j=1}^n \overline{W^u(q_j)}
$$
is the union of the image of a finite number of loops $\gamma_1, \ldots ,\gamma_i$ based at $x_0$. 
Moreover, since the coefficients in front of the index zero critical points not equal to $x_0$
in the sum for $\partial_1^G\left(\sum_{j=1}^n g_j q_j\right)$ are all zero, we must have
$$
\partial_1^G  \left( \sum_{j=1}^n g_j q_j \right) = (g_1 - (\gamma_1)_\ast(g_1)) x_0 + \cdots +
(g_i - (\gamma_i)_\ast(g_i)) x_0
$$
for some $g_1,\ldots , g_i \in G_{x_0}.$

\proofend

\begin{remark}
The preceding theorem shows that $H_0((C_\ast(f;G),\partial_\ast^G))$ is independent of
the Morse-Smale pair $(f,\mathsf{g})$ and depends only on the representation
$\pi_1(M,x_0) \rightarrow \text{Aut}(G_{x_0})$, i.e. the isomorphism class of $G$.
Also, note that the fact that every path from $x_0 \in Cr_0(f)$ to $p_0\in Cr_0(f)$ is 
homotopic rel endpoints to a path that lies in the $1$-skeleton of $f$ follows from
the Cellular Approximation Theorem, cf. Theorem IV.11.4 of \cite{BreTop}, and the results
about unstable manifolds and CW-structures discussed in Section \ref{unstableCW}.
\end{remark}


\subsection{Morse Eilenberg Theorem}\label{EMTheorem}
In this subsection we prove a Morse theoretic version of Eilenberg's theorem relating the
homology with local coefficients of a space to the equivariant homology of its universal 
cover. Before we state the theorem we first discuss the Morse-Smale-Witten chain complex on 
a covering space $\widetilde{M}$ of a closed smooth Riemannian manifold $M$
when the coordinate charts, the metric, and the Morse function are all pulled back to the
covering space. We show that under these assumptions the homology of the (untwisted) Morse-Smale-Witten
chain complex is isomorphic to the singular homology of the covering space, even if the covering
space is not compact.

\subsubsection{The Morse-Smale-Witten chain complex on a covering space}
Let $(M,\mathsf{g})$ be a closed smooth Riemannian manifold of dimension $m < \infty$, 
and let $f:M \rightarrow \mathbb{R}$ be a smooth Morse-Smale function on $M$.
Let $\pi:\widetilde{M}\rightarrow M$ be a covering space of $M$ with the coordinate charts
on $M$ pulled back to $\widetilde{M}$ so that $\pi$ is a local diffeomorphism. 
Let $\widetilde{\mathsf{g}}$ be the pullback of $\mathsf{g}$ to $\widetilde{M}$, and let 
$\tilde{f} = f \circ \pi$. 
$$
\xymatrix{
\tilde{f}: (\widetilde{M},\widetilde{\mathsf{g}})  \ar[r]^-{\pi} & (M,\mathsf{g})  \ar[r]^-{f} & \mathbb{R}
}
$$
Nondegeneracy and Morse-Smale transversality are local conditions, and hence
$\tilde{f}$ is a smooth Morse-Smale function on $(\widetilde{M},\widetilde{\mathsf{g}})$.
Moreover, for all $k=0,\ldots, m$ the set of critical points of $\tilde{f}$ of index $k$ is 
$$
Cr_k(\tilde{f}) = \{\tilde{q} \in \widetilde{M} |\ \pi(\tilde{q}) = q \text{ where } q \in Cr_k(f)\},
$$
because the Hessian of $\tilde{f}$ at $\tilde{q}$ is the pullback of the Hessian of $f$ at $q=\pi(\tilde{q})$.

\smallskip 
If $\widetilde{M}$ is compact, then the Morse-Smale-Witten chain complex $(C_\ast(\tilde{f}),
\tilde{\partial}_\ast)$ of $\tilde{f}$ is well-defined, and its homology is isomorphic to the 
singular homology $H_\ast(\widetilde{M};\mathbb{Z})$ by the Morse Homology Theorem, cf. 
Theorem 7.4 of \cite{BanLec}. 

\begin{example}[The universal cover of a real projective space]\label{projectivespacecover}
Consider the Morse-Smale-Witten chain complex of the function $f:\mathbb{R}P^2 \rightarrow 
\mathbb{R}$ described in Example \ref{projectivespace}. The function $f$ lifts to a function 
$\tilde{f}:S^2 \rightarrow \mathbb{R}$ on the universal cover $S^2$ of $\mathbb{R}P^2$, and we can 
pull back the orientations of the local unstable manifolds of $f$ to those of $\tilde{f}:S^2\rightarrow 
\mathbb{R}$. The resulting phase diagram of $\tilde{f}$ and the signs associated to the gradient flow
lines between critical points of relative index one are shown in the following diagram.

\begin{figure}[h]
\includegraphics{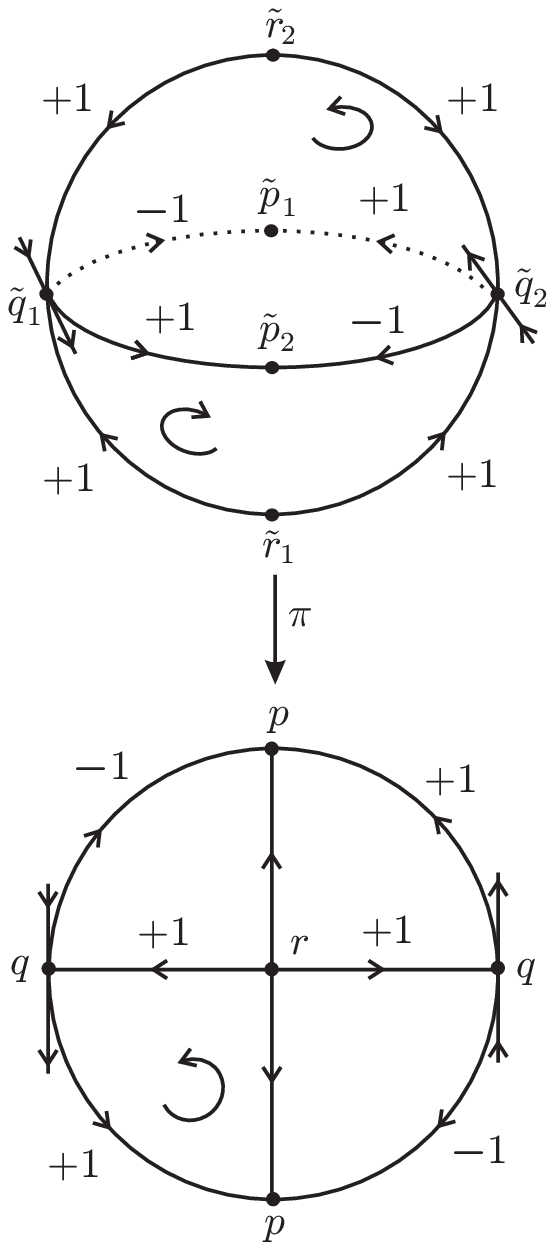}
\end{figure}

\noindent
We have $Cr_2(\tilde{f}) = \{ \tilde{r}_1,\tilde{r}_2 \}$ with
$$
\tilde{\partial}_2(\tilde{r}_1) = \tilde{\partial}_2(\tilde{r}_2) = 
\tilde{q}_1 + \tilde{q}_2,
$$
and $Cr_1(\tilde{f}) = \{ \tilde{q}_1,\tilde{q}_2 \}$ with
$$
\tilde{\partial}_1(\tilde{q}_1) = \tilde{p}_2 - \tilde{p}_1 \quad \text{ and } \quad
\tilde{\partial}_1(\tilde{q}_2) = \tilde{p}_1 - \tilde{p}_2,
$$
where $Cr_0(\tilde{f}) = \{\tilde{p}_1, \tilde{p}_2 \}$.  Therefore,
\begin{eqnarray*}
H_2((C_\ast(\tilde{f}),\tilde{\partial}_\ast)) &  = &  <\tilde{r}_1 - \tilde{r}_2>\  \approx 
\ \mathbb{Z}\\
H_1((C_\ast(\tilde{f}),\tilde{\partial}_\ast)) &  = &  <\tilde{q}_1 + \tilde{q}_2> / 
<\tilde{q}_1 + \tilde{q_2}>\  \approx\  0\\
H_0((C_\ast(\tilde{f}),\tilde{\partial}_\ast)) &  = &  <\tilde{p}_2 , \tilde{p}_1> / 
<\tilde{p}_1 - \tilde{p}_2>\  \approx \mathbb{Z}
\end{eqnarray*}
as expected.
\end{example}

\smallskip
We will now show that even if the covering space $\widetilde{M}$ is not compact, then 
$\tilde{\partial}_\ast$ is still a well-defined boundary operator and $H_k((C_\ast(\tilde{f}),
\tilde{\partial}_\ast)) \approx H_k(\widetilde{M};\mathbb{Z})$ for all $k=0,\ldots , m$. 
To see this, first recall that for all $k=0,\ldots, m$, the group $C_k(\tilde{f})$ is the free abelian
group generated by the critical points of index $k$. Hence, $C_k(\tilde{f})$ consists of sums of the form
$$
\sum_{\tilde{q} \in Cr_k(\tilde{f})} n_{\tilde{q}} \tilde{q}, 
$$
where all but finitely many of the coefficients $n_{\tilde{q}} \in \mathbb{Z}$ are zero. So, even if
$Cr_k(\tilde{f})$ is infinite, $C_k(\tilde{f})$ consists of finite sums. However, in order to
know that $\tilde{\partial}_\ast$ is well-defined we need to know that if 
$\tilde{q}\in Cr_k(\tilde{f})$ and $\tilde{p} \in Cr_{k-1}(\tilde{f})$, then the number of gradient
flow lines from $\tilde{q}$ to $\tilde{p}$ is finite, for all $1 \leq k \leq m$. To see this, note 
that a gradient flow line of $(\tilde{f}, \widetilde{\mathsf{g}})$ from $\tilde{q}$ to $\tilde{p}$ projects
to a gradient flow line of $(f,\mathsf{g})$ from $q = \pi(\tilde{q})$ to $p = \pi(\tilde{p})$. Moreover,
every gradient flow line of $(f,\mathsf{g})$ from $q$ to $p$ lifts to a unique gradient flow line of 
$(\tilde{f},\widetilde{\mathsf{g}})$ starting at $\tilde{q}$ by the path lifting property of 
$\pi:\widetilde{M} \rightarrow M$. Therefore, the number of gradient flow lines from 
$\tilde{q}$ to $\tilde{p}$ is less than or equal to the number of gradient flow lines from $q$ to $p$, 
which is finite since $M$ is compact and $(f,\mathsf{g})$ satisfies the Morse-Smale transversality
condition, cf. Corollary 6.29 of \cite{BanLec}.

\begin{example}[The universal cover of a circle]
Consider the height function on $S^1$ with its universal cover $\mathbb{R}$, where the arrows
in the following diagram indicate the direction of minus the gradient flow.

\begin{figure}[h]
\includegraphics{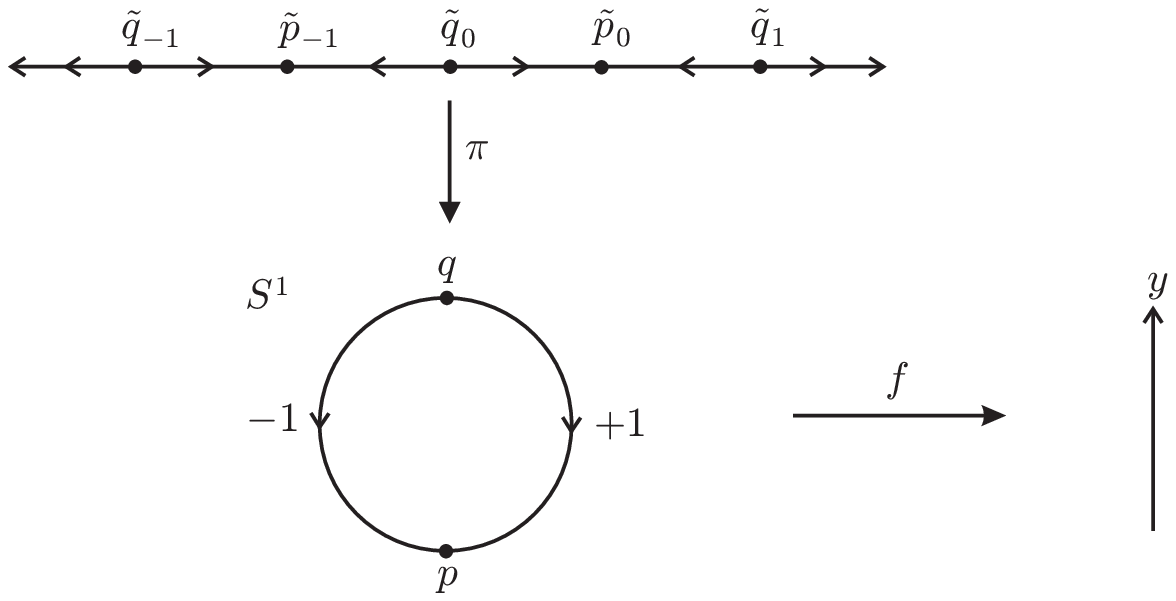}
\end{figure}
Above the index $1$ critical point $q$ of $f$ there are an infinite number of index
$1$ critical points $\{\ldots, \tilde{q}_{-1}, \tilde{q}_0, \tilde{q}_1,\ldots\}$ of $\tilde{f} 
= f \circ \pi$. Similarly, there are an infinite number of index $0$ critical points
$\{\ldots, \tilde{p}_{-1}, \tilde{p}_0, \tilde{p}_1,\ldots\}$ above the index $0$ critical
point $p$. Note that there are two gradient flow lines from $q$ to $p$, but when these gradient
flow lines are lifted to the universal cover starting at a critical point $\tilde{q}_j$
the gradient flow lines in the universal cover end at different critical points 
$\tilde{p}_{j-1}$ and $\tilde{p}_j$. So, for any $j, j' \in \mathbb{Z}$ the number of gradient
flow lines from $\tilde{q}_j$ to $\tilde{p}_{j'}$ is strictly less than the number of gradient 
flow lines from $q$ to $p$. 

If we pull back the orientations of $W^u(q)$ and $W^u(p)$ to the universal cover,
then the signs associated to the gradient flow lines in the universal cover are the same
as the signs associated to the gradient flow lines that they correspond to under $\pi$.  Thus, for
any $j \in \mathbb{Z}$ we have
$$
\tilde{\partial}_1 (\tilde{q}_j) = \tilde{p}_j - \tilde{p}_{j-1},
$$
which implies that $H_0((C_\ast(\tilde{f}),\tilde{\partial}_\ast)) \approx \mathbb{Z}$.
It also implies that for any element in $C_1(\tilde{f})$ we have
$$
\tilde{\partial}_1 \left( \sum_{j \in \mathbb{Z}} n_j \tilde{q}_j \right) \neq 0,
$$
since only finitely many of the $n_j$ are nonzero. That is, there is a largest $j$ such that
$n_j \neq 0$, and thus $n_j \tilde{p}_j \neq 0$ in the sum $\sum_{j \in \mathbb{Z}} n_j 
\tilde{\partial}_1(\tilde{q}_j )$. This implies that $H_1((C_\ast(\tilde{f}),\tilde{\partial}_\ast)) \approx 0$.
\end{example}

\begin{theorem}[Morse Homology Theorem for a Covering Space]\label{MorHomCover}
Let $\pi:\widetilde{M} \rightarrow M$ be a covering space of a closed smooth Riemannian
manifold $(M,\mathsf{g})$ of finite dimension $m$. Assume that the smooth structure on
$\widetilde{M}$ is given by pulling back the coordinate charts on the base $M$,
the metric $\widetilde{\mathsf{g}}$ on $\widetilde{M}$ is the pullback of the metric $\mathsf{g}$ on $M$, 
and $\tilde{f}:\widetilde{M} \rightarrow \mathbb{R}$ is the pullback of a function $f:M \rightarrow \mathbb{R}$.
If $(f,\mathsf{g})$ is a Morse-Smale pair on $M$, then $(\tilde{f},\widetilde{\mathsf{g}})$ is a
Morse-Smale pair on $(\widetilde{M},\widetilde{\mathsf{g}})$, the Morse-Smale-Witten chain complex
$(C_\ast(\tilde{f}), \tilde{\partial}_\ast)$ is well-defined, and its homology is isomorphic to the singular
homology $H_\ast(\widetilde{M}; \mathbb{Z})$.
\end{theorem}

\proofstart
As discussed above, the pair $(\tilde{f},\widetilde{\mathsf{g}})$ is Morse-Smale since $\pi$ is a 
local isometry, and $\tilde{\partial}_\ast$ is well-defined because of the path lifting property. 
The proof that $H_k((C_\ast(\tilde{f}), \tilde{\partial}_\ast))  \approx H_k(\widetilde{M};
\mathbb{Z})$ for all $k=0,\ldots , m$ is essentially the same as the proof of The Morse Homology 
Theorem for finite dimensional compact smooth manifolds given in Chapter 7 of \cite{BanLec}.

More explicitly, a filtration of $M$
$$
\emptyset = N_{-1} \subseteq N_0 \subseteq N_1 \subseteq \cdots \subseteq N_m = M
$$
such that $(N_k,N_{j-1})$ is a cofibered index pair for 
$$
W(k,j)\  = \bigcup_{j \leq \lambda_p \leq \lambda_q \leq k} W(q,p) \subseteq M
$$
for all $0 \leq j \leq k \leq m$, where $W(q,p) = W^u(q) \cap W^s(q) \subset M$, 
pulls back under $\pi:\widetilde{M} \rightarrow M$ to a filtration of $\widetilde{M}$
$$
\emptyset = \widetilde{N}_{-1} \subseteq \widetilde{N}_0 \subseteq \widetilde{N}_1 \subseteq
\cdots \subseteq \widetilde{N}_m = \widetilde{M}
$$
such that $(\widetilde{N}_k,\widetilde{N}_{j-1})$ is a cofibered index pair for 
$$
\widetilde{W}(k,j)\  = \bigcup_{j \leq \lambda_{\tilde{p}} \leq \lambda_{\tilde{q}} \leq k} 
\widetilde{W}(\tilde{q},\tilde{p}) \subseteq \widetilde{M}
$$
for all $0 \leq j \leq k \leq m$, where $\widetilde{W}(\tilde{q} ,\tilde{p}) = W^u(\tilde{q})
\cap W^s(\tilde{p}) \subset \widetilde{M}$. One can then show that the following diagram commutes for all 
$k=1,\ldots , m$,
$$
\xymatrix{
C_k(\tilde{f}) \ar@{<->}[d]^\approx  \ar[r]^-{\tilde{\partial}_k} & C_{k-1}(\tilde{f}) 
\ar@{<->}[d]^\approx\\
H_k(\widetilde{N}_k,\widetilde{N}_{k-1} ) \ar[r]^-{\tilde{\delta}_\ast} & 
H_{k-1}(\widetilde{N}_{k-1},\widetilde{N}_{k-2})
}
$$
where $\tilde{\delta}_\ast$ is the connecting homomorphism of  the triple $(\widetilde{N}_k,\widetilde{N}_{k-1},
\widetilde{N}_{k-2})$. The proof is essentially the same as the argument used to prove Lemma 7.21 of
\cite{BanLec}, which is local in the sense that it works on an isolating neighborhood of 
$\widetilde{W}(\tilde{q},\tilde{p}) \cup \{\tilde{q},\tilde{p} \}$, where $\tilde{q} \in C_k(\tilde{f})$
and $\tilde{p} \in C_{k-1}(\tilde{f})$. Also, note that although \cite{BanLec} assumes that $M$ is orientable,
which determines an orientation on the stable manifolds, these orientations are not used in the proof of 
Lemma 7.21 since it is the orientations of the unstable manifolds that determine generators for
$H_k(\widetilde{N}_k,\widetilde{N}_{k-1})$ and $H_{k-1}(\widetilde{N}_{k-1},\widetilde{N}_{k-2})$.  
The rest of the proof is identical to the proof of Theorem 7.4 of \cite{BanLec}.
\proofend

\smallskip\noindent
{\bf Note:} A similar result is stated in Section 6.4.3 of \cite{PajCir}, where Pajitnov
refers to Morse complex on the universal cover as the ``universal Morse complex''. Pajitnov
also proved that the universal Morse complex is simply homotopy equivalent to the chain complex
of a triangulation, cf. Theorem A.5 of \cite{PazOnt}.


\subsubsection{The Morse Eilenberg Theorem}
Let $f:M \rightarrow \mathbb{R}$ be a smooth Morse function. Pick a basepoint $x_0 \in Cr_0(f)$ for $M$ 
and a basepoint $\tilde{x}_0 \in \pi^{-1}(x_0)$ for the universal cover $\pi:\widetilde{M} \rightarrow M$.
Assume that the smooth structure on $\widetilde{M}$ is given by pulling back the charts on $M$ so that
$\pi$ is a local diffeomorphism, and let $\tilde{f} = f \circ \pi$. 

The action of $\pi_1(M,x_0)$ on the universal cover by deck transformations
$$
\pi_1(M,x_0) \times \widetilde{M} \rightarrow \widetilde{M}
$$
restricts to an action on $Cr_k(\tilde{f})$ because
$$
Cr_k(\tilde{f}) = \bigcup_{q \in Cr_k(f)} \pi^{-1}(q),
$$
and hence there is left action 
$$
\pi_1(M,x_0) \times C_k(\tilde{f}) \rightarrow C_k(\tilde{f})
$$
on the free abelian group generated by the critical points of $\tilde{f}$ of index $k$.
If $G$ is a bundle of abelian groups over $M$, then there is also a left action
$$
\pi_1(M,x_0) \times G_{x_0} \rightarrow G_{x_0}
$$
given by $[\gamma] \cdot g = \gamma_\ast(g)$ for all $[\gamma] \in \pi_1(M,x_0)$ and
$g \in G_{x_0}$, which is converted to a right action by
$$
g \cdot [\gamma]  \stackrel{\text{def}}{=} [\gamma]^{-1} \cdot g =   \gamma^{-1}_\ast (g).
$$

Now consider $G_{x_0} \otimes_{\mathbb{Z}} C_k(\tilde{f})$, and let $K_k(\tilde{f};G_{x_0})$ be 
subgroup generated by elements of the form 
$$
g \cdot[\gamma] \otimes \tilde{q} - g \otimes [\gamma] \cdot \tilde{q},
$$
so that 
$$
G_{x_0} \otimes_{\pi_1} C_k(\tilde{f}) = \left( G_{x_0} \otimes_{\mathbb{Z}} C_k(\tilde{f})  \right)
/K_k(\tilde{f};G_{x_0}),
$$
where the first tensor product is over $\pi_1(M,x_0)$. 
Pick a metric $\mathsf{g}$ on $M$ such that $(f,\mathsf{g})$ is a Morse-Smale pair and let 
$\widetilde{\mathsf{g}} = \pi^\ast(\mathsf{g})$. If we orient the unstable manifolds of 
$(\tilde{f}, \widetilde{\mathsf{g}})$ by pulling back
the orientations of the unstable manifolds of $(f,\mathsf{g})$, then the (untwisted) Morse-Smale-Witten 
boundary operator $\partial_k:C_k(f) \rightarrow C_{k-1}(f)$ induces a boundary operator 
on $C_\ast(\tilde{f})$. Hence, it induces a boundary operator
$$
\tilde{\partial}_k: G_{x_0} \otimes_{\mathbb{Z}} C_k(\tilde{f}) \rightarrow
G_{x_0} \otimes_{\mathbb{Z}} C_{k-1}(\tilde{f}).
$$
This boundary operator commutes with the action of $\pi_1(M,x_0)$ on $\widetilde{M}$, and thus it maps
$K_k(\tilde{f};G_{x_0})$ to $K_{k-1}(\tilde{f};G_{x_0})$. Hence, there is an induced boundary operator
$$
\bar{\partial}_k: G_{x_0} \otimes_{\pi_1} C_k(\tilde{f}) \rightarrow
G_{x_0} \otimes_{\pi_1} C_{k-1}(\tilde{f}).
$$

\begin{theorem}[Morse Eilenberg]\label{EilenbergMorse}
Let $G$ be a bundle of abelian groups over a closed smooth Riemannian manifold
$(M,\mathsf{g})$ of dimension $m < \infty$, and let $f:M \rightarrow \mathbb{R}$ be
a smooth Morse-Smale function on $M$. Then,
$$
H_k((G_{x_0} \otimes_{\pi_1} C_\ast(\tilde{f}), \bar{\partial}_\ast) )  \approx
H_k((C_\ast(f;G),\partial^G_\ast) )
$$
for all $k=0,\ldots, m$.
\end{theorem}

\proofstart
Since $\widetilde{M}$ is simply connected, for any $\tilde{q} \in Cr_k(\tilde{f})$
there is a unique homotopy class of paths rel endpoints from $\tilde{q}$ to the basepoint 
$\tilde{x}_0\in \widetilde{M}$, where $0 \leq k \leq m$. So, if  $\tilde{\gamma}_{\tilde{q}}$ is 
any path in $\widetilde{M}$ from $\tilde{q}$ to $\tilde{x}_0$, then $\gamma_{\tilde{q}} 
\equiv \pi \circ \tilde{\gamma}_{\tilde{q}}$ determines a well-defined homotopy class of paths 
rel endpoints in $M$ from $q \equiv \pi(\tilde{q}) \in Cr_k(f)$ to the basepoint 
$x_0=\pi(\tilde{x}_0)$ of $M$. Thus, there is a homomorphism
$$
\tilde{\Phi}_k:G_{x_0} \otimes_{\mathbb{Z}} C_k(\tilde{f}) \rightarrow C_k(f;G)
$$ 
defined on a generator $g \otimes \tilde{q}$ by 
$$
\tilde{\Phi}_k(g \otimes \tilde{q}) \stackrel{\text{def}}{=} (\gamma_{\tilde{q}})_\ast(g) q,
$$
where $g \in G_{x_0}$ and $\tilde{q} \in Cr_k(\tilde{f})$.

We claim that $K_k(\tilde{f};G_{x_0})$ is the kernel of $\tilde{\Phi}_k$. To see this, take any generator 
$$
g \cdot[\gamma] \otimes \tilde{q} - g \otimes [\gamma] \cdot \tilde{q} \in K_k(\tilde{f};G_{x_0}),
$$ 
where $g \in G_{x_0}$, $\tilde{q} \in Cr_k(\tilde{f})$, and $[\gamma]\in \pi_1(M,x_0)$, and note
that $[\gamma]\cdot \tilde{\gamma}_{\tilde{q}}$ is a path from $[\gamma]\cdot \tilde{q}$ to
$[\gamma]\cdot \tilde{x}_0$. Lifting the path $\gamma$ to a path starting at $\tilde{x}_0$ gives
a path from $\tilde{x}_0$ to $[\gamma]\cdot \tilde{x}_0$ whose inverse is a path 
$\tilde{\gamma}_{[\gamma] \cdot \tilde{x}_0}$ from $[\gamma]\cdot \tilde{x}_0$ to $\tilde{x}_0$.
The concatenation $([\gamma]\cdot \tilde{\gamma}_{\tilde{q}})\tilde{\gamma}_{[\gamma] \cdot \tilde{x}_0}$
is a path from $[\gamma]\cdot \tilde{q}$ to $\tilde{x}_0$.  Therefore,
\begin{eqnarray*}
\tilde{\Phi}_k \left(g \cdot[\gamma] \otimes \tilde{q} - g \otimes [\gamma] \cdot \tilde{q} \right)
& = & \tilde{\Phi}_k \left(g \cdot[\gamma] \otimes \tilde{q}\right) - 
      \tilde{\Phi}\left(g \otimes [\gamma] \cdot \tilde{q}\right)\\
& = & (\gamma_{\tilde{q}})_\ast\left(\gamma^{-1}_\ast (g)\right) q - 
      \left((\gamma_{\tilde{q}})_\ast \circ \gamma^{-1}_\ast \right)(g) q\\
& = & 0,
\end{eqnarray*}
since $\pi  \circ \left( [\gamma]\cdot \tilde{\gamma}_{\tilde{q}} \right) = \gamma_{\tilde{q}}$ and 
$\pi \circ \tilde{\gamma}_{[\gamma]\cdot \tilde{x}_0} = \gamma^{-1}$.

\begin{figure}[h]
\includegraphics{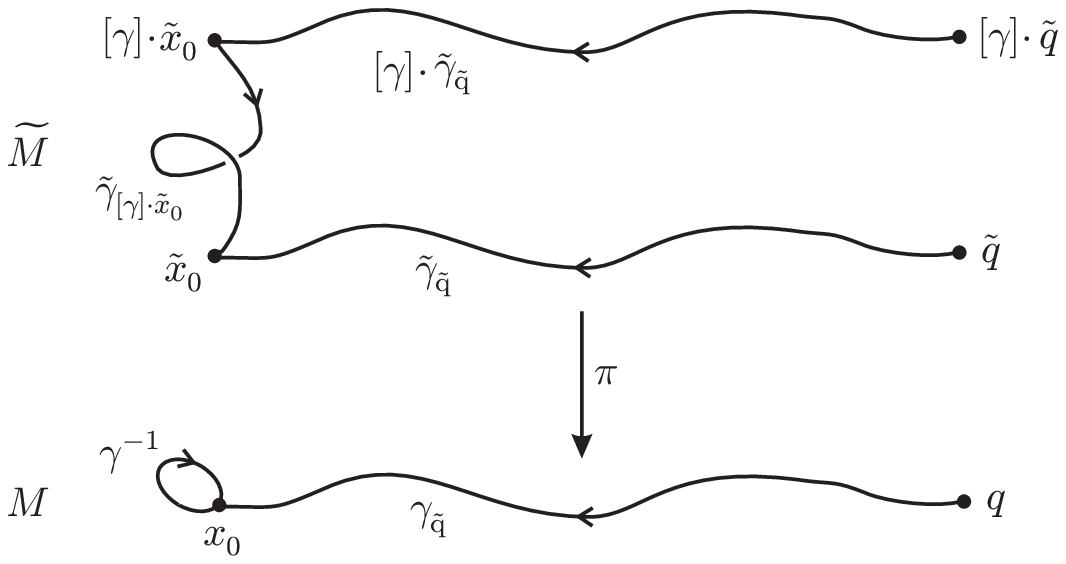}
\end{figure}

Moreover, if $\tilde{\Phi}_k$ sends some element of $G_{x_0} \otimes_{\mathbb{Z}} C_k(\tilde{f})$ 
to $0\in C_k(f;G)$, then the coefficients in front of each element of $Cr_k(f)$ in the image must be
zero. Thus, we may assume that the element that $\tilde{\Phi}_k$ maps to zero is of the form 
$\sum_{i=1}^n g_i \otimes \tilde{q}_i$, where $\pi(\tilde{q}_i) = q \in Cr_k(f)$ for all $i=1,\ldots, n$.
Then for all $i=2,\ldots ,n$, there exist $[\gamma_i] \in \pi_1(M,x_0)$ such that 
$\tilde{q}_i = [\gamma_i]\cdot \tilde{q}_1$. Hence,
\begin{eqnarray*}
0 & = & \tilde{\Phi}_k\left(g_1\otimes \tilde{q}_1 + \sum_{i=2}^n g_i \otimes [\gamma_i] \cdot 
        \tilde{q}_1\right)\\
  & = & (\gamma_{\tilde{q}_1})_\ast(g_1)q + \sum_{i=2}^n \left((\gamma_{\tilde{q}_1})_\ast \circ 
        (\gamma^{-1}_i)_\ast\right)(g_i) q\\
  & = & (\gamma_{\tilde{q}_1})_\ast \left(g_1 + \sum_{i=2}^n (\gamma^{-1}_i)_\ast(g_i) \right) q,
\end{eqnarray*}
which implies that 
$$
g_1 + \sum_{i=2}^n (\gamma^{-1}_i)_\ast(g_i) = 0,
$$
since $(\gamma_{\tilde{q}_1})_\ast$ is an isomorphism.  Therefore,
\begin{eqnarray*}
\sum_{i=1}^n g_i \otimes \tilde{q}_i & = & g_1 \otimes \tilde{q}_1 + \sum_{i=2}^n 
                                           g_i \otimes [\gamma_i] \cdot \tilde{q}_1\\
& = & -\sum_{i=2}^n (\gamma^{-1}_i)_\ast(g_i) \otimes \tilde{q}_1 +
       \sum_{i=2}^n g_i \otimes [\gamma_i] \cdot \tilde{q}_1\\
& = & \left. \left. \sum_{i=2}^n \right( g_i \otimes [\gamma_i]\cdot \tilde{q}_1 - g_i \cdot [\gamma_i] 
      \otimes \tilde{q}_1 \right) \\
& \in & K_k(\tilde{f};G_{x_0}).
\end{eqnarray*}

We have shown that $K_k(\tilde{f};G_{x_0}) = \text{ker }\tilde{\Phi}_k$. To see that $\tilde{\Phi}_k$
is surjective note that for any elementary chain $gq \in C_k(f;G)$ we can pick a path $\gamma^{-1}_q$ 
from $x_0$ to $q$ and lift it to a path in $\widetilde{M}$ starting at $\tilde{x}_0$. The end of the 
lifted path will be a critical point $\tilde{q}$ with $\pi(\tilde{q}) = q$, and $\tilde{\Phi}_k
\left((\gamma^{-1}_q)_\ast (g) \otimes \tilde{q}\right) = 
(\gamma_{\tilde{q}})_\ast ((\gamma^{-1}_q)_\ast (g)) q = gq$. Therefore, we have an induced isomorphism
$$
\Phi_k:\left(G_{x_0} \otimes_\mathbb{Z} C_k(\tilde{f})\right)/K_k(\tilde{f};G_{x_0}) \rightarrow C_k(f;G)
$$
for all $k=0,\ldots ,m$.
$$
\xymatrix{
G_{x_0} \otimes_{\mathbb{Z}} C_k(\tilde{f}) \ar[d] \ar[drr]^{\tilde{\Phi}_k} & & \\
G_{x_0} \otimes_{\pi_1} C_k(\tilde{f}) \ar@{<->}[rr]^-{\Phi_k}_-{\approx} & & C_k(f;G)
}
$$

It remains to show that $\Phi_k$ is a chain map for all $0 \leq k \leq m$. To do this, 
we will show that $\tilde{\Phi}_k$ is a chain map. It then follows that $\Phi_k$ is a chain
map because the boundary operator $\bar{\partial}_\ast$ on $G_{x_0} 
\otimes_{\pi_1} C_\ast(\tilde{f})$ is induced from the boundary operator 
$\tilde{\partial}_\ast$ on $G_{x_0} \otimes_{\mathbb{Z}} C_\ast(\tilde{f})$.
So, let $g \otimes \tilde{q}$ be a generator of $G_{x_0} \otimes_{\mathbb{Z}} C_k(\tilde{f})$, 
and note that 
\begin{eqnarray*}
\partial^G_k(\tilde{\Phi}_k(g \otimes \tilde{q})) 
& = & \partial^G_k((\gamma_{\tilde{q}})_\ast(g) q)\\
& = & \sum_{p \in Cr_{k-1}(f)} \sum_{\nu \in \mathcal{M}(q,p)} \epsilon(\nu)
\gamma^{\nu}_\ast((\gamma_{\tilde{q}})_\ast(g))p\\ 
\end{eqnarray*}
where $\gamma^{\nu}:[0,1] \rightarrow M$ is any continuous path from $p$ to $q$ whose
image coincides with the image of $\nu \in \mathcal{M}(q,p)$ and $\epsilon(\nu) = \pm 1$ 
is the sign determined by the orientation on $\mathcal{M}(q,p)$. 
On the other hand,
\begin{eqnarray*}
\tilde{\Phi}_{k-1}\left(\tilde{\partial}_k(g \otimes \tilde{q})\right) 
& = & \tilde{\Phi}_{k-1} \left( g \otimes \sum_{\tilde{p} \in Cr_{k-1}(\tilde{f})} 
      \sum_{\tilde{\nu} \in \mathcal{M}(\tilde{q},\tilde{p})} \epsilon(\tilde{\nu}) 
      \tilde{p} \right)\\
& = & \sum_{\tilde{p} \in Cr_{k-1}(\tilde{f})} \sum_{\tilde{\nu} \in \mathcal{M}(\tilde{q},\tilde{p})}
      \epsilon(\tilde{\nu})\ \tilde{\Phi}_{k-1} \left(g \otimes \tilde{p}\right)\\
& = & \sum_{\tilde{p} \in Cr_{k-1}(\tilde{f})} \sum_{\tilde{\nu} \in \mathcal{M}(\tilde{q},\tilde{p})}
      \epsilon(\tilde{\nu})(\gamma_{\tilde{p}})_\ast(g) p\\
& = & \sum_{\tilde{p} \in Cr_{k-1}(\tilde{f})} \sum_{\tilde{\nu} \in \mathcal{M}(\tilde{q},\tilde{p})}
      \epsilon(\tilde{\nu})(\pi \circ \gamma^{\tilde{\nu}}\tilde{\gamma}_{\tilde{q}})_\ast(g)) p\\
& = & \sum_{\tilde{p} \in Cr_{k-1}(\tilde{f})} \sum_{\tilde{\nu} \in \mathcal{M}(\tilde{q},\tilde{p})}
      \epsilon(\tilde{\nu})(\pi \circ \gamma^{\tilde{\nu}})_\ast((\gamma_{\tilde{q}})_\ast(g)) p,\\
\end{eqnarray*}
where $\gamma^{\tilde{\nu}}:[0,1] \rightarrow \widetilde{M}$ is any continuous path from $\tilde{p}$ to 
$\tilde{q}$ whose image coincides with the image of $\tilde{\nu} \in 
\mathcal{M}(\tilde{q},\tilde{p})$, $\epsilon(\tilde{\nu}) = \pm 1$ is the sign determined 
by the orientation on $\mathcal{M}(\tilde{q},\tilde{p})$, and $\gamma^{\tilde{\nu}}
\tilde{\gamma}_{\tilde{q}}$ (the concatenation of two paths in $\widetilde{M}$) is a path from
$\tilde{p}$ to $\tilde{x}_0$. 

Finally, note that we can pick the paths $\gamma^{\tilde{\nu}}$ in the sum for 
$\tilde{\Phi}_{k-1}\left(\tilde{\partial}_k(g \otimes \tilde{q})\right)$ to be lifts of the paths 
$\gamma^\nu$ in the sum for $\partial^G_k(\tilde{\Phi}_k(g \otimes \tilde{q}))$. Moreover, 
there is a bijection between the collection of paths $\gamma^\nu$ in the first sum and the collection
of paths $\gamma^{\tilde{\nu}}$ in the second sum with $\epsilon(\nu) = \epsilon(\tilde{\nu})$ 
if $\gamma^{\tilde{\nu}}$ is a lift of $\gamma^\nu$, since the orientations of the unstable manifolds of 
$(\tilde{f},\widetilde{\mathsf{g}})$ were determined by pulling back the orientations of the 
unstable manifolds of $(f,\mathsf{g})$. Therefore, $\partial^G_k(\tilde{\Phi}_k(g \otimes 
\tilde{q})) = \tilde{\Phi}_{k-1}\left(\tilde{\partial}_k(g \otimes \tilde{q})\right)$.

\proofend

\begin{remark}
Note that the chain complex $(G_{x_0} \otimes_{\pi_1} C_\ast(\tilde{f}), \bar{\partial}_\ast)$
was defined using a basepoint $x_0 \in M$, whereas the definition of $(C_\ast(f;G),\partial^G_\ast)$
does not involve $x_0$. Thus, the Morse Eilenberg Theorem (Theorem \ref{EilenbergMorse}) implies that the
isomorphism class of $H_k((G_{x_0} \otimes_{\pi_1} C_\ast(\tilde{f}), \bar{\partial}_\ast))$
does not depend on the basepoint $x_0 \in M$ for all $k=0,\ldots, m$.

On the other hand, the boundary operator of the twisted Morse-Smale-Witten chain complex
$(C_\ast(f;G),\partial^G_\ast)$ was defined using the homomorphisms that define the bundle of abelian
groups $G$ (Definition \ref{twistedboundary}), whereas the boundary operator for the chain complex
$(G_{x_0} \otimes_{\pi_1} C_\ast(\tilde{f}), \bar{\partial}_\ast)$ was defined using
the representation $\pi_1(M,x_0) \times G_{x_0} \rightarrow G_{x_0}$ induced by the homomorphisms that
define $G$. Thus, Theorem \ref{EilenbergMorse} shows that the homology of the twisted Morse-Smale-Witten 
chain complex $(C_\ast(f;G),\partial^G_\ast)$ only depends on the isomorphism class of $G$ 
(see Theorem \ref{actionsystem}). We will give an independent proof of this fact in the next
section, where we also show that the homology does not depend on the Morse-Smale pair $(f,\mathsf{g})$.
\end{remark}


\section{Homology determined by the isomorphism class of $G$}\label{continuation}

The main goal of this section is to prove an invariance theorem 
(Theorem \ref{homologyindependence}) which shows that on a closed finite dimensional
smooth Riemannian manifold $(M,\mathsf{g})$ the homology of the twisted
Morse-Smale-Witten chain complex $(C_\ast(f;G),\partial_\ast^G)$ is independent of the
Morse-Smale pair $(f,\mathsf{g})$ and depends only on the isomorphism class of the bundle of
abelian groups $G$.

\subsection{A chain map}

Let $\rho: \overline{\mathbb{R}} \rightarrow [-1,1]$ be a smooth strictly increasing
function with $\rho(-\infty) = - 1$, $\rho(+\infty)= 1$, and $\lim_{t \rightarrow\pm \infty}
\rho'(t) = 0$, and let $(f_1,\mathsf{g}_1)$ and $(f_2,\mathsf{g}_2)$ be smooth Morse-Smale pairs on
$M$. By adding a constant to $f_1$ we may assume that $\inf f_1 > \sup f_2$.
Let $F_{21}:\mathbb{R} \times M \rightarrow \mathbb{R}$ be a smooth function that is 
strictly decreasing in its first component such that for some large $T \gg 0$ we have
$$
F_{21}(t,x) = \left\{ \begin{array}{llc}
f_1(x) - \rho(t) & \text{ if } & t < -T \\
h_t(x)   & \text{ if } & -T\leq t \leq T \\
f_2(x) - \rho(t) & \text{ if } & t > T
\end{array}\right.
$$
where $h_t(x)$ is an approximation to $\frac{1}{2T}(T-t)(f_1(x) - \rho(t)) + 
\frac{1}{2T}(T+t)(f_2(x) - \rho(t))$ with $\frac{d}{dt} h_t(x) < 0$ that makes $F_{21}$ smooth.
Extend $F_{21}$ to a smooth function on $\overline{\mathbb{R}} \times M$ by defining 
$F_{21}(-\infty,x) = f_1(x) + 1$ and $F_{21}(+\infty,x) = f_2(x) - 1$. Note that the critical 
points of $F_{21}$ are all on the boundary of $\overline{\mathbb{R}}\times M$ since 
$\frac{d}{dt} F_{21}(t,x) < 0$ for all $t \in \mathbb{R}$.  Pick a Riemannian metric $\mathsf{g}$
on $\overline{\mathbb{R}}\times M$ such that
$$
\mathsf{g}(t,x) = \left\{
\begin{array}{ll}
\mathsf{g}_1(x) + dt^2 & \text{if } t < -T\\
\mathsf{g}_2(x) + dt^2 & \text{if } t > T
\end{array}\right.
$$
for all $x \in M$.
(Compare with Lemma 1.17 of \cite{CorRig}.)

Let $\varphi_\tau:\overline{\mathbb{R}} \times M \rightarrow \overline{\mathbb{R}}\times M$ denote the
flow associated to the negative gradient $-\nabla F_{21}$ with respect to the metric $\mathsf{g}$. Let 
$q_1 \in Cr(f_1)$, $q_2 \in Cr(f_2)$, and define
\begin{eqnarray*}
W^u_F(q_1) & = & \{(t,x) \in \overline{\mathbb{R}} \times M| \lim_{\tau \rightarrow -\infty}
\varphi_\tau(t,x) = (-\infty,q_1) \}\\
W^s_F(q_2) & = & \{(t,x) \in \overline{\mathbb{R}} \times M| \lim_{\tau \rightarrow +\infty}
\varphi_\tau(t,x) = (+\infty,q_2) \}\\
W_F(q_1,q_2) & = & W^u_F(q_1) \cap W^s_F(q_2) \subset \mathbb{R}\times M.
\end{eqnarray*}
Since both $(f_1,\mathsf{g}_1)$ and $(f_2,\mathsf{g}_2)$ are Morse-Smale pairs we can perturb either 
$$
\text{(1) the approximation } h_t(x) \text{ or } \text{(2) the Riemannian metric } \mathsf{g}
$$
so that $W^u_F(q_1) \pitchfork W^s_F(q_2)$ for all $q_1 \in Cr(f_1)$ and $q_2 \in Cr(f_2)$.
(Compare with Proposition 1.12 and Lemma 1.13 of \cite{CorRig}.)

Choosing orientations for the unstable manifolds $W^u(q_1)$ and $W^u(q_2)$ then determines
an orientation on $W_F(q_1,q_2)$ via the short exact sequence
$$
\xymatrix{
T_\ast W_F(q_1,q_2) \ar@{^{(}->}[r] & T_\ast W_F^u(q_1) |_{W_F(q_1,q_2)} \ar[r]
& \nu_\ast(W_F(q_1,q_2), W_F^u(q_1))|_{W_F(q_1,q_2)} \ar[r] & 0
}
$$
where the fibers of the normal bundle are canonically isomorphic to $T_{q_2}W^u_{f_2}(q_2)$ via
the flow of $-\nabla F_{21}$ and $W_F^u(q_1)$ is oriented as follows. 
For $\tau << 0$, $\varphi_\tau(t,x)$ will be in the region where $t < -T$.
In that region the tangent space to $W^u_F(q_1)$ has the product orientation $\mathbb{R} \times
T_\ast W^u_{f_1}(q_1)$, and this orientation then determines an orientation on
$T_{(t,x)} W_F^u(q_1)$ via $\varphi_{-\tau}$. 

Taking a quotient by the action of $\mathbb{R}$ 
given by the negative gradient flow then gives an oriented smooth manifold
$$
\mathcal{M}_F(q_1,q_2) = (W^u_F(q_1) \cap W^s_F(q_2))/\mathbb{R}
$$
of dimension $\lambda_{q_1} - \lambda_{q_2}$, called the moduli space of time dependent
gradient flow lines. The orientation on $\mathcal{M}_F(q_1,q_2)$ is chosen by identifying
the space with a level set and putting $-\nabla F_{21}$ first, analogous to the way the orientation
on $\mathcal{M}(q,p)$ was chosen.


\medskip
Denote the pullback of a bundle of abelian groups $G$ over $M$ under the projection 
$\pi: \overline{\mathbb{R}} \times M \rightarrow M$ by $G^\ast$. Thus, $G^\ast$ associates
to every point $(t,x) \in \overline{\mathbb{R}} \times M$ the abelian group $G_x$ and to every
continuous path $\gamma:[0,1] \rightarrow \overline{\mathbb{R}} \times M$ the 
homomorphism $\gamma_\ast:G_{\pi(\gamma(1))} \rightarrow G_{\pi(\gamma(0))}$
determined by the path $\pi\circ \gamma$.

\begin{definition}[Chain map]\label{chainmapdef}
Let $(f_1,\mathsf{g}_1)$ and $(f_2,\mathsf{g}_2)$ be smooth Morse-Smale pairs on $M$, and let 
$G_1$ and $G_2$ be bundles of abelian groups over $M$.  Assume that $G_1$ and $G_2$ are
isomorphic, so there exists a family of isomorphisms $\Phi:G_1 \rightarrow G_2$ making
the diagram in Definition \ref{isobundles} commute.  Let $(F_{21})_\Box:(C_\ast(f_1;G_1),
\partial_\ast^{G_1}) \rightarrow (C_\ast(f_2;G_2), \partial_\ast^{G_2})$ be the
linear map defined on an elementary chain $gq_1 \in C_k(f_1;G_1)$ by
$$
(F_{21})_\Box(gq_1)\ \ = (-1)^{\lambda_{q_1}} \sum_{q_2 \in Cr_k(f_2)} \sum_{\nu_{F} 
\in \mathcal{M}_F(q_1,q_2)} \epsilon(\nu_{F}) (\gamma_F)_\ast(\Phi(g)) q_2
$$
where $\gamma_{F}:[0,1] \rightarrow \overline{\mathbb{R}}\times M$ is any continuous
path from $(+\infty,q_2)$ to $(-\infty,q_1)$ such that the image $\gamma_F((0,1))$coincides
with $Im(\nu_{F})$, $\Phi(g) \in (G_2)_{q_1} = (G^\ast_2)_{(-\infty,q_1)}$, 
$\epsilon(\nu_{F}) = \pm 1$ is the sign determined by the orientation on 
$\mathcal{M}_F(q_1,q_2)$, and $\lambda_{q_1}$ denotes the index of $q_1$ as
a critical point of $f_1:M \rightarrow \mathbb{R}$.
\end{definition}

\smallskip\noindent
Note: If $G_1=G_2 = \mathbb{Z}$ and $\Phi$ is the identity, then $(F_{21})_\Box(q_1)$
is $(-1)^{\lambda_{q_1}}$ times the continuation map defined in Section 4.1.3 of \cite{SchMor}.

\medskip
The above sum over $\mathcal{M}_F(q_1,q_2)$ is finite because of the following result,
cf. Section 2.4.3 of \cite{SchMor}.

\begin{theorem}[Compactification]\label{timedepcompactification}
For any $q_1 \in Cr(f_1)$ and $q_2 \in Cr(f_2)$ the moduli space $\mathcal{M}_F(q_1,q_2)$ has
a compactification $\overline{\mathcal{M}}_F(q_1,q_2)$, consisting of all the piecewise gradient
flow lines from $(-\infty,q_1)$ to $(+\infty,q_2)$, including both time dependent and time
independent gradient flow lines.
\end{theorem}
We will now show how this structure on the $1$-dimensional compactified moduli spaces of
time dependent gradient flow lines implies that $(F_{21})_\Box$ is a chain map. In order to
distinguish between the time independent gradient flow lines and the time dependent gradient
flow lines, in the following we will denote the (zero dimensional) moduli spaces of time independent
gradient flow lines of $f_1:\{-\infty\} \times M \rightarrow \mathbb{R}$ by $\mathcal{M}_{f_1}(q_1,p_1)$,
the time independent moduli spaces of gradient flow lines of $f_2:\{+\infty\} \times M \rightarrow
\mathbb{R}$ by $\mathcal{M}_{f_2}(q_2,p_2)$, and the time dependent gradient flow lines 
of $F_{21}: \overline{\mathbb{R}} \times M \rightarrow \mathbb{R}$ by $\mathcal{M}_F(q_1,q_2)$
or $\mathcal{M}_F(p_1,p_2)$. 

If $\lambda_{q_1} - \lambda_{p_2} = 1$, then $\overline{\mathcal{M}}_F(q_1,p_2)$ is
a $1$-dimensional smooth manifold with boundary, and
$$
\partial \overline{\mathcal{M}}_F(q_1,p_2) = \left(\bigcup_{p_1} \mathcal{M}_{f_1}(q_1,p_1)
\times \mathcal{M}_F(p_1,p_2) \right) \coprod \left(\bigcup_{q_2} \mathcal{M}_F(q_1,q_2) \times 
\mathcal{M}_{f_2}(q_2,p_2)\right).
$$ 
So, there are three possibilities for a non-empty boundary of a path component 
$\overline{\mathcal{M}}_F(q_1,p_2;[\nu_F])$.

\smallskip\noindent
1) We have $\partial \overline{\mathcal{M}}_F(q_1,p_2;[\nu_F]) = \{(\nu_1,\nu_F^p)$, 
$(\tilde{\nu}_1,\tilde{\nu}_F^p)\}$, where $\nu_1 \in \mathcal{M}_{f_1}(q_1,p_1)$,
$\nu_F^p \in \mathcal{M}_F(p_1,p_2)$, $\tilde{\nu}_1 \in \mathcal{M}_{f_1}(q_1,\tilde{p}_1)$,
and $\tilde{\nu}_F^p \in \mathcal{M}_F(\tilde{p}_1,p_2)$.
\begin{figure}[h]
\includegraphics{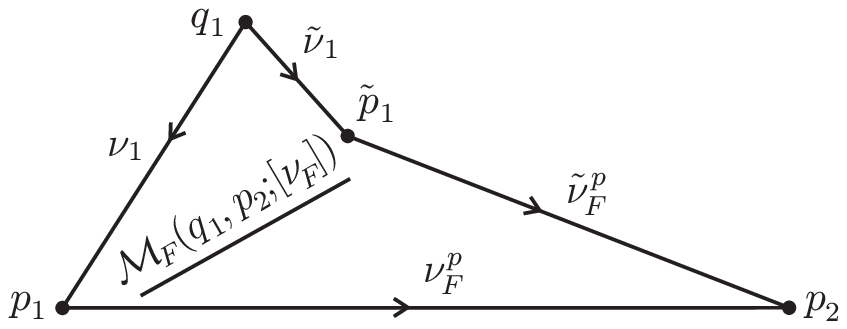}
\end{figure} 

\noindent
2) We have $\partial \overline{\mathcal{M}}_F(q_1,p_2;[\nu_F]) = \{(\nu_1,\nu_F^p)$, $(\nu_F^q,\nu_2)\}$, 
where $\nu_1 \in \mathcal{M}_{f_1}(q_1,p_1)$, $\nu_F^p \in \mathcal{M}_F(p_1,p_2)$, $\nu_F^q \in 
\mathcal{M}_F(q_1,q_2)$, $\nu_2 \in \mathcal{M}_{f_2}(q_2,p_2)$.
\begin{figure}[h]
\includegraphics{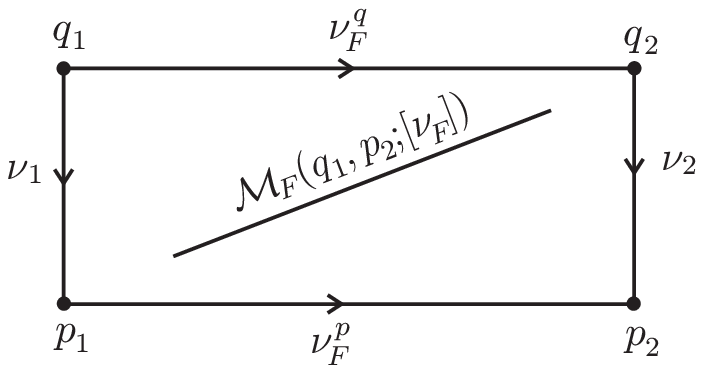}
\end{figure}

\noindent
3) We have $\partial \overline{\mathcal{M}}_F(q_1,p_2;[\nu_F]) = \{(\nu_F^q,\nu_2)$, 
$(\tilde{\nu}_F^q,\tilde{\nu}_2)\}$, where $\nu_F^q \in \mathcal{M}_F(q_1,q_2)$, 
$\nu_2 \in \mathcal{M}_{f_2}(q_2,p_2)$, $\tilde{\nu}_F^q \in \mathcal{M}_F(q_1,\tilde{q}_2)$,
$\tilde{\nu}_2 \in \mathcal{M}_{f_2}(\tilde{q}_2,p_2)$.
\begin{figure}[h]
\includegraphics{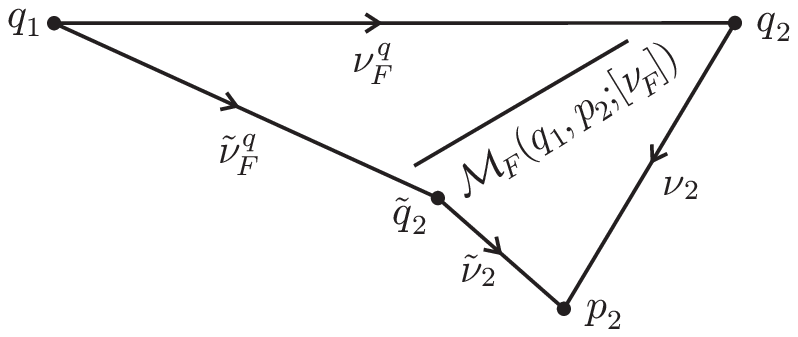}
\end{figure} 

\noindent
In all three cases we have orientation formulas analogous to those in Lemma 
\ref{componentboundary}.
\begin{enumerate}
\item $\partial \overline{\mathcal{M}}_F(q_1,p_2;[\nu_F]) = 
(\tilde{\nu}_1,\tilde{\nu}_F^p) - (\nu_1,\nu_F^p)$ and  
$\epsilon(\tilde{\nu}_1)\epsilon(\tilde{\nu}_F^p) + \epsilon(\nu_1)\epsilon(\nu_F^p) = 0$.
\item $\partial \overline{\mathcal{M}}_F(q_1,p_2;[\nu_F]) = 
(\nu_F^q,\nu_2) - (\nu_1,\nu_F^p)$ and 
$\epsilon(\nu_F^q)\epsilon(\nu_2) + \epsilon(\nu_1)\epsilon(\nu_F^p) = 0$.
\item $\partial \overline{\mathcal{M}}_F(q_1,p_2;[\nu_F]) = 
(\nu^q_F,\nu_2) - (\tilde{\nu}^q_F,\tilde{\nu}_2)$ and  
$\epsilon(\nu_F^q)\epsilon(\nu_2)  + \epsilon(\tilde{\nu}_F^q)\epsilon(\tilde{\nu}_2) = 0$.
\end{enumerate}
These formulas can be verified by first checking them for specific orientations chosen
on $W^u(q_1)$, $W^u(q_2)$, $W^u(p_1)$, and $W^u(p_2)$ and then noting that changing an
orientation changes a complimentary pair of signs.

\begin{proposition}\label{chainmap}
The map $(F_{21})_\Box:(C_\ast(f_1;G_1),\partial_\ast^{G_1}) \rightarrow 
(C_\ast(f_2;G_2), \partial_\ast^{G_2})$ is a chain map. That is, the map preserves 
degree and $\partial_\ast^{G_2} \circ (F_{21})_\Box = (F_{21})_\Box \circ 
\partial_\ast^{G_1}$.
\end{proposition}

\proofstart
Let $gq_1 \in C_k(f_1;G_1)$ be an elementary chain for some $k=1,\ldots, m$ where 
$m = \text{dim }M$. We have $\partial_k^{G_2}((F_{21})_\Box(gq_1))$
$$
\begin{array}{lcl}
& = & \displaystyle \partial_k^{G_2}\left((-1)^{\lambda_{q_1}} \sum_{q_2 \in Cr_k(f_2)}
      \sum_{\nu_{F}^q \in \mathcal{M}_F(q_1,q_2)} \epsilon(\nu_{F}^q)
      (\gamma_{F}^q)_\ast(\Phi(g)) q_2 \right) \\
& = & \displaystyle (-1)^{k} \sum_{\stackrel{q_2 \in Cr_k(f_2)}{\nu_{F}^q \in 
      \mathcal{M}_F(q_1,q_2)}} \epsilon(\nu_{F}^q) \partial_k^{G_2}
      \left((\gamma_{F}^q)_\ast(\Phi(g)) q_2\right) \\
& = & \displaystyle (-1)^{k} \sum_{\stackrel{q_2 \in Cr_k(f_2)}{\nu_{F}^q \in 
      \mathcal{M}_F(q_1,q_2)}} \epsilon(\nu_{F}^q) \sum_{\stackrel{p_2 \in 
      Cr_{k-1}(f_2)}{\nu_2 \in \mathcal{M}_{f_2}(q_2,p_2)}} \epsilon(\nu_2) 
      (\gamma_2)_\ast\left((\gamma_{F}^q)_\ast(\Phi(g))\right) p_2 \\
& = & \displaystyle (-1)^{k}\sum_{\stackrel{q_2 \in Cr_k(f_2)}{\nu_{F}^q \in 
      \mathcal{M}_F(q_1,q_2)}} \sum_{\stackrel{p_2 \in Cr_{k-1}(f_2)}{\nu_2 \in 
      \mathcal{M}_{f_2}(q_2,p_2)}} \epsilon(\nu_{F}^q) \epsilon(\nu_2) 
      (\gamma_2(\pi\circ\gamma_{F}^q))_\ast(\Phi(g)) p_2
\end{array}
$$
where $\gamma_{F}^q$ parameterizes $\nu_{F}^q$ from $(+\infty,q_2)$ to 
$(-\infty,q_1)$ and $\gamma_2$ parameterizes $\nu_2$ from $p_2$ to $q_2$.
Whereas, $(F_{21})_\Box(\partial_k^{G_1}(gq_1))$
$$
\begin{array}{lcl}
& = & \displaystyle (F_{21})_\Box\left(\sum_{p_1 \in Cr_{k-1}(f_1)} \sum_{\nu_1 \in 
      \mathcal{M}_{f_1}(q_1,p_1)} \epsilon(\nu_1)(\gamma_1)_\ast(g) p_1 \right) \\
& = & \displaystyle \sum_{\stackrel{p_1 \in Cr_{k-1}(f_1)}{\nu_1 \in \mathcal{M}_{f_1}(q_1,p_1)}}
      \epsilon(\nu_1) (F_{21})_\Box\left((\gamma_1)_\ast(g)p_1\right) \\
& = & \displaystyle (-1)^{\lambda_{p_1}} 
      \sum_{\stackrel{p_1 \in Cr_{k-1}(f_1)}{\nu_1 \in \mathcal{M}_{f_1}(q_1,p_1)}}
      \epsilon(\nu_1)
      \sum_{\stackrel{p_2 \in Cr_{k-1}(f_2)}{\nu_{F}^p\in \mathcal{M}_F(p_1,p_2)}}
      \epsilon(\nu_{F}^p)(\gamma_{F}^p)_\ast\left(\Phi((\gamma_1)_\ast(g))\right) p_2 \\
& = & \displaystyle (-1)^{k-1}
      \sum_{\stackrel{p_1 \in Cr_{k-1}(f_1)}{\nu_1 \in \mathcal{M}_{f_1}(q_1,p_1)}} 
      \sum_{\stackrel{p_2 \in Cr_{k-1}(f_2)}{\nu_{F}^p\in \mathcal{M}_F(p_1,p_2)}} 
      \epsilon(\nu_1) \epsilon(\nu_{F}^p) 
     ((\pi\circ\gamma_{F}^p)\gamma_1)_\ast(\Phi(g))p_2 \\
\end{array}
$$
where $\gamma_{F}^p$ parameterizes $\nu_{F}^p$ from $(+\infty,p_2)$ to $(-\infty,p_1)$ 
and $\gamma_1$ parameterizes $\nu_1$ from $p_1$ to $q_1$.

Now recall that the terms in the above sums for $\partial_k^{G_2}((F_{21})_\Box(gq_1))$
and $(F_{21})_\Box(\partial_k^{G_1}(gq_1))$ correspond to the boundary points of the
path components of $\overline{\mathcal{M}}_F(q_1,p_2)$. Considering the three
cases discussed before the statement of the proposition, we observe the following.
The terms corresponding to boundary points in case 1) cancel each other out
in the sum for $(F_{21})_\Box(\partial_k^{G_1}(gq_1))$ because the homomorphism 
$((\pi\circ\gamma_{F}^p)\gamma_1)_\ast$ is constant on each path component and
$\epsilon(\tilde{\nu}_1)\epsilon(\tilde{\nu}_F^p) + \epsilon(\nu_1)\epsilon(\nu_F^p) = 0$.
Similarly, the terms corresponding to boundary points in case 3) cancel each other out
in the sum for $\partial_k^{G_2}((F_{21})_\Box(gq_1))$ because the homomorphism 
$(\gamma_2(\pi\circ\gamma_{F}^q))_\ast$ is constant on each path component and
$\epsilon(\nu_F^q)\epsilon(\nu_2)  + \epsilon(\tilde{\nu}_F^q)\epsilon(\tilde{\nu}_2) = 0$.
Finally, for the boundary points in case 2) we have 
\begin{eqnarray*}
((\pi\circ\gamma_{F}^p)\gamma_1)_\ast(\Phi(g))p_2 & = &
 (\gamma_2(\pi\circ\gamma_{F}^q))_\ast(\Phi(g)) p_2\\
\epsilon(\nu_1) \epsilon(\nu_{F}^p) & = & - \epsilon(\nu_{F}^q) \epsilon(\nu_2)
\end{eqnarray*}
for each pair of endpoints in the same path component $[(\nu_1,\nu_F^p)] = [(\nu_F^q,\nu_2)]$.
Therefore, $\partial_k((F_{21})_\Box(q_1)) = (F_{21})_\Box(\partial_k(q_1))$.

\proofend

\begin{corollary}
The map $(F_{21})_\Box$ induces a homomorphism in homology
$$
(F_{21})_\ast:H_k(C_\ast(f_1;G_1),\partial_\ast^{G_1}) \rightarrow 
H_k(C_\ast(f_2;G_2), \partial_\ast^{G_2})
$$
for all $k=0,\ldots ,m$.
\end{corollary}


\subsection{A chain homotopy}

Assume that we have four Morse-Smale pairs $(f_j,\mathsf{g}_j)$ and
four isomorphic bundles of abelian groups $G_j$, where $j=1,2,3,4$.
For $j=2,3,4$ there are families of isomorphisms $\Phi_{j1}:G_1 \rightarrow G_j$
making the diagram in Definition \ref{isobundles} commute, and if we define 
$\Phi_{42} = \Phi_{41} \circ \Phi^{-1}_{21}$, $\Phi_{43} = \Phi_{41} \circ \Phi^{-1}_{31}$,
$\Phi_{ij} = \Phi_{ji}^{-1}$ when $i<j$, and $\Phi_{jj} = id$ for all $j$, then for all 
$i,j=1,2,3,4$ we have a family of isomorphisms $\Phi_{ji}:G_i \rightarrow G_j$ making 
the diagram in Definition \ref{isobundles} commute such that $\Phi_{42} \circ \Phi_{21} 
= \Phi_{41}$ and $\Phi_{43} \circ \Phi_{31} = \Phi_{41}$. Also, for all $i,j=1,2,3,4$ 
there are chain maps
$$
(F_{ji})_\Box:(C_\ast(f_i;G_i),\partial_\ast^{G_i}) \rightarrow 
(C_\ast(f_j;G_j), \partial_\ast^{G_j})
$$ 
from Definition \ref{chainmapdef} and Proposition \ref{chainmap}, defined using the family of
isomorphisms $\Phi_{ji}$ and functions $F_{ji}:\overline{\mathbb{R}} \times M  \rightarrow \mathbb{R}$. 

We will show that under these assumptions the moduli spaces of gradient flow lines of
a smooth function $H:\overline{\mathbb{R}} \times \overline{\mathbb{R}} \times M \rightarrow \mathbb{R}$
meeting certain transversality requirements can be used to construct a chain homotopy between the
chain maps $(F_{42})_\Box \circ (F_{21})_\Box$ and $(F_{43})_\Box \circ (F_{31})_\Box$.
By adding constants to $f_1, f_2$, and $f_3$ we may assume that $\inf f_1 > \sup f_2$, 
$\inf f_1 > \sup f_3$, $\inf f_2 > \sup f_4$, and $\inf f_3 > \sup f_4$.

Following \cite{BanMor} and \cite{WebThe} we let $H: \overline{\mathbb{R}} \times \overline{\mathbb{R}}
\times M \rightarrow  \mathbb{R}$ be a smooth function that is strictly decreasing
in its first two components such that for some large $T \gg 0$ we have
$$
H(s,t,x) = \left\{ \begin{array}{lllll}
f_1(x) - \rho(s) - \rho(t) & \text{ if } & s < -T & \text{ and } & t < -T\\
f_2(x) - \rho(s) - \rho(t) & \text{ if } & s > T  & \text{ and } & t < -T\\
f_3(x) - \rho(s) - \rho(t) & \text{ if } & s < -T & \text{ and } & t > T\\
f_4(x) - \rho(s) - \rho(t) & \text{ if } & s > T  & \text{ and } & t > T
\end{array}\right.
$$
where $\rho: \overline{\mathbb{R}} \rightarrow [-1,1]$ is a smooth strictly increasing
function with $\rho(-\infty) = - 1$, $\rho(+\infty)= 1$, and 
$\lim_{t \rightarrow\pm \infty} \rho'(t) = 0$.
The critical points of $H$ are all on the boundary of $\overline{\mathbb{R}} \times \overline{\mathbb{R}}
\times M$ since $H$ is strictly decreasing in its first two components.  Pick a Riemannian metric 
$\mathsf{g}$ on $\overline{\mathbb{R}} \times \overline{\mathbb{R}}  \times M$ such that 
$$
\mathsf{g}(s,t,x) = \left\{
\begin{array}{lllll}
\mathsf{g}_1(x) + ds^2 + dt^2 & \text{ if } & s < -T & \text{ and } & t < -T\\
\mathsf{g}_2(x) + ds^2 + dt^2 & \text{ if } & s > T  & \text{ and } & t < -T\\
\mathsf{g}_3(x) + ds^2 + dt^2 & \text{ if } & s < -T & \text{ and } & t > T\\
\mathsf{g}_4(x) + ds^2 + dt^2 & \text{ if } & s > T  & \text{ and } & t > T
\end{array}\right.
$$
for all $x \in M$.
The negative gradient flow of $H:\overline{\mathbb{R}} \times \overline{\mathbb{R}} \times M \rightarrow \mathbb{R}$
with respect to $\mathsf{g}$ can be pictured as follows, where $M$ is in the vertical direction.

\begin{center}
\includegraphics{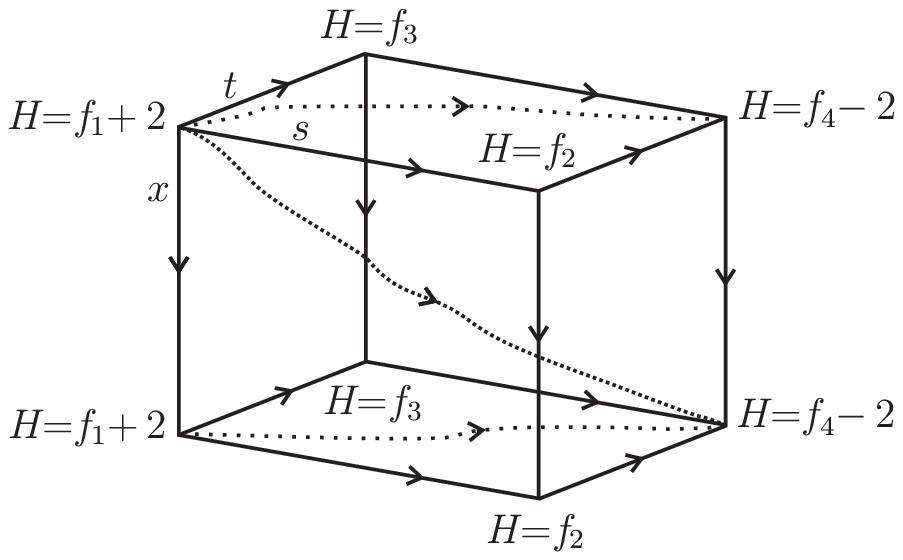}
\end{center}

Let $\varphi_\tau:\overline{\mathbb{R}} \times \overline{\mathbb{R}} \times M \rightarrow  
\overline{\mathbb{R}} \times \overline{\mathbb{R}} \times M$ be the flow associated to $-\nabla H$, the 
negative gradient of $H$ with respect to $\mathsf{g}$, and for $p_1 \in Cr(f_1)$ and $q_j \in Cr(f_j)$ for some
$j=2,3,4$ define
\begin{eqnarray*}
W^u_H(p_1) & = & \{(s,t,x) \in \overline{\mathbb{R}} \times \overline{\mathbb{R}} \times M |
\lim_{\tau \rightarrow -\infty} \varphi_\tau(s,t,x) = (-\infty,-\infty,p_1) \}\\
W^u_H(q_2) & = & \{(s,t,x) \in \overline{\mathbb{R}} \times \overline{\mathbb{R}} \times M |
\lim_{\tau \rightarrow -\infty} \varphi_\tau(s,t,x) = (+\infty,-\infty,q_2) \}\\
W^s_H(q_2) & = & \{(s,t,x) \in \overline{\mathbb{R}} \times \overline{\mathbb{R}} \times M |
\lim_{\tau \rightarrow +\infty} \varphi_\tau(s,t,x) = (+\infty,-\infty,q_2) \}\\
W^u_H(q_3) & = & \{(s,t,x) \in \overline{\mathbb{R}} \times \overline{\mathbb{R}} \times M |
\lim_{\tau \rightarrow -\infty} \varphi_\tau(s,t,x) = (-\infty,+\infty,q_3) \}\\
W^s_H(q_3) & = & \{(s,t,x) \in \overline{\mathbb{R}} \times \overline{\mathbb{R}} \times M |
\lim_{\tau \rightarrow +\infty} \varphi_\tau(s,t,x) = (-\infty,+\infty,q_3) \}\\
W^s_H(q_4) & = & \{(s,t,x) \in \overline{\mathbb{R}} \times \overline{\mathbb{R}} \times M |
\lim_{\tau \rightarrow +\infty} \varphi_\tau(s,t,x) = (+\infty,+\infty,q_4) \}\\
W_H(p_1,q_j) & = & W^u_H(p_1) \cap W^s_H(q_j)  \subset \overline{\mathbb{R}} \times \overline{\mathbb{R}} \times M \\
W_H(q_j,q_4) & = & W^u_H(q_j) \cap W^s_H(q_4)  \subset \overline{\mathbb{R}} \times \overline{\mathbb{R}} \times M.
\end{eqnarray*}
Since $f_j$ satisfies the Morse-Smale transversality condition with respect to $\mathsf{g}$
for all $j=1,2,3,4$, we can choose $H$ so that
\begin{itemize}
\item  $W^u_H(p_1) \pitchfork W^s_H(q_j)$ for all $p_1 \in Cr(f_1)$ and $q_j \in Cr(f_j)$ for $j=2,3,4$,
\item  $W^u_H(q_j) \pitchfork W^s_H(q_4)$ for all $q_j \in Cr(f_j)$ and $q_4 \in Cr(f_4)$ for $j=2,3$.
\end{itemize}
These conditions imply that $H$ restricted to each face of $\overline{\mathbb{R}} \times \overline{\mathbb{R}}
\times M$ defines a chain map. That is, if we define $F_{21}(s,x) = H(s,-\infty,x)$,
$F_{42}(t,x) = H(+\infty,t,x)$, $F_{31}(t,x) = H(-\infty,t,x)$, and $F_{43}(s,x) = H(s,+\infty,x)$,
then the chain maps $(F_{21})_\Box$, $(F_{42})_\Box$, $(F_{31})_\Box$, $(F_{43})_\Box$ from 
Definition \ref{chainmapdef} and Proposition \ref{chainmap} are defined. Conversely, given functions
$F_{21}$, $F_{42}$, $F_{31}$, $F_{43}$ on the faces of $\overline{\mathbb{R}} \times \overline{\mathbb{R}}
\times M$ that satisfy the transversality conditions needed to define the chain maps in Definition 
\ref{chainmapdef}, we can choose an $H$ meeting the above transversality conditions that agrees with
these functions on the faces of $\overline{\mathbb{R}} \times \overline{\mathbb{R}} \times M$.

Choosing orientations for the unstable manifolds $W^u_{f_1}(p_1)$ and $W^u_{f_4}(q_4)$ then determines an 
orientation on $W_H(p_1,q_4)$ via the short exact sequence
$$
\xymatrix{
T_\ast W_H(p_1,q_4) \ar@{^{(}->}[r] & T_\ast W_H^u(p_1) |_{W_H(p_1,q_4)} \ar[r]
& \nu_\ast(W_H(p_1,q_4), W_H^u(p_1))|_{W_H(p_1,q_4)} \ar[r] & 0
}
$$
where the fibers of the normal bundle are canonically isomorphic to $T_{q_4}W^u_{f_4}(q_4)$ via
the flow of $-\nabla H$ and $W_H^u(q_1)$ is oriented as follows. 
For $\tau << 0$, $\varphi_\tau(s,t,x)$ will be in the region where $s,t < -T$.
In that region the tangent space to $W^u_H(q_1)$ has the product orientation 
$\mathbb{R} \times \mathbb{R} \times T_\ast W^u_{f_1}(q_1)$, and this orientation then determines an
orientation on $T_{(s,t,x)} W_H^u(q_1)$ via $\varphi_{-\tau}$. 
Taking a quotient by the action of $\mathbb{R}$ given by the negative gradient flow then gives an
oriented smooth manifold
$$
\mathcal{M}_H(p_1,q_4) = (W^u_H(p_1) \cap W^s_H(q_4))/\mathbb{R},
$$
of dimension $\lambda_{p_1} - \lambda_{q_4} + 1$, where the orientation on
$\mathcal{M}_H(p_1,q_4)$ is chosen by putting $-\nabla H$ first.

\begin{definition}[Chain homotopy]\label{chainhomotopydef}
Let $f_j:M \rightarrow \mathbb{R}$ be smooth Morse-Smale functions on $M$ and let $G_j$ be
isomorphic bundles of abelian groups over $M$, for $j=1,2,3,4$. Let $G_4^\ast$ denote the
pullback of $G_4$ under the projection $\pi:\overline{\mathbb{R}} \times \overline{\mathbb{R}} \times M 
\rightarrow M$, and let $\Phi_{41}:G_1 \rightarrow G_4$ be a family of isomorphisms making the diagram in
Definition \ref{isobundles} commute.  Assume that $H: \overline{\mathbb{R}} \times \overline{\mathbb{R}}
\times M \rightarrow \mathbb{R}$ and the metric $\mathsf{g}$ on $\overline{\mathbb{R}} \times 
\overline{\mathbb{R}} \times M$ satisfy the conditions listed above, and define 
$H_\Box:(C_\ast(f_1;G_1),\partial_\ast^{G_1}) \rightarrow (C_\ast(f_4;G_4), \partial_\ast^{G_4})$
to be the linear map given on an elementary chain $gp_1 \in C_k(f_1;G_1)$ by 
$$
H_\Box(gp_1)\ \ = \sum_{q_4 \in Cr_{k+1}(f_4)} \sum_{\nu_H 
\in \mathcal{M}_H(p_1,q_4)} \epsilon(\nu_H) (\gamma_H)_\ast(\Phi_{41}(g))q_4
$$
where $\gamma_H:[0,1] \rightarrow \overline{\mathbb{R}} \times \overline{\mathbb{R}} \times M$
is any continuous path from $(+\infty,+\infty,q_4)$ to $(-\infty,-\infty,p_1)$
such that the image $\gamma_H((0,1))$ coincides with  $Im(\nu_H)$, $\Phi_{41}(g) \in
(G_4)_{p_1} = (G_4^\ast)_{(-\infty,-\infty,p_1)}$, $\epsilon(\nu_H) = \pm 1$ is the
sign determined by the orientation on $\mathcal{M}_H(p_1,q_4)$, and $\lambda_{p_1}$ denotes
the index of $p_1$ as a critical point of $f_1:M \rightarrow \mathbb{R}$.  
\end{definition}

There are compactification results for the moduli spaces of $H$ similar to those stated in the
previous subsection, cf. Section 2.4.4 of \cite{SchMor}. Hence, the moduli space 
$\mathcal{M}_H(p_1,q_4)$ has a compactification $\overline{\mathcal{M}}_H(p_1,q_4)$, 
consisting of all the piecewise gradient flow lines of $H$ from $(-\infty,-\infty,p_1)$ to
$(+\infty,+\infty,q_4)$. This implies that the above sum over $\mathcal{M}_H(p_1,q_4)$ is finite.

\begin{proposition}\label{chainhomotopy}
Let $f_j:M \rightarrow \mathbb{R}$ be smooth Morse-Smale functions on $M$ and let $G_j$ be
isomorphic bundles of abelian groups over $M$, for $j=1,2,3,4$. 
The map $H_\Box:(C_\ast(f_1;G_1),\partial_\ast^{G_1}) \rightarrow 
(C_\ast(f_4;G_4), \partial_\ast^{G_4})$ is a chain homotopy between 
$(F_{42})_\Box \circ (F_{21})_\Box$ and $(F_{43})_\Box \circ (F_{31})_\Box$, i.e.
$$
(F_{43})_\Box \circ (F_{31})_\Box - (F_{42})_\Box \circ (F_{21})_\Box
= \partial_{k+1}^{G_4} H_\Box + H_\Box \partial_k^{G_1}
$$
for all $k=0,\ldots, m$.
\end{proposition}

\proofstart
Choose families of isomorphisms $\Phi_{ji}:G_i \rightarrow G_j$ for $i,j=1,2,3,4$ making 
the diagram in Definition \ref{isobundles} commute such that $\Phi_{42} \circ \Phi_{21} =
\Phi_{43} \circ \Phi_{31} = \Phi_{41}$. Let $q_1 \in Cr_{k}(f_1)$ and $q_4 \in Cr_{k}(f_4)$
for some $k=0,\ldots, m$, where $m = \text{dim }M$, and consider the compactified moduli space 
$\overline{\mathcal{M}}_H(q_1,q_4)$, where $H:\overline{\mathbb{R}} \times 
\overline{\mathbb{R}} \times M \rightarrow \mathbb{R}$ is a smooth function satisfying the
transversality conditions listed before the proposition such that $H(s,-\infty,x) = F_{21}(s,x)$,
$H(+\infty,t,x) = F_{42}(t,x)$, $H(-\infty,t,x) = F_{31}(t,x)$, and $H(s,+\infty,x) = F_{43}(s,x)$.

Since $\lambda_{q_1} = \lambda_{q_4}$, this compactified moduli space is an oriented compact
smooth manifold with boundary of dimension one, and hence it is diffeomorphic to a disjoint 
union of oriented intervals. The endpoints of these intervals correspond to broken gradient flow
lines of $H:\overline{\mathbb{R}} \times \overline{\mathbb{R}} \times M \rightarrow \mathbb{R}$
from $(-\infty,-\infty,q_1)$ to $(+\infty,+\infty,q_4)$, i.e. elements of the following
zero dimensional spaces.
\begin{eqnarray}
\mathcal{M}_{H}(q_1,q_2) \times \mathcal{M}_{H}(q_2,q_4) & \text{ some } &
q_2 \in Cr_{k}(f_2)\\
\mathcal{M}_{H}(q_1,q_3) \times \mathcal{M}_{H}(q_3,q_4) & \text{ some } &
q_3 \in Cr_{k}(f_3)\\
\mathcal{M}_{f_1}(q_1,p_1) \times \mathcal{M}_H(p_1,q_4) & \text{ some } &
p_1 \in Cr_{k-1}(f_1)\\
\mathcal{M}_H(q_1,r_4) \times \mathcal{M}_{f_4}(r_4,q_4) & \text{ some } &
r_4 \in Cr_{k+1}(f_4)
\end{eqnarray}

\begin{center}
\includegraphics{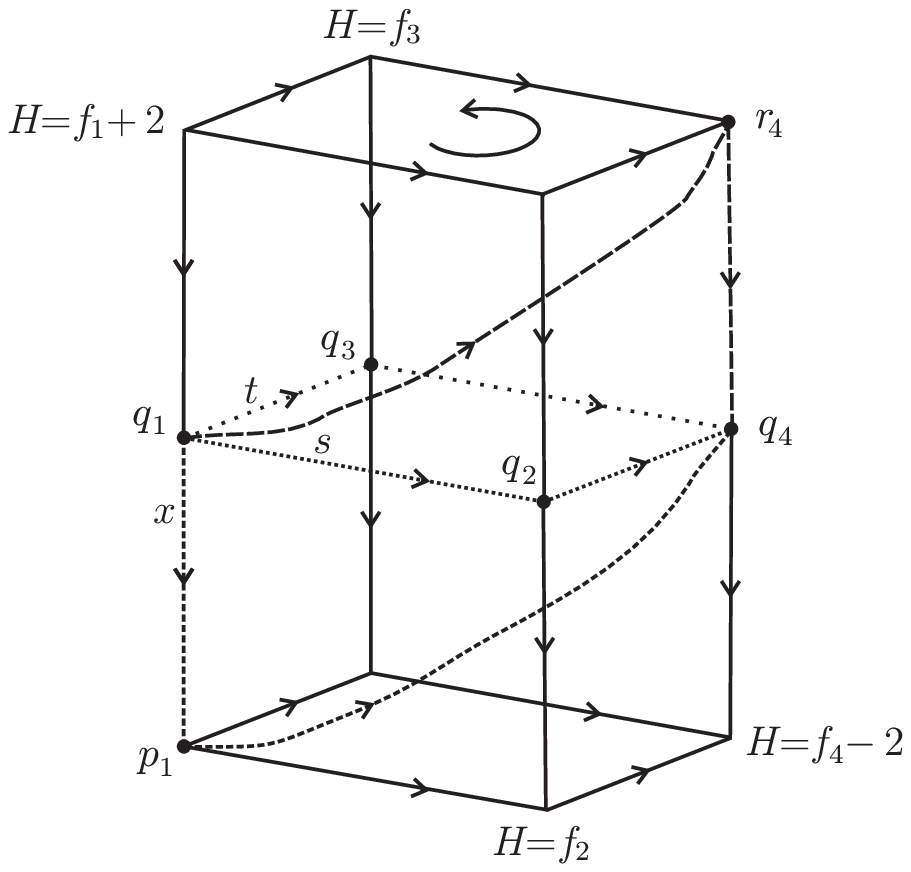}
\end{center}

Now let $g_1 \in (G_1)_{q_1}$ and note that up to sign we have the following.
\begin{enumerate}
\item $((F_{42})_\Box \circ (F_{21})_\Box)(g_1q_1)$ is defined by summing over elements of type $(1)$
\item $((F_{43})_\Box \circ (F_{31})_\Box)(g_1q_1)$ is defined by summing over elements of type $(2)$
\item $(H_\Box \partial_k^{G_1})(g_1q_1)$ is defined by summing over elements of type $(3)$
\item $(\partial_{k+1}^{G_4} H_\Box)(g_1q_1)$ is defined by summing over  elements of type $(4)$
\end{enumerate}
However, these sums are all defined using signs determined by orientations, and the orientations 
on the moduli spaces of $H$ aren't necessarily the same as the orientations on the moduli spaces 
of $F_{ji}$ for $(i,j)=(1,2)$, $(1,3)$, $(2,4)$, $(3,4)$.

\smallskip
Choose orientations for $W^u_{f_j}(q_j)$ for $j=1,2,3,4$. These orientations then determine
orientations on $W^u_H(q_i)$ and $W^u_{F_{ji}}(q_i)$ for $(i,j)=(1,2)$, $(1,3)$, $(2,4)$, $(3,4)$
using the product orientation and the gradient flow. Note that $W^u_H(q_2) = W^u_{F_{42}}(q_2)$,
$W^u_H(q_3) = W^u_{F_{43}}(q_3)$, and $W^u_H(q_4) = W^u_{f_4}(q_4)$ as oriented manifolds,
and we have the following dimensions.
\begin{eqnarray*}
\text{dim }W^u_H(q_1) & = & k+2\\
\text{dim }W^u_H(q_2) = \text{dim }W^u_H(q_3) & = & k+1\\
\text{dim }W^u_H(q_4) & = & k.
\end{eqnarray*}
The spaces $W_H(q_i,q_j)$ are oriented using the short exact sequences
$$
\xymatrix{
T_\ast W_H(q_i,q_j) \ar@{^{(}->}[r] & T_\ast W_H^u(q_i) |_{W_H(q_i,q_j)} \ar[r]
& \nu_\ast(W_H(q_i,q_j), W_H^u(q_i))|_{W_H(q_i,q_j)} \ar[r] & 0
}
$$
where the fibers of the normal bundle are canonically isomorphic to $T_{q_j}W^u_H(q_j)$, 
and $\mathcal{M}_H(q_i,q_j) = W_H(q_i,q_j)/\mathbb{R}$ is oriented by putting $-\nabla H$ first.
The spaces $W_{F_{ji}}(q_i,q_j)$ are oriented similarly using the short
exact sequences
$$
\xymatrix{
T_\ast W_{F_{ji}}(q_i,q_j) \ar@{^{(}->}[r] & T_\ast W_{F_{ji}}^u(q_i) |_{W_{F_{ji}}(q_i,q_j)} \ar[r]
& \nu_\ast(W_{F_{ji}}(q_i,q_j), W_{F_{ji}}^u(q_i)) |_{W_{F_{ji}}(q_i,q_j)} \ar[r] & 0
}
$$
where the fibers of the normal bundle are canonically isomorphic to $T_{q_j}W^u_{f_j}(q_j)$,
and $\mathcal{M}_{F_{ji}}(q_i,q_j) = W_{F_{ji}}(q_i,q_j)/\mathbb{R}$ is oriented by putting $-\nabla F_{ji}$ first. 

\smallskip
Since $W^u_H(q_2) = W^u_{F_{42}}(q_2)$ and $W^u_H(q_4) = W^u_{f_4}(q_4)$ as oriented manifolds,
the above orientation conventions imply that
$$
W_H(q_2,q_4) = W_{F_{42}}(q_2,q_4)
$$
as oriented manifold. Similarly, $W^u_H(q_3) = W^u_{F_{43}}(q_3)$ and $W^u_H(q_4) = W^u_{f_4}(q_4)$
imply
$$
W_H(q_3,q_4) = W_{F_{43}}(q_3,q_4).
$$

Now, let $y = (s,-\infty,x)\in W^u_H(q_1)$ be in the region where $H(s,t,x) = f_1(x) - 
\rho(s) - \rho(t)$, and let $(v_1,v_2,w_1, w_2, \ldots, w_k)$ be a positive basis for 
$T_y W^u_H(q_1) = \mathbb{R} \times \mathbb{R} \times T_x W^u_{f_1}(q_1)$. 
Then $(v_1,w_1,w_2,\ldots, w_k)$ is a positive basis for $T_yW^u_{F_{21}}(q_1) = \mathbb{R} \times
T_x W^u_{f_1}(q_1)$. If $(w_1',w_2',\ldots, w_k')$ is a positive basis for $T_{q_2}W^u_{f_2}(q_2)$,
then $(v_2, w_1',w_2',\ldots, w_k')$ is a positive basis for $T_{q_2}W^u_H(q_2)$, where we have
identified $v_2 \in T_yW^u_H(q_1)$ with its image under the map induced by the gradient flow
in the $s$ direction. The gradient flow in the $s$ direction transports the basis 
$(w_1',w_2',\ldots, w_k')$ along a gradient flow line $\gamma \in W_{F_{21}}(q_1,q_2)$ to a 
basis for $\nu_y(W_{F_{21}}(q_1,q_2), W^u_{F_{21}}(q_1))$, 
which we still denote by $(w_1',w_2',\ldots, w_k')$. With this same identification, 
$(v_2,w_1',w_2',\ldots , w_k')$ is then a basis for $\nu_y(W_H(q_1,q_2), W^u_H(q_1))$.
Thus, $\gamma \in W_{F_{21}}(q_1,q_2)$ is oriented in the $v_1$ direction if 
$(w_1,w_2,\ldots, w_k)$ and $(w_1',w_2',\ldots , w_k')$ determine the same orientation on
$T_x W^u_{f_1}(q_1)$, since $(v_1, w_1',w_2',\ldots , w_k')$ will be a positive basis for
$T_yW^u_{F_{21}}(q_1)$, and $\gamma$ is oriented in the $-v_1$ direction if these orientations
are opposite, because $(-v_1, w_1',w_2',\ldots , w_k')$ would then be a positive basis for
$T_yW^u_{F_{21}}(q_1)$. Since the same is true for $\gamma$ viewed as an element of $W_H(q_1,q_2)$,
we have shown that
$$
W_H(q_1,q_2) = W_{F_{21}}(q_1,q_2)
$$
as oriented manifolds.

Finally, let $y = (-\infty,t,x)\in W^u_H(q_1)$ be in the region where $H(s,t,x) = f_1(x) - 
\rho(s) - \rho(t)$, and let $(v_1,v_2,w_1, w_2, \ldots, w_k)$ be a positive basis for 
$T_y W^u_H(q_1) = \mathbb{R} \times \mathbb{R} \times T_x W^u_{f_1}(q_1)$. 
Then $(v_2,w_1,w_2,\ldots, w_k)$ is a positive basis for $T_yW^u_{F_{31}}(q_1) = \mathbb{R} \times
T_x W^u_{f_1}(q_1)$. If $(w_1',w_2',\ldots, w_k')$ is a positive basis for $T_{q_3}W^u_{f_3}(q_3)$,
then $(v_1, w_1',w_2',\ldots, w_k')$ is a positive basis for $T_{q_3}W^u_H(q_3)$, where we have 
identified $v_1 \in T_yW^u_H(q_1)$ with its image under the map induced by the gradient flow in 
the $t$ direction. The gradient flow in the $t$ direction transports the basis 
$(w_1',w_2',\ldots, w_k')$ along a gradient flow line $\gamma \in W_{F_{31}}(q_1,q_3)$ to a 
basis for $\nu_y(W_{F_{31}}(q_1,q_3), W^u_{F_{31}}(q_1))$, which we still denote by 
$(w_1',w_2',\ldots, w_k')$. With this same identification, 
$(v_1,w_1',w_2',\ldots , w_k')$ is then a basis for $\nu_y(W_H(q_1,q_3), W^u_H(q_1))$.
Thus, $\gamma \in W_{F_{31}}(q_1,q_3)$ is oriented in the $v_2$ direction if 
$(w_1,w_2,\ldots, w_k)$ and $(w_1',w_2',\ldots , w_k')$ determine the same orientation on
$T_x W^u_{f_1}(q_1)$, and $\gamma$ is oriented in the $-v_2$ direction if these orientations
are opposite.  On the other hand, $\gamma$ viewed as an element of $W_H(q_1,q_3)$ is oriented
in the $-v_2$ direction if $(w_1,w_2,\ldots, w_k)$ and $(w_1',w_2',\ldots , w_k')$ determine 
the same orientation on $T_x W^u_{f_1}(q_1)$, since $(-v_2,v_1, w_1',w_2',\ldots , w_k')$ will be
a positive basis for $T_yW^u_H(q_1)$, and $\gamma$ is oriented in the $v_2$ direction
if these orientations are opposite, because $(v_2,v_1,w_1',w_2',\ldots , w_k')$ would then
be a positive basis for $T_yW^u_H(q_1)$. Hence,
$$
W_H(q_1,q_3) = - W_{F_{31}}(q_1,q_3)
$$
as oriented manifolds.

This shows that we have the following orientation relations between moduli spaces of broken
flow lines.
\begin{eqnarray*}
\mathcal{M}_H(q_1,q_2) \times \mathcal{M}_H(q_2,q_4) & = & + \mathcal{M}_{F_{21}}(q_1,q_2) \times 
\mathcal{M}_{F_{43}}(q_2,q_4)\\
\mathcal{M}_H(q_1,q_3) \times \mathcal{M}_H(q_3,q_4) & = & - \mathcal{M}_{F_{31}}(q_1,q_3) \times 
\mathcal{M}_{F_{43}}(q_3,q_4) 
\end{eqnarray*}
To complete the proof, for $(i,j)=(1,2)$, $(1,3)$, $(2,4)$, $(3,4)$ let $(H_{ji})_\Box$ be the 
chain map defined using the same sum that defines $(F_{ji})_\Box$ in Definition \ref{chainmapdef},
but with the orientations from the moduli spaces $\mathcal{M}_{F_{ji}}(q_i,q_j)$ replaced with
those of $\mathcal{M}_H(q_i,q_j)$. Since $\overline{\mathcal{M}}_H(q_1,q_4)$ is an oriented compact
smooth one dimensional manifold with a coherent orientation on its boundary, summing over the moduli
spaces of broken flow lines that make up the boundary of $\overline{\mathcal{M}}_H(q_1,q_4)$ gives
$$
(H_{43})_\Box \circ (H_{31})_\Box + (H_{42})_\Box \circ (H_{21})_\Box +
\partial_{k+1}^{G_4} H_\Box + H_\Box \partial_k^{G_1} = 0.
$$  
The above relations concerning the orientations of the moduli spaces of $H$ versus those
of $F_{ji}$ for $(i,j)=(1,2)$, $(1,3)$, $(2,4)$, $(3,4)$ imply that
$(H_{42})_\Box \circ (H_{21})_\Box = (F_{42})_\Box \circ (F_{21})_\Box$ and 
$(H_{43})_\Box \circ (H_{31})_\Box = -(F_{43})_\Box \circ (F_{31})_\Box$.  Therefore,
$$
(F_{43})_\Box \circ (F_{31})_\Box - (F_{42})_\Box \circ (F_{21})_\Box
= \partial_{k+1}^{G_4} H_\Box + H_\Box \partial_k^{G_1}.
$$
\proofend


\subsection{An invariance theorem}
We now prove the main theorem in this section.

\begin{lemma}\label{identity}
Let $f_1:M \rightarrow \mathbb{R}$ be a Morse-Smale function on $(M,\mathsf{g}_1)$
and $\Phi: G_1 \rightarrow G_2$ a family of isomorphisms between two bundles of abelian
groups $G_1$ and $G_2$ over $M$. For all $k=0,\ldots ,m$ the chain map 
$$
(F_{21})_\Box:C_k(f_1;G_1) \rightarrow C_k(f_2;G_2)
$$
from Definition \ref{chainmapdef} is an isomorphism if we take $f_1=f_2$, $\mathsf{g}$ 
the product metric $\mathsf{g}_1 + dt^2$ on $\overline{\mathbb{R}} \times M$, and 
$F_{21}(t,x) = f_1(x) - \rho(t)$.  Moreover, on an elementary chain
$gq \in C_k(f_1;G_1)$ we have $(F_{21})_\Box(gq) = (-1)^k \Phi(g)q$.
Thus, $(F_{12})_\Box \circ (F_{21})_\Box = id$ if $(F_{12})_\Box$ is defined
using the product metric, $F_{12}(t,x) = f_2(x) - \rho(t)$, and
$\Phi^{-1}:G_2 \rightarrow G_1$. 
\end{lemma}

\proofstart
Note that $-\nabla F_{21} = (\rho'(t),-\nabla f_1)$ where $\rho'(t) > 0$ for all $t$ since
$\mathsf{g}$ is the product metric. Recalling that $f_1=f_2$ decreases along its gradient flow lines,
we see that for any elementary chain $gq \in C_k(f_1;G_1)$ we have
$$
\begin{array}{lcl}
\displaystyle (F_{21})_\Box(gq)\ \ 
& = & \displaystyle (-1)^{\lambda_{q}} \sum_{q_2 \in Cr_k(f)} 
      \sum_{\nu_F \in \mathcal{M}_F(q,q_2)} \epsilon(\nu_F) (\gamma_F)_\ast(\Phi(g))q_2\\
& = & \displaystyle (-1)^k \epsilon(\nu_F) (\gamma_F)_\ast(\Phi(g))q\\
& = & (-1)^k \Phi(g)q,
\end{array}
$$
where $(\gamma_F)_\ast = id$ since $\pi \circ \gamma_F \equiv q$ and 
$\epsilon(\nu_F) = +1$ because for $\tau << 0$ the tangent bundle of
$W^u_F(q)$ is oriented as $\mathbb{R} \times T_\ast W^u_{f_1}(q)$ and $W_F(q,q)$ is oriented
by the relation
$$
T_\ast W_F(q,q) \oplus T_q W^u_{f_2}(q) \approx T_\ast W^u_F(q)|_{W_F(q,q)}.
$$
\proofend

\begin{corollary}\label{canonical}
For any two Morse-Smale pairs $(f_1,\mathsf{g}_1)$ and $(f_2,\mathsf{g}_2)$ and a family of
isomorphisms $\Phi:G_1 \rightarrow G_2$ between bundles of abelian groups $G_1$ and $G_2$
over $M$, the time dependent gradient flow lines from $f_1$ to $f_2$ determine a canonical
homomorphism 
$$
(F_{21})_\ast:H_\ast((C_\ast(f_1;G_1),\partial_\ast^{G_1})) \rightarrow 
H_\ast((C_\ast(f_2;G_2), \partial_\ast^{G_2})), 
$$
i.e. the map $(F_{21})_\ast$ is independent of the choices made in the definition of 
$(F_{21})_\Box$.
\end{corollary}

\smallskip\noindent
Proof: Let $f_2 = f_3 = f_4$, $\mathsf{g}_2=\mathsf{g}_3=\mathsf{g}_4$, and $G_2=G_3=G_4$ in 
Theorem \ref{chainhomotopy}, and let $F_{21}:\mathbb{R} \times M \rightarrow \mathbb{R}$ 
and $F_{31}:\mathbb{R} \times M\rightarrow \mathbb{R}$ be two functions that define
time dependent gradient flow lines from $f_1$ to $f_2 = f_3 = f_4$ with respect to
metrics $\mathsf{g}_1$ and $\mathsf{g}_2$ on $\overline{\mathbb{R}} \times M$ that are equal to
$\mathsf{g}_1 + dt^2$ and $\mathsf{g}_2 + dt^2$ near the respective ends of 
$\overline{\mathbb{R}} \times M$. Theorem \ref{chainhomotopy} implies
$$
(F_{43})_\ast \circ (F_{31})_\ast = (F_{42})_\ast \circ (F_{21})_\ast
$$
and Lemma \ref{identity} implies that $(F_{42})^{-1}_\ast = (F_{24})_\ast$
satisfies $(F_{42})^{-1}_\ast \circ (F_{43})_\ast = id$, if we define 
$(F_{24})_\Box = (F_{34})_\Box$ using the choices in Lemma \ref{identity}.
Thus, $(F_{31})_\ast = (F_{21})_\ast$.
\proofend

%
%

\begin{theorem}[Invariance Theorem]\label{homologyindependence}
Let $(M,\mathsf{g})$ be a closed finite dimensional smooth Riemannian manifold, 
and let $G$ be a bundle of abelian groups over $M$. Then the homology of the twisted
Morse-Smale-Witten chain complex $(C_\ast(f;G), \partial_\ast^{G})$ is independent
of the Morse-Smale pair $(f,\mathsf{g})$ and depends only on the isomorphism class of the
bundle of abelian groups $G$.
\end{theorem}

\proofstart
Let $(f_1, \mathsf{g}_1)$ and $(f_2, \mathsf{g}_2)$ be Morse-Smale pairs on $M$, and let 
$\Phi:G_1 \rightarrow G_2$ be a family of isomorphisms between bundles of abelian groups 
$G_1$ and $G_2$ over $M$. Letting $f_1=f_3=f_4$, $\mathsf{g}_1=\mathsf{g}_3=\mathsf{g}_4$, and 
$G_1 = G_3 = G_4$ in Theorem \ref{chainhomotopy} we have
$$
(F_{43})_\ast \circ (F_{31})_\ast = (F_{42})_\ast \circ (F_{21})_\ast
$$
where $(F_{43})_\ast \circ (F_{31})_\ast = id$ by Lemma \ref{identity} and
Corollary \ref{canonical}. Therefore, $(F_{12})_\ast \circ (F_{21})_\ast = id$.
Similarly, $(F_{21})_\ast \circ (F_{12})_\ast = id$.
\proofend

\noindent
Combining the preceding theorem with Claim \ref{etaisomorphic} we have the following.

\begin{corollary}\label{etadeRhamclass}
Let $\eta \in \Omega^1_{\text{cl}}(M,\mathbb{R})$ be a closed one form on a Riemannian
manifold $(M,\mathsf{g})$. Then the homology of the $\eta$-twisted Morse-Smale-Witten chain
complex $(C_\ast(f) \otimes \mathbb{R},\partial^\eta_\ast)$ from Definition \ref{etatwisted}
is independent of the Morse-Smale pair $(f,\mathsf{g})$ and depends only on the de Rham 
cohomology class of $\eta$.
\end{corollary}


\section{Singular and CW-Homology with local coefficients}\label{singularCW}

The main goal of this section is to prove the following theorem.

\begin{theorem}[Twisted Morse Homology Theorem]\label{twistedsame}
Let $f:M \rightarrow \mathbb{R}$ be a smooth Morse-Smale function on a closed
finite dimensional smooth Riemannian manifold $(M,\mathsf{g})$, and let $G$ be a bundle of
abelian groups over $M$. The homology of the Morse-Smale-Witten chain complex
with coefficients in $G$ is isomorphic to the singular homology of $M$ with coefficients in $G$, i.e.
$$
H_k((C_\ast(f;G), \partial_\ast^G)) \approx H_k(M;G)
$$
for all $k=0,\ldots ,m$.
\end{theorem}

We will prove this theorem by observing that if the unstable manifolds of the Morse-Smale
function $f:M \rightarrow \mathbb{R}$ determine a \textbf{regular} CW-structure on
$(M,\mathsf{g})$, then the Morse-Smale-Witten chain complex $(C_\ast(f;G), \partial_\ast^G)$
and Steenrod's CW-chain complex with coefficients in the bundle of abelian groups $G$ coincide
(Lemma \ref{boundarysame}). We will then show how to construct a Morse-Smale pair $(f,\mathsf{g})$
on $M$ whose unstable manifolds determine a regular CW-structure on $M$ 
(Theorem \ref{regularMorse}) and apply the above Invariance Theorem (Theorem
\ref{homologyindependence}). This proves the Twisted Morse Homology Theorem
because, for a regular CW-complex, the homology of Steenrod's CW-chain complex with
coefficients in $G$ is isomorphic to the singular homology with coefficients in $G$
(Lemma \ref{CWsingular}).


\subsection{Singular homology with local coefficients}\label{singularlocal}
For the convenience of the reader we now recall the definition of singular homology with
coefficients in a bundle of abelian groups.  For more details see Chapter VI of \cite{WhiEle}.
\nocite{FarTop}

\smallskip
Let $G$ be a bundle of abelian groups over a topological space $X$. Let $\Delta^k$ denote
the standard $k$-simplex with vertices $e_0,\ldots ,e_k$, and let $C_k(X;G)$ be the set of
all functions $c$ such that the following conditions hold.
\begin{enumerate}
\item For every singular $k$-simplex $u:\Delta^k \rightarrow X$, $c(u) \in G_{u(e_0)}$
      is defined.
\item The set of singular simplices $u$ such that $c(u) \neq 0$ is finite.
\end{enumerate}

\noindent
Elements of the abelian group $C_k(X;G)$ are called {\bf singular $k$-chains with
coefficients in $G$}, and every $ c \in C_k(X;G)$ can be represented as a finite sum
$$
c = \sum_{i=1}^n c(u_i)\cdot u_i
$$
where $u_1,\ldots , u_n$ are the singular simplices such that $c(u_i) \neq 0$ and
$c(u_i) \in G_{u_i(e_0)}$ for all $i=1,\ldots , n$. The \textbf{relative singular $k$-chains
with coefficients in $G$} for a subspace $A \subseteq X$, denoted $C_k(X,A;G)$, are
defined similarly.

\begin{definition}\label{twistedsingular}
The \textbf{singular boundary operator with coefficients in $G$} is defined to be the 
homomorphism $\partial_k:C_k(X;G) \rightarrow C_{k-1}(X;G)$ given on an elementary chain 
$c = g\cdot u$ by
$$
\partial_k(g\cdot u) = (\gamma_u)_\ast (g) \cdot u \circ F_0 + \sum_{i=1}^k (-1)^i g \cdot u \circ F_i
$$
where $(\gamma_u)_\ast:G_{u(e_0)} \rightarrow G_{u(e_1)}$ is the homomorphism
associated to the path $\gamma_u(t) = u((1-t) e_1 + t e_0)$ from $u(e_1)$ to $u(e_0)$
and $F_i:\Delta^{k-1} \hookrightarrow \Delta^k$ is the inclusion onto the face
opposite $e_i$ for all $i=0,\ldots , k-1$.  The pair $(C_\ast(X;G),\partial_\ast)$ is
a chain complex, and its homology groups $H_\ast(X;G)$ are called the {\bf homology
groups of $X$ with coefficients in the bundle $G$}.  The relative homology groups
$H_\ast(X,A;G)$ for a subspace $A \subseteq X$ are defined similarly.
\end{definition}

Eilenberg showed that singular homology with local coefficients is closely related to
equivariant homology. Suppose that $(X,x_0)$ is a connected topological space and $G_0$
is an abelian group on which $\pi_1(X,x_0)$ acts.  Since $\pi_1(X,x_0)$ also acts on the universal cover
$\widetilde{X}$ and this action commutes with the boundary operator on singular chains in
$\widetilde{X}$, there is a chain complex $(G_0 \otimes_{\pi_1} C_\ast(\widetilde{X}), \bar{\partial}_\ast)$,
where the tensor product is taken over $\pi_1(X,x_0)$ and the boundary operator $\bar{\partial}_\ast$
is induced from the boundary operator on the singular chains in $\widetilde{X}$.  The homology
groups of this complex are the equivariant homology groups $E_\ast(\widetilde{X};G_0)$, cf.
Section \ref{EMTheorem}.  

Let $G$ be a bundle of abelian groups on $X$ in the isomorphism class determined by the action of 
$\pi_1(X,x_0)$ on $G_0$ (Theorem \ref{actionsystem}). We have the following.

\begin{theorem}[Eilenberg]\label{Eilenberg}
For all $k$, $H_k(X;G)$ is isomorphic to $E_k(\widetilde{X};G_0)$.
\end{theorem}

\noindent
For a proof of the preceding theorem see Section VI.3 of \cite{WhiEle}.


\subsection{Regular CW-complexes}
We now recall the definition of a regular CW-complex and we prove that, for a regular
CW-complex, Steenrod's CW-boundary operator with coefficients in a bundle of abelian
groups $G$ (Definition \ref{twistedregular}) is induced from the singular boundary 
operator with coefficients in $G$ (Lemma \ref{CWsingular}). The reader should note that
Steenrod's CW-boundary operator may not be well-defined for a general CW-complex.
So, Steenrod's CW-chain complex with local coefficients is only defined for some
restricted class of CW-complexes; although a CW-chain complex with local coefficients
for a general CW-complex can be defined using a boundary operator that is induced from
a connecting homomorphism.

\begin{definition}\label{regularCW}
A CW-complex $X$ is \textbf{regular} if every closed $k$-cell $e^k$, with $k > 0$, is 
homeomorphic to $\Delta^k$.
\end{definition}

\begin{remark}
Not every CW-complex is regular. For instance, the CW-structure on $S^1$
given by the Morse function in Example \ref{circle} is not regular. In that example the
closed 1-cell $e^1$ is homeomorphic to $S^1$ instead of $\Delta^1$.
Similarly, the CW-structure on the 2-torus with 4 cells is not regular,
cf. Example 7.11 of \cite{BanLec}.
\end{remark}

\smallskip
Regular cell complexes satisfy the following properties, which are not necessarily 
satisfied by nonregular CW-complexes. For more details see Section IX.6 of \cite{MasABa}
or Section II.6 of \cite{WhiEle}.

\begin{enumerate}
\item If $e^k$ is a $k$-cell, then its \textbf{boundary} $\dot{e}^k = e^k - 
\text{int}(e^k)$ is the union of finitely many $(k-1)$-cells. (Note that the ``boundary'' in
this context may be different than the topological boundary of $e^k$.)
\item If $j < k$ and $e^j$ and $e^k$ are cells such that $e^j \cap \dot{e}^k \neq \emptyset$,
then $e^j \subset \dot{e}^k$.
\item For any $k$-cell $e^k$ of $X$ with $k \geq 0$, $e^k$ and $\dot{e}^k$ are 
the underlying spaces of sub-complexes of $X$.
\item If $e^k$ and $e^{k+2}$ are cells such that $e^k$ is a face of $e^{k+2}$, then
there are exactly two $(k+1)$-cells $e^{k+1}$ such that $e^k$ is a proper face of $e^{k+1}$
and $e^{k+1}$ is a proper face of $e^{k+2}$, i.e. $e^k < e^{k+1} < e^{k+2}$.
\item The incidence number $[e^k:e^{k-1}]$ is $\pm 1$ if $e^{k-1} < e^k$ and zero otherwise. 
\end{enumerate}

\smallskip
For any system of local coefficients $G$ over a CW-complex $X$ (not necessarily regular)
the triple $(X^{(k-2)},X^{(k-1)},X^{(k)})$, where $X^{(k)}$ denotes the $k$-skeleton of $X$,
determines a connecting homomorphism 
$$
\xymatrix{
H_k(X^{(k)},X^{(k-1)};G) \ar[r]^-{\delta_k} & H_{k-1}(X^{(k-1)};G)
}
$$
that can be composed with the map $H_{k-1}(X^{(k-1)};G) \stackrel{j_\ast}{\rightarrow} 
H_{k-1}(X^{(k-1)},X^{(k-2)};G)$ induced from the inclusion $j:X^{(k-2)} \hookrightarrow X^{(k-1)}$
to give a map
$$
\xymatrix{
H_k(X^{(k)},X^{(k-1)};G) \ar[r]^-{\tilde{\partial}_k} & H_{k-1}(X^{(k-1)},X^{(k-2)};G).
}
$$
The above map satisfies $\tilde{\partial}_{k-1}\circ \tilde{\partial}_k = 0$, and the
homology groups of the chain complex with boundary operator $\tilde{\partial}_k$ and 
$k^{\text{th}}$-chain group $H_k(X^{(k)},X^{(k-1)};G)$ are isomorphic to the singular
homology groups of $X$ with coefficients in the bundle $G$ from Definition
\ref{twistedsingular}, cf. Theorems VI.2.4 and VI.4.4 of \cite{WhiEle}.

\smallskip
Now if $X$ is regular, then the homology group $H_k(X^{(k)},X^{(k-1)};G)$ can be
represented as a direct sum of the images of the maps induced by the characteristic maps
$f_\sigma:(\Delta^k,\dot{\Delta}^k) \rightarrow (e^k_\sigma,\dot{e}^k_\sigma) \subseteq
(X^{(k)},X^{(k-1)})$ of the $k$-cells $e^k_\sigma$ of $X$ as follows. For each $k$-cell 
$e^k_\sigma$ we choose a basepoint $x(e^k_\sigma)$, which determines an isomorphism
$$
\bigoplus_\sigma (f_\sigma)_\ast: \bigoplus_\sigma H_k(\Delta^k,\dot{\Delta}^k;G_{x(e^k_\sigma)})
\stackrel{\approx}{\longrightarrow} H_k(X^{(k)},X^{(k-1)};G),
$$
cf. Theorem VI.4.1 of \cite{WhiEle}. Note that the definition of the induced map
$(f_\sigma)_\ast$ requires both a map of spaces $f_\sigma:(\Delta^k,\dot{\Delta}^k)
\rightarrow (X^{(k)},X^{(k-1})$ and a homomorphism $\gamma_\ast:G_{x(e^k_\sigma)}
\rightarrow f_\sigma^\ast(G)$, cf. Section VI.2 of \cite{WhiEle}. We take the homomorphism
$\gamma_\ast$ to be the one defined by restricting the local coefficient system $G$ to
the simply connected space $e_\sigma^k$, i.e. for any point $x \in \Delta^k$ there is
a unique homotopy class of paths rel endpoints from $f_\sigma(x)$ to $x(e_\sigma^k)$ and
hence a well-defined homomorphism $G_{x(e^k_\sigma)} \rightarrow G_{f_\sigma(x)}$.
Since $G_{x(e^k_\sigma)} \approx H_k(\Delta^k,\dot{\Delta}^k;G_{x(e^k_\sigma)})$,
we can use the above isomorphisms to identify
\begin{eqnarray*}
CW_k(X;G) & \stackrel{\text{def}}{=} & \left.\left\{\sum_\sigma g e^k_\sigma\right|\ g 
            \in G_{x(e^k_\sigma)}\right\}\\
          & = & \bigoplus_\sigma G_{x(e^k_\sigma)}\\
          & = & \bigoplus_\sigma H_k(\Delta^k,\dot{\Delta}^k;G_{x(e^k_\sigma)})\\
          & \stackrel{\oplus_\sigma (f_\sigma)_\ast}{=} & H_k(X^{(k)},X^{(k-1)};G)
\end{eqnarray*}

\begin{definition}\label{twistedregular}
\textbf{Steenrod's CW-boundary operator with coefficients in $G$} is defined to be the
homomorphism $\partial_k:CW_k(X;G) \rightarrow CW_{k-1}(X;G)$ given on an elementary
chain $g e^k$ by 
$$
\partial_k(g e^k) = \sum_{e^{k-1} < e^k} [e^k:e^{k-1}] (\gamma_{e^{k-1} e^k})_\ast(g) e^{k-1},
$$
where $(\gamma_{e^{k-1} e^k})_\ast: G_{x(e^k)} \rightarrow G_{x(e^{k-1})}$ denotes the
isomorphism determined by any path from $x(e^{k-1})$ to $x(e^k)$ contained in the closure 
of $e^k$.
We will call the pair $(CW_\ast(X;G), \partial_\ast)$ \textbf{Steenrod's CW-chain complex
with coefficients in the bundle $G$}.
\end{definition}

\noindent
Note that $(\gamma_{e^{k-1} e^k})_\ast$ does not depend on the path from 
$x(e^{k-1})$ to $x(e^k)$ since $e^k \approx \Delta^k$ is simply connected.
Moreover, one can show directly that $\partial_{k-1} \circ \partial_k = 0$, and the
homology of Steenrod's CW-chain complex for a regular CW-complex $X$ with coefficients 
in $G$ is independent of the choice of basepoints $x(e^k)$; see \cite{CooHom} or 
Section 31.2 of \cite{SteThe} for more details.

\begin{remark}
If the CW-complex $X$ is not regular, then the homomorphism 
$(\gamma_{e^{k-1} e^k})_\ast$ may not be well-defined. For instance, the above formula for
$\partial_1(g e^1)$ is not well-defined for the CW-structure given by the unstable
manifolds of the height function on $S^1$ in Example \ref{circle}. Similarly, the above
formula is not well-defined for the CW-structure on the 2-torus with 4 cells, 
cf. Example 7.11 of \cite{BanLec}.
\end{remark}

\begin{lemma}\label{CWsingular}
If $X$ is a regular CW-complex and $G$ is a bundle of abelian groups over $X$, then the
singular boundary operator with coefficients in $G$ from Definition \ref{twistedsingular}
induces Steenrod's CW-boundary operator with coefficients in $G$ from 
Definition \ref{twistedregular}. That is, the following diagram commutes.
$$
\xymatrix{
CW_k(X;G) \ar[r]^-{\partial_k} \ar@{<->}[d] & CW_{k-1}(X;G)  \ar@{<->}[d]\\
H_k(X^{(k)},X^{(k-1)};G) \ar[r]^-{\tilde{\partial}_k} & H_{k-1}(X^{(k-1)},X^{(k-2)};G)
}
$$
Thus, the homology of Steenrod's CW-chain complex $(CW_\ast(X;G),\partial_\ast)$ is
isomorphic to the singular homology of $X$ with coefficients in the bundle $G$.
\end{lemma}

\proofstart
The connecting homomorphism $\delta_k$ is natural, cf. Theorem VI.2.4 of \cite{WhiEle}.
Hence, the following diagram commutes
$$
\xymatrix{
H_k(\Delta^k,\dot{\Delta}^k;G_{x(e_\sigma^k)}) \ar[r]^-{\bar{\delta}_k} \ar[d]^{(f_\sigma)_\ast}
  & H_{k-1}(\dot{\Delta}^k;G_{x(e_\sigma^k)}) \ar[d]^{(f_{\partial\sigma})_\ast}
  & \\
H_k(X^{(k)},X^{(k-1)};G) \ar[r]^-{\delta_k} & H_{k-1}(X^{(k-1)};G) \ar[r]^-{j_\ast} & 
  H_{k-1}(X^{(k-1)},X^{(k-2)};G)\\
& & \bigoplus_\tau H_{k-1}\left(\Delta^{k-1},\dot{\Delta}^{k-1};G_{x(e_\tau^{k-1})}\right)
    \ar[u]_{\oplus_\tau (f_\tau)_\ast}^\approx
}
$$
where $f_\sigma$ and $f_\tau$ are characteristic maps, $f_{\partial \sigma} = f_\sigma |_{\dot{\Delta}^k}$,
the homomorphisms $G_{x(e_\sigma^k)} \rightarrow f_\sigma^\ast G$ and $G_{x(e_\sigma^k)}
\rightarrow f_{\partial\sigma}^\ast G$ are defined by restricting the local coefficient system to the
simply connected space $e_\sigma^k$, the homomorphism $G_{x(e^{k-1}_\tau)} \rightarrow f_\tau^\ast G$ 
is defined by restricting the local system to the simply connected space $e^{k-1}_\tau$, and
$j:(X^{(k-1)},\emptyset) \rightarrow (X^{(k-1)},X^{(k-2)})$ is the inclusion.
Moreover, we have fixed isomorphisms 
$$
G_{x(e_\sigma^k)} = H_k(\Delta^k,\dot{\Delta}^k;G_{x(e_\sigma^k)}) \quad \mbox{and} \quad
\bigoplus_\tau H_{k-1}\left(\Delta^{k-1},\dot{\Delta}^{k-1};G_{x(e_\tau^{k-1})}\right) = 
CW_{k-1}(X;G).
$$
The statement of the lemma claims that these isomorphisms and the map given by tracing
the above diagram from the upper left to the lower right agree with $\partial_k$
from Definition \ref{twistedregular}.

If we let $X^{(k-1)}/X^{(k-2)}$ be the space obtained by identifying $X^{(k-2)}$ to a point
when $k > 1$ or $X^{(0)}$ union a disjoint basepoint $\ast$ when $k=1$, then 
$X^{(k-1)}/X^{(k-2)}$ is a bouquet of $(k-1)$-spheres $S^{k-1}_\tau$ unioned at the 
basepoint $\ast$. For each $\tau$ let $p_\tau$ be the composition $X^{(k-1)} \rightarrow
X^{(k-1)}/X^{(k-2)} \rightarrow S^{k-1}_\tau$, where the second map collapses every sphere
other than $S^{k-1}_\tau$ to the basepoint. Then $p_\tau:(X^{(k-1)},X^{(k-2)}) \rightarrow
(S^{k-1}_\tau,\ast) $ is a map of pairs, and $p_\tau \circ f_{\tau'}: \Delta^{k-1} \rightarrow
S^{k-1}_\tau$ is a constant map to the basepoint if $\tau \neq \tau'$ and the map that
identifies the boundary of $\Delta^{k-1}$ to the basepoint when $\tau = \tau'$. Identifying 
$H_{k-1}(S_\tau^{k-1},\ast; G_{x(e^{k-1}_\tau)}) = H_{k-1}(\Delta^{k-1}, 
\dot{\Delta}^{k-1};G_{x(e^{k-1}_\tau)})$ we see that
$$
\bigoplus_\tau (p_\tau)_\ast: H_{k-1}(X^{(k-1)},X^{(k-2)};G) \stackrel{\approx}{\longrightarrow} 
\bigoplus_\tau H_{k-1}(\Delta^{k-1},\dot{\Delta}^{k-1};G_{x(e^{k-1}_\tau)}) 
$$
is the inverse of $\bigoplus_\tau (f_\tau)_\ast$, cf. Proposition 2.14 of \cite{BanLec}.

Now, the characteristic maps $f_\sigma:\Delta^k \rightarrow e^k_\sigma$ and
$f_\tau:\Delta^{k-1} \rightarrow e^{k-1}_\tau$ determine orientations on the cells
$e^k_\sigma$ and $e^{k-1}_\tau$, and the sign $[e^k_\sigma:e^{k-1}_\tau] = \pm 1$ keeps
track of the compatibility of these orientations.  That is, if $e^{k-1}_\tau < e^k_\sigma$,
then $[e^k_\sigma:e^{k-1}_\tau]$ is $+1$ or $-1$ depending on whether the following
isomorphism preserves or reverses the orientations determined by $f_\sigma$ and
$f_\tau$ 
$$
\xymatrix{
\mathbb{Z} \approx H_k(e^k_\sigma,\dot{e}^k_\sigma) \ar[r]^-{\bar{\delta}_k} & H_{k-1}(\dot{e}^k_\sigma) 
\ar[r]^-{j_\ast} & H_{k-1}(\dot{e}^k_\sigma, \overline{\dot{e}^k_\sigma - e^{k-1}_\tau}) 
\ar[r]^-{\approx} & H_{k-1}(e^{k-1}_\tau, \dot{e}^{k-1}_\tau) \approx \mathbb{Z},
}
$$
where the last isomorphism is given by excision, cf. Section II.6 of \cite{WhiEle}.
Moreover, the exact sequence
$$
\xymatrix{
H_k(\Delta^k;G_{x(e_\sigma^k)}) \ar[r] & H_k(\Delta^k,\dot{\Delta}^k;G_{x(e_\sigma^k)})
\ar[r]^-{\bar{\delta}_k} & H_{k-1}(\dot{\Delta}^k;G_{x(e_\sigma^k)}) \ar[r] & 
H_{k-1}(\Delta^k;G_{x(e_\sigma^k)})
}
$$
shows that $\bar{\delta}_k$ is an isomorphism when $k > 1$ and injective
when $k=1$, where the orientation of $\dot{\Delta}^k$ is determined by the
orientation of $\Delta^k$.  Thus, by the Universal Coefficient Theorem
we have
\begin{eqnarray*}
H_k(\Delta^k,\dot{\Delta}^k;G_{x(e^k_\sigma)}) & \approx & \mathbb{Z} \otimes_\mathbb{Z}
G_{x(e^k_\sigma)} \approx G_{x(e^k_\sigma)} \quad \text{ for all } k \geq 0, \\
H_{k-1}(\dot{\Delta}^k;G_{x(e^k_\sigma)}) & \approx &  
\left\{
\begin{array}{ll}
G_{x(e^k_\sigma)}                          & \text{ if } k > 1\\
G_{x(e^k_\sigma)} \oplus G_{x(e^k_\sigma)} & \text{ if } k = 1,
\end{array}
\right.
\end{eqnarray*}
and
$$
\bar{\delta}_k(g) = 
\left\{
\begin{array}{ll}
g      & \text{ if } k \geq 1\\
(g,-g) & \text{ if } k = 1,
\end{array}
\right.
$$
cf. Example 2.2 of \cite{BanLec}.  Since $p_\tau \circ j \circ f_{\partial \sigma}:
\dot{\Delta}^{k-1} \rightarrow S^{k-1}_\tau$ is a map of degree $[e^k_\sigma:e^{k-1}_\tau]$,
this shows that for any $g \in G_{x(e^k_\sigma)}$ we have
\begin{eqnarray*}
(p_\tau)_\ast(j_\ast(\delta_k((f_\sigma)_\ast(g)))) 
  & = & (p_\tau)_\ast(j_\ast((f_{\partial\sigma})_\ast\bar{\delta}_k((g))))\\
  & = & [e^k_\sigma:e^{k-1}_\tau] (\gamma_{e^k_\sigma e^{k-1}_\tau})_\ast(g) \in 
        G_{x(e^{k-1}_\tau)}.
\end{eqnarray*}
\proofend


\subsection{Unstable manifolds and regular CW-structures}\label{unstableCW}

The unstable manifolds of the Morse-Smale functions in Examples \ref{circle} and
\ref{projectivespace} do not determine regular CW-structures on $S^1$ or $\mathbb{R}P^2$.
This can be seen directly from Definition \ref{regularCW} or by noting that, in both cases,
there are cells $e^k$ and $e^{k-1}$ with $e^{k-1} < e^k$ where the incidence number
$[e^k:e^{k-1}] \neq  \pm 1$. However, we will show in Theorem \ref{regularMorse} that 
on any closed finite dimensional smooth manifold $M$ it is always possible to find a Morse
function $f:M \rightarrow \mathbb{R}$ and a Riemannian metric $\mathsf{g}$ such that $f$ satisfies
the Morse-Smale transversality condition with respect to $\mathsf{g}$ and the unstable manifolds of
$(f,\mathsf{g})$ determine a regular CW-structure on $M$.

\begin{example}[A deformed circle]\label{deformedcircle}
Consider the unit circle $M = S^1$ and the Morse-Smale function $f:S^1 \rightarrow \mathbb{R}$
given by the height function in the following picture, where the arrows indicate the
orientations of the unstable manifolds.  It's clear that the unstable manifolds of
$f$ determine a regular CW-structure on $S^1$ since $\overline{W^u(q_1)} = W^u(q_1) \cup
\{p_1,p_2\} \approx \Delta^1$ and $\overline{W^u(q_2)} = W^u(q_2) \cup \{p_1,p_2\}
\approx \Delta^1$. 
\begin{figure}[h]
\includegraphics{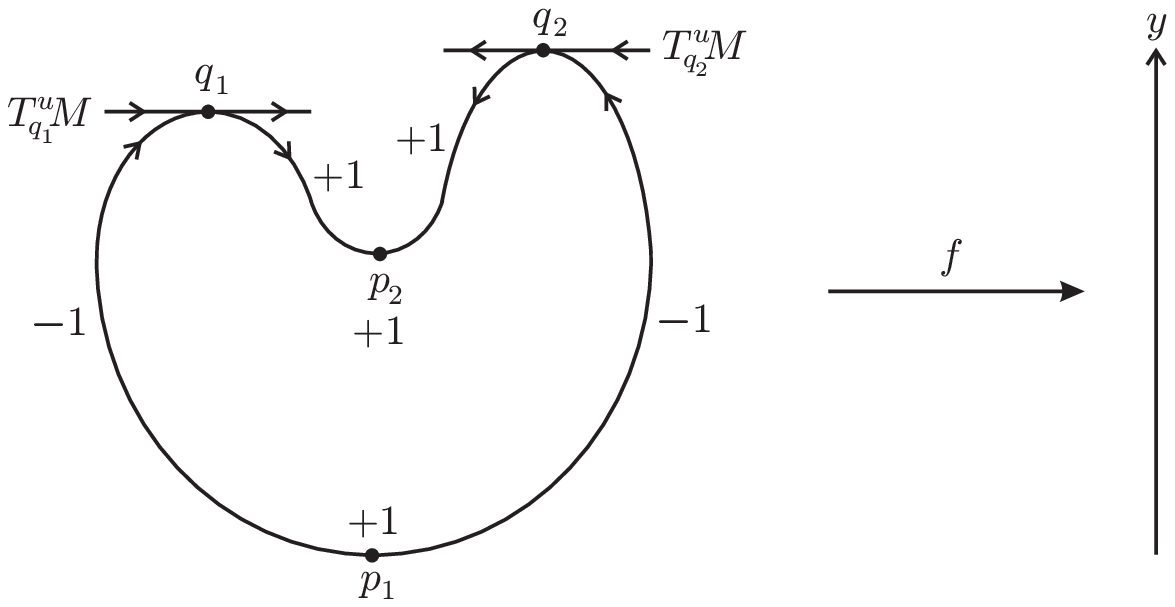}
\end{figure}

The (non-twisted) Morse-Smale-Witten complex of $f$ is 
$$
\xymatrix{
0 \ar[r] & C_1(f) \ar[r]^-{\partial_1}\ar@{<->}[d]^{\approx} & C_0(f) \ar@{<->}[d]^{\approx}
\ar[r] & 0\\
0 \ar[r] & <q_1,q_2> \ar[r]^-{\partial_1} & <p_1,p_2> \ar[r] & 0
}
$$
with $\partial_1(q_1) = \partial_1(q_2) = p_2-p_1$.   
The homology group $H_1(C_\ast(f),\partial_\ast) \approx \mathbb{Z}$ is generated by
the class $[q_1-q_2]$, and the image of $\partial_1$ is the group $<p_2-p_1>$.
Hence, $H_0(C_\ast(f),\partial_\ast) \approx <p_1,p_2>/<p_2-p_1> \approx \mathbb{Z}$.

\smallskip
Now, let $G$ be a bundle of abelian groups over $S^1$, and let $\gamma_{p_iq_j}:[0,1]
\rightarrow M$ be a path from $p_i$ to $q_j$ whose image coincides with the image of a 
gradient flow line $\nu_{q_jp_i}$, for $i,j=1,2$. The twisted Morse-Smale-Witten
boundary operator is determined by
\begin{eqnarray*}
\partial^G_1(g_1q_1) & = & (\gamma_{p_2q_1})_\ast (g_1) p_2 - (\gamma_{p_1q_1})_\ast (g_1) p_1
\quad \text{ for all } g_1 \in G_{q_1}\\
\partial^G_1(g_2q_2) & = & (\gamma_{p_2q_2})_\ast (g_2) p_2 - (\gamma_{p_1q_2})_\ast (g_2) p_1
\quad \text{ for all } g_2 \in G_{q_2}.
\end{eqnarray*}
Note that in this example it is not immediately obvious that the homology of the twisted
chain complex $(C_\ast(f;G);\partial^G_\ast)$ depends only on the isomorphism class of $G$
(Theorem \ref{homologyindependence}).  On the other hand, since the unstable manifolds of
the Morse-Smale function determine a regular CW-structure $X^{(0)} \subseteq X^{(1)}$ 
on $S^1$, we can compare the twisted Morse-Small-Witten chain complex 
$(C_\ast(f;G);\partial^G_\ast)$ to Steenrod's twisted CW-chain complex
$(CW_\ast(X;G), \partial_\ast)$.

If we choose the critical point $q_j = x(e^1_j)$ as the basepoint for the cell $e^1_j =
\overline{W^u(q_j)}$ for $j=1,2$, then the homomorphism $(\gamma_{e^0_ie^1_j})_\ast:G_{q_j}
\rightarrow G_{p_i}$ from Definition \ref{twistedregular} coincides with
$(\gamma_{p_iq_j})_\ast$ for $i,j=1,2$.  Moreover, one can check directly that the signs
$\epsilon(\nu_{q_jp_i})$ associated to the gradient flow lines agree with the corresponding
incidence numbers $[e^1_j:e^0_i]$.  Hence, Steenrod's CW-boundary operator with coefficients
in $G$ (Definition \ref{twistedregular}) is the same as the Morse-Smale-Witten boundary
operator with coefficients in $G$ (Definition \ref{twistedboundary}) once we identify the
corresponding generators. Therefore, $H_k((C_\ast(f;G), \partial^G_\ast)) \approx 
H_\ast((CW_\ast(X;G), \partial_\ast)) \approx H_k(S^1;G)$ for all $k$.
\end{example}

\begin{lemma}\label{boundarysame}
Let $f:M \rightarrow \mathbb{R}$ be a Morse-Smale function on a closed finite dimensional
smooth Riemannian manifold $(M,\mathsf{g})$. Assume that the unstable manifolds of $(f,\mathsf{g})$ 
determine a regular CW-structure $X$ on $M$, and for each closed $k$-cell $e^k_q$
with $W^u(q) \subset e^k_q$ choose the basepoint $x(e^k_q)$ to be the critical point
$q\in Cr_k(f)$. Let $G$ be a bundle of abelian groups over $M$, and identify the generator 
$e^k_q$ in the CW-chain complex $(CW_\ast(X;G),\partial_\ast)$ with the generator $q$ 
in the Morse-Smale-Witten chain complex $(C_\ast(f;G),\partial_\ast^G)$.

Under these identifications, Steenrod's CW-boundary operator with coefficients in $G$
from  Definition \ref{twistedregular} coincides with the twisted Morse-Smale-Witten boundary
operator from Definition \ref{twistedboundary}.  Thus,
$H_k((CW_\ast(X;G),\partial_\ast)) \approx H_k((C_\ast(f;G),\partial_\ast^G))$ for all
$k = 0,\ldots ,m$.
\end{lemma}

\proofstart
The assumption that the CW-complex is regular implies that any two gradient
flow lines from $q \in Cr_k(f)$ to $p \in Cr_{k-1}(f)$ are homotopic rel endpoints. 
Hence, if $\gamma^{\nu}$ is any continuous path from $p$ to $q$ whose images coincide with
$\nu \in \mathcal{M}(q,p)$, then the isomorphism $(\gamma^{\nu})_\ast:G_q \rightarrow G_p$
is independent of $\nu$.  That is, in the notation of Definition \ref{twistedregular}, 
we have $(\gamma^{\nu})_\ast = (\gamma_{e^{k-1}_pe^k_q})_\ast$.
Thus, for the twisted Morse-Smale-Witten boundary operator from Definition
\ref{twistedboundary} we have for any $g \in G_q$
\begin{eqnarray*}
\partial_k^G(g q) & = & \sum_{p \in Cr_{k-1}(f)} \sum_{\nu \in \mathcal{M}(q,p)}
                             \epsilon(\nu) \gamma_\ast^\nu(g) p\\
    & = & \sum_{p \in Cr_{k-1}(f)} \left(\sum_{\nu \in \mathcal{M}(q,p)}
          \epsilon(\nu) \right) (\gamma_{e^{k-1}_pe^k_q})_\ast(g) p\\
    & = & \sum_{e_p^{k-1} < e_q^k} \#\mathcal{M}(q,p) (\gamma_{e^{k-1}_pe^k_q})_\ast(g) p,
\end{eqnarray*}
where the last equality follows from the fact that
$$
\overline{W^u(q)} = \bigcup_{q \succeq p} W^u(p),
$$
cf. Corollary 6.27 of \cite{BanLec}, and
$$
 \#\mathcal{M}(q,p) \stackrel{\text{def}}{=} \sum_{\nu \in \mathcal{M}(q,p)} \epsilon(\nu).
$$

Now, if the Riemannian metric $\mathsf{g}$ is locally trivial with respect to the Morse charts of
$f$, i.e. if one can choose the Morse charts to be isometries with respect to the
standard Euclidean metric on $\mathbb{R}^m$ (see Definition 2.16 of \cite{QinOnm}), then
by Theorem 3.9 of \cite{QinOnm} for any $q \in Cr_k(f)$ and $p \in Cr_{k-1}(f)$ we have
$\#\mathcal{M}(q,p) = [e^k_q:e^{k-1}_p]$ where $[e^k_q:e^{k-1}_p]$ is the incidence number 
used to define the CW-boundary operator. The same result was shown to hold for a general
Morse-Smale metric in Theorem 9.3 of \cite{QinAna}. This proves the lemma since Steenrod's
CW-boundary operator from Definition \ref{twistedregular} is given by
$$
\partial_k(g e^k_q) = \sum_{e^{k-1}_p < e^k_q} [e^k_q:e^{k-1}_p] 
(\gamma_{e^{k-1}_pe^k_q})_\ast(g) e^{k-1}_p.
$$
\proofend

\begin{remark}\label{signdegree}
The proof of the previous lemma relied on the identity
$$
\#\mathcal{M}(q,p) = [e^k_q:e^{k-1}_p].
$$
As noted in the proof, this identity was explicitly proved by Qin in \cite{QinOnm} under the
assumption that the metric $\mathsf{g}$ is locally trivial with respect to the Morse charts of
$f$ and then later extended to all metrics by Qin in \cite{QinAna}. It was also proved by
Audin and Damian for gradient-like vector fields in Appendix 4.9 of \cite{AudMor}.
The proof of the identity relies on the manifolds with corners structure on the compactified
moduli space $\overline{\mathcal{M}}(q,p)$, which has been studied by several authors, including Qin,
under various assumptions on the Riemannian metric or the form of the gradient vector
field in the Morse charts, cf. \cite{AudMor}, \cite{BFKOnt}, \cite{BurOnt}, \cite{LatExi}, \cite{QinOnm},
\cite{QinAna}.

In Theorem \ref{regularMorse} we show that it is always possible to find a Morse-Smale pair
$(f,\mathsf{g})$ on $M$ such that there are Morse charts of $f$ that are isometries with respect
to the standard Euclidean metric on $\mathbb{R}^m$ around every critical point and the
unstable manifolds of $(f,\mathsf{g})$ determine a regular CW-structure on $M$. Moreover,
the above identity is easy to establish directly for the function constructed in
Theorem \ref{regularMorse}.  Thus, our proof of Theorem \ref{twistedsame} is independent
of the more general results proved in the above referenced papers.
\end{remark}


\subsection{A Morse-Smale function that determines a regular CW-structure}

It was proved by Laudenbach in the appendix to \cite{BisAne} that if $M$ is a 
finite dimensional closed smooth manifold and $f:M \rightarrow \mathbb{R}$ is Morse-Smale
with respect to a Riemannian metric $\mathsf{g}$ such that the gradient vector field is
``Special Morse'', i.e. if the gradient vector field of $f$ in every Morse chart is the
gradient of $f$ with respect to the standard Euclidean metric on $\mathbb{R}^m$, then
the unstable manifolds of $(f,\mathsf{g})$ determine a CW-structure on $M$. A similar result
was announced by Burghelea, Friedlander, and Kappeler in the Epilogue to \cite{BFKOnt}
for gradient-like vector fields and proved by Audin and Damian for gradient-like vector 
fields in Appendix 4.9 of \cite{AudMor}.  Qin also proved that the unstable manifolds determine 
a CW-structure under the assumption that the Riemannian metric is locally trivial with respect
to the Morse charts of $f$ \cite{QinOnm} and then later extended his proof to general Riemannian 
metrics \cite{QinAna}.

However, none of these results address the question of whether or not the CW-structure
determined by the unstable manifolds is \textbf{regular}.
Examples \ref{circle} and \ref{projectivespace} and the standard height function on the
torus $T^2$, cf. Example 7.11 of \cite{BanLec}, show that the CW-structure determined by
the unstable manifolds of a Morse-Smale pair $(f,\mathsf{g})$ might not be regular.
However, the following theorem shows that it is always possible to find at least one
Morse-Smale pair $(f,\mathsf{g})$ on $M$ such that CW-structure determined by the unstable
manifolds of $(f,\mathsf{g})$ is regular.

\begin{theorem}\label{regularMorse}
On any closed finite dimensional smooth manifold $M$ there exists a smooth Morse-Smale 
pair $(f,\mathsf{g})$ such that the unstable manifolds of $(f,\mathsf{g})$ determine a regular
CW-structure on $M$. Moreover, the Riemannian metric $\mathsf{g}$ can be chosen so that 
there are Morse charts of $f$ around every critical point that are isometries with respect
to the standard Euclidean metric on $\mathbb{R}^m$.
\end{theorem}

\proofstart
Following an idea dual to a method used in one of Thurston's unpublished manuscripts,
cf. Chapter III of \cite{BanSur}, we choose a triangulation of $M$ fine enough so that
every $m$-simplex is contained in a coordinate chart, where $m < \infty$ is the dimension
of the manifold.
Let $n_k$ denote the number of $k$-simplices in the triangulation for all $k=0,\ldots, m$,
let $D_j^k$ denote the $j^\text{th}$ $k$-simplex for all $j=1,\ldots ,n_k$, 
let $p^0_j = D^0_j$ for all $j=1,\ldots, n_0$, and for each $k>0$ pick a basepoint 
$p^k_j \in D^k_j$ in the interior of each $k$-simplex $D^k_j$ for all $j=1,\ldots , n_k$.
Pick any smooth Riemannian metric $\widetilde{\mathsf{g}}$ on $M$, and let 
$\rho:\mathbb{R} \rightarrow \mathbb{R}$ be a smooth bump function such that
\begin{itemize}
\item $0 < \rho(x)  \leq 1$ if $|x| < 1$,
\item $\rho(-x) = \rho(x)$ for all $x$,
\item $\rho(x) = 0$ if $|x| \geq 1$,
\item $\rho(x) = 1$ if $|x| \leq \frac{1}{2}$,
\item $\rho\ '(x) \geq 0$ if $x \leq 0$ and $\rho\ '(x) \leq 0$ if $x \geq 0$,
\item $\rho^{(n)}(\pm 1) = 0$ for all $n\geq 0$.
\end{itemize}
(See the problems at the end of Chapter 6 in \cite{BanLec} for a specific example.)
\begin{figure}[h]
\includegraphics{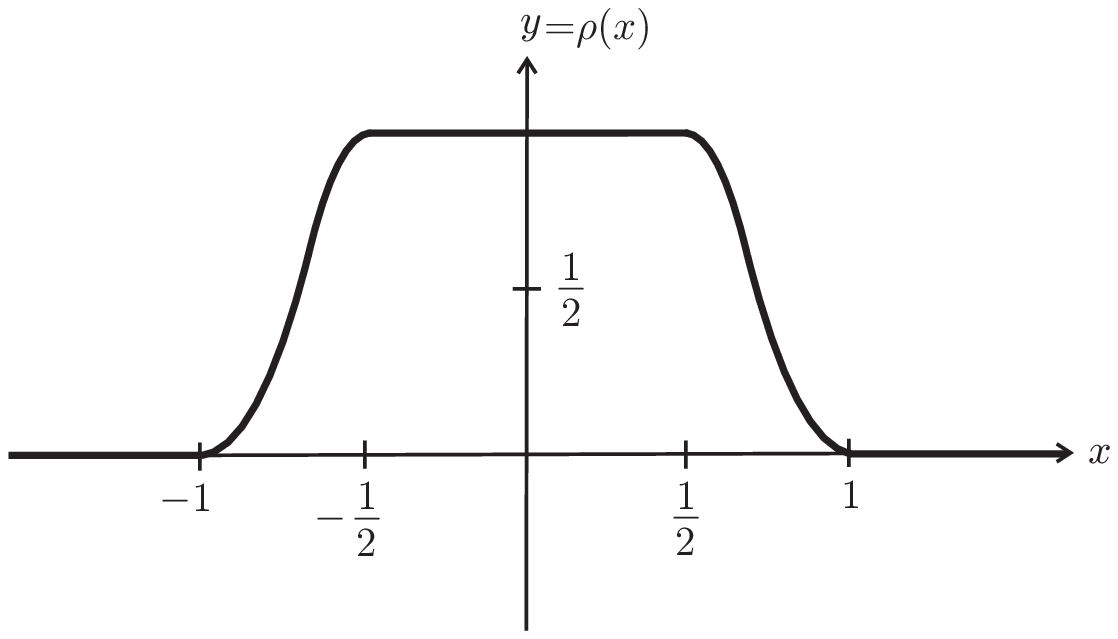}
\end{figure}

\smallskip
Around every $0$-simplex $D^0_j=p_j^0$ choose a small open coordinate neighborhood $T^0_{j,1}$
centered at $p^0_j$ such that $\overline{T^0_{i,1}} \cap \overline{T^0_{j,1}} = \emptyset$ 
if  $i \neq j$. By making a change of coordinates on $T^0_{j,1}$ we will assume
that in the local coordinates $(x_1,\ldots ,x_m)$ we have
$$
T^0_{j,1} = \left.\left\{(x_1,\ldots ,x_m)\right| x_1^2 + \cdots + x_m^2 < 1\right\},
$$
where $(0,\ldots ,0)$ corresponds to the point $p^0_j$. Inside each coordinate
neighborhood we have smaller coordinate neighborhoods $T^0_{j,s} \subseteq T^0_{j,1}$, 
where $s \in (0,1]$.
$$
T^0_{j,s} = \left\{(x_1,\ldots ,x_m) \in T^0_{j,1}\left|\ x_1^2 + \cdots + x_m^2 < s^2
\right\}\right.\\
$$

Modify $\widetilde{\mathsf{g}}$ on $T^0_{j,\frac{1}{2}}$ to a smooth metric $\mathsf{g}^0$ that is
the pullback of the standard metric on $\mathbb{R}^m$ under the coordinate chart, 
cf. Remark 6.31 of \cite{BanLec}, and define a function
$$
f_j^0(x_1,\ldots , x_m) = \rho(x_1^2 + \cdots +x_m^2)(-1+x_1^2 + \cdots + x_m^2)
$$
using the local coordinates on $T^0_{j,1}$.
Note that
$$
f_j^0(x_1,\ldots , x_m)  =  -1 + x_1^2 + \cdots + x_m^2  \text{ for all } 
                         (x_1,\ldots, x_m) \in T^0_{j,\frac{1}{2}}
$$
and for all $(x_1,\ldots ,x_m) \in T^0_{j,1}$ we have
\begin{eqnarray*}
\frac{\partial f_j^0}{\partial x_i}(x_1,\ldots ,x_m) & < & 0 \text { if } \quad x_i < 0\\
\frac{\partial f_j^0}{\partial x_i}(x_1,\ldots ,x_m) & > & 0 \text { if } \quad x_i > 0
\end{eqnarray*}
for all $i=1,\ldots ,m$. Moreover, all the higher order partial derivatives of $f_j^0$ 
are zero on the boundary of $T^0_{j,1}$, and hence we can extend $f_j^0$ to a smooth
function on $M$ by setting it equal to zero outside of $T^0_{j,1}$.
\begin{figure}[h]
\includegraphics{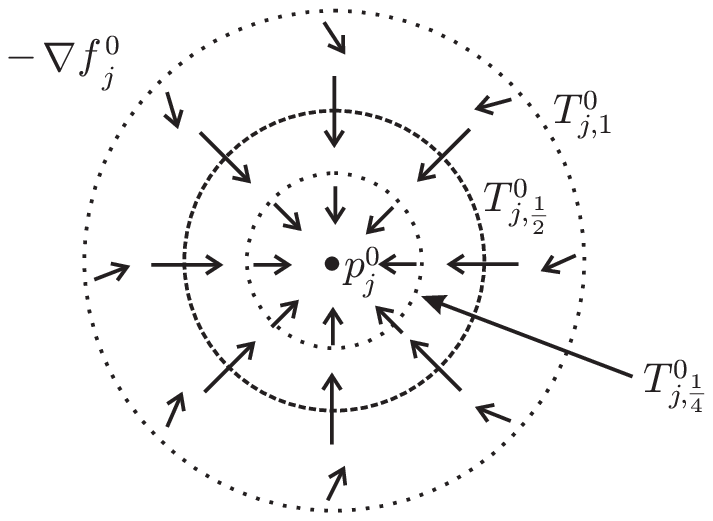}
\end{figure}

\noindent
For all $s \in (0,1]$ define
$$
T^0_s = \bigcup_{j=1}^{n_0} T^0_{j,s} \quad \text{and} \quad f^0 = \sum_{j=1}^{n_0} f_j^0. 
$$
The function $f^0:M \rightarrow \mathbb{R}$ is a smooth function that is Morse-Smale with
respect to the metric $\mathsf{g}^0$ on $T^0_{\frac{1}{2}}$, it has a critical point $p^0_j$ 
of index $0$ at each $0$-simplex $D^0_j$ in the triangulation, and it is equal to
zero on $M - T_1^0$.


\smallskip
Next, for every $1$-simplex $D^1_j$ in the triangulation choose a small open tubular
neighborhood $T^1_{j,1,1}$ of $D^1_j - \overline{T^0_{\frac{1}{4}}}$ with coordinates
$(x_1,\ldots, x_m)$ centered at the basepoint $p^1_j \in D_j^1$, where $x_1$ is 
the coordinate along the simplex $D^1_j$ and $(x_2,\ldots,x_m)$ are coordinates 
normal to $D^1_j$. Choose the tubular neighborhoods small enough so that
$\overline{T^1_{i,1,1}} \cap \overline{T^1_{j,1,1}} = \emptyset$ if $i \neq j$. By making a
change of coordinates on $T^1_{j,1,1}$ we will assume that in the local coordinates 
$(x_1,\ldots , x_m)$ we have
$$
T^1_{j,1,1} = \left.\left\{(x_1,\ldots ,x_m)\right| x_1^2 < 1 \text{ and } x_2^2 + \cdots + 
x_m^2 < 1 \right\},
$$
where $(0,\ldots ,0)$ corresponds to $p^1_j$. Inside each coordinate neighborhood
we have smaller open tubular neighborhoods $T^1_{j,r,s} \subseteq T^1_{j,1,1}$,
where $r,s \in (0,1]$.
$$
T^1_{j,r,s} = \left\{(x_1,\ldots ,x_m) \in T^1_{j,1,1}\left|\ x_1^2 < r^2
\text{ and } \ x_2^2 + \cdots + x_m^2 < s^2 \right\}\right.
$$
By making another change of coordinates on $T^1_{j,1,1}$ we will assume that 
the boundary of the $1$-disk
$$
D^1_{j,\frac{1}{2}}(x_1) = \left.\left\{(x_1,\ldots , x_m) \in T^1_{j,1,1} 
\right| x_1^2 < \left(\frac{1}{2}\right)^2 \text{ and } x_2=\cdots =x_m=0\right\}, 
$$
i.e. $(\pm\frac{1}{2},0,\ldots ,0)$, corresponds to points in $\overline{T^0_{\frac{1}{2}}}$
on or near its boundary and the fibers of $T^1_{j,\frac{1}{2},\frac{1}{2}}$ above
$(\pm\frac{1}{2},0,\ldots ,0)$ are contained in
$\overline{T^0_{\frac{1}{2}}}$.
\begin{figure}[h]
\includegraphics{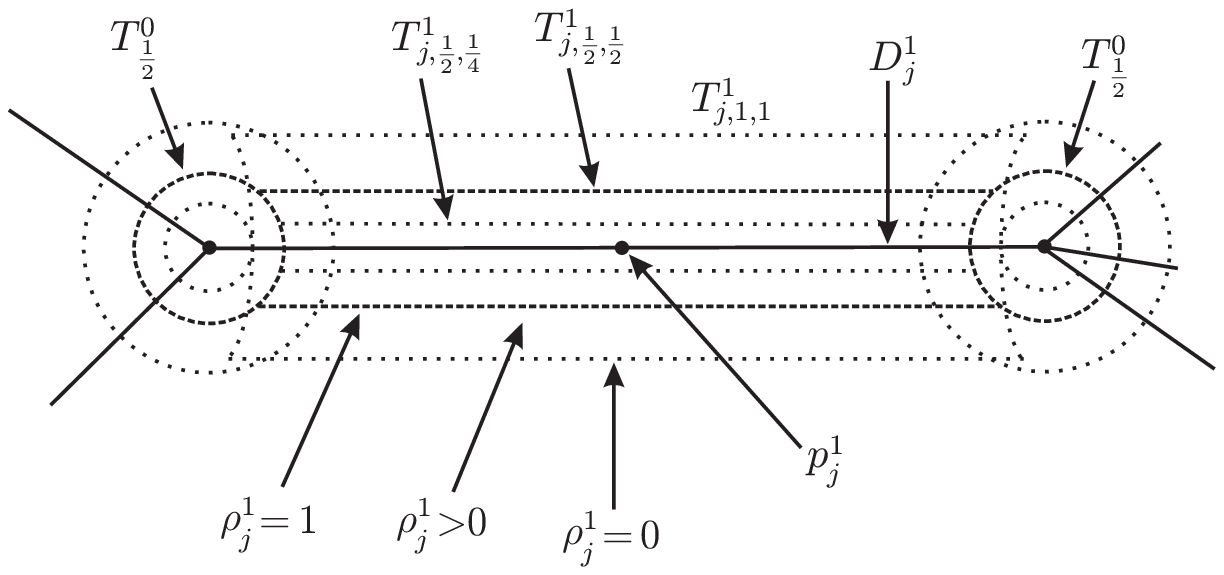}
\end{figure}

Modify the metric $\mathsf{g}^0$ to a smooth metric $\mathsf{g}^1$ that is the pullback of the standard
metric on $\mathbb{R}^m$ under the coordinate chart on $T^1_{j,\frac{1}{2},\frac{1}{2}}$
and define a function
$$
f_j^1(x_1,\ldots ,x_m) = \rho(x_1^2)\rho(x_2^2 + \cdots  + x_m^2)(1-x_1^2 + x_2^2 + 
\cdots + x_m^2)
$$
using the local coordinates on $T^1_{j,1,1}$. We have
$$
f_j^1(x_1,\ldots , x_m)  =  1-x_1^2 + x_2^2 + \cdots + x_m^2 \text{ for all } 
                            (x_1,\ldots, x_m) \in T^1_{j,\frac{1}{2},\frac{1}{2}},
$$
and for all $(x_1,\ldots ,x_m) \in T^1_{j,1,1}$, 
\begin{eqnarray*}
\frac{\partial f_j^1}{\partial x_1}(x_1,\ldots ,x_m) & > & 0 \text { if }  x_1 < 0\\
\frac{\partial f_j^1}{\partial x_1}(x_1,\ldots ,x_m) & < & 0 \text { if }  x_1 > 0.
\end{eqnarray*}
Hence, $p_j^1$ is the only critical point of $f_j^1$ in $T^1_{j,1,\frac{1}{2}}$.
Moreover, all the higher order partial derivatives of $f_j^1$ are zero on the boundary of
$T^1_{j,1,1}$, and hence we can extend $f_j^1$ to a smooth function on $M$ by setting it
equal to $0$ outside of $T^1_{j,1,1}$. 
For all $r,s \in (0,1]$ define
$$
T^1_{r,s} = \bigcup_{j=1}^{n_1} T^1_{j,r,s}, \quad T^{\leq 1}_{r,s} = T^0_{s} 
\cup T^1_{r,s}, \quad f^1 = \sum_{j=1}^{n_1} f_j^1 \quad \text{and} \quad f^{\leq 1} = f^0 + f^1.
$$
The function $f^{\leq 1}$ is a smooth function on $M$ that is Morse-Smale with respect to
the metric $\mathsf{g}^1$ on $T^{\leq 1}_{\frac{1}{2},\frac{1}{2}}$. It has a critical point
$p^0_j$ of index $0$ at each $0$-simplex $D^0_j$ in the triangulation, and the closure of
the unstable manifold $W^u(p_j^1)$ is homeomorphic to $D^1_j \approx \Delta^1$
for each critical point $p^1_j$ of index $1$.  Moreover, there are Morse charts 
for $f^{\leq 1}$ that are isometries with respect to $\mathsf{g}^1$ and the standard Euclidean 
metric on $\mathbb{R}^m$ around every critical point of index $0$ and $1$.

\smallskip
Note that for all $j=1,\ldots , n_1$ the function $f^{\leq 1}$ has a single
critical point $p_j^1 \in T^1_{j,1,\frac{1}{2}}$ and $3^{m-1}-1$ critical points 
$T^1_{j,1,1} - T^1_{j,1,\frac{1}{2}}$, corresponding to the $3^{m-1}$ solutions to the
system of equations
$$
\frac{\partial f_j^1}{\partial x_i}(0,x_2,\ldots ,x_m) = 0
$$ 
for all $i=2,\ldots , m$. We will eliminate the critical points in $T^1_{j,1,1} - 
T^1_{j,1,\frac{1}{2}}$ by adding in functions $f^j$ for $j=2,\ldots, m$ whose gradients
dominate the gradient of $f^{\leq 1}$ on this region.


\smallskip
For every $2$-simplex $D^2_j$ in the triangulation choose a small open tubular
neighborhood $T^2_{j,1,1}$ of 
$$
D^2_j - \overline{T^{\leq 1}_{1,\frac{1}{4}}} 
$$
with coordinates $(x_1,\ldots, x_m)$ centered at the basepoint point $p^2_j \in D^2_j$,
where $(x_1,x_2)$ are coordinates on the simplex $D^2_j$ and $(x_3,\ldots,x_m)$ are
coordinates normal to $D^2_j$. Choose the tubular neighborhoods small enough so that
$\overline{T^2_{i,1,1}} \cap \overline{T^2_{j,1,1}} = \emptyset$ if $i \neq j$. 
By making a change of coordinates on $T^2_{j,1,1}$ we will assume that in
the local coordinates $(x_1,\ldots , x_m)$ we have
$$
T^2_{j,1,1} = \left.\left\{(x_1,\ldots ,x_m)\right| x_1^2 + x_2^2 < 1 \text{ and } 
x_3^2 + \cdots + x_m^2 < 1 \right\},
$$
where $(0,\ldots ,0)$ corresponds to $p^2_j$. Inside each coordinate neighborhood
we have smaller open tubular neighborhoods $T^2_{j,r,s} \subseteq T^2_{j,1,1}$
and $2$-disks $D^2_{j,r}(x_1,x_2)$, where $r,s \in (0,1]$.
\begin{eqnarray*}
T^1_{j,r,s} & = & \left\{(x_1,\ldots ,x_m) \in T^2_{j,1,1}\left|\ x_1^2+x_2^2 < r^2
\text{ and } \ x_3^2 + \cdots + x_m^2 < s^2 \right\}\right.\\
D^2_{j,r}(x_1,x_2) & = & \left.\left\{(x_1,\ldots , x_m) \in T^2_{j,1,1}
\right| x_1^2+x_2^2 < r^2 \text{ and } x_3=\cdots =x_m=0\right\} 
\end{eqnarray*}
By making another change of coordinates on $T^2_{j,1,1}$ we will assume that 
the boundary of the $2$-disk $D^2_{j,\frac{1}{2}}(x_1,x_2)$
corresponds to points in $\overline{T^{\leq 1}_{\frac{1}{2},\frac{1}{2}}}$ on or near
its boundary and the fibers of $T^2_{j,\frac{1}{2},\frac{1}{2}}$ above the boundary points
of $D^2_{j,\frac{1}{2}}(x_1,x_2)$ are contained in $\overline{T^{\leq 1}_{\frac{1}{2},
\frac{1}{2}}}$.
\begin{figure}[h]
\includegraphics{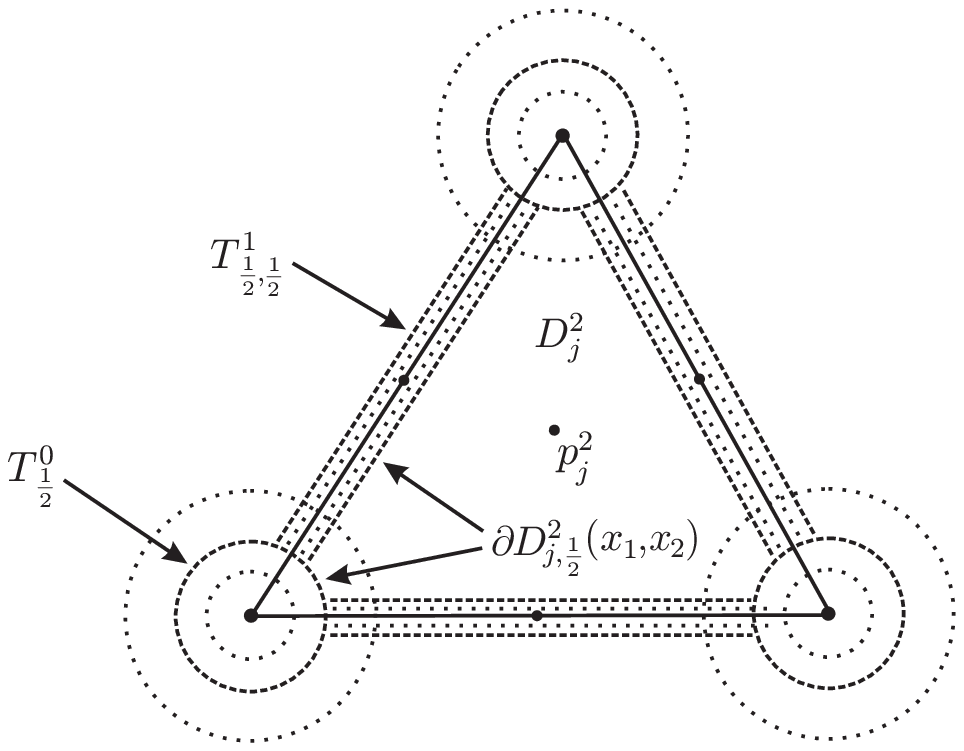}
\end{figure}

Modify the metric $\mathsf{g}^1$ to a smooth metric $\mathsf{g}^2$ that is the pullback of the standard
metric on $\mathbb{R}^m$ under the coordinate chart on $T^2_{j,\frac{1}{2},\frac{1}{2}}$
and define a function
$$
f_j^2(x_1,\ldots ,x_m) = \rho(x_1^2+x_2^2)\rho(x_3^2 + \cdots  + x_m^2)(C_j^2-x_1^2 - x_2^2
+ x_3^2 + \cdots + x_m^2)
$$
using the local coordinates on $T_{j,1,1}^2$, where $C^2_j$ is chosen large enough so that
$f^{\leq 1} + f_j^2$ does not have any critical points in the fibers of $T^2_{j,1,1}$ 
above $D^2_{j,1}(x_1,x_2) - \{p_j^2\}$.  The constant $C^2_j$ exists because
for all $(x_1,\ldots ,x_m) \in T^2_{j,1,1}$ we have
\begin{eqnarray*}
\frac{\partial f_j^2}{\partial x_i}(x_1,\ldots ,x_m) & > & 0 \text { if }  x_i < 0\\
\frac{\partial f_j^2}{\partial x_i}(x_1,\ldots ,x_m) & < & 0 \text { if }  x_i > 0
\end{eqnarray*}
for $i=1,2$, and the magnitude of these partial derivatives at any point in the fibers
of $T^2_{j,1,1}$ above $D^2_{j,1}(x_1,x_2) - \{p_j^2\}$ can be made arbitrarily
large by increasing $C^2_j$. Note that
$$
f_j^2(x_1,\ldots , x_m)  =  C_2-x_1^2 - x_2^2 + x_3^2 + \cdots + x_m^2 \text{ for all } 
                            (x_1,\ldots, x_m) \in T^2_{j,\frac{1}{2},\frac{1}{2}},
$$
and all the higher order partial derivatives of $f_j^2$ are zero on the boundary of 
$T^2_{j,1,1}$. Hence, we can extend $f_j^2$ to a smooth function on $M$ by setting it
equal to $0$ outside of $T^2_{j,1,1}$.  For all $r,s \in (0,1]$ define
$$
T^2_{r,s} = \bigcup_{j=1}^{n_2} T^2_{j,r,s}, \quad T^{\leq 2}_{r,s} = T^{\leq 1}_{r,s} 
\cup T^2_{r,s}, \quad f^2 = \sum_{j=1}^{n_2} f^2_j \quad \text{and} \quad f^{\leq 2} =
f^{\leq 1} + f^2.
$$
The function $f^{\leq 2}$ is a smooth function on $M$ that is Morse-Smale with respect to
the metric $\mathsf{g}^2$ on $T^{\leq 2}_{\frac{1}{2},\frac{1}{2}}$. It has a critical point
$p^0_j$ of index $0$ at each $0$-simplex $D^0_j$ in the triangulation, the closure of
the unstable manifold $W^u(p_j^1)$ is homeomorphic to $D^1_j \approx \Delta^1$
for each critical point $p^1_j$ of index $1$, and the closure of the unstable manifold
$W^u(p_j^2)$ is homeomorphic to $D^2_j \approx \Delta^2$ for each critical point $p^2_j$
of index $2$.  Moreover, there are Morse charts for $f^{\leq 2}$ that are isometries with 
respect to $\mathsf{g}^2$ and the standard Euclidean metric on $\mathbb{R}^m$ around every
critical point $p_j^i$ in $T^{\leq 2}_{\frac{1}{2},\frac{1}{2}}$.


\smallskip
Continuing in this fashion, for every $k=0, \ldots ,m-1$ there are open tubular 
neighborhoods
$$
T^{\leq k}_{r,s} = T^0_s \cup T^1_{r,s} \cup \cdots \cup T^{k}_{r,s}, 
$$
where $r,s \in (0,1]$, a smooth metric $\mathsf{g}^{k}$, and a smooth function 
$f^{\leq k}:M \rightarrow \mathbb{R}$ that satisfy the following.
\begin{itemize}
\item The function $f^{\leq k}$ is Morse-Smale with respect to $\mathsf{g}^{k}$ on
$T^{\leq k}_{\frac{1}{2},\frac{1}{2}}$.
\item The basepoints $p_j^i$ for $i=0,\ldots ,k$ and $j=1,\ldots, n_i$
are the only critical points of $f^{\leq k}$ in $T^{\leq k}_{\frac{1}{2},\frac{1}{2}}$,
and there is a Morse chart for $f^{\leq k}$ around each $p_j^i$ that is an isometry with
respect to $\mathsf{g}^{k}$ and the standard Euclidean metric on $\mathbb{R}^m$.
\item For all $i=0,\ldots ,k$ and $j=1,\ldots, n_i$, the unstable manifold $W^u(p_j^i)$
of the critical point $p_j^i$ is homeomorphic to the simplex $D^i_j \approx \Delta^i$.
\end{itemize}
Finally, around every critical point $p_j^m$ we have the open neighborhood
$$
T^m_{j,1} = D^m_j - \overline{T^{\leq m-1}_{1,\frac{1}{4}}}
$$
with coordinates $(x_1,\ldots, x_m)$ centered at the basepoint $p^m_j \in D^m_j$.
By making a change of coordinates on $T^m_{j,1}$ we will assume that in the
local coordinates we have
$$
T^m_{j,1} = \left.\left\{(x_1,\ldots ,x_m)\right| x_1^2 + \cdots + x_m^2 < 1\right\},
$$
where $(0,\ldots ,0)$ corresponds to the point $p^m_j$.  Inside each coordinate
neighborhood we have smaller coordinate neighborhoods $T^m_{j,r} \subseteq T^m_{j,1}$, 
where $r \in (0,1]$.
$$
T^m_{j,r} = \left\{(x_1,\ldots ,x_m) \in T^m_{j,1}\left|\ x_1^2 + \cdots + x_m^2 < r^2
\right\}\right.\\
$$

Modify $\mathsf{g}^{m-1}$ on $T^m_{j,\frac{1}{2}}$ to a smooth metric $\mathsf{g}^m = \mathsf{g}$
that is the pullback of the standard metric on $\mathbb{R}^m$ under the coordinate chart, 
and define a function
$$
f_j^m(x_1,\ldots , x_m) = \rho(x_1^2 + \cdots +x_m^2)(C^m_j-x_1^2 - \cdots - x_m^2)
$$
using the local coordinates on $T^m_{j,1}$, where $C^m_j$ is chosen large enough so that
$f^{\leq m-1} + f_j^m$ does not have any critical points in $T^m_{j,1} - \{p^m_j\}$.
Note that
$$
f_j^m(x_1,\ldots , x_m)  =  C^m_j - x_1^2 - \cdots - x_m^2  \text{ for all } 
                         (x_1,\ldots, x_m) \in T^0_{j,\frac{1}{2}}
$$
and for all $(x_1,\ldots ,x_m) \in T^m_{j,1}$ we have
\begin{eqnarray*}
\frac{\partial f_j^m}{\partial x_i}(x_1,\ldots ,x_m) & > & 0 \text { if } \quad x_i < 0\\
\frac{\partial f_j^m}{\partial x_i}(x_1,\ldots ,x_m) & < & 0 \text { if } \quad x_i > 0
\end{eqnarray*}
for all $i=1,\ldots ,m$. Moreover, all the higher order partial derivatives of $f_j^m$ 
are zero on the boundary of $T^m_{j,1}$, and hence we can extend $f_j^m$ to a smooth
function on $M$ by setting it equal to zero outside of $T^m_{j,1}$.
Define 
$$
f^m = \sum_{j=1}^{n_m} f^m_j \quad \text{and} \quad f = f^{\leq m} = f^{\leq m-1} + f^m.
$$

The pair $(f,\mathsf{g})$ is a smooth Morse-Smale pair on $M$ whose unstable manifolds
are homeomorphic to the simplices in the chosen triangulation of $M$.  Hence, the
unstable manifolds of $(f,\mathsf{g})$ determine a regular a regular CW-structure on $M$.
Moreover, $\mathsf{g}$ is the pullback of the Euclidean metric on $\mathbb{R}^m$ under a Morse
chart around every critical point of $f$.

\proofend

\begin{remark}\label{CombMorse}
The smooth Morse-Smale pair $(f,\mathsf{g})$ in Theorem \ref{regularMorse} was constructed
by starting with a sufficiently fine triangulation of $M$ and then defining $f$ and
$\mathsf{g}$ so that the unstable manifolds of $(f,\mathsf{g})$ mimic the chosen triangulation.
Hence, it seems likely that Theorem \ref{regularMorse} will have applications to combinatorial
Morse theory. For instance, Theorem 3.1 of \cite{GalCom} proves that given a generic
Morse function $f:M \rightarrow \mathbb{R}$ on a smooth closed oriented manifold $M$ there 
exists a $C^1$-triangulation $T$ of $M$ and a combinatorial Morse vector field $V$ on $T$ that
realizes the smooth Morse-Smale-Witten chain complex.  Theorem \ref{regularMorse} should be
useful for proving a converse to this theorem.
\end{remark}

\begin{remark}\label{signdegreespecial}
To verify the identity $\#\mathcal{M}(q,p) = [e_q^k:e_p^{k-1}]$ directly for the function
constructed in Theorem \ref{regularMorse}, first note that if $\lambda_q - \lambda_p = 1$,
then there is exactly one gradient flow line from $q$ to $p$ if $p \in \overline{W^u(q)}$
and no gradient flow lines from $q$ to $p$ otherwise.

If there are no gradient flow lines from $q$ to $p$, then $e_p^{k-1} \not< e_q^k$ and 
the identity is trivial. If there is one gradient flow line from $q$ to $p$, then the
incidence number $[e_q^k:e_p^{k-1}] = \pm 1$, since the CW-complex is regular and 
$e_p^{k-1} < e_q^k$. So, the identity reduces to a claim concerning compatible orientations
in this case.  

The sign $\#\mathcal{M}(q,p)$ is determined by choosing orientations on $W^u(q)$, $W^u(p)$,
orienting $W(q,p)$ via the short exact sequence 
$$
\xymatrix{
0 \ar[r] & T_\ast W(q,p) \ar@{^{(}->}[r] & T_\ast W^u(q) |_{W(q,p)} \ar[r] & 
\nu_\ast(W(q,p), W^u(q))|_{W(q,p)} \ar[r] & 0,
}
$$
where the fibers of the normal bundle are canonically isomorphic to $T_p W^u(p)$ via the
gradient flow, and then setting $\#\mathcal{M}(q,p) = \pm 1$ depending on whether or not
the resulting orientation on $W(q,p)$ agrees with the orientation given by $-\nabla f$.
In other words, the vector field $-\nabla f$ determines an outward pointing normal vector
on the boundary of the closed $k$-disk $\overline{W^u(q)}$ at $p\in \partial 
\overline{W^u(q)}$, and $\#\mathcal{M}(q,p) = \pm 1$ depending on whether or not
this outward pointing normal vector followed by a positive basis for $T^u_pW^u(p)$ gives
a positive basis for $T_p\overline{W^u(q)}$. Since $\partial \overline{W^u(q)} \approx
S^{k-1}$ is oriented similarly using an outward pointing normal vector (see Section
\ref{pathcomponents}), the incidence number $[e_q^k:e_p^{k-1}]$ is also $\pm 1$ following
the same rule, cf. Section 2.3 of \cite{BanLec} for the definition of incidence numbers.
\end{remark}

\begin{example}[A regular CW-structure on a real projective space]\label{projectivetriangulated}
The following diagram illustrates the gradient flow of the Morse-Smale function constructed
in Theorem \ref{regularMorse} for a minimal triangulation of $\mathbb{R}P^2$ with
ten $2$-simplices. The triangulation in the diagram is one of two possible irreducible
triangulations of $\mathbb{R}P^2$, cf. Figure 3 in \cite{BarGen}. As in Example 
\ref{projectivespace}, diametrically opposed points on the boundary of the closed disk are
identified.  

\begin{figure}[h]\label{RP2triangulation}
\includegraphics{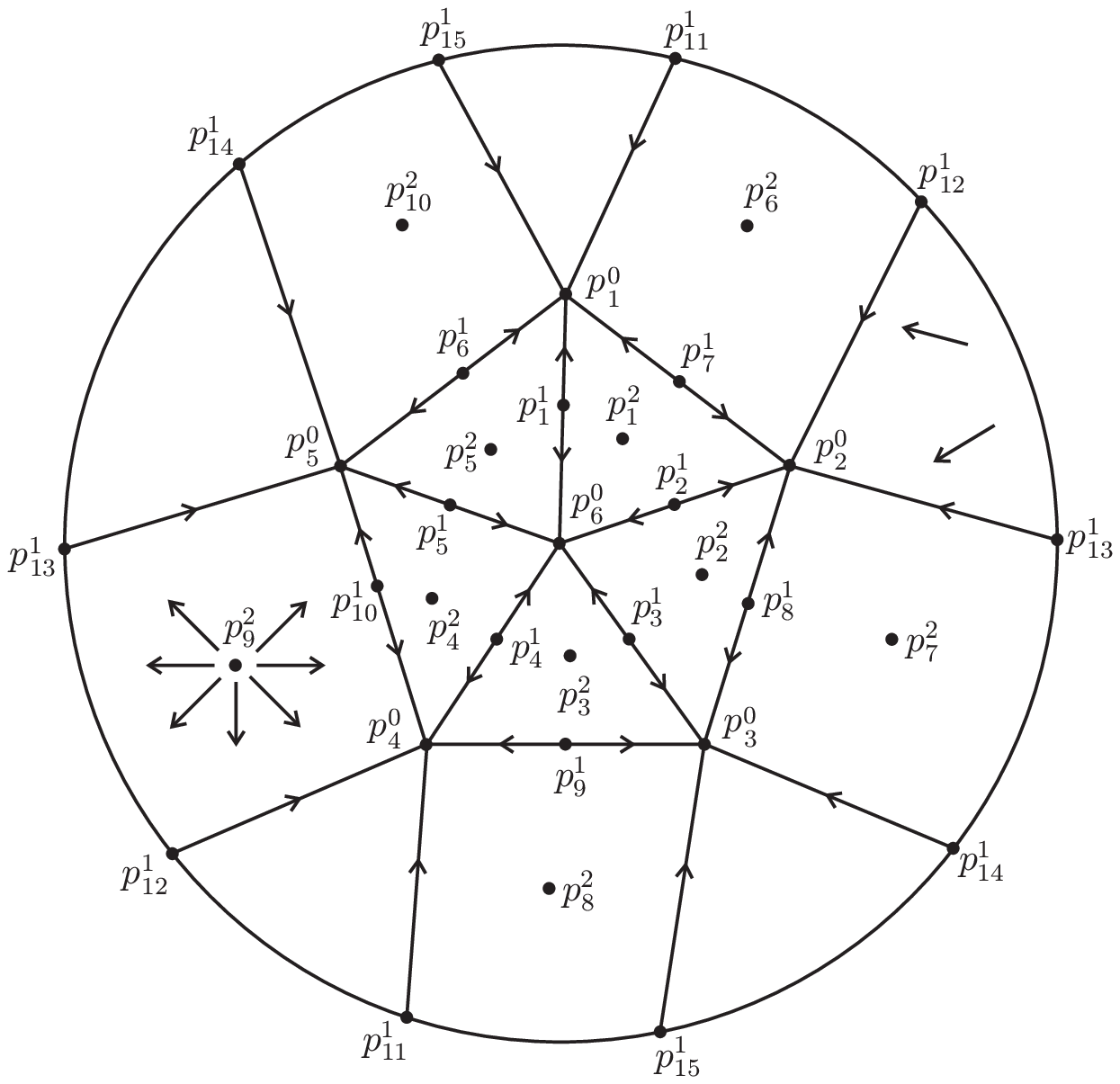}

A regular CW-structure determined by unstable manifolds on $\mathbb{R}P^2$
\end{figure}

Note that the five 2-simplices in the center of the diagram can be combined to create a
single $2$-cell. Moreover, the proof of Theorem \ref{regularMorse} shows that the
resulting regular CW-structure with six $2$-cells is determined by the unstable manifolds
of a Morse-Smale function on $\mathbb{R}P^2$.
\end{example}

\smallskip\noindent
{\bf Proof of Theorem \ref{twistedsame}:}
By the Invariance Theorem (Theorem \ref{homologyindependence}) and Theorem 
\ref{regularMorse} we may assume that $(f,\mathsf{g})$ is a Morse-Smale pair whose unstable
manifolds determine a regular CW-structure $X$ on $M$. Thus, by Lemmas \ref{CWsingular} 
and \ref{boundarysame} we have
$$
H_k(M;G) \approx H_k((CW_\ast(X;G),\partial_\ast)) \approx 
H_k((C_\ast(f;G),\partial_\ast^G))
$$
for all $k = 0, \ldots , m$.
\proofend


\subsection{Local coefficient systems of $R$-modules and the Euler number}

The Twisted Morse Homology Theorem (Theorem \ref{twistedsame}) was proved under the
assumption that $G$ is a bundle of abelian groups over $M$, cf. Definition \ref{bundleofgroups}.
However, some local coefficient systems have additional algebraic structures that may be of interest,
cf. Remark \ref{modulebundle}. For instance, the fiber of the bundle might be a ring, a module over
a ring, a vector space, or a field. In each case, the isomorphisms defining the local system
should be isomorphisms in the appropriate category.

In this section we will consider bundles $\mathcal{L}$ of free $R$-modules over $M$ for some 
commutative ring $R$ with unit. That is, let $M$ be a topological space, and assume that for 
every point $x \in M$ we have a free $R$-module $\mathcal{L}_x$, and for every path 
$\gamma:[0,1] \rightarrow M$ we have an $R$-module isomorphism $\gamma_\ast:\mathcal{L}_{\gamma(1)}
\rightarrow \mathcal{L}_{\gamma(0)}$ that satisfies the conditions listed in 
Definition \ref{bundleofgroups}. This additional algebraic structure allows us to consider 
the \textbf{rank} of the fiber $\mathcal{L}_x$ as a module over $R$, cf. Theorem 5.1 of \cite{RomAdv}.

In general, a submodule of a free module might not be free. So, we will assume 
that $R$ is a principal ideal domain. This assumption guarantees that every submodule of a
free $R$-module $C$ is itself free, and the rank of the submodule is less than or
equal to the rank of $C$, cf. Theorem 6.5 of \cite{RomAdv}. In general, the \text{rank} of
a finitely generated $R$-module $C$ is defined to be the rank of $C/C_{\text{tor}}$, i.e. 
the rank of the module modulo its torsion submodule, cf. Theorem 6.8 of \cite{RomAdv}.

\begin{lemma}\label{ranklemma}
Let $R$ be a principal ideal domain, assume that $C$ and $D$ are finitely generated free 
$R$-modules, and let $\partial:C \rightarrow D$ be an $R$-module homomorphism. Then,
$$
\text{rank }C = \text{rank ker }\partial + \text{rank im } \partial.
$$
Moreover, if $C' \subseteq C$, then
$$
\text{rank }C = \text{rank }C' + \text{rank } C/C'.
$$
\end{lemma}

\proofstart
For any commutative ring with unit, Theorem 5.6 of \cite{RomAdv} implies that
$$
C \approx \text{ker }\partial \oplus \text{im } \partial.
$$
Since $R$ is a principal ideal domain and $C$ is free and finitely generated,
$\text{ker }\partial$ is free and finitely generated, cf. Theorem 6.5 of \cite{RomAdv}.
Similarly, $\text{im } \partial$ is free and finitely generated. Therefore,
$$
\text{rank }C = \text{rank ker }\partial + \text{rank im } \partial.
$$

Now consider the quotient map $C \rightarrow C/C'$. Since $C$ is finitely generated, the
quotient module $C/C'$ is also finitely generated. Therefore,
$$
\frac{C/C'}{(C/C')_{\text{tor}}}
$$
is finitely generated and free since $R$ is a principal ideal domain, 
cf. Theorem 6.8 of \cite{RomAdv}.  Moreover, 
$$
\text{rank ker} \left(C \rightarrow C/C' \right) = 
\text{rank ker} \left(C \rightarrow \frac{C/C'}{(C/C')_{\text{tor}}} \right),
$$
because for any basis $\{c_1,\ldots ,c_j\}$ for the kernel of the homomorphism on the right
there are elements $r_1,\ldots, r_j \in R$ such that $\{r_1c_1,\ldots ,r_jc_j\}$ is a basis for 
the kernel of the homomorphism on the left. Thus, applying the first assertion of the lemma to
the $R$-module homomorphism
$$
C \rightarrow \frac{C/C'}{(C/C')_{\text{tor}}}
$$
yields the second assertion.
\proofend

\begin{theorem}[Euler-Poincar\'e Theorem]\label{EulerPoincare}
Let $R$ be a principal ideal domain, and let 
$$
\xymatrix{
0 \ar[r]^-{\partial_{m+1}} & C_m \ar[r]^-{\partial_m} & C_{m-1} \ar[r]^-{\partial_{m-1}} & \cdots 
\ar[r]^{\partial_2} & C_1 \ar[r]^-{\partial_1} & C_{0} \ar[r]^-{\partial_0} & 0
}
$$
be a finite chain complex of finitely generated free $R$-modules, i.e. assume that $C_k$ is
a finitely generated free $R$-module and $\partial_k$ is an $R$-module homomorphism for all $k$. Then
$$
\sum_{k=0}^m (-1)^k \text{rank }C_k = \sum_{k=0}^m (-1)^k \text{rank }H_k((C_\ast,\partial_\ast)).
$$
\end{theorem}

\proofstart
By Lemma \ref{ranklemma} we have
$$
\text{rank }C_k = \text{rank ker } \partial_k + \text{rank im }\partial_k
$$
and
$$
\text{rank }H_k((C_\ast,\partial_\ast)) = \text{rank ker }\partial_k - \text{rank im }\partial_{k+1} 
$$
for all $k$. Hence, the two sums are equal since  $\text{im }\partial_{m+1}$ and
$\text{im }\partial_0$ are both $0$.
\proofend

\begin{definition}
Let $X$ be a finite regular CW-complex of dimension $m$, and let $\mathcal{L}$ be a bundle of 
finitely generated free $R$-modules over $X$, where $R$ is a principal ideal domain. The 
\textbf{$\mathcal{L}$-twisted Euler number} $\mathcal{X}_{\mathcal{L}}(X)$ is defined to be
$$
\mathcal{X}_{\mathcal{L}}(X) = \sum_{k=0}^m (-1)^k \text{rank }H_k(X;\mathcal{L}),
$$
where $\text{rank }H_k(X;\mathcal{L})$ denotes the rank of the $k^{\text{th}}$ singular homology 
group of $X$ with coefficients in the bundle $\mathcal{L}$ as a module over $R$.
\end{definition}

\begin{theorem}[Invariance of the Twisted Euler Number]\label{Eulerinvariant}
If $X$ is a finite regular CW-complex and $\mathcal{L}$ is a bundle of finitely generated free
$R$-modules over $X$, where $R$ is a principal ideal domain, then the $\mathcal{L}$-twisted
Euler number is well defined, and 
$$
\mathcal{X}_{\mathcal{L}}(X) = \sum_{k=0}^m (-1)^k \text{rank }H_k(X;\mathcal{L}) = 
\sum_{k=0}^m (-1)^k \text{rank }H_k(X;\mathcal{L}_{x_0}) = \mathcal{X}_{\mathcal{L}_{x_0}}(X),
$$
where $x_0\in X$ is any basepoint and $m$ is the dimension of $X$. That is, the 
$\mathcal{L}$-twisted Euler number is the same as the (untwisted) Euler number for homology 
with coefficients in the fiber $\mathcal{L}_{x_0}$. 
\end{theorem}

\proofstart
Since $X$ is a finite CW-complex and $\mathcal{L}$ is a bundle of finitely generated
free $R$-modules, the CW-chain group
$$
CW_k(X;\mathcal{L}) = \bigoplus_{e^k_\sigma} \mathcal{L}_{x(e^k_\sigma)}
$$
is a finitely generated free $R$-module for all $k$, where the sum runs over the $k$-cells
$e^k_\sigma$ in $X$. By Lemma \ref{CWsingular}, the singular boundary operator with coefficients in
$\mathcal{L}$ induces Steenrod's CW-boundary operator with coefficients in $\mathcal{L}$
on a regular CW-complex. Moreover, it's clear that Steenrod's boundary operator is an $R$-module
homomorphism when the isomorphisms that $\mathcal{L}$ associates to paths are $R$-module
homomorphisms, cf. Definition \ref{twistedregular}. Similarly, the singular boundary
operator with coefficients in $\mathcal{L}$ is an $R$-module homomorphism, cf. 
Definition \ref{twistedsingular}. Therefore, the $R$-module homology of Steenrod's CW-chain 
complex with coefficients in $\mathcal{L}$ is isomorphic to the $R$-module singular homology 
of $M$ with coefficients in $\mathcal{L}$, i.e.
$$
H_k((CW_\ast(X;\mathcal{L}), \partial_\ast)) \approx H_k(X;\mathcal{L})
$$
for all $k=0,\ldots ,m$. 
This shows that $H_k(X;\mathcal{L})$ is a finitely generated $R$-module for all $k$, and 
therefore $\mathcal{X}_{\mathcal{L}}(X)$ is well-defined.

The invariance of the $\mathcal{L}$-twisted Euler number follows from the Euler-Poincar\'e Theorem
(Theorem \ref{EulerPoincare}). That is, 
$$
CW_k(X;\mathcal{L}) = \bigoplus_{e^k_\sigma} \mathcal{L}_{x(e^k_\sigma)} \approx \ 
\bigoplus_{e^k_\sigma} \mathcal{L}_{x_0} = CW_k(X;\mathcal{L}_{x_0})
$$
for all $k=0,\ldots, m$, and therefore
$$
\mathcal{X}_{\mathcal{L}}(X) = \sum_{k=0}^m (-1)^k \text{rank }CW_k(X;\mathcal{L}) = 
\sum_{k=0}^m (-1)^k \text{rank }CW_k(X;\mathcal{L}_{x_0}) = \mathcal{X}_{\mathcal{L}_{x_0}}(X).
$$
\proofend

\begin{corollary}\label{Eulerclassical}
If $X$ is a finite regular CW-complex and $\mathcal{L}$ is a bundle of finitely generated free
$R$-modules of rank one over $X$, where $R$ is a principal ideal domain, 
then
$$
\mathcal{X}_{\mathcal{L}}(X) = \mathcal{X}(X),
$$
the classical Euler number of $X$.
\end{corollary}

\proofstart
By the Invariance of the Twisted Euler Number (Theorem \ref{Eulerinvariant}), the 
Euler-Poincar\'e Theorem (Theorem \ref{EulerPoincare}), and the assumption that the fiber 
$\mathcal{L}_{x_0}$ has rank one we have 
$$
\mathcal{X}_{\mathcal{L}}(X) = \mathcal{X}_{\mathcal{L}_{x_0}}(X) = 
\sum_{k=0}^m (-1)^k \text{rank }CW_k(X;\mathcal{L}_{x_0}) = \sum_{k=0}^m (-1)^k |X^{(k)}|,
$$
where $|X^{(k)}|$ denotes the number of $k$-cells in $X$.
The rightmost sum is $\mathcal{X}(X)$.
\proofend

Every closed finite dimensional smooth manifold has a regular CW-structure, cf.
Theorem \ref{regularMorse}. Thus, the previous corollary and the classical Morse inequalities, 
cf. Theorem 3.33 of \cite{BanLec}, yield the following.

\begin{corollary}\label{twistedmMorseequality}
Let $f:M \rightarrow \mathbb{R}$ be a smooth Morse-Smale function on a closed
finite dimensional smooth Riemannian manifold $(M,\mathsf{g})$, and let $\mathcal{L}$ be
a bundle of rank one $R$-modules over $M$, where $R$ is a principal ideal domain.
Then
$$
\mathcal{X}_{\mathcal{L}}(M) = \sum_{k=0}^m (-1)^k \nu_k,
$$
where $\nu_k$ denotes the number of critical points of $f$ of index $k$.
\end{corollary}

\begin{remark}
The rest of the Morse inequalities, cf. Theorem 3.33 of \cite{BanLec},
also hold for the twisted Betti numbers
$$
b_k(\mathcal{L}) \stackrel{\text{def}}{=} \text{rank } H_k(M;\mathcal{L})
$$
when $\mathcal{L}$ is a bundle of finitely generated free $R$-modules of rank one on a 
closed finite dimensional smooth manifold $M$. This follows from The Twisted Morse Homology
Theorem (Theorem \ref{twistedsame}), The Euler-Poincar\'e Theorem (Theorem \ref{EulerPoincare}),
and the proof of Theorem 3.33 in \cite{BanLec}.
\end{remark}


\section{Twisted Morse cohomology and Lichnerowicz cohomology}\label{CoLich}

Each of the chain complexes defined in the previous sections has an associated
cochain complex, e.g. singular cochain complexes with coefficients in a bundle of abelian
groups $G$ \cite{SpaSin}, Steenrod's CW-cochain complexes with coefficients in $G$
(for regular CW-complexes) \cite[Section 31]{SteThe}, and Morse-Smale-Witten cochain
complexes with coefficients in $G$. In this section we will prove that the $\eta$-twisted
Morse cohomology groups are isomorphic to the Lichnerowicz cohomology groups 
associated to $-\eta$.

\subsection{Twisted Morse cohomology}

\begin{definition}[Twisted Morse-Smale-Witten cochain complex]\label{twistedcochain}
Let $f:M \rightarrow \mathbb{R}$ be a smooth Morse-Smale function on a closed 
smooth Riemannian manifold $(M,\mathsf{g})$ of dimension $m < \infty$. Fix orientations on the
unstable manifolds of $(f,\mathsf{g})$, and let $G$ be a bundle of abelian groups over $M$.
For any $k=0,\ldots, m$, a \textbf{Morse-Smale-Witten $k$-cochain with coefficients in $G$}
is defined to be a function $\theta$ that assigns to each critical point $p\in Cr_k(f)$ an
element $\theta(p) \in G_p$. The \textbf{$k^\text{th}$ Morse-Smale-Witten cochain group}
is the collection of $k$-cochains, where the group operation is pointwise application
of the group operation in $G_p$. Hence, 
$$
C^k(f;G)\ \approx \bigoplus_{p \in Cr_k(f)} G_p.
$$
The \textbf{Morse-Smale-Witten cochain complex with coefficients in $G$} is
the chain complex $(C^\ast(f;G),\delta_\ast^G)$ where
$\delta^G_k:C^k(f;G) \rightarrow C^{k+1}(f;G)$ is defined on a 
$k$-cochain $\theta \in C^k(f;G)$ by
$$
(\delta_k^G\theta)(q)\ \ = \sum_{p \in Cr_{k}(f)} \sum_{\nu \in \mathcal{M}(q,p)}
\epsilon(\nu) (\gamma_\nu)_\ast(\theta(p)) \in G_q, 
$$
for any critical point $q \in Cr_{k+1}(f)$, where $\gamma_\nu:[0,1] \rightarrow M$
is any continuous path from $q$ to $p$ whose image coincides with the image of 
$\nu \in \mathcal{M}(q,p)$ and $\epsilon(\nu) = \pm 1$ is the sign determined by the
orientation on $\mathcal{M}(q,p)$. 
\end{definition}

\begin{remark}
As in Definition \ref{twistedboundary}, $-1 \cdot g$ denotes the inverse of $g \in G_q$.
The proof that $(\delta_\ast^G)^2=0$ is similar to the proof of Lemma \ref{boundarysquared}.
\end{remark}

\begin{proposition}\label{isobundlescoho}
If $G_1$ and $G_2$ are isomorphic bundles of abelian groups over $M$, then
$H_k((C^\ast(f;G_1),\delta_\ast^{G_1})) \approx H_k((C^\ast(f;G_2),\delta_\ast^{G_2}))$
for all $k=0, \ldots, m$.
\end{proposition}

\proofstart
By Definition \ref{isobundles} there is a family of isomorphisms $\Phi:G_1 \rightarrow G_2$
that commute with the homomorphisms associated to a path. Thus, for every $k=0, \ldots ,m$
there is an induced isomorphism $\Phi:C^k(f;G_1) \rightarrow C^k(f;G_2)$ given by mapping
$\theta \in C^k(f;G_1)$ to $\Phi \circ \theta \in C^k(f;G_2)$ and $\theta \in C^k(f;G_2)$
to $\Phi^{-1} \circ \theta \in C^k(f;G_1)$. Moreover, $\Phi:(C^\ast(f;G_1),\delta_\ast^{G_1})
\rightarrow (C^\ast(f;G_2),\delta_\ast^{G_2})$ is a chain map, because for all
$\theta \in C^k(f;G_1)$ and $q \in Cr_{k+1}(f)$ we have
\begin{eqnarray*}
(\delta^{G_2}_k (\Phi \circ \theta))(q) & = & 
  \sum_{p \in Cr_{k}(f)} \sum_{\nu \in \mathcal{M}(q,p)}\epsilon(\nu) 
  (\gamma_\nu)_\ast((\Phi\circ\theta)(p))\\
& = & \Phi\left(\sum_{p \in Cr_{k}(f)} \sum_{\nu \in \mathcal{M}(q,p)}\epsilon(\nu) 
      (\gamma_\nu)_\ast(\theta(p))\right)\\
& = & \Phi((\delta^{G_1}_k \theta)(q)) \in (G_2)_q.
\end{eqnarray*}
A similar computation shows that $\Phi^{-1}:(C^\ast(f;G_2),\delta_\ast^{G_2}) \rightarrow
(C^\ast(f;G_1),\delta_\ast^{G_1})$ is a chain map, and hence $\Phi$ is a chain
equivalence.
\proofend

Let $G$ be a bundle of abelian groups over a closed smooth Riemannian manifold
$(M,\mathsf{g})$ of dimension $m < \infty$, let $f:M \rightarrow \mathbb{R}$ be
a smooth Morse-Smale function on $M$, and choose a basepoint $x_0 \in Cr_0(f)$ of $M$.
Recall that if $M$ is connected, then  the isomorphism class of $G$ is determined by a representation
$$
\pi_1(M,x_0) \times G_{x_0} \rightarrow G_{x_0},
$$
and let $H_{x_0}^\ast \subset G_{x_0}$ be the subgroup defined by
$$
H_{x_0}^\ast = \{g \in G_{x_0} |\ \gamma_\ast(g) = g \text{ for all } [\gamma] \in \pi_1(M,x_0) \}.
$$
The following theorem gives the $0$-dimensional twisted Morse cohomology group of $M$ in
terms of the above action, cf. Theorem VI.3.2 of \cite{WhiEle}.

\begin{theorem}\label{H^0}
If $M$ is connected, then the $0$-dimensional twisted Morse cohomology group of $M$
is isomorphic to $H_{x_0}^\ast$, i.e.
$$
H_0((C^\ast(f;G),\delta_\ast^G)) \approx H_{x_0}^\ast.
$$
\end{theorem}

\proofstart
Let $\theta$ be a Morse-Smale-Witten $0$-cochain with coefficients in $G$, i.e. a function that assigns to each
$p \in Cr_0(f)$ an element $\theta(p) \in G_p$.  If $\delta^G_0\theta = 0$, then for all $q \in Cr_1(f)$ we have
$$
(\delta^G_0 \theta)(q) = (\gamma_{\nu_+})_\ast(\theta(p_+)) - (\gamma_{\nu_-})_\ast(\theta(p_-)) 
= 0 \in G_q,
$$
where $\gamma_{\nu_+}$ is a path from $q$ to the positive end $p_+$ of $\overline{W^u(q)}$ and 
$\gamma_{\nu_-}$ is a path from $q$ to the negative end $p_-$ of $\overline{W^u(q)}$. In other words,
the value of $\theta(p_+) \in G_{p_+}$ is determined by the value $\theta(p_-) \in G_{p_-}$ via
the homomorphism associated to a path in $\overline{W^u(q)}$ from $p_+$ to $p_-$.

Since $M$ is connected, every critical point $p \in Cr_0(f)$ can be connected to the basepoint 
$x_0 \in Cr_0(f)$ by a path contained in the 1-skeleton of $f$ (see the proof of Theorem \ref{H_0}).
Hence, if $\theta \in \text{ker } \delta^G_0$, then $\theta(p)$ is determined by $\theta(x_0)$ 
for all $p \in Cr_0(f)$. Moreover, every loop based at $x_0$ can be homotoped through loops based at 
$x_0$ to a loop based at $x_0$ contained in the $1$-skeleton of $f$, and hence $\gamma_\ast(\theta(x_0)) = 
\theta(x_0) \in G_{x_0}$ for all $[\gamma] \in \pi_1(M,x_0)$. 

This shows that mapping $\theta \in \text{ker }\delta^G_0$ to $\theta(x_0)$ gives
an injective homomorphism $\text{ker }\delta^G_0 \rightarrow H_0^\ast$.  To see that this
map is surjective, simply note that if $g \in G_{x_0}$ satisfies $\gamma_\ast(g) = g$ for all 
$[\gamma] \in \pi_1(M,x_0)$, then there is a unique Morse-Smale-Witten 0-cochain $\theta \in \text{ker }
\delta_0^G$ such that $\theta(x_0) = g$.

\proofend

\medskip
In this section we are mainly interested in cohomology with coefficients in a flat line 
bundle $e^\eta$ determined by a closed $1$-form $\eta$.  Note that in this case,
$e^\eta$ is a bundle of rings over $M$ with $e^\eta_q \approx \mathbb{R}$ for all $q\in M$.

\begin{definition}[$\eta$-Twisted Morse-Smale-Witten cochain complex]\label{etacochain}
Let $f:M \rightarrow \mathbb{R}$ be a smooth Morse-Smale function on a closed finite
dimensional smooth Riemannian manifold $(M,\mathsf{g})$. Fix orientations on the
unstable manifolds of $(f,\mathsf{g})$, and let $\eta \in \Omega^1_{cl}(M,\mathbb{R})$.
The \textbf{$\eta$-twisted Morse-Smale-Witten cochain complex} is the chain complex
$(C^\ast(f;e^\eta),\delta_\ast^\eta)$, where $\delta^\eta_k:C^k(f;e^\eta)
\rightarrow C^{k+1}(f;e^\eta)$ is defined on a $k$-cochain $\theta \in C^k(f;e^\eta)$ by
$$
(\delta_k^\eta\theta)(q)\ \ = \sum_{p \in Cr_{k}(f)} \sum_{\nu \in \mathcal{M}(q,p)}
\epsilon(\nu) \exp{\left(\int_{\overline{\mathbb{R}}} (\gamma^\nu)^\ast(\eta)
\right)}\theta(p) \in e^\eta_q, 
$$
for any critical point $q \in Cr_{k+1}(f)$, where $\gamma^\nu:\overline{\mathbb{R}}
\rightarrow M$ is any continuous path from $p$ to $q$ whose image coincides with the image of 
$\nu \in \mathcal{M}(q,p)$ and $\epsilon(\nu) = \pm 1$ is the sign determined by the
orientation on $\mathcal{M}(q,p)$.
\end{definition}

As a consequence of Claim \ref{etaisomorphic} and Proposition \ref{isobundlescoho} we have the
following.

\begin{corollary}\label{cohomologyetainvariant}
If $\eta_1, \eta_2 \in \Omega^1_{cl}(M,\mathbb{R})$ are in the same de Rham cohomology class,
then $H_k((C^\ast(f;e^{\eta_1}), \delta^{\eta_1}_\ast)) \approx H_k((C^\ast(f;e^{\eta_2}),
\delta^{\eta_2}_\ast))$
for all $k=0,\ldots , m$.
\end{corollary}


\subsection{Lichnerowicz cohomology and LCS manifolds}\label{LichLCS}
A closed $1$-form $\eta \in \Omega^1_{cl}(M,\mathbb{R})$ on a finite dimensional
smooth manifold $M$ can be used to deform the differential of the de Rham cochain
complex as follows. For any $k$-form $\xi\in \Omega^k(M,\mathbb{R}$) define
$$
d_\eta \xi = d\xi + \eta \wedge \xi.
$$
It is easy to verify that $d_\eta \circ d_\eta = 0$, and hence $d_\eta$ defines
a cochain complex
$$
\xymatrix{
\Omega^0(M,\mathbb{R}) \ar[r]^-{d_\eta} & \Omega^1(M,\mathbb{R}) \ar[r]^-{d_\eta} & 
\Omega^2(M,\mathbb{R}) \ar[r]^-{d_\eta} & \cdots
}
$$
called the \textbf{Lichnerowicz cochain complex}.  The homology of this complex is
the \textbf{Lichnerowicz cohomology}, denoted by $H^\ast_\eta(M)$; it is 
sometimes referred to as the \textbf{adapted cohomology} of the pair
$(M,\eta)$ \cite{BanSom} \cite{GuiGeo} \cite{OrnMor} \cite{VaiRem}.

\begin{remark}
The differential $d_\eta$ in the Lichnerowicz cochain complex can be viewed
as a generalization of the Witten deformation \cite{WitSup} to closed $1$-forms.
That is, if $\eta = dh$ is exact, then $d_\eta \xi = e^{-h} d (e^h\xi)$ for all 
$\xi\in \Omega^\ast(M,\mathbb{R}$)
\end{remark}

\smallskip
Lichnerowicz cohomology is an invariant used to study locally conformal symplectic
manifolds and locally conformal (almost) K\"ahler manifolds \cite{VaiOnL} \cite{VaiLoc}.

\begin{definition}\label{LCSDef}
A \textbf{locally conformal symplectic (LCS)} form $\Omega$ on a finite dimensional 
smooth manifold $M$ is a smooth nondegenerate $2$-form such that there exists an open cover 
$\mathcal{U} = \{U_i\}_{i \in I}$ of $M$ and smooth positive functions $\lambda_i > 0$
on each $U_i$ such that $\lambda_i \Omega|_{U_i}$ is a symplectic form on $U_i$, i.e. 
$\lambda_i \Omega|_{U_i}$ is closed.

A \textbf{locally conformal almost K\"ahler manifold} is a triple $(M,J,\mathsf{g})$, where
$M$ is a finite dimensional smooth manifold, $J$ is a smooth almost complex structure on $M$,
and $\mathsf{g}$ is a Hermitian Riemannian metric on $M$, such that the $2$-form defined by
$\Omega(X,Y) = \mathsf{g}(X,JY)$ is LCS.

Two LCS forms $\Omega$ and $\Omega'$ on $M$ are said to be \textbf{conformally equivalent}
if and only if there exists a smooth positive function $h > 0$ such that $\Omega' = h \Omega$.
\end{definition}

The connection between the conformal equivalence class of an LCS form $\Omega$ and the 
de Rham cohomology class of an associated closed $1$-form $\eta \in
\Omega^1_{cl}(M,\mathbb{R})$ is given by the following two propositions \cite{LeeAki}
\cite{VaiOnL}.

\begin{proposition}\label{LCSLee}
Let $(M,\Omega)$ be an LCS manifold with an open cover $\mathcal{U}= \{U_i\}_{i \in I}$
and associated smooth positive functions $\lambda_i > 0$ such that $\lambda_i \Omega|_{U_i}$
is a symplectic form.  Then the forms $\{d (\ln \lambda_i)\}_{i \in I}$
fit together to give a smooth closed $1$-form $\eta$ such that
$$
d \Omega = -\eta \wedge \Omega,
$$
and $\eta$ is uniquely determined by the nondegenerate $2$-form $\Omega$.

Conversely, if $\Omega$ is a nondegenerate $2$-form on a smooth manifold $M$ such that
$d \Omega = -\eta \wedge \Omega$ for some closed $1$-form $\eta$, then $\Omega$ is LCS. 
\end{proposition}

\proofstart
Assume $\Omega_i = \lambda_i \Omega|_{U_i}$ is closed. Then
$d \lambda_i \wedge \Omega|_{U_i} + \lambda_i d\Omega|_{U_i} = 0$, and hence
$$
d\Omega|_{U_i} = - d(\ln \lambda_i) \wedge \Omega|_{U_i}.
$$
Moreover, if $U_i \cap U_j \neq \emptyset$ we have
$$
d \Omega|_{U_i \cap U_j} = - d(\ln \lambda_i) \wedge \Omega|_{U_i\cap U_j}
= - d(\ln \lambda_j) \wedge \Omega|_{U_i\cap U_j},
$$
which implies that $d(\ln\lambda_i) = d(\ln \lambda_j)$ since $\Omega$ is
nondegenerate.  Thus, the forms $\{d (\ln \lambda_i)\}_{i \in I}$ fit together to give
a closed $1$-form $\eta$ such that $d \Omega = -\eta \wedge \Omega$.

Conversely, if there exists some $\eta \in \Omega^1_{cl}(M,\mathbb{R})$ such that
$d \Omega = -\eta \wedge \Omega$, then on any contractible open neighborhood $U_i$ we have
$\eta|_{U_i} = df_i$ for some function $f_i$, and setting $\lambda_i = e^{f_i}$ we have
$\eta|_{U_i} = d(\ln \lambda_i)$.  Moreover,
\begin{eqnarray*}
d(\lambda_i \Omega|_{U_i}) 
& = & d\lambda_i \wedge \Omega|_{U_i} + \lambda_i\,  d\Omega|_{U_i}\\
& = & d\lambda_i \wedge \Omega|_{U_i} + \lambda_i (-\eta|_{U_i} \wedge \Omega|_{U_i})\\
& = & d\lambda_i \wedge \Omega|_{U_i} - \lambda_i (d(\ln \lambda_i) \wedge \Omega|_{U_i})\\
& = & \left(d\lambda_i - \lambda_i (d\lambda_i/\lambda_i)\right) \wedge \Omega|_{U_i}\\
& = & 0.
\end{eqnarray*}
\proofend

\begin{definition}
The smooth closed $1$-form $\eta$ in the preceding proposition is called the \textbf{Lee form}
associated to the LCS $2$-form $\Omega$.
\end{definition}

\begin{proposition}\label{LCSConformal}
If $\Omega$ is an LCS form on a finite dimensional smooth manifold $M$ with associated
Lee form $\eta$ and $\Omega' = h \Omega$ for some smooth positive function $h > 0$, 
i.e. $\Omega'$ is conformally equivalent to $\Omega$, then the Lee form associated to
$\Omega'$ is $\eta - d (\ln h)$. Thus, the de Rham cohomology class of the Lee form 
$[\eta] \in H^\ast_{\text{dR}}(M;\mathbb{R})$ is an invariant of the conformal class
of $\Omega$.
\end{proposition}
\proofstart
\begin{eqnarray*}
d \Omega' & = & dh \wedge \Omega + h\, d \Omega\\
& = & dh \wedge \Omega + h (-\eta \wedge \Omega)\\
& = & \frac{dh}{h} \wedge h \Omega - \eta \wedge h\, \Omega\\
& = & -(\eta - d(\ln h)) \wedge \Omega'
\end{eqnarray*}
\proofend

\begin{corollary}\label{etaconformalinv}
Let $(M,\Omega)$ be a closed, smooth, finite dimensional LCS manifold with Lee form 
$\eta \in \Omega^1_{cl}(M,\mathbb{R})$. Then the $\eta$-twisted Morse homology groups
$H_\ast((C_\ast(f) \otimes \mathbb{R},\partial^\eta_\ast))$ and the $\eta$-twisted Morse
cohomology groups $H_\ast((C^\ast(f;e^\eta),\delta^\eta_\ast))$ are invariants of the conformal 
class of $\Omega$.
\end{corollary}

\proofstart
For homology this follows from Corollary  \ref{etadeRhamclass}. For cohomology
this follows from Corollary \ref{cohomologyetainvariant}.
\proofend

\smallskip
We now prove the following generalization of Proposition 4.4 of \cite{LeoOnt}, which
shows that the Lichnerowicz cohomology groups $H^\ast_\eta(M)$ are also an invariant of
the conformal class of the LCS manifold $(M,\Omega)$.

\begin{theorem}\label{Lichcoclass}
For any finite dimensional smooth manifold $M$, the Lichnerowicz cohomology groups
$H_\eta^\ast(M)$ depend only on the cohomology class $[\eta] \in H^\ast_{\text{dR}}
(M;\mathbb{R})$. In particular, if $\eta$ is exact then the Lichnerowicz cohomology
groups are isomorphic to the de Rham cohomology groups, i.e. $H_\eta^k(M) \approx 
H^k_{\text{dR}}(M;\mathbb{R})$ for all $k=0,\ldots , m$.
\end{theorem}

\proofstart
First note that every smooth exact $1$-form $df$ can be written as $dh/h$ for some
smooth positive function $h>0$ by setting $h=e^f$, because $df = d (\ln h) = dh/h$.  So, if
$\eta$ and $\eta'$ are closed $1$-forms on $M$ with $[\eta] = [\eta'] \in 
H^1_{\text{dR}}(M;\mathbb{R})$, then there exists a smooth positive function 
$h:M \rightarrow \mathbb{R}$ such that $\eta' = \eta + dh/h$.

Now, define $\phi:\Omega^\ast(M,\mathbb{R}) \rightarrow \Omega^\ast(M,\mathbb{R})$ 
by $\phi(\xi) = \xi/h$ and $\psi:\Omega^\ast(M,\mathbb{R}) \rightarrow 
\Omega^\ast(M,\mathbb{R})$ by $\psi(\xi) = h\xi$. Clearly, $\phi$ and $\psi$
are isomorphisms on $\Omega^k(M,\mathbb{R})$ for all $k=0,\ldots, m$ with 
$\phi \circ \psi = id = \psi \circ \phi$.
Moreover, $\phi:(\Omega^\ast(M,\mathbb{R}), d_{\eta}) \rightarrow 
(\Omega^\ast(M,\mathbb{R}), d_{\eta'})$ is a chain map, because for all $\xi \in
\Omega^\ast(M,\mathbb{R})$ we have
\begin{eqnarray*}
d_{\eta'} (\phi(\xi)) & = & d \left( \frac{1}{h} \xi \right) + \left( \eta + \frac{dh}{h}\right)
                            \wedge \frac{1}{h}\xi\\
& = & -\frac{1}{h^2} dh \wedge \xi + \frac{1}{h} d\xi + \eta \wedge \frac{1}{h} \xi +
      \frac{1}{h^2} dh \wedge \xi \\
& = & \frac{1}{h} \left( d\xi + \eta \wedge \xi \right)\\
& = & \phi(d_{\eta} \xi).
\end{eqnarray*}

\noindent
Similarly, $\psi:(\Omega^\ast(M,\mathbb{R}), d_{\eta'}) \rightarrow
(\Omega^\ast(M,\mathbb{R}), d_\eta)$ is a chain map, because for all $\xi \in
\Omega^\ast(M,\mathbb{R})$ we have
\begin{eqnarray*}
d_\eta(\psi(\xi)) & = & d(h \xi) + \eta \wedge h\xi\\
& = & dh \wedge \xi + h d\xi + \eta \wedge h\xi \\
& = & h \left(d\xi + \eta \wedge \xi + \frac{dh}{h} \wedge \xi \right)\\
& = & \psi(d_{\eta'} \xi).
\end{eqnarray*}

\noindent
Thus, $\phi$ is a chain equivalence with inverse $\psi$, and
$$
\phi_\ast: H_k(\Omega^\ast(M,\mathbb{R}), d_\eta) \rightarrow 
H_k(\Omega^\ast(M,\mathbb{R}), d_{\eta'})
$$
is an isomorphism for all $k=0,\ldots , m$.
\proofend



\subsection{Mapping differential forms to Morse-Smale-Witten cochains}\label{deRhamSec}
Let $f:M \rightarrow \mathbb{R}$ be a smooth Morse-Smale function on a closed smooth
Riemannian manifold $(M,\mathsf{g})$ of finite dimension $m$. Fix orientations on the 
unstable manifolds of $(f,\mathsf{g})$, and assume that the unstable manifolds determine a regular
CW-structure on $M$. In particular, for any $p \in Cr(f)$ the closure of its unstable
manifold $U_p \stackrel{\text{def}}{=} \overline{W^u(p)}$ is simply connected.

\begin{definition}\label{cochainmap}
Fix any $\eta \in \Omega^1_{\text{cl}} (M,\mathbb{R})$ and note that for any
$p \in Cr(f)$ the restriction $-\eta|_{U_p} = d(\ln h)$ for some smooth positive function
$h:U_p \rightarrow \mathbb{R}$, since $-\eta|_{U_p}$ is exact.
For any $\xi \in \Omega^k(M,\mathbb{R})$, define $\theta_\xi(p) = \xi (p)$ if $k=0$ or
$$
\theta_\xi(p) =  \frac{1}{h(p)} \int_{U_p} h \xi \in e^\eta_p
$$
if $k=1,\ldots, m$. Note that this definition is independent of the choice of $h$ because
$-\eta|_{U_p} = d(\ln \tilde{h}) = d(\ln h)$ implies $\tilde{h} = Ch$ for some $C \in \mathbb{R}$.
Thus, for any $k=0,\ldots , m$ a $k$-form $\xi$ determines a well-defined Morse-Smale-Witten 
$k$-cochain $\theta_\xi$, and $F(\xi) = \theta_\xi$ defines a linear map 
$F:\Omega^k(M,\mathbb{R}) \rightarrow C^k(f;e^\eta)$.
\end{definition}

\begin{remark}
The map $F:\Omega^k(M,\mathbb{R}) \rightarrow C^k(f;e^\eta)$ defined above is similar to
the one used to prove the de Rham Theorem, cf. Section V.5 of \cite{BreTop}, as well as the map
defined in Section 3.3 of \cite{AusMor}.
\end{remark}

\begin{proposition}\label{Fchain}
$F:(\Omega^\ast(M,\mathbb{R}),d_{-\eta}) \rightarrow (C^\ast(f;e^\eta), \delta^\eta_\ast)$ is a
chain map. That is, the map preserves degree and $F \circ d_{-\eta} = \delta^\eta_\ast \circ F$.
\end{proposition}

\proofstart
Pick any $q \in Cr_{k+1}(f)$, let $-\eta|_{U_q} = d(\ln h)$ for some smooth 
positive function $h$ on $U_q = \overline{W^u(q)} \approx \Delta^{k+1}$, and note that for any
$\xi \in \Omega^k(M,\mathbb{R})$ we have
$$
d_{-\eta} \xi = d\xi + \frac{dh}{h} \wedge \xi = \frac{1}{h} d(h\xi)
$$
on $U_q$.  Moreover, the orientations on the unstable manifolds and the signs 
$\epsilon(\nu) = \pm 1$ satisfy the relation 
$$
\partial\overline{W^u(q)}\ = \bigcup_{p \in Cr_k(f)} \bigcup_{\nu \in \mathcal{M}(q,p)}
\epsilon(\nu) \overline{W^u(p)}
$$
as oriented manifolds, cf. Remark \ref{signdegreespecial}. Hence,
\begin{eqnarray*}
(F \circ d_{-\eta}(\xi))(q) = \theta_{d_{-\eta}\xi}(q) 
& = & \frac{1}{h(q)} \int_{U_q} h d_{-\eta} \xi\\
& = & \frac{1}{h(q)} \int_{U_q} d(h\xi)\\
& = & \frac{1}{h(q)} \int_{\partial U_q} h\xi\\
& = & \frac{1}{h(q)} \sum_{p \in Cr_k(f)} \sum_{\nu \in \mathcal{M}(q,p)} 
      \epsilon(\nu) \int_{U_p} h \xi\\
& = & \frac{1}{h(q)} \sum_{p \in Cr_k(f)} \sum_{\nu \in \mathcal{M}(q,p)} 
      \epsilon(\nu) h(p)\theta_\xi(p)\\
& = & \sum_{p \in Cr_k(f)} \sum_{\nu \in \mathcal{M}(q,p)} \epsilon(\nu)
      e^{\ln h(p) - \ln h(q)} \theta_\xi(p)\\
& = & \sum_{p \in Cr_k(f)} \sum_{\nu \in \mathcal{M}(q,p)} \epsilon(\nu)
      \exp\left(\int_{\overline{\mathbb{R}}} (\gamma^\nu)^\ast(-d(\ln h)) \right)
      \theta_\xi(p)\\
& = & \sum_{p \in Cr_k(f)} \sum_{\nu \in \mathcal{M}(q,p)} \epsilon(\nu)
      \exp\left(\int_{\overline{\mathbb{R}}} (\gamma^\nu)^\ast(\eta) \right)
      \theta_\xi(p)\\
& = & (\delta^\eta_k \circ F(\xi))(q),
\end{eqnarray*}
were $\gamma^\nu$ is any parameterization of $\nu$ from $p$ to $q$.

\proofend

To prove that the chain map $F$ induces an isomorphism between the Lichnerowicz groups
$H^\ast_{-\eta}(M)$ and the $\eta$-twisted Morse cohomology groups of $(M,f,\mathsf{g},\eta)$,
we will use the following local version of the $\eta$-twisted Morse-Smale-Witten 
cochain complex.

\begin{definition}[Local $\eta$-twisted Morse-Smale-Witten cochain complex]\label{localetacochain}
For any open set $U \subseteq M$ and for all $k=0,\ldots ,m$ define $C^k_U(f;e^\eta) \subseteq
C^k(f;e^\eta)$ to be the $k$-cochains $\theta$ such that $\theta(p) = 0$ for all 
$p \in Cr_k(f) - U$, and for any $\theta \in C^k(f;e^\eta)$ define $\theta|_U \in C^k_U(f;e^\eta)$
to be 
$$
\theta|_U(p) = \left\{
\begin{array}{ll}
\theta(p) & \text{ if }p \in Cr_k(f) \cap U\\
    0     & \text{ if }p \in Cr_k(f) - U.
\end{array}\right.
$$
Define $\delta^\eta_k|_U:C^k_U(f;e^\eta) \rightarrow C^{k+1}_U(f;e^\eta)$ by restricting the sum
defining $\delta^\eta_k$ in Definition \ref{etacochain} to the critical points $p\in Cr_k(f) \cap U$
when evaluating $\delta^\eta_k|_U \theta$ on a critical point $q \in Cr_{k+1}(f) \cap U$.

We will call an open subset $U \subseteq M$ \textbf{unstably closed} if for every critical point
$p \in Cr(f) \cap U$ we have $\overline{W^u(p)} \subseteq U$. We will call a set $U \subseteq M$
\textbf{unstably 0-connected} if for every $p_1,p_2 \in Cr_0(f) \cap U$ there exists
$q_1, q_2, \ldots , q_n \in Cr_1(f) \cap U$ such that $\bigcup_{i=1}^n \overline{W^u(q_i)}$
is the image of a continuous path from $p_1$ to $p_2$.
\end{definition}

\begin{remark}
If $U \subseteq M$ is unstably closed, then for all $p,q\in Cr(f) \cap U$
such that $\lambda_q - \lambda_p$ is $1$ or $2$ we have $Im(\overline{\mathcal{M}(q,p)}) \subset U$,
cf. Corollary 6.27 of \cite{BanLec}. Thus, $(\delta^\eta_\ast|_U)^2 = 0$ whenever $U$ is unstably
closed, cf. Lemmas \ref{componentboundary} and  \ref{boundarysquared}, and 
$(C^\ast_U(f;e^\eta),\delta_\ast^\eta|_U)$ is a cochain complex. The additional assumption that $U$ is
unstably 0-connected implies that $H_0((C^\ast_U(f;e^\eta),\delta_\ast^\eta|_U)) \approx \mathbb{R}$.
\end{remark}

\begin{lemma}[$\eta$-Twisted Morse Poincar\'e Lemma]\label{MSWPoincare}
Let $f:M \rightarrow \mathbb{R}$ be a smooth Morse-Smale function on a closed finite dimensional
smooth Riemannian manifold $(M,\mathsf{g})$ such that the unstable manifolds of $(f,\mathsf{g})$ determine
a regular CW-structure on $M$. If $U\subseteq M$ is an open contractible set that is unstably
closed and unstably 0-connected, then
$$
F_\ast: H^k_{-\eta}(U) \rightarrow H_k((C_U^\ast(f;e^\eta),\delta^\eta_\ast|_U))
$$
is an isomorphism for all $k=0,\ldots, m$.  Moreover, $H^0_{-\eta}(U) \approx \mathbb{R}$
and $H^k_{-\eta}(U) = 0$ for $k =1, \ldots, m$.
\end{lemma}

\proofstart
Since $U$ is contractible and $\eta$ is closed, $H^k_{-\eta}(U)\approx H^k_{dR}(U;\mathbb{R})$
for all $k$ by Theorem \ref{Lichcoclass}, and the second assertion follows from the
usual Poincar\'e Lemma.  

\smallskip
Now, for every $k=0,\ldots, m$ pick an open neighborhood $N_k$ of $Cr_k(f)\cap U$ consisting of
a union of contractible neighborhoods $N_k^p$ around each critical point $p$ of index
$k$ in $U$, small enough so that no component of $N_k$ contains more than one critical point and
$N_j \cap N_k = \emptyset$ if $j \neq k$.  For all $k=1,\ldots, m$ let $\xi_{k-1}$ be a smooth
$(k-1)$-form with support in $N_k$ such that 
$$
\frac{1}{h(p)} \int_{U_p} h d_{-\eta}\xi_{k-1} = 1 \in e_p^\eta
$$  
for all $p \in Cr_k(f) \cap U$, where $U_p = \overline{W^u(p)}$ and $-\eta|_U = d( \ln h)$.

For all $k=0,\ldots, m$ and every $\theta_k \in C^k(f;e^\eta)$ let $\overline{\theta_k}$
be the extension of $\theta_k$ to $N_k$ that is constant on each component $N_k^p$ of $N_k$, 
i.e. $\overline{\theta_k}|_{N_k^p} \equiv \theta_k(p)$, and define a linear map 
$G:C^k(f;e^\eta) \rightarrow \Omega^k(M,\mathbb{R})$ by
$$
G(\theta_k) = \overline{\theta_k} d_{-\eta} \xi_{k-1} + \overline{\delta^\eta_k \theta_k} \xi_k,
$$ 
where we define $\xi_{-1} = \xi_{m} = 0$.
Note that $G(\theta_k)$ is a smooth $k$-form, even though $\overline{\theta_k}$ and
$\overline{\delta^\eta_k\theta_k}$ are discontinuous, because $d_{-\eta}\xi_{k-1}$ is smooth
with support in $N_k$ and $\xi_k$ is smooth with support in $N_{k+1}$.
In other words, for every $k=1,\ldots, m$, $N_k$ is a disjoint union of contractible 
open sets $N_k = \bigcup_{p \in Cr_k(f)} N_k^p$, and we can define $\xi_{k-1}^p$ to 
be the smooth $(k-1)$-form that is equal to $\xi_{k-1}$ on the component $N_k^p$ and zero
outside of $N_k^p$.  Then
$$
G(\theta_k) = \sum_{p \in Cr_k(f)} \theta_k(p) d_{-\eta} \xi_{k-1}^p\ +
              \sum_{q \in Cr_{k+1}(f)} (\delta^\eta_k \theta_k)(q) \xi_k^q,
$$
which is a smooth $k$-form on $M$.

For all $k=0,\ldots, m$ and any $\theta_k \in C^k(f;e^\eta)$ we have
\begin{eqnarray*}
d_{-\eta}(G(\theta_k)) & = & d_{-\eta} \left( \overline{\theta_k} d_{-\eta} \xi_{k-1} 
                                    + \overline{\delta^\eta_k \theta_k} \xi_k\right)\\
& = & \overline{\theta_k} d_{-\eta}^2 \xi_{k-1} + \overline{\delta^\eta_k \theta_k} 
      d_{-\eta} \xi_k\\
& = & \overline{\delta^\eta_k \theta_k} d_{-\eta} \xi_k + \overline{\delta^\eta_{k+1}
      \delta^\eta_k \theta_k} \xi_{k+1}\\
& = & G(\delta^\eta_k \theta_k).
\end{eqnarray*}
Thus, $d_{-\eta}  \circ G = G \circ \delta^\eta_\ast$, i.e. $G:(C^\ast(f;e^\eta),
\delta_\ast^\eta) \rightarrow (\Omega^\ast(M,\mathbb{R}),d_{-\eta})$ is a chain map. Moreover,
for all $k=1, \ldots, m$ we have for any $\theta_k \in C^k(f;e^\eta)$
\begin{eqnarray*}
((F \circ G)(\theta_k))(p) & = & F(\overline{\theta_k} d_{-\eta} \xi_{k-1} + \overline{\delta_k^\eta \theta_k} \xi_k)(p)\\
& = & \frac{1}{h(p)} \int_{U_p} h \overline{\theta_k}d_{-\eta}\xi_{k-1}\\
& = & \theta_k(p) \frac{1}{h(p)}\int_{U_p}  h d_{-\eta}\xi_{k-1}\\
& = & \theta_k(p) \in e^\eta_p,
\end{eqnarray*}
for all $p \in Cr_k(f)$, and hence $(F \circ G)_\ast: H_k((C_U^\ast(f;e^\eta),\delta^\eta_\ast|_U))
\rightarrow H_k((C_U^\ast(f;e^\eta),\delta^\eta_\ast|_U))$ is the identity for all $k=1,\ldots , m$.
This implies that $F_\ast:H_{-\eta}^k(M) \rightarrow H_k((C^\ast_U(f;e^\eta), \delta^\eta_\ast|_U))$
is surjective for $k=1,\ldots , m$, and hence $F_\ast$ is the trivial isomorphism for all 
$k=1, \ldots, m$, since $H_{-\eta}^k(M) = 0$.

\smallskip
Now, if $\xi \in \Omega^0(U)$, then
$$
d_{-\eta} \xi = d \xi + \frac{dh}{h} \wedge \xi = \frac{1}{h} d(h\xi) 
$$
implies that $\xi \in \text{ker }d_{-\eta}$ if and only if $\xi = C/h$ for some $C \in \mathbb{R}$,
where $-\eta|_U = d(\ln h)$. Moreover, Proposition \ref{Fchain} and the assumption that $U$ is
unstably closed implies that if $\xi \in \text{ker }d_{-\eta}$, then $\theta_\xi = F(\xi)$
satisfies $(\delta^\eta_0 \theta_\xi)(q) = 0$ for all $q \in Cr_1(f) \cap U$, 
i.e. $F:\text{ker }d_{-\eta} \rightarrow \text{ker }\delta_0^\eta$. Note that $F(C/h) =
C/h|$, where $h|$ denotes the restriction of $h:U \rightarrow \mathbb{R}$ to $Cr_0(f)\cap U$.

Define a map $G:\text{ker }\delta_0^\eta \rightarrow \text{ker }d_{-\eta}$ by $G(\theta_0) = 
\theta_0(p)h(p)/h$ for any $p \in Cr_0(f) \cap U$. To see that this definition is independent of
$p \in Cr_0(f) \cap U$, first note that if $p_+, p_- \in Cr_0(f) \cap U$ are the boundary points
of an unstable manifold $W^u(q)$ with $q \in Cr_1(f) \cap U$, then $\mathcal{M}(q,p_+)$ consists
of a single element $\nu_{p_+}$ and $\mathcal{M}(q,p_-)$ consists of a single element $\nu_{p_-}$.
Moreover,  $\epsilon(\nu_{p_+}) = - \epsilon(\nu_{p_-})$, and we can label the points so that
$\epsilon(\nu_{p_+}) = +1$. Thus, 
\begin{eqnarray*}
0 & = & (\delta^\eta_0 \theta_0)(q)\\
   & = & \exp \left( \int_{\overline{\mathbb{R}}} (\gamma^{\nu_{p_+}})^\ast ( \eta)  \right) \theta_0(p_+) - 
         \exp \left( \int_{\overline{\mathbb{R}}} (\gamma^{\nu_{p_- }})^\ast  (\eta)  \right) \theta_0(p_-)\\
   & = & e^{\ln h(p_+) - \ln h(q)} \theta_0(p_+) - e^{\ln h(p_-) - \ln h(q)} \theta_0(p_-)\\
   & = & \frac{h(p_+)}{h(q)} \theta_0(p_+) - \frac{h(p_-)}{h(q)} \theta_0(p_-),
\end{eqnarray*}
where $\gamma^{\nu_{p_+}}$ and $\gamma^{\nu_{p_-}}$ are paths that end at $q$.  Therefore,
$h(p_+)\theta_0(p_+) = h(p_-) \theta_0(p_-)$. The fact that $G(\theta_0) = \theta_0(p)h(p)/h$
does not depend on the choice of $p \in Cr_0(f) \cap U$ now follows inductively using the assumption
that $U$ is unstably 0-connected.

So, for any $\theta_0 \in \text{ker }\delta_0^\eta$ we have $(F\circ G)(\theta_0) = 
F\left(\theta_0(p)h(p)/h\right) = \theta_0(p)h(p)/h| = \theta_0$, since $\theta_0(p)h(p)$
does not depend on $p \in Cr_0(f) \cap U$.  Moreover, for any $C/h \in \text{ker } d_{-\eta}$
we have $(G\circ F)(C/h)  =  G(C/h|) = \frac{C}{h(p)} h(p)/h = C/h$, and hence $F \circ G = 
G \circ F = id$. Therefore, $G:\text{ker }\delta_0^\eta \rightarrow \text{ker }d_{-\eta}$ is a 
well-defined group homomorphism with inverse $F$, and $F_\ast:H_{-\eta}^0(M) \rightarrow
H_0((C^\ast_U(f;e^\eta), \delta^\eta_\ast|_U))$ is an isomorphism.
\proofend

\begin{theorem}[$\eta$-Twisted Morse de Rham Theorem]\label{twistedDeRham}
Let $f:M \rightarrow \mathbb{R}$ be a smooth Morse-Smale function on a closed finite
dimensional smooth Riemannian manifold $(M,\mathsf{g})$.  Fix orientations on the unstable 
manifolds of $(f,\mathsf{g})$ and assume that the unstable manifolds determine a regular
CW-structure on $M$. For any $\eta \in \Omega^1_{\text{cl}} (M,\mathbb{R})$, the 
$\eta$-twisted Morse cohomology groups are isomorphic to the Lichnerowicz cohomology
groups defined by $-\eta$, i.e.
$$
H_k((C^\ast(f;e^\eta),\delta_\ast^\eta)) \approx H^k_{-\eta}(M)
$$
for all $k=0,\ldots ,m$.
\end{theorem}

\proofstart
Recall that for any two open sets $U,V$ in $M$ there is short exact
sequence of cochain complexes
$$
\xymatrix{
0 \ar[r] & \Omega^k(U \cup V) \ar[r] & \Omega^k(U) \oplus \Omega^k(V) \ar[r] & 
\Omega^k(U \cap V) \ar[r] & 0,
}
$$
where the second map is defined by $\xi \mapsto (\xi|_{U},-\xi|_{V})$ and the third map
is given by $(\xi,\xi') \mapsto \xi|_{U \cap V} + \xi'|_{U \cap V}$. Assuming that $U$ and $V$
are unstably closed, we can define a similar short exact
sequence of cochain complexes
$$
\xymatrix{
0 \ar[r] & C^k_{U \cup V}(f;e^\eta) \ar[r] & C^k_{U}(f;e^\eta) \oplus C^k_{V}(f;e^\eta)
\ar[r] & C^k_{U \cap V}(f;e^\eta) \ar[r] & 0,
}
$$
where the second map is defined by $\theta \mapsto (\theta|_{U},-\theta|_{V})$ and the 
third map is given by $(\theta,\theta') \mapsto \theta|_{U \cap V} + \theta'|_{U \cap V}$.

Moreover, the chain map $F:(\Omega^\ast(M,\mathbb{R}),d_{-\eta}) \rightarrow (C^\ast(f;e^\eta),
\delta^\eta_\ast)$ induces chain maps that fit into the commutative diagram
$$
\xymatrix{
0 \ar[r] & \Omega^k(U \cup V) \ar[r] \ar[d]^F & \Omega^k(U) \oplus \Omega^k(V) \ar[r]
\ar[d]^{F \oplus F} & \Omega^k(U \cap V) \ar[r] \ar[d]^F & 0 \\
0 \ar[r] & C^k_{U \cup V}(f;e^\eta) \ar[r] & C^k_{U}(f;e^\eta) \oplus C^k_{V}(f;e^\eta)
\ar[r] & C^k_{U \cap V}(f;e^\eta) \ar[r] & 0,
}
$$
which in turn induces the following ladder commutative diagram with exact rows.
$$
\xymatrix{
\cdots \ar[r] & H_{-\eta}^k(U \cup V) \ar[r] \ar[d]^{F_\ast} & H_{-\eta}^k(U) \oplus H_{-\eta}^k(V) \ar[r]
\ar[d]^{F_\ast \oplus F_\ast} & H_{-\eta}^k(U \cap V) \ar[r] \ar[d]^{F_\ast} & \cdots \\
\cdots \ar[r] & H^k_{U \cup V}(f;e^\eta) \ar[r] & H^k_{U}(f;e^\eta) \oplus H^k_{V}(f;e^\eta)
\ar[r] & H^k_{U \cap V}(f;e^\eta) \ar[r] & \cdots
}
$$
The $5$-lemma applied to this diagram implies that if $F_\ast$ is an isomorphism on
$U$, $V$, and, $U \cap V$ for $k = i-1$ and $i$, then $F_\ast$ is an isomorphism on $U \cup V$
when $k=i$.

\smallskip
Now, note that both the union and intersection of unstably closed sets is unstably closed, and $F_\ast$
is an isomorphism for all $k$ on every open unstably closed set with components that are contractible
and unstably 0-connected by Lemma \ref{MSWPoincare}. Thus, if $V_1$ and $V_2$ are open and unstably
closed, and $V_1$, $V_2$, and $V_1 \cap V_2$ have components that are contractible and unstably 0-connected,
then the above ladder commutative diagram implies that $F_\ast$ is an isomorphism for all $k$ on
$V_1 \cup V_2$. Continuing by induction, assume that we have a finite collection $\{V_j\}_{j=1}^n$ of
open unstably closed sets with components that are contractible and unstably 0-connected
such that 
\begin{itemize}
\item $\left(\bigcup_{j=1}^{n-1} V_j\right) \cap V_n$ has components that are contractible and unstably 0-connected
\item $F_\ast$ is an isomorphism for all $k$ on $\bigcup_{j=1}^{n-1} V_j$.
\end{itemize}
Then taking $U = \bigcup_{j=1}^{n-1} V_j$ and $V = V_n$, Lemma \ref{MSWPoincare} and the above
ladder commutative diagram imply that $F_\ast$ is an isomorphism for all $k$ on 
$\bigcup_{j=1}^n V_j$. This shows that $F_\ast$ is an isomorphism for all $k$ on
any unstably closed set $U$ that can be written as a finite union $U = V_1 \cup \cdots \cup V_n$ of
open unstably closed sets $\{V_j\}_{j=1}^n$ with components that are contractible and unstably
0-connected such that
\begin{enumerate}
\item[$(\ast)$] $\left(\bigcup_{j=1}^{i-1} V_j\right) \cap V_i$  has components that are contractible
                and unstably 0-connected for all $i=2, \ldots, n$.
\end{enumerate} 

\smallskip
To complete the proof, take a finite open cover $\{U_j\}_{j=1}^n$ of contractible sets
that are unstably closed and unstably 0-connected consisting of a small thickening of
$\overline{W^u(p_j)}$ for each critical point $p_j$ of index $m = \text{dim }M$. 
If $U_1 \cap U_2 \neq \emptyset$, then this intersection can be written as a finite union of open
unstably closed sets $\{V_j\}_{j=1}^n$ with components that are contractible and unstably 0-connected satisfying
condition $(\ast)$ because the unstable manifolds of $(f,\mathsf{g})$ determine a regular CW-structure on $M$.
Thus, Lemma \ref{MSWPoincare} and the above ladder commutative diagram imply that $F_\ast$ is an isomorphism
for all $k$ on $U_1 \cup U_2$. 

Now assume that we have shown that $F_\ast$ is an isomorphism for all $k$ on $U_1 \cup \cdots \cup U_{i-1}$ 
for some $i = 2, \ldots, n$. The set $\left( \bigcup_{j=1}^{i-1} U_i \right) \cap U_i$ can be written as
a finite union of open unstably closed sets $\{V_j\}_{j=1}^n$ with components that are contractible
and unstably 0-connected satisfying condition $(\ast)$ since the unstable manifolds of $(f,\mathsf{g})$
determine a regular CW-structure. Therefore, Lemma \ref{MSWPoincare} and the above ladder commutative
diagram imply that $F_\ast$ is an isomorphism for all $k$ on $U_1 \cup \cdots \cup U_i$ and hence on
$U_1 \cup \cdots \cup U_n = M$. 

\proofend

\begin{remark}
The above proof is modeled on the proof of the de Rham Theorem found in Section V.9 of \cite{BreTop}.
\end{remark}

\begin{corollary}
The homology groups of the $\eta$-twisted Morse-Smale-Witten cochain complex $(C^\ast(f;e^\eta), 
\delta^\eta_\ast)$ do not depend on the pair $(f,\mathsf{g})$, as long as the unstable manifolds of
$(f,\mathsf{g})$ determine a regular CW-structure on the Riemannian manifold $(M,\mathsf{g})$.
\end{corollary}

\begin{remark}
The preceding corollary can also be proved using the methods found in Section \ref{continuation}.
Using that approach, we could remove the assumption that the unstable manifolds of $(f,\mathsf{g})$
determine a regular CW-structure.
\end{remark}

\begin{corollary}[Invariance of the $\eta$-twisted Euler number]\label{coEulerindex}
Let $M$ be a closed finite dimensional smooth manifold and $\eta \in \Omega^1_{\text{cl}}(M,\mathbb{R})$
a closed 1-form on $M$. The $\eta$-twisted Lichnerowicz Euler number
$$
\mathcal{X}^\eta(M) \stackrel{\text{def}}{=} \ \sum_{k=0}^{m} (-1)^k \text{dim }H^k_{\eta}(M)
$$
is well-defined and equal to the classical Euler number $\mathcal{X}(M)$.
\end{corollary}

\proofstart
Note that $(\Omega^\ast(M,\mathbb{R}),d_{-\eta})$ is a complex of $\mathbb{R}$-modules which is
not finitely generated, whereas $(C^\ast(f;e^\eta),\delta_\ast^\eta)$ is a complex 
of finitely generated free $\mathbb{R}$-modules. Moreover, the map $F:\Omega^k(M,\mathbb{R})
\rightarrow C^k(f;e^\eta)$ from Definition \ref{cochainmap} is an $\mathbb{R}$-module homomorphism.
Therefore, the $\eta$-Twisted Morse de Rham Theorem (Theorem \ref{twistedDeRham}) implies that
$H_{-\eta}^k(M)$ is a finite dimensional real vector space for all $k$ and all 
$\eta \in \Omega^1_{\text{cl}}(M,\mathbb{R})$. This implies that the $\eta$-twisted Lichnerowicz
Euler number
$$
\mathcal{X}^\eta(M) = \sum_{k=0}^{m} (-1)^k \text{dim }H^k_{\eta}(M)
= \sum_{k=0}^{m} (-1)^k \text{dim }H_k((C^\ast(f;e^{-\eta}),\delta_\ast^{-\eta}))
$$
is well-defined. The fact that $\mathcal{X}^\eta(M) = \mathcal{X}(M)$ follows from
the $\eta$-Twisted Morse de Rham Theorem (Theorem \ref{twistedDeRham})
and the Euler-Poincar\'e Theorem (Theorem \ref{EulerPoincare}), 
cf. Theorem \ref{Eulerinvariant}.
\proofend

\begin{remark}\label{AnalyticEuler}
Corollary \ref{coEulerindex} can also be proved for the complex 
$(\Omega^\ast(M,\mathbb{R}),d_{\eta})$ using the index theory of elliptic operators. 
That is, in Section 6 of their classic paper \cite{AtiIndIII} Atiyah and Singer note that
on any closed finite dimensional smooth manifold $M$ the de Rham complex of complex valued differential
forms is elliptic, and its Euler characteristic is the usual Euler number. They then consider the
operator $d + d^\ast$, where $d^\ast$ denotes the adjoint of $d$ with respect to a Riemannian
metric on $M$. Restricting $d + d^\ast$ to the even dimensional complex valued differential forms gives
an elliptic differential operator of order 1 whose index is the Euler number of the 
de Rham complex, cf. Section III.4.D of \cite{BooTop} and Proposition 19.1.16 of \cite{HorTheIII}.
Perturbing $d$ to $d_\eta$ is a lower order perturbation, and hence the index of 
$d + d^\ast$ does not change under the perturbation.
Therefore,
$$
\mathcal{X}(M) = \text{index} (d + d^\ast) = \text{index} (d_\eta + (d_\eta)^\ast) = \mathcal{X}^\eta(M).
$$
\end{remark}


\subsection{Relationship to sheaf cohomology}\label{sheafSec}
Proposition \ref{Fchain} shows that the map $F$ given in Definition \ref{cochainmap} is a chain map
from the Lichnerowicz cochain complex $(\Omega^\ast(M,\mathbb{R}), d_{-\eta})$ to the $\eta$-twisted
Morse-Smale-Witten cochain complex $(C^\ast(f;e^\eta),\delta^\eta_\ast)$, and Theorem \ref{twistedDeRham}
shows that this chain map induces an isomorphism of the corresponding cohomology groups. A result due
to Vaisman shows that the Lichnerowicz cohomology groups $H^\ast_{-\eta}(M)$ are also
isomorphic to the cohomology of $M$ with coefficients in a sheaf \cite{VaiRem}.

\smallskip
Vaisman defines $\mathscr{F}_\eta(M)$ to be the sheaf of germs of differentiable functions
$f:M \rightarrow \mathbb{R}$ such that 
$$
d_{-\eta} f = df - \eta f = 0.
$$
Letting $\mathscr{A}^k(M)$ be the sheaf of germs of differentiable $k$-forms on $M$ for all
$k=0, \ldots , m$, he shows that 
$$
\xymatrix{
0 \ar[r] & \mathscr{F}_\eta(M) \ar@{^{(}->}[r] & \mathscr{A}^0(M) \ar[r]^{d_{-\eta}} & 
\mathscr{A}^1(M) \ar[r]^{d_{-\eta}} & \cdots
}
$$
is a fine resolution of $\mathscr{F}_\eta(M)$ by proving that $d_{-\eta}$ satisfies a
Poincar\'e Lemma. This proves that $H^k(M;\mathscr{F}_\eta(M))$, the $k^{th}$-cohomology group of $M$
with coefficients in the sheaf $\mathscr{F}_\eta(M)$, is isomorphic to the Lichnerowicz
cohomology group $H^k_{-\eta}(M)$ for all $k=0, \ldots, m$ (Proposition 3.1 of \cite{VaiRem}).
Applying Theorem \ref{twistedDeRham}, we see that the $\eta$-twisted Morse cohomology group
$H_k((C^\ast(f;e^\eta),\delta_\ast^\eta))$ is also isomorphic to $H^k(M;\mathscr{F}_\eta(M))$ 
for all $k=0, \ldots , m$.

\smallskip  
Now, if $X$ is a paracompact Hausdorff space that is locally path connected and semilocally 
$1$-connected, then there is a bijection between equivalence classes of local coefficient 
systems on $X$ and equivalence classes of locally constant sheaves on $X$, cf. Proposition 2.5.1
of \cite{DimShe} or Exercises 6.F of \cite{SpaAlg}. The above cohomology isomorphisms suggest that
the equivalence class of the sheaf $\mathscr{F}_\eta(M)$ should correspond to the equivalence 
class of the local system $e^\eta$ under this bijection. To prove this, first note that the
fact that $\mathscr{F}_\eta(M)$ is locally constant is implicit in Section 3 of \cite{VaiRem};
we reformulate this explicitly in the following lemma.

\begin{lemma}[\cite{VaiRem}]
If $M$ is a finite dimensional smooth manifold and $U$ is an open connected subset of $M$
such that $\eta|_U$ is exact, then
$$
\mathscr{F}_\eta(M)|_U \approx \mathbb{R}|_U.
$$
Moreover, if $\eta|_U$ is not exact, then the only section of $\mathscr{F}_\eta(M)|_U$ is the zero section.
\end{lemma}

\proofstart
Suppose that $f:U \rightarrow \mathbb{R}$ is a differentiable function such that
$$
d_{-\eta} f = df - \eta|_U f = 0.
$$
If $\eta|_U$ is exact, then $-\eta|_U = dh/h$ for some positive smooth function $h:U \rightarrow
\mathbb{R}$ and
$$
d_{-\eta} f = df +\frac{dh}{h} f = \frac{1}{h} d(hf) = 0.
$$
Therefore, $f = C/h$ for some $C \in \mathbb{R}$. 

If $f(x) \neq 0$ for all $x \in U$, then the first equation shows that $\eta|_U = df/f = d (\ln f)$,
i.e. $\eta|_U$ is exact. Thus, if $\eta|_U$ is not exact, then $f(x_0) = 0$ for some $x_0 \in U$, which
implies that $f \equiv 0$ on any open neighborhood $V$ of $x_0$ where $\eta|_V$ is exact. For any other
point $x_1 \in U$ we can cover a path from $x_0$ to $x_1$ with finitely many open sets $\{V_j\}_{j=1}^n$
such that $\eta|_{V_j}$ is exact for all $j=1, \ldots, n$ and conclude that $f(x_1) = 0$. Therefore, 
the only section of $\mathscr{F}_\eta(M)|_U$ is the zero section when $\eta|_U$ is not exact.
\proofend

Following Exercise 6.F.1 of \cite{SpaAlg}, we now recall that the local system $e^\eta$ on $M$ determines
a presheaf $\mathscr{F}_{e^\eta}$ on $M$ as follows: for any open set $U \subseteq M$, let 
$\mathscr{F}_{e^\eta}(U)$ be the set of all functions $f$ assigning to each $x \in U$ an element 
$f(x) \in \mathbb{R}_x$ with the property that for any continuous path $\gamma:[0,1] \rightarrow U$ 
we have 
$$
f(\gamma(1)) = e^{\int_0^1 \gamma^\ast(\eta)}f(\gamma(0)).
$$
According to Exercise 6.F.1 of \cite{SpaAlg}, the presheaf $\mathscr{F}_{e^\eta}$ satisfies the two
additional axioms necessary to be a sheaf, cf. Section 6.7.5 of \cite{SpaAlg}. Moreover, if $\eta|_U$
is exact, then $-\eta|_U = d(\ln h)$ for some positive smooth function $h:U \rightarrow \mathbb{R}$
and we have
$$
f(\gamma(1)) = \frac{h(0)f(\gamma(0))}{h(1)}.
$$
That is, $f = C/h$ for $C =h(0)f(\gamma(0))$ on any connected open set $U$ such that $-\eta|_U = d(\ln h)$,
and  hence $\mathscr{F}_{e^\eta}(U) \approx \mathbb{R}|_U \approx \mathscr{F}_\eta(M)|_U$. 
If $\eta|_U$ is not exact, then $e^{\int_0^1 \gamma^\ast(\eta)}$ can have different values on
different paths in $U$, and hence $\mathscr{F}_{e^\eta}(U) \approx \{0\}$. This shows that
$\mathscr{F}_{e^\eta}$ is isomorphic to $\mathscr{F}_\eta(M)$, cf. Proposition 5.8 of \cite{WarFou}. 

Continuing with Exercise 6.F.7 of \cite{SpaAlg}, we see that the \v{C}ech cohomology of $M$ with 
coefficients in the sheaf $\mathscr{F}_\eta(M)$ is isomorphic to the singular cohomology of $M$
with coefficients in the local system $e^\eta$, i.e
$$
H^k(M;\mathscr{F}_\eta(M)) \approx H^k(M;e^\eta)
$$
for all $k=0, \ldots, m$. To summarize, we have the following isomorphisms between
$\eta$-twisted Morse cohomology, Lichnerowicz cohomology, sheaf cohomology, and singular
cohomology with coefficients in the local system $e^\eta$.

\begin{theorem}\label{sheafiso}
Let $f:M \rightarrow \mathbb{R}$ be a smooth Morse-Smale function on a closed finite
dimensional smooth Riemannian manifold $(M,\mathsf{g})$, and assume that the unstable manifolds 
determine a regular CW-structure on $M$. For any $\eta \in \Omega^1_{\text{cl}} (M,\mathbb{R})$,
$$
H_k((C^\ast(f;e^\eta),\delta_\ast^\eta)) \approx H^k_{-\eta}(M) \approx H^k(M;\mathscr{F}_\eta(M))
\approx H^k(M;e^\eta)
$$
for all $k=0,\ldots ,m$.
\end{theorem} 

\begin{remark}\label{removeregular}
The arguments in Section \ref{continuation} can be used to remove the assumption that the 
unstable manifolds of $(f,\mathsf{g})$ determine a regular CW-structure, and the techniques used in 
Section \ref{singularCW} can be used to prove directly that the $\eta$-twisted Morse cohomology 
groups are isomorphic to the singular cohomology groups with coefficients in $e^\eta$. 
\end{remark}


\section{Applications and computations}\label{App}

In this section we collect some applications of homology and cohomology with local coefficients 
and examine them within the context of twisted Morse homology and cohomology.


\subsection{Parallel $1$-forms and Lichnerowicz cohomology computations}\label{ParLich}

Let $\eta$ be a closed 1-form on a finite dimensional smooth Riemannian manifold
$(M,\mathsf{g})$ of dimension $m$. For any $\xi \in \Omega^k(M,\mathbb{R})$ let
\begin{eqnarray*}
d_\eta \xi & = & d\xi + \eta\wedge \xi\\
\delta \xi & = & (-1)^{mk+m+1} \star(d(\star\xi))\\
i_U(\xi)   & = & (-1)^{mk+m} \star(d(\star\xi) + \eta \wedge \star \xi)\\
\delta_\eta\xi & = & (\delta + i_U) \xi,
\end{eqnarray*}
where $\star$ is the Hodge star operator. Using the Hodge decomposition 
$$
\Omega^k(M,\mathbb{R}) = \mathcal{H}^k_\eta(M) \oplus d_\eta(\Omega^{k-1}(M)) 
\oplus \delta_\eta(\Omega^{k+1}(M))
$$
where
$$
\mathcal{H}^k_\eta(M) = \{\xi \in \Omega^k(M,\mathbb{R})|\ d_\eta\xi = 0 \text{ and }
\delta_\eta \xi = 0 \}
$$
Le\'on, L\'opez, Marrero, and Padr\'on proved the following.
\begin{theorem}[Theorem 4.5 of \cite{LeoOnt}]\label{Leo4.5}
Let $M$ be a compact differentiable manifold and $\eta$ a closed 1-form on $M$ with $\eta \neq 0$.
Supposed that $\mathsf{g}$ is a Riemannian metric on $M$ such that $\eta$ is parallel with respect
to $\mathsf{g}$. Then, the cohomology $H^\ast_\eta(M)$ is trivial.
\end{theorem}

Combining the $\eta$-Twisted Morse de Rham Theorem (Theorem \ref{twistedDeRham})
with this result gives the following.

\begin{theorem}\label{parallelMorse}
Let $f:M \rightarrow \mathbb{R}$ be a smooth Morse-Smale function on a closed finite
dimensional smooth Riemannian manifold $(M,\mathsf{g})$. If $\eta$ is a nonzero closed 1-form
on $M$ that is parallel with respect to $\mathsf{g}$, then 
$$
H_k((C^\ast(f;e^\eta),\delta_\ast^\eta)) \approx 0
$$
for all $k$.
\end{theorem}

\proofstart
The 1-form $\eta$ is parallel with respect to $\mathsf{g}$ if any only if $-\eta$ is parallel 
with respect to $\mathsf{g}$. Thus, Theorem \ref{twistedDeRham}, Remark \ref{removeregular},
and Theorem \ref{Leo4.5}, imply that
$$
H_k((C^\ast(f;e^\eta),\delta_\ast^\eta)) \approx H^k_{-\eta}(M) = 0
$$
for all $k$.
\proofend

\begin{corollary}[Parallel $1$-Form Obstruction]\label{notparallel}
Let $f:M \rightarrow \mathbb{R}$ be a smooth Morse-Smale function on a closed finite
dimensional smooth Riemannian manifold $M$, and assume there exists a nonzero closed
1-form $\eta$ on $M$ such that $H_k((C^\ast(f;e^\eta),\delta_\ast^\eta)) \neq 0$ for some
$k$. Then for any nonzero closed 1-form $\zeta$ on $M$ such that $[\zeta] = [\eta] \in 
H^1(M;\mathbb{R})$ the 1-form $\zeta$ is not parallel with respect to any Riemannian metric
on $M$.
\end{corollary}

\proofstart
By Theorem \ref{Lichcoclass}, Theorem \ref{twistedDeRham}, and Remark \ref{removeregular}
we see that 
$$
H_k((C^\ast(f;e^\eta),\delta_\ast^\eta)) \neq 0 \quad \Rightarrow \quad
H_k((C^\ast(f;e^\zeta),\delta_\ast^\zeta)) \neq 0.
$$
By Theorem \ref{twistedDeRham} and Remark \ref{removeregular} these homology groups are independent
of the Riemannian metric used to define them, and hence the result follows from 
Theorem \ref{parallelMorse}.
\proofend

\begin{corollary}[Euler Number Parallel $1$-Form Obstruction]\label{EulPar}
If  $M$ is a closed finite dimensional smooth manifold with nonzero Euler number, then there are no nonzero 
closed parallel $1$-forms on $M$ with respect to any Riemannian metric on $M$.
\end{corollary}

\proofstart
Corollary \ref{coEulerindex} implies that for any nonzero  closed 1-form $\eta \in \Omega_{\text{cl}}(M,\mathbb{R})$, 
$\mathcal{X}^\eta(M) = \mathcal{X}(M) \neq 0$. Therefore, $H^\ast_\eta(M) \neq 0$, and Theorem \ref{Leo4.5} implies
that  $\eta$ is not parallel with respect to any Riemannian metric on $M$.
\proofend

\begin{example}[The Lichnerowicz cohomology of a cylinder]\label{cylinder}
Let $(N,\mathsf{g})$ be a compact smooth Riemannian manifold and let $M = N \times S^1$. 
Define $\eta = \pi_2^\ast(d\theta)$ to be the 1-form on $M$ given by 
pulling back the closed 1-form $d\theta$ from Example \ref{circle} under 
the projection $\pi_2:M \times S^1 \rightarrow S^1$ onto the second component. 
The 1-form $\eta$ is parallel on $M$ with respect to the Riemannian metric
$$
\widetilde{\mathsf{g}}  =  \pi_1^\ast(\mathsf{g}) + \eta \otimes \eta,
$$
where $\pi_1:M \times S^1 \rightarrow M$ denotes projection onto the first factor. 
Therefore, by Theorem \ref{parallelMorse} we see that for any Morse-Smale 
function $f:M \rightarrow \mathbb{R}$ the $\eta$-twisted Morse cohomology groups
$$ 
H_k((C^\ast(f;e^\eta),\delta_\ast^\eta)) \approx 0
$$
for all $k$.

\smallskip
Now assume in addition that $(N,\mathsf{g},\lambda)$ is a contact manifold, and
define
$$
\Omega = \pi_1^\ast(d\lambda) + \pi_2^\ast(d\theta) \wedge \pi_1^\ast(\lambda).
$$
Then $\Omega$ is a non-degenerate 2-form on $M$ such that
\begin{eqnarray*}
d\Omega & = & \pi_1^\ast(d^2 \lambda) + \pi_2^\ast(d(d\theta)) \wedge \pi_1^\ast(\lambda) -
\pi_2^\ast(d\theta) \wedge \pi_1^\ast(d\lambda) \\
& = & -\pi_2^\ast(d\theta) \wedge \pi_1^\ast(d\lambda)\\
& = & -\pi_2^\ast(d\theta) \wedge \Omega.
\end{eqnarray*}
Hence, $M = N \times S^1$ is a locally conformal symplectic (LCS) manifold 
with Lee form $\eta = \pi_2^\ast(d\theta)$, cf. Proposition \ref{LCSLee}.
Since $\eta$ is parallel with respect to the Riemannian metric $\widetilde{\mathsf{g}}$ on $M$,
this gives a class of LCS manifolds $(M,\Omega)$ with Lee form $\eta$ such that the
$\eta$-twisted Morse cohomology groups $H_k((C^\ast(f;e^\eta),\delta^\eta_\ast))$ are zero
for all $k$, cf. Example 4.10(2) of \cite{LeoOnt}.
\end{example}

\begin{remark}
In the previous example the $\eta$-twisted Euler number
$$
\mathcal{X}^\eta(M) = \sum_{k=0}^{m} (-1)^k \text{dim }H^k_{\eta}(M) = 0.
$$
Thus, Corollary \ref{coEulerindex} implies that the untwisted Euler number $\mathcal{X}(M)$ is also zero.
The fact that $\mathcal{X}(N \times S^1) = 0$ can also be seen directly from the K\"unneth Theorem,
cf. Theorem VI.3.2 of \cite{BreTop}, or the fact that $\mathcal{X}(N \times S^1) = \mathcal{X}(N)
\mathcal{X}(S^1)$. 
\end{remark}

\begin{example}[The Lichnerowicz cohomology of a surface of genus 2]\label{genus2Lich}
Consider a surface $S$ of genus 2 as the connected sum of two tori. If we view each tori as
a square whose opposite sides are identified, then the surface of genus 2 can be viewed as
an octagon where every other edge on each half of the octagon is identified, cf. Section I.5 of
\cite{MasABa}. By applying the construction used in the proof of Theorem \ref{regularMorse} to the octagon 
representing $S$ we can construct a Morse-Smale pair $(f,\mathsf{g})$ with one critical point
of index 0, four critical points of index 1, and one critical point of index 2. The critical point 
$p^0$ of index 0 will correspond to the vertices of the octagon, which are all identified, the four
critical points $p^1_1,p^1_2,p^1_3,p^1_4$, of index 1 will be on the edges of the octagon, which 
are identified in pairs, and the critical point $p^2$ of index 2 will be on the interior of the octagon.   
\begin{figure}[h]
\includegraphics{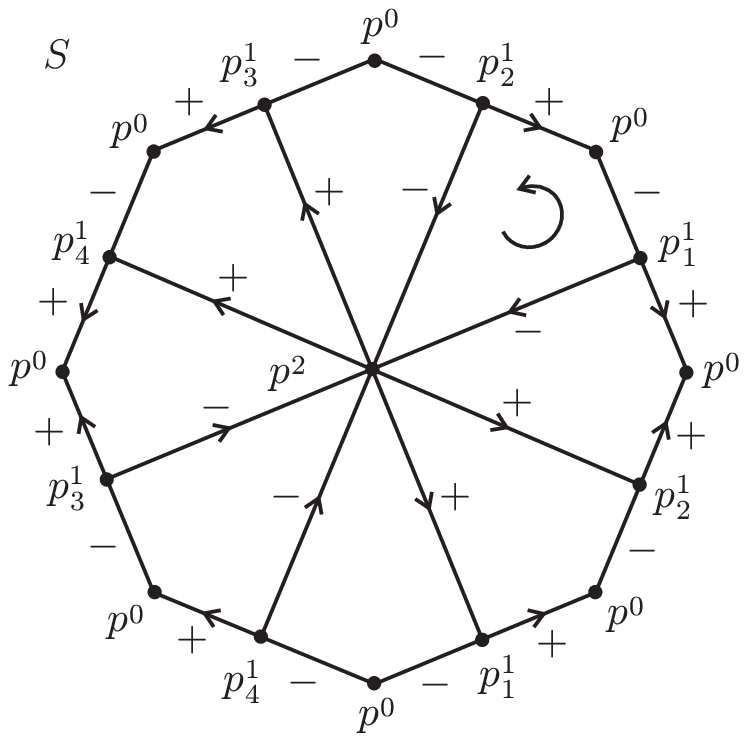}
\end{figure}

As shown in the diagram, there are exactly two gradient flow lines between each pair of critical
points of relative index 1, and they have opposite signs. Hence, the (untwisted) coboundary operators
$\delta_\ast$ are all zero in the Morse-Smale-Witten cochain complex.
$$
\xymatrix{
0 \ar[r] & C^0(f;\mathbb{R}) \ar[r]^-{\delta_0} & C^1(f;\mathbb{R}) \ar[r]^-{\delta_1} & 
C^2(f;\mathbb{R}) \ar[r]^-{\delta_2} & 0\\ 
}
$$
Therefore,
$$
H_k((C^\ast(f;\mathbb{R}),\delta_\ast)) \approx 
\left\{
\begin{array}{ll}
\mathbb{R}                                                       & \text{if } k=0\\
\mathbb{R} \oplus \mathbb{R} \oplus \mathbb{R} \oplus \mathbb{R} & \text{if } k=1\\
\mathbb{R}                                                       & \text{if } k=2,
\end{array}
\right.
$$
and Theorem \ref{Lichcoclass}, Theorem \ref{twistedDeRham}, and Remark \ref{removeregular}
imply that $H^1_{\text{dR}}(S;\mathbb{R}) \approx \mathbb{R}^4$. 
Moreover, there are smooth paths $\gamma_i$, such that $\text{Im}(\gamma_i) = \overline{W^u(p^1_i)}$
for $i=1,2,3,4$, that represent a basis for $H_1(S;\mathbb{R})$. Hence, 
$H^1_{\text{dR}}(S;\mathbb{R})$ has a basis represented by closed 1-forms 
$\eta_j\in \Omega_{\text{cl}}(M,\mathbb{R})$ for $j=1,2,3,4$ such that
$$
\int_{\gamma_i} \eta_j = 
\left\{ \begin{array}{cl}
1 & \text{if } i=j\\
0 & \text{if } i \neq j,
\end{array}\right.
$$
since the identification of $H^1_{\text{dR}}(S;\mathbb{R})$ with $\text{Hom}(H_1(S;\mathbb{R}),
\mathbb{R})$ is given by integrating a representative of the de Rham cohomology
class over a smooth representative of the homology class, cf. Sections V.5 and V.7 of \cite{BreTop}.

\smallskip
To compute $H_k((C^\ast(f;e^{\eta_j}),\delta_\ast^{\eta_j})) \approx H_{-\eta_j}^k(S)$ we need to
consider the integrals of $\eta_j$ over the gradient flow lines from $p^2$ to the critical points
of index 1. To address this we apply Stokes' Theorem for smooth singular chains, 
cf. Theorem 8.4 of \cite{SpiDif1} or Theorem 4.7 of \cite{WarFou}. 
For instance, consider the integral of $\eta_1$ over the boundary of a 
smooth singular 2-cube $\sigma:[0,1] \times [0,1] \rightarrow S$ whose four faces map injectively onto
the following smoothly embedded loops and closed intervals in the counterclockwise direction
shown in the diagram.
\begin{enumerate}
\item The loop $\overline{W^u(p_1^1)}$.
\item The closure of the image of the positive gradient flow line from $p_2^1$ to $p^0$.
\item The loop $\overline{W(p^2,p_2^1)}$.
\item The closure of the image of the positive gradient flow line from $p_2^1$ to $p^0$.
\end{enumerate}
\begin{figure}[h]
\includegraphics{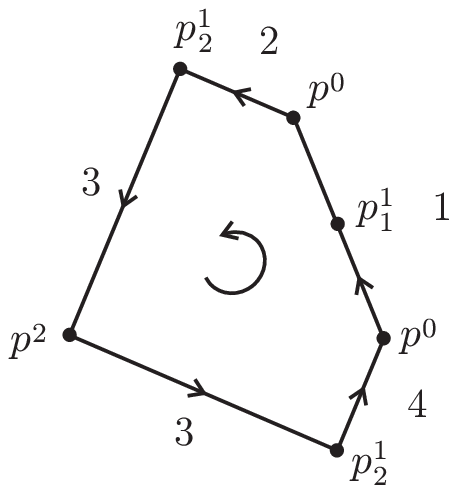}
\end{figure}
Since
$$
\int_{\partial \sigma} \eta_1 = \int_\sigma d\eta_1 = 0,
$$
and the integrals over the closure of the image of the positive gradient flow line from $p_2^1$ 
to $p^0$ cancel out, we see that the integrals of $\eta_1$ over the loops $\overline{W^u(p_1^1)}$ 
and $\overline{W(p^2,p^1_2)}$ are opposite each other. Similarly, the integrals of $\eta_1$ over 
the loops $\overline{W(p^2,p^1_1)}$, $\overline{W(p^2,p^1_3)}$, and $\overline{W(p^2,p^1_4)}$ 
are all zero, since the integrals of $\eta_1$ over the loops $\overline{W^u(p^1_2)}$,
$\overline{W^u(p^1_4)}$, and $\overline{W^u(p^1_3)}$ are zero.
This gives enough information to compute $H_k((C^\ast(f;e^{\eta_1}),\delta_\ast^{\eta_1}))
\approx H_{-\eta_1}^k(S)$ for all $k$.
$$
\xymatrix{
0 \ar[r] & C^0(f;e^{\eta_1}) \ar[r]^-{\delta^{\eta_1}_0} & C^1(f;e^{\eta_1}) 
\ar[r]^-{\delta^{\eta_1}_1} & C^2(f;e^{\eta_1}) \ar[r]^-{\delta^{\eta_1}_2} & 0\\ 
}
$$

To see this, let $\theta \in C^0(f;e^{\eta_1})$, and recall from Definition \ref{etacochain} that
$$
(\delta_0^{\eta_1}\theta)(p^1_i)\ \ = \sum_{\nu \in \mathcal{M}(p^1_i,p^0)}
\epsilon(\nu) \exp{\left(\int_{\overline{\mathbb{R}}} (\gamma^\nu)^\ast(\eta_1)
\right)}\theta(p^0) \in e^{\eta_1}_{p^1_i}, 
$$
for any critical point $p^1_i \in Cr_1(f)$, where $\gamma^\nu:\overline{\mathbb{R}}
\rightarrow M$ is any continuous path from $p^0$ to $p^1_i$ whose image coincides with the
image of $\nu \in \mathcal{M}(p^1_i,p^0)$ and $\epsilon(\nu) = \pm 1$ is the sign determined by the
orientation on $\mathcal{M}(p^1_i,p^0)$. For $i \neq 1$, the integral of $\eta_1$ over the loop
$\overline{W^u(p^1_i)}$ is zero, and hence the two integrals of $\eta_1$ in the above sum are the same.
Thus, $(\delta_0^{\eta_1}\theta)(p^1_i) = 0$ for $i \neq 1$, since the two gradient flow lines 
from $p^1_i$ to $p^0$ have opposite signs. However, $(\delta_0^{\eta_1}\theta)(p^1_1) \neq 0$ 
for $\theta \neq 0$ because the two integrals have different values. Specifically, if the integral 
of $\eta_1$ along the positive gradient flow line is $a$, then the integral of $\eta_1$ along the
negative gradient flow line will be $a \pm 1$, where the sign is determined by whatever orientation
is chosen for the generator $\gamma_1$ of $H_1(S;\mathbb{R})$. Therefore, 
$$
(\delta_0^{\eta_1}\theta)(p^1_1) = (e^a - e^{a\pm 1}) \theta(p^0) \neq 0
$$
if $\theta(p^0) \neq 0$, and
$$
H_0((C^\ast(f;e^{\eta_1}),\delta_\ast^{\eta_1})) \approx 0.
$$
Moreover, if $\theta_j$ is a basis for $C^1(f;e^{\eta_1})$ such that
$$
\theta_j(p^1_i) = 
\left\{ \begin{array}{cl}
1 & \text{if } i=j\\
0 & \text{if } i \neq j
\end{array}\right.
$$
for $j=1,2,3,4$, then the above argument shows that $\theta_1$ generates
$\text{Im }\delta_0^{\eta_1}$.

Now consider $\delta_1^{\eta_1}(\theta_j)$ for $j=1,2,3,4$. Using the notation from 
Definition \ref{etacochain} we have
\begin{eqnarray*}
(\delta_1^{\eta_1}\theta_j)(p^2)\ \ 
& = & \sum_{p^1_i \in Cr_{1}(f)} \sum_{\nu \in \mathcal{M}(p^2,p^1_i)}
\epsilon(\nu) \exp{\left(\int_{\overline{\mathbb{R}}} (\gamma^\nu)^\ast(\eta_1)
\right)}\theta_j(p^1_i)\\ 
& = & \sum_{\nu \in \mathcal{M}(p^2,p^1_j)} \epsilon(\nu) 
\exp{\left(\int_{\overline{\mathbb{R}}} (\gamma^\nu)^\ast(\eta_1)
\right)} 1 \in e^{\eta_1}_{p^2}.
\end{eqnarray*} 
Since the integrals of $\eta_1$ over the loops $\overline{W(p^2,p^1_1)}$, 
$\overline{W(p^2,p^1_3)}$, and $\overline{W(p^2,p^1_4)}$ are zero and
the two gradient flow lines from $p^2$ to $p^1_i$ have opposite signs, we see that 
$\theta_1,\theta_3, \theta_4 \in \text{ker }\delta_1^{\eta_1}$. However, the integral of $\eta_1$ 
over the loop $\overline{W(p^2,p^1_2)}$ is $\pm 1$, where the sign is determined by whatever 
orientation is chosen for the generator $\gamma_1$ of $H_1(S;\mathbb{R})$. Therefore, 
$\delta_1^{\eta_1}(\theta_2) \neq 0$, and we have
$$
H_{-\eta_1}^k(S) \approx H_k((C^\ast(f;e^{\eta_1}),\delta_\ast^{\eta_1})) \approx 
\left\{
\begin{array}{ll}
0                             & \text{if } k=0\\
\mathbb{R} \oplus \mathbb{R}  & \text{if } k=1\\
0                             & \text{if } k=2.
\end{array}
\right.
$$

It's clear that the Lichnerowicz cohomology $H_{-\eta_j}^k(S) \approx H_k((C^\ast(f;e^{\eta_j}),
\delta_\ast^{\eta_j}))$ is the same for $j=2,3,4$. The Lichnerowicz cohomology $H_{-\eta}^k(S)$
for a general closed $1$-form $\eta \in \Omega_{\text{cl}}^1(S,\mathbb{R})$ can be computed
using the invariance of the $\eta$-twisted Euler number (Corollary \ref{coEulerindex}), i.e.
$\mathcal{X}^\eta(S) = \mathcal{X}(S) = -2$, and the fact that $\{[\eta_1],[\eta_2],[\eta_3],[\eta_4]\}$
is a basis for $H^1_{\text{dR}}(S;\mathbb{R})$. That is, Corollary \ref{coEulerindex} and the placement
of the 6 generators of $(C^\ast(f;e^{\eta}),\delta_\ast^{\eta})$ imply that there are only 4 possibilities
for $H_{-\eta}^\ast(S)$. We claim for any $\eta \in \Omega_{\text{cl}}^1(S,\mathbb{R})$ with 
$[\eta] \neq 0$, $H_{-\eta}^\ast(S)$ is the same as $H_{-\eta_1}^\ast(S)$.

To see this, let $\eta = a_1 \eta_1 + a_2 \eta_2 + a_3 \eta_3 + a_4 \eta_4$ for some 
$a_1,a_2,a_3,a_4 \in \mathbb{R}$ and assume that $a_j \neq 0$ for some $j$, i.e. $\eta$ is not exact.
Letting $\gamma^+$ and $\gamma^-$ be paths from $p^0$ to $p^1_j$ parameterizing the elements of
$\mathcal{M}(p^1_j,p^0)$ with signs $+$ and $-$ respectively we compute as follows. For any 
$\theta \in C^0(f;e^\eta)$, 
\begin{eqnarray*}
(\delta_0^{\eta}\theta)(p^1_j)\ \ & = &\sum_{\nu \in \mathcal{M}(p^1_j,p^0)}
\epsilon(\nu) \exp{\left(\int_{\overline{\mathbb{R}}} (\gamma^\nu)^\ast(\eta)
\right)}\theta(p^0) \\
& = & \left( e^{\int_{\gamma^+} \eta} - e^{\int_{\gamma^-} \eta} \right) \theta(p^0)\\
& = & \left( e^{a_1\int_{\gamma^+} \eta_1}e^{a_2\int_{\gamma^+}\eta_2}e^{a_3\int_{\gamma^+} \eta_3}
e^{a_4\int_{\gamma^+}\eta_4} - 
e^{a_1\int_{\gamma^-} \eta_1}e^{a_2\int_{\gamma^-}\eta_2}e^{a_3\int_{\gamma^-} \eta_3}
e^{a_4\int_{\gamma^-}\eta_4} \right) \theta(p^0)\\
& = & e^Ae^Be^C(e^{a_ja} - e^{a_j(a \pm 1)})\theta(p^0) 
\end{eqnarray*}
for some $a,A,B,C \in \mathbb{R}$, because
$$
\int_{\gamma^-} \eta_i\ = \int_{\gamma^+} \eta_i\ \text{ if }i \neq j \quad \text{and} \quad
\int_{\gamma^-} \eta_j = a \pm 1\ \text{ if } \int_{\gamma^+} \eta_j = a.
$$
This shows that $H_0((C^\ast(f;e^{\eta}),\delta_\ast^{\eta})) \approx 0$.

Now consider the basis element $\theta_l \in C^1(f;e^\eta)$, where $\overline{W(p^2,p^1_l)}$ is on 
the opposite side of the singular 2-cube from $\overline{W^u(p^1_j)}$, and let $\gamma^+$ and
$\gamma^-$ be paths from $p^1_l$ to $p^2$ parameterizing the elements of $\mathcal{M}(p^2,p^1_l)$
with signs $+$ and $-$ respectively.
\begin{figure}[h]
\includegraphics{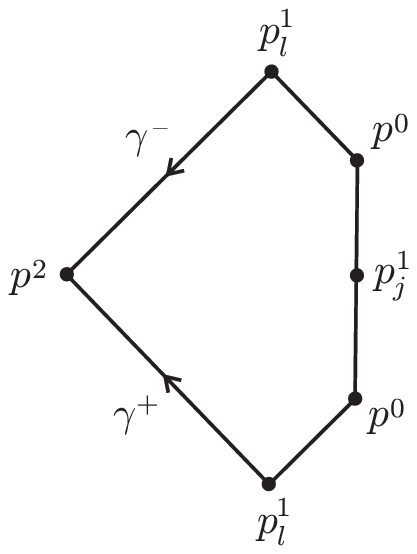}
\end{figure}

\begin{eqnarray*}
(\delta_1^{\eta}\theta_l)(p^2)\ \ 
& = & \sum_{p^1_i \in Cr_{1}(f)} \sum_{\nu \in \mathcal{M}(p^2,p^1_i)}
\epsilon(\nu) \exp{\left(\int_{\overline{\mathbb{R}}} (\gamma^\nu)^\ast(\eta)
\right)}\theta_l(p^1_i)\\ 
& = & \sum_{\nu \in \mathcal{M}(p^2,p^1_l)} \epsilon(\nu) 
\exp{\left(\int_{\overline{\mathbb{R}}} (\gamma^\nu)^\ast(\eta)
\right)} 1\\
& = & \left( e^{a_1\int_{\gamma^+} \eta_1}e^{a_2\int_{\gamma^+}\eta_2}e^{a_3\int_{\gamma^+} \eta_3}
e^{a_4\int_{\gamma^+}\eta_4} - 
e^{a_1\int_{\gamma^-} \eta_1}e^{a_2\int_{\gamma^-}\eta_2}e^{a_3\int_{\gamma^-} \eta_3}
e^{a_4\int_{\gamma^-}\eta_4} \right)\\
& = & e^Ae^Be^C(e^{a_ja} - e^{a_j(a \pm 1)}) 
\end{eqnarray*}
for some $a,A,B,C \in \mathbb{R}$, because
$$
\int_{\gamma^-} \eta_i\ = \int_{\gamma^+} \eta_i\ \text{ if }i \neq j \quad \text{and} \quad
\int_{\gamma^-} \eta_j = a \pm 1\ \text{ if } \int_{\gamma^+} \eta_j = a.
$$
Thus $\delta_1^{\eta}\theta_l \neq 0$ when $a_j \neq 0$, and
$H_2((C^\ast(f;e^{\eta}),\delta_\ast^{\eta})) \approx 0$.
Therefore, 
$$
H_{-\eta}^k(S) \approx H_k((C^\ast(f;e^{\eta}),\delta_\ast^{\eta})) \approx 
\left\{
\begin{array}{ll}
0                             & \text{if } k=0\\
\mathbb{R} \oplus \mathbb{R}  & \text{if } k=1\\
0                             & \text{if } k=2
\end{array}
\right.
$$
for any $\eta \in \Omega_{\text{cl}}^1(S,\mathbb{R})$ with $[\eta] \neq 0$,
because $\mathcal{X}^\eta(S) = -2$.
\end{example}


\subsection{H-spaces}\label{Hspace}

Recall that an H-space is a topological space $X$ together with a continuous map $m:X \times X 
\rightarrow X$ and an element $e \in X$ such that $m(e,\cdot):X \rightarrow X$ and $m(\cdot,e):X 
\rightarrow X$ are homotopic to the identity through maps that preserve $e$. The map $m$ is called 
\textbf{multiplication} and the element $e$ is called the \textbf{homotopy unit}. This definition can
be weakened by removing the assumption that the homotopies preserve $e$ or strengthened by requiring that
$e$ be a strict identity, and all three definitions are equivalent for CW-complexes, 
cf. Section 3.C of \cite{HatAlg}. Standard examples of H-spaces include the based loop space 
of a topological space, where the multiplication map is given by concatenation,
and topological groups, i.e. topological spaces with a group structure such that both the multiplication 
and inverse maps are continuous. For both of these examples the multiplication map $m$ is associative.

We will only consider H-spaces where $m$ is associative up to homotopy because we are interested in
the {\it Pontryagin ring}. For any commutative ring $R$ with unit the multiplication map $m$ on an 
H-space $X$ induces what is known as the {\bf Pontryagin product} on the homology:
$$
\xymatrix{
H_\ast(X;R)  \otimes H_\ast(X;R) \ar[r]^-{\times} & H_\ast(X \times X; R) \ar[r]^-{m_\ast} & H_\ast(X;R),
}
$$
where the first map is the homology cross product.  The assumption that $m$ is associative up to 
homotopy implies that the Pontryagin product is associative, and hence $H_\ast(X;R)$ is a graded
ring with unit known as the {\bf Pontryagin ring}, cf. Section VII.3 of \cite{DolLec} or Section III.7
of \cite{WhiEle}.

\smallskip
Recently, Albers, Frauenfelder, and Oancea extended the Pontryagin product to homology with local 
coefficients in a system of rank one $R$-modules $\mathcal{L}$ on a path connected space $X$, where $R$
is a commutative ring with unit \cite{AlbLoc}. That is, $\mathcal{L}$ is a bundle of abelian groups such
that each fiber has the structure of a free rank one $R$-module and the homomorphisms $\gamma_\ast$ are all 
$R$-module isomorphisms, i.e. $\gamma_\ast$ is given by multiplication by an element of the group of
invertible elements $R^\times \subset R$. Any bundle of abelian groups with fiber $\mathbb{Z}$ is an
example of a local coefficient system of rank one $\mathbb{Z}$-modules, and $e^\eta$, where $\eta$ is
a closed smooth $1$-form on a smooth manifold, is an example of a local coefficient system of rank one
$\mathbb{R}$-modules (see Remark \ref{modulebundle}).

\begin{proposition}[Proposition 1 \cite{AlbLoc}]\label{pontryagin}
Let $R$ be a commutative ring with unit and $X$ a path-connected associative H-space.
For any local coefficient system of rank one $R$-modules $\mathcal{L}$ on $X$, the multiplication
$m:X \times X \rightarrow X$ induces a unital ring structure on $H_\ast(X;\mathcal{L})$. 
The unit is represented by the class of a point.
\end{proposition}

\smallskip\noindent
The ring structure is defined at the chain level as the composition
$$
\xymatrix{
C_\ast(X;\mathcal{L})  \otimes C_\ast(X;\mathcal{L}) \ar[r]^-{B} & 
C_\ast(X \times X; pr_1^\ast\mathcal{L} \otimes_R pr_2^\ast\mathcal{L}) \ar[r]^-{m_\ast} & 
C_\ast(X;\mathcal{L}),
}
$$
where $B$ is the Eilenberg-MacLane shuffle map with local coefficients in $\mathcal{L}$,
$pr_i: X \times X \rightarrow X$ is a projection map for $i=1,2$, and $m_\ast$ is induced by the
multiplication $m$. The fact that $\mathcal{L}$ is a local system of \textit{rank one} $R$ modules is central to
the proof of Proposition \ref{pontryagin}, since this implies that there is an isomorphism
$$
pr_1^\ast\mathcal{L} \otimes_R pr_2^\ast\mathcal{L} \cong m^\ast \mathcal{L}.
$$
This fact is proved in Lemma 1 of \cite{AlbLoc} using the natural isomorphism $R \otimes_R R \cong
R$ given by multiplication in $R$.

\smallskip
Albers, Frauenfelder, and Oancea used Proposition \ref{pontryagin} to prove a vanishing result for
the homology of an H-space with local coefficients in certain systems of rank one $R$-modules 
(Proposition 3 of \cite{AlbLoc}).

\begin{proposition}\label{singularvanish}
Let $R$ be a commutative ring with unit and $\mathcal{L}$ a local coefficient system of rank one 
$R$-modules on a path connected associative H-space $X$ with homotopy unit $e \in X$. Assume
there exists $[\gamma] \in \pi_1(X,e)$ such that
\begin{enumerate}
\item $\gamma_\ast(s) = g s$ for some $g \in R^\times_e $ with $g \neq 1$, i.e. $\gamma_\ast \neq id$, and
\item $g^{-1} - 1 \in R^\times_e$, i.e. $g^{-1} - 1$ is invertible.
\end{enumerate}
Then $H_\ast(X;\mathcal{L}) = 0$.  
\end{proposition}

\proofstart 
If $e \in X$ is the homotopy unit, then Proposition \ref{pontryagin} says that 
$[1\cdot e] \in H_0(X;\mathcal{L})$ is the unit in the Pontryagin ring $H_\ast(X;\mathcal{L})$,
where $1 \in R_e$ (see Section \ref{singularlocal}). Thus, it suffices to show that (1) and (2)
imply that $[1 \cdot e] = [0] \in H_\ast(X;\mathcal{L})$.

If $\gamma:[0,1] \rightarrow X$ is a loop based at $e$, then $\gamma_\ast(s) = g s$
for some $g \in R^\times_e$, and for any $s \in R_e$ we have $s \cdot \gamma\in C_1(X;\mathcal{L})$
with
$$
\partial_1(s \cdot \gamma) = \gamma^{-1}_\ast(s) \cdot e - s \cdot e = 
\left( g^{-1}s - s \right)\cdot e
$$
(see Definition \ref{twistedsingular}).
In particular, for $s=1 \in R_e$ we have $\partial_1(1 \cdot \gamma) = \left( g^{-1} - 1 \right)\cdot e$.
Hence, if $[\gamma] \in \pi_1(M,e)$ satisfies conditions (1) and (2), we can set 
$s = \left(g^{-1} - 1\right)^{-1}$ in the above equation for $\partial_1(s \cdot \gamma)$ and we have
\begin{eqnarray*}
\partial_1\left((g^{-1} - 1)^{-1} \cdot \gamma \right) 
& = & \left(g^{-1} (g^{-1} - 1)^{-1} - (g^{-1} - 1)^{-1}\right) \cdot e\\
& = & \left((g^{-1} - 1) (g^{-1} - 1)^{-1}\right) \cdot e\\
& = & 1 \cdot e.
\end{eqnarray*}
Therefore, $[1 \cdot e] = [0] \in H_\ast(X;\mathcal{L})$.
\proofend

\smallskip\noindent
\textbf{Note:} If $R = \mathbb{Z}$, then conditions (1) and (2) cannot hold simultaneously because
$\gamma_\ast(s) = \pm s$ when $R_e = \mathbb{Z}$ and $-2 \not\in \mathbb{Z}^\times$. Similarly,
if $R = \mathbb{Z}_2$, then condition (1) cannot hold.
However, if $R$ is a field and $\mathcal{L}$ is not simple, then there always exists some 
$[\gamma] \in \pi_1(X,e)$ such that (1) and (2) hold. This observation leads to the following
obstruction to the existence of an associative H-space structure.

\begin{corollary}\label{singularHspaceobstruction}
If $X$ is a path connected topological space and there exists a local coefficient system $\mathcal{L}$
of rank one vector spaces on $X$, i.e. a line bundle over some field, such that
\begin{enumerate}
\item $\mathcal{L}$ is not simple, i.e. $\mathcal{L}$ is not isomorphic to a constant bundle, and
\item $H_\ast(X;\mathcal{L}) \neq 0$,
\end{enumerate}
then $X$ is not an associative H-space.
\end{corollary}

\smallskip
We now interpret and apply the above results within the context of twisted Morse homology.

\begin{theorem}\label{HspaceMorse}
Let $f:M \rightarrow \mathbb{R}$ be a smooth Morse-Smale function on a closed
path connected finite dimensional smooth Riemannian manifold $(M,\mathsf{g})$, and let $\mathcal{L}$ be a
local coefficient system of rank one $R$-modules on $M$, where $R$ is a commutative ring with
unit.  Assume that $M$ is an associative H-space with homotopy unit $e \in M$, and 
there exists $[\gamma] \in \pi_1(X,e)$ such that
\begin{enumerate}
\item $\gamma_\ast(s) = g s$ for some $g \in R^\times_e $ with $g \neq 1$, i.e. $\gamma_\ast \neq id$, and
\item $g^{-1} - 1 \in R^\times_e$, i.e. $g^{-1} - 1$ is invertible.
\end{enumerate}
Then $H_k((C_\ast(f;\mathcal{L}),\partial_\ast^\mathcal{L})) = 0$ for all $k=0,\ldots ,m$.
\end{theorem}

\proofstart
Note that if we ``forget'' the product structure, then $\mathcal{L}$ is a bundle of abelian groups.
So, the twisted Morse-Smale-Witten chain complex $(C_\ast(f;\mathcal{L}),\partial_\ast^{\mathcal{L}})$
is defined, and by the Twisted Morse Homology Theorem (Theorem \ref{twistedsame}) the homology groups of
$(C_\ast(f;\mathcal{L}),\partial_\ast^{\mathcal{L}})$ are isomorphic to the singular homology groups of $M$
with coefficients in the bundle of abelian groups $\mathcal{L}$. Moreover, the singular homology
groups of $M$ with coefficients in the bundle of abelian groups obtained by ``forgetting'' the product
structure on $\mathcal{L}$ are isomorphic to the singular homology groups of $M$ with coefficients
in the $R$-module $\mathcal{L}$. The result now follows from Proposition \ref{singularvanish}.
\proofend

\begin{corollary}[Associative H-space Obstruction]\label{Hspaceobstruction}
Let $f:M \rightarrow \mathbb{R}$ be a smooth Morse-Smale function on a closed
path connected finite dimensional smooth Riemannian manifold $(M,\mathsf{g})$.  If there exists
a local coefficient system $\mathcal{L}$ of rank one vector spaces on $M$, i.e. a line bundle
over some field, such that
\begin{enumerate}
\item $\mathcal{L}$ is not simple, i.e. $\mathcal{L}$ is not isomorphic to a constant bundle, and
\item $H_k((C_\ast(f;\mathcal{L}),\partial_\ast^{\mathcal{L}})) \neq 0$ for some $k$,
\end{enumerate}
then $M$ is not an associative H-space.
\end{corollary}

\begin{corollary}[Euler Number Associative H-space Obstruction]\label{Eulerobstruction}
If $M$ is a closed path connected finite dimensional smooth manifold with
$H^1_{\text{dR}}(M;\mathbb{R}) \neq 0$ and $\mathcal{X}(M) \neq 0$, then
$M$ is not an associative H-space.
\end{corollary}

\proofstart
Pick a closed 1-form $\eta\in \Omega_{\text{cl}}(M,\mathbb{R})$ such that the de Rham cohomology
class $[\eta] \neq 0$. Then $e^\eta$ is a non-simple local coefficient system of rank one vector spaces
on $M$. Moreover, by Theorem \ref{regularMorse} and Corollary \ref{Eulerclassical} we have that
$H_\ast(M;e^\eta) \neq 0$, since $\mathcal{X}_{e^\eta}(M) = \mathcal{X}(M) \neq 0$. The result now
follows from Corollary \ref{singularHspaceobstruction}.
\proofend

\smallskip\noindent
\textbf{Note:} Since both based loop spaces and topological groups are associative
H-spaces, Corollary \ref{Hspaceobstruction} and Corollary \ref{Eulerobstruction} give 
obstructions to closed finite dimensional smooth manifolds being based loop spaces or 
topological groups, cf. \cite{BauFin}.

\begin{example}[The twisted homology of a circle]\label{circleHspace}
The manifold $S^1$ is a topological group, and hence an associative H-space.  Thus, Theorem 
\ref{HspaceMorse} implies that the twisted Morse homology of $S^1$ vanishes for any non-simple
local coefficient system $\mathcal{L}$ of rank one vector spaces on $S^1$.
\begin{figure}[h]
\includegraphics{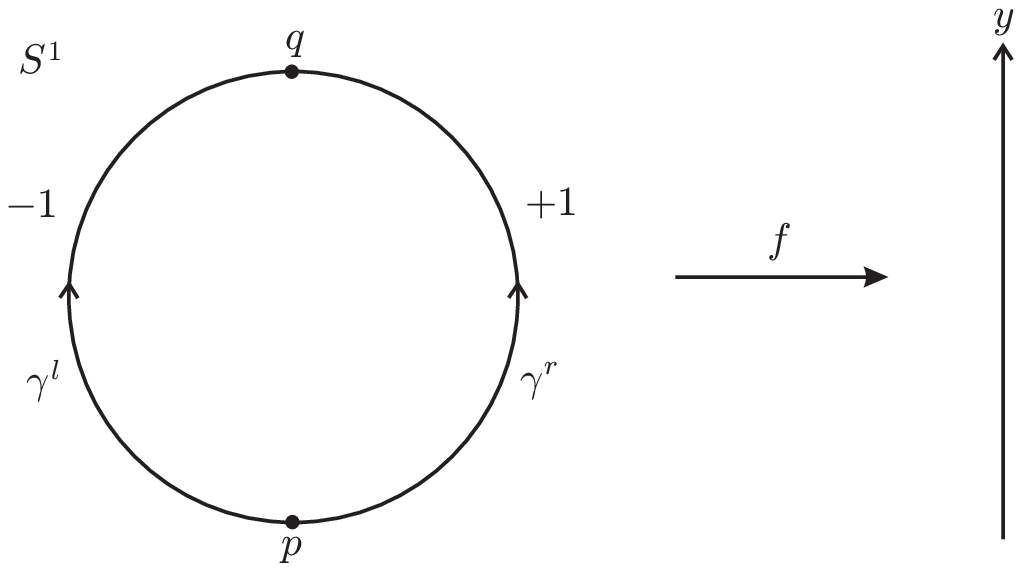}
\end{figure}

To see this directly, note that the isomorphism class of $\mathcal{L}$ is determined by
the invertible linear map $\gamma_\ast:\mathbb{F} \rightarrow \mathbb{F}$ associated to the clockwise
generator $[\gamma] \in \pi_1(S^1,q)$, i.e. an element $a \in \mathbb{F}^\ast$ where 
$\gamma_\ast(x) = ax$ and $\mathbb{F}$ is the ground field for the vector space.
Choosing the constant path at $q$ and the path from $q$ to $p$ on the right half of $S^1$ to define
a local system associated to the representation, we have
\begin{eqnarray*}
\partial_1^{\mathcal{L}}(xq) & = & \left(\gamma^r_\ast(x) - \gamma^l_\ast(x)\right) p\\
                      & = & \left(x-ax\right)p.
\end{eqnarray*}
Therefore, $H_k((C_\ast(f;\mathcal{L}),\partial_\ast^{\mathcal{L}})) \approx  0$ for all $k$ whenever
$a \neq 1$.
\end{example}

\begin{example}[The twisted homology of a torus]\label{torusgroup}
The torus $T^n = (S^1)^n$ is a topological group. Hence, for any Morse-Smale pair $(f,\mathsf{g})$
and any local coefficient system $\mathcal{L}$ of rank one $R$-modules such that 
there exists a $[\gamma] \in \pi_1(T^n,e)$ with
$$
\gamma_\ast(s) = g s \text{ for some }g \in R^\times_e \text{ where } g^{-1} - 1 \text{ is invertible},
$$
$H_k((C_\ast(f;\mathcal{L}),\partial_\ast^{\mathcal{L}})) = 0$ for all $k=0,\ldots ,m$ by
Theorem \ref{HspaceMorse}.

For the local coefficient system of rank one vector spaces $e^\eta$ this can be seen directly 
using a computation similar to the one we present for the system of rank 
one Novikov modules on $T^2$ in Example \ref{Novtorus}. Note that this result is consistent with 
Example \ref{cylinder}; although Example \ref{cylinder} was about twisted cohomology rather 
than twisted homology.
\end{example}

\begin{example}[A surface of genus $2$ is not an associative H-space]\label{genus2notH}
In Example \ref{genus2Lich} we saw that the twisted Morse cohomology 
$H_{-\eta_1}^\ast(S) \approx H_k((C^\ast(f;e^{\eta_1}),\delta_\ast^{\eta_1}))$
is nonzero, where $e^{\eta_1}$ is a non-trivial local coefficient system of real
line bundles on a surface $S$ of genus $2$. This suggests that the twisted Morse
homology with coefficients in $e^{\eta_1}$ should also be nonzero. However, the standard
Universal Coefficient Theorems, which compare homology and cohomology do not apply
to homology and cohomology with local coefficients, cf. Section V.7 of \cite{BreTop}.

So, consider the $\eta_1$-twisted Morse-Smale-Witten chain complex
$$
\xymatrix{
0 \ar[r]^-{\partial^{\eta_1}_3} & C_2(f;e^{\eta_1}) \ar[r]^-{\partial^{\eta_1}_2} & C_1(f;e^{\eta_1}) 
\ar[r]^-{\partial^{\eta_1}_1} & C_0(f;e^{\eta_1}) \ar[r] & 0\\ 
}
$$
for the function and metric studied in Example \ref{genus2Lich}. Using the notation
from Example \ref{genus2Lich}, $C_2(f;e^{\eta_1})$ consists of real multiples of $p^2$, 
$C_1(f;e^{\eta_1})$ consists of real linear combinations of $p^1_1$, $p^1_2$, $p^1_3$,
and $p^1_4$, and $C_0(f;e^{\eta_1})$ consists of real multiples of $p^0$. Moreover,
$\partial_2^{\eta_1}(p^2)$ is a nonzero multiple of $p^1_2$, the kernel of
$\partial_1^{\eta_1}$ consists of real linear combinations of $p^1_2$, $p^1_3$, and $p^1_4$,
and $\partial_1^{\eta_1}(p^1_1)$ is a nonzero multiple of $p^0$.
Therefore,
$$
H_k((C_\ast(f;e^{\eta_1}),\partial_\ast^{\eta_1})) \approx 
\left\{
\begin{array}{ll}
0                             & \text{if } k=0\\
\mathbb{R} \oplus \mathbb{R}  & \text{if } k=1\\
0                             & \text{if } k=2,
\end{array}
\right.
$$
and Corollary \ref{Hspaceobstruction} implies that a surface of genus 2 is not an associative H-space.
Note that this result also follows from Corollary \ref{Eulerobstruction}.
\end{example}

\begin{example}[$\mathbb{R}P^{2n}$ is not an associative H-space]\label{RP2notH}
Consider $M = \mathbb{R}P^2$ with the Morse-Smale function and orientations from Example
\ref{projectivespace}.
\begin{figure}[h]
\includegraphics{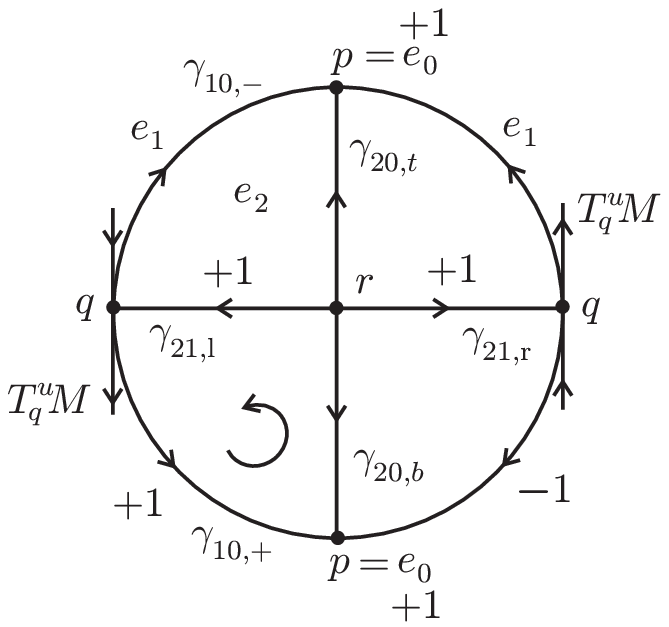}
\end{figure}
Since $H^1(\mathbb{R}P^2;\mathbb{R})  =  0$, the local coefficient system $e^\eta$ is simple
for every $\eta \in \Omega_{\text{cl}}^1(\mathbb{R}P^2,\mathbb{R})$. However, we can
construct a non-simple local coefficient system $\mathcal{L}$ of rank one $\mathbb{R}$ vector spaces
on $\mathbb{R}P^2$ using homomorphisms analogous to those defined in Example \ref{projectivespace}, where
the non-trivial element $[\gamma] \in \pi_1(\mathbb{R}P^2,r)$ corresponds to the linear map
$\gamma_\ast:\mathbb{R} \rightarrow \mathbb{R}$ given by $\gamma_\ast(x) = -x$.
\begin{eqnarray*}
(\gamma_{12,r})_\ast(x) & = & +x\\
(\gamma_{12,l})_\ast(x) & = & -x\\
(\gamma_{01,+})_\ast(x) & = & +x\\
(\gamma_{01,-})_\ast(x) & = & -x\\
\end{eqnarray*}
With this non-simple local system we have
$$
H_2((C_\ast(f;\mathcal{L}),\partial_\ast^{\mathcal{L}})) \approx \mathbb{R}, \quad
H_1((C_\ast(f;\mathcal{L}),\partial_\ast^{\mathcal{L}})) \approx 0, \quad
H_0((C_\ast(f;\mathcal{L}),\partial_\ast^{\mathcal{L}}))  \approx 0.
$$
Therefore, Corollary \ref{Hspaceobstruction} implies that $\mathbb{R}P^2$ is not an associative H-space.

\medskip
A similar argument can be used to show that $\mathbb{R}P^n$ is not an associative H-space when
$n$ is even, but the obstruction gives no information when $n$ is odd. That is, the function
$\tilde{f}:S^n \rightarrow \mathbb{R}$ defined by 
$$
\tilde{f}(x_1,\ldots, x_n) = \sum_{j=2}^n (j-1)x_j^2
$$
is a Morse-Smale function with respect to the standard Riemannian metric on $S^n$, and it satisfies
$\tilde{f}(-x_1,\ldots, -x_n) = \tilde{f}(x_1,\ldots, x_n)$. So, $\tilde{f}$ induces a Morse-Smale
function $f:\mathbb{R}P^n \rightarrow \mathbb{R}$, which has one critical point $p_k$ of index $k$
for all $k=0,\ldots , n$. There are exactly two gradient flow lines between any two critical points of
$p_k$ and $p_{k-1}$ of relative index one, and the signs associated to these two gradient flow lines
will be the same when $k$ is even and opposite when $k$ is odd. 

We can define a non-simple local coefficient system $\mathcal{L}$ of rank one $\mathbb{R}$ vector spaces
on $\mathbb{R}P^n$ as above, where for all $k=1, \ldots ,n$ we have $(\gamma_1)_\ast(x) = x$ for one
of the gradient flow lines from $p_k$ to $p_{k-1}$ and $(\gamma_2)_\ast(x) = -x$ for the other gradient
flow line from $p_k$ to $p_{k-1}$. Then the Morse-Smale-Witten chain complex with coefficients in 
$\mathcal{L}$ is 
$$
\xymatrix{
\mathbb{R} \ar[r]^{\partial_n} & \mathbb{R} \ar[r]^{\partial_{n-1}} & \cdots \ar[r]^{\partial_2} 
& \mathbb{R} \ar[r]^{\partial_1}  & \mathbb{R} \ar[r] & 0,
}
$$
where $\partial_k(x) = 0$ when $k$ is even and $\partial_k(x) = 2x$ when $k$ is odd.
Therefore, 
$$
H_k((C_\ast(f;\mathcal{L}),\partial_\ast^{\mathcal{L}})) \approx 0
$$
for all $k$ when $n$ is odd, but 
$$
H_n((C_\ast(f;\mathcal{L}),\partial_\ast^{\mathcal{L}})) \approx \mathbb{R}
$$
when $n$ is even. Thus, Corollary \ref{Hspaceobstruction} implies that $\mathbb{R}P^{2n}$ is not
an associative H-space.

\smallskip\noindent
{\bf Note:} The real projective space $\mathbb{R}P^n$ is an H-space only if $n=$ 1, 3, or 7, cf. Example
3C.3 of \cite{HatAlg}.
\end{example}


\subsection{Novikov homology}\label{NovHomology}

The construction of the (twisted or untwisted) Morse-Smale-Witten chain complex only required an exact 
1-form $df\in \Omega^1(M,\mathbb{R})$, rather than a function $f:M \rightarrow \mathbb{R}$. That is, 
the critical points of $f$ are the zeros of $df$, and the gradient flow lines of $f$ are determined by 
$df$ and the Riemannian metric on $M$. This observation naturally leads one to ask if it is possible
to construct a chain complex analogous to the Morse-Smale-Witten chain complex using a closed Morse 1-form
$\zeta \in \Omega^1_{\text{cl}}(M,\mathbb{R})$ in place of the exact Morse 1-form $df$, or at least obtain
inequalities for $\zeta$ analogous to the Morse inequalities for $df$. These questions were first addressed
by S.P. Novikov \cite{NovMul} \cite{NovThe}, and later by several other authors, including 
M. Farber \cite{FarMor} \cite{FarTop}, M. Farber and A. Ranicki \cite{FarThe}, A. Pajitnov \cite{PajCir}
\cite{PazOnt}, A. Ranicki \cite{RanThe}, and D. Sch\"utz \cite{SchOne}, \cite{SchGeo}. In this
section we discuss how the homology groups associated to a closed 1-form, i.e. Novikov homology,
can be viewed as a special case of twisted Morse homology. 

\subsubsection{A covering space associated to a 1-form}
Novikov's approach to constructing a chain complex associated to a closed 1-form $\zeta$ on a 
connected smooth manifold $M$ involves pulling back $\zeta$ to a covering space 
$\widetilde{M}_\xi$ of $M$ where the pullback is exact. More specifically, if we pick 
a basepoint $x_0$ for $M$, then a closed 1-form $\zeta \in \Omega^1_{\text{cl}}(M,\mathbb{R})$
defines a homomorphism of periods $\text{Per}_\xi:\pi_1(M, x_0) \rightarrow \mathbb{R}$ given by
$$
\text{Per}_\xi([\gamma]) = \int_\gamma \zeta,
$$ 
which only depends on the de Rham cohomology class $\xi = [\zeta] \in H^1(M;\mathbb{R})$. 
The kernel $\Delta_\xi \stackrel{\text{def}}{=} \text{ker Per}_\xi$ is a subgroup of $\pi_1(M,x_0)$, 
and hence it determines a covering space $\pi:\widetilde{M}_\xi \rightarrow M$ such that 
$\pi_{\#}: \pi_1(\widetilde{M}_\xi,\tilde{x}_0) \rightarrow\Delta_\xi \subseteq \pi_1(M,x_0)$
is an isomorphism, where $\tilde{x}_0 \in \widetilde{M}_\xi$ is any basepoint satisfying 
$\pi(\tilde{x}_0) = x_0$. In fact, $\pi_1(M,x_0)$ acts on the universal cover of $M$
via deck transformations $\pi_1(M,x_0) \times \widetilde{M} \rightarrow \widetilde{M}$ and 
$\widetilde{M}_\xi = \widetilde{M}/\Delta_\xi$, cf. Theorem III.8.1 of \cite{BreTop}. 
Moreover, $\widetilde{M}_\xi$ is a regular covering space of $M$
because $\text{ker Per}_\xi$  is a normal subgroup of $\pi_1(M,x_0)$, 
and the pullback $\tilde{\zeta} = \pi^\ast(\zeta)$ is exact on $\widetilde{M}_\xi$ 
because the function $\tilde{f}:\widetilde{M}_\xi \rightarrow \mathbb{R}$ given by
$$
\tilde{f}(\tilde{x}) = \int_{\tilde{\gamma}} \tilde{\zeta},
$$
where $\tilde{\gamma}$ is any path from $\tilde{x}_0$ to $\tilde{x}$, is well-defined and satisfies
$d\tilde{f} = \tilde{\zeta}$.

\smallskip\noindent
{\bf Note:} If $\zeta$ is an {\bf integer} valued 1-form, i.e. if $\xi = [\zeta] \in H^1(M;\mathbb{Z})$,
then there is a well-defined circle valued function $f:M \rightarrow S^1$ given by
$$
f(x) = e^{2\pi i \tilde{f}(\tilde{x})},
$$
where $\tilde{x}$ is any element of the fiber $\pi^{-1}(x)$. Moreover, every smooth circle valued function
$f:M \rightarrow S^1$ determines a corresponding smooth integer valued closed 1-form by pulling back 
$\frac{1}{2\pi} d\theta$ (see Example \ref{circle}), where $[\frac{1}{2\pi} d\theta]$ is a generator for 
$H^1(S^1;\mathbb{Z})$. Thus, circle valued Morse theory can be viewed as a special case of 
Novikov's theory for closed 1-forms, cf. Lemma 2.1 of \cite{FarTop}, Section III of \cite{NovMul}, 
or Section 2.1.4 of \cite{PajCir}. 

\begin{definition}
The \textbf{rank} of a cohomology class $\xi \in H^1(M;\mathbb{R})$ (or a closed 1-form)
is defined to be the rank of the group of periods of $\xi$, i.e. the rank of the finitely generated abelian
group $\Gamma_\xi \stackrel{\text{def}}{=} \text{Im } \text{Per}_\xi \subseteq \mathbb{R}$. 
\end{definition}

Since $H_1(M;\mathbb{Z})$ is isomorphic to $\pi_1(M,x_0)$ modulo its commutator subgroup and $\mathbb{R}$
is commutative, the period homomorphism $\text{Per}_\xi$ factors through $H_1(M;\mathbb{Z})$. That is,
there is a homomorphism $\overline{\text{Per}}_\xi:H_1(M;\mathbb{Z}) \rightarrow \mathbb{R}$ that makes
the following diagram commute.
$$
\xymatrix{
\pi_1(M,x_0) \ar[rr]^-{\text{Per}_\xi} \ar[dr] & & \mathbb{R}\\ 
& H_1(M;\mathbb{Z}) \ar[ur]_-{\overline{\text{Per}}_\xi}
}
$$
Hence, the rank of any cohomology class $\xi \in H^1(M;\mathbb{R})$  is bounded above by the
rank of $H_1(M;\mathbb{Z})$, i.e. the number of elements in a basis for $H_1(M;\mathbb{Z})$ modulo
its torsion subgroup.

\begin{remark}\label{rankonedense}
If $\xi$ is a \textbf{rank one} cohomology class, then we can divide $\xi$ by the 
smallest positive element in $\Gamma_\xi$ to get an integral cohomology class. Thus, rank one cohomology
classes are real multiples of integral cohomology classes. Any rational cohomology class 
$\xi \in H^1(M;\mathbb{Q})$ is of rank one, because $\text{Im }\text{Per}_\xi$ is generated by a 
finite set of rational numbers $\frac{a_1}{b_1},\frac{a_2}{b_2}, \ldots ,\frac{a_r}{b_r}$, which is 
a subgroup of the cyclic group generated by $\frac{1}{b_1b_2\cdots b_r}$. Thus, rank one cohomology
classes are dense in $H^1(M;\mathbb{R})$, cf. Theorem 1.44 or Corollary 2.2  of \cite{FarTop}. 
Moreover, the Novikov homology of a general cohomology class $\xi = [\zeta] \in H^1(M;\mathbb{R})$ 
can be constructed by perturbing a Morse form $\zeta$ to a rank one Morse form that agrees with $\zeta$ 
in a neighborhood of its critical points, cf. Section 4 of \cite{SchGeo}.
\end{remark} 

We now observe that the group of periods $\Gamma_\xi$ is the quotient of $\pi_1(M,x_0)$ by the normal
subgroup $\Delta_\xi$, i.e.
$$
\pi_1(M,x_0)/\Delta_\xi \approx \Gamma_\xi,
$$
since $\Delta_\xi$ and $\Gamma_\xi$ are the kernel and image of the homomorphism
$\text{Per}_\xi:\pi_1(M,x_0) \rightarrow \mathbb{R}$. Moreover, $\pi_1(M,x_0)/\Delta_\xi$
is isomorphic to the deck transformation group of $\widetilde{M}_\xi$, cf. Corollary III.6.9 
of \cite{BreTop} or Proposition 1.39 of \cite{HatAlg}, and hence the deck transformation group
$\text{Aut}(\pi)$ of $\pi:\widetilde{M}_\xi \rightarrow M$ is isomorphic to $\Gamma_\xi$. The
homomorphism $\text{Aut}(\pi) \rightarrow \Gamma_\xi$ can be described as follows. 
Starting with a deck transformation $D\in \text{Aut}(\pi)$, take any path
$\tilde{\gamma}$ from the basepoint $\tilde{x}_0 \in \widetilde{M}_\xi$ to $D \tilde{x}_0$. 
The image of $\tilde{\gamma}$ under the projection $\pi: \widetilde{M}_\xi \rightarrow M$ maps
to a loop $\gamma$ in $M$ based at $x_0$, and $\text{Per}_\xi([\gamma]) \in \Gamma_\xi$, which
independent of the choice of the path $\tilde{\gamma}$.

\smallskip
This leads to the following, cf. Section 14.6 of \cite{SteThe}.

\begin{claim}\label{periodfiber}
The fiber and the deck transformation group of the regular covering space 
$\pi:\widetilde{M}_\xi \rightarrow M$ can be identified with the group of periods $\Gamma_\xi$.
$$
\xymatrix{
\Gamma_\xi \ar[r] & \widetilde{M}_\xi \ar[d]^\pi & \\
                                    & \widetilde{M}_\xi/\Gamma_\xi \ar@{<->}[r]^-\approx & M.
}
$$
\end{claim}

\proofstart
Recall that $\widetilde{M}_\xi$ is a regular cover because $\pi_\#(\pi_1(\widetilde{M}_\xi, \tilde{x}_0))
= \Delta_\xi$ is the kernel of $\text{Per}_\xi$, and hence a normal subgroup of $\pi_1(M,x_0)$.
Thus, the group of deck transformations of $\widetilde{M}_\xi$ acts transitively on the fibers 
of $\pi:\widetilde{M}_\xi \rightarrow M$.  Moreover, a deck transformation is determined by where
it maps the basepoint of $\widetilde{M}_\xi$, cf. Lemma III.4.4 of \cite{BreTop}, and hence there is
a bijection between the group of deck transformations and the fiber containing the basepoint.
The fact that $\Gamma_\xi$ is isomorphic to the deck transformation group of $\widetilde{M}_\xi$
completes the proof.
\proofend

\begin{corollary}\label{rank1cyclic}
If $\xi \in H^1(M;\mathbb{R})$ is a rank one cohomology class, then $\pi:\widetilde{M}_\xi \rightarrow M$ 
is a cyclic covering space. 
$$
\xymatrix{
\mathbb{Z} \ar[r] & \widetilde{M}_\xi \ar[d]^\pi \\
                                  & M
}
$$
\end{corollary}


\subsubsection{Novikov rings}
If $M$ is compact, then every Morse form $\zeta\in \Omega_{\text{cl}}^1(M,\mathbb{R})$ will have finitely
many zeros. However, the pullback $\tilde{\zeta}$ to $\widetilde{M}_\xi$ of any closed Morse 1-form $\zeta$
with $[\zeta] = \xi \neq 0$ will have an infinite number of zeros, i.e. the Morse function 
$\tilde{f}:\widetilde{M}_\xi \rightarrow \mathbb{R}$ will have infinitely many critical points,
since the fiber $\Gamma_\xi$ of $\widetilde{M}_\xi$ is infinite. 
To account for this Novikov introduced a ring $\text{Nov}(\Gamma_\xi)$, where $\Gamma_\xi$ is the 
group of periods of $\zeta$. 

\begin{definition}\label{Novikovring}
Let $\Gamma \subseteq \mathbb{R}$ be an additive subgroup. The {\bf Novikov ring} $\text{Nov}(\Gamma)$
is defined to be the ring of all formal power series of the form
$$
\sum_{\gamma \in \Gamma} n_\gamma t^\gamma
$$
with integer coefficients $n_\gamma \in \mathbb{Z}$, such that at most countably many $n_\gamma \neq 0$,
and 
$$
\text{for any } c \in \mathbb{R} \text{ the set } \{\gamma \in \Gamma |\ n_\gamma \neq 0 \text{ and }
\gamma > c\} \text{ is finite.}
$$ 
\end{definition}

\smallskip\noindent
{\bf Note:} $\text{Nov}(\Gamma)$ is a commutative ring consisting of ``half infinite'' formal series
with integer coefficients and exponents in $\Gamma\subseteq \mathbb{R}$. That is, every element 
$g \in \text{Nov}(\Gamma)$ can be written as $g = \sum_{i=1}^\infty n_i t^{\gamma_i}$, where 
$n_i \in \mathbb{Z}$, $\gamma_i \in \Gamma$, and $\gamma_1 > \gamma_2 > \gamma_3 > \cdots$. 
Addition and multiplication in $\text{Nov}(\Gamma)$ are defined formally, cf. Section 1.2.1 of \cite{FarTop}.
An element of $\text{Nov}(\Gamma)$ is invertible if and only if its top coefficient is invertible in 
$\mathbb{Z}$, and $\text{Nov}(\Gamma)$ is a principal ideal domain, cf. Lemmas 1.9 and 1.10 of
\cite{FarTop} or Theorems 4.1 and 4.2 of \cite{HofFlo}. The ring $\text{Nov}(\mathbb{R})$ is denoted
by $\text{Nov}$.

\smallskip\noindent
{\bf Note:} If $\xi \in H^1(M;\mathbb{R})$ is a rank one cohomology class, then 
$\Gamma_\xi \approx \mathbb{Z}$ and $\text{Nov}(\Gamma_\xi)$ is isomorphic to the Novikov ring 
$\text{Nov}(\mathbb{Z})$.  That is, the ring with elements of the form
$$
g = \sum_{i \in \mathbb{Z}} n_i t^i,
$$
where only finitely many $n_i$ with $ i > 0$ are nonzero.


\subsubsection{A local coefficient system of rank one $\text{Nov}$-modules}
If $\gamma$ is a closed loop based at $x_0$, then $[\gamma]\in H_1(M;\mathbb{R})$ and
$$
<\xi,[\gamma]> = \int_\gamma \zeta,
$$
where $\xi = [\zeta] \in H^1(M;\mathbb{R})$. Moreover, $t^{<\xi,[\gamma]>} \in \text{Nov}(\Gamma_\xi)$
is invertible with inverse $t^{<\xi,[\gamma^{-1}]>} \in \text{Nov}(\Gamma_\xi)$. Thus, the 
cohomology class $\xi \in H^1(M;\mathbb{R})$ determines a representation
$$
\pi_1(M,x_0)  \times \text{Nov}(\Gamma_\xi) \rightarrow \text{Nov}(\Gamma_\xi)
$$
defined by $(\gamma,g) \mapsto t^{<\xi,[\gamma]>}g$, and hence a bundle $\mathcal{L}(\Gamma_\xi)$ 
of rank one $\text{Nov}(\Gamma_\xi)$-modules on $M$ (defined up to isomorphism).

\smallskip
If we extend the above representation to a representation on the Novikov ring $\text{Nov} 
\stackrel{\text{def}}{=} \text{Nov}(\mathbb{R})$, then we can give the following explicit description 
of a bundle $\mathcal{L}_\zeta$ in the isomorphism class $\mathcal{L}_\xi$
of bundles of rank one $\text{Nov}$-modules determined by $\xi = [\zeta] \in H^1(M;\mathbb{R})$. 

\begin{definition}\label{zetasystem}
For any $\zeta \in \Omega^1_{\text{cl}}(M,\mathbb{R})$, the local coefficient system
$\mathcal{L}_\zeta$ of rank one $\text{Nov}$-modules is defined as follows. The fiber at every
point is $\text{Nov}$, and for every path $\gamma:[0,1] \rightarrow M$, the isomorphism 
$\gamma_\ast:\text{Nov}_{\gamma(1)} \rightarrow \text{Nov}_{\gamma(0)}$ is defined to be 
multiplication by $t^{\int_0^1 \gamma^\ast(\zeta)}$, i.e.
$$
g \mapsto t^{\int_0^1 \gamma^\ast(\zeta)}\ g.
$$
\end{definition}

\begin{claim}\label{zetarep}
Let $\zeta \in \Omega^1_{\text{cl}}(M,\mathbb{R})$ with $[\zeta] = \xi \in H^1(M;\mathbb{R})$.
The local coefficient system $\mathcal{L}_\zeta$ is in the isomorphism class $\mathcal{L}_\xi$
of local coefficient systems determined by the representation
$$
\pi_1(M,x_0)  \times \text{Nov} \rightarrow \text{Nov}
$$
given by $(\gamma,g) \mapsto t^{<\xi,[\gamma]>}g$. 
In particular, if $\zeta_1, \zeta_2 \in \Omega^1_{\text{cl}}(M,\mathbb{R})$ represent the 
same de Rham cohomology class, then $\mathcal{L}_{\zeta_1}$ is isomorphic to $\mathcal{L}_{\zeta_2}$.
\end{claim}

\proofstart
The local coefficient system associated to the representation is defined by assigning $\text{Nov}$
as the fiber at every point in $M$ and fixing a homotopy class of paths rel endpoints from a basepoint 
$x_0 \in M$ to every point in $M$. Suppose that $\gamma_0$ represents the chosen homotopy class of 
paths from $x_0$ to $\gamma(0)$ and $\gamma_1$ represents the chosen homotopy class of paths from $x_0$
to $\gamma(1)$. The map $\gamma_\ast:\text{Nov}_{\gamma(1)} \rightarrow 
\text{Nov}_{\gamma(0)}$ is defined to be the homomorphism that the representation associates to
the concatenated loop $\gamma_0\gamma\gamma_1^{-1}$. (For more details see the proof of Theorem VI.1.12 
in \cite{WhiEle}.) 
To complete the proof, note that the following diagram commutes.
$$
\xymatrix{
\text{Nov}_{\gamma(1)} \ar[rrr]^-{\times t^{\int_0^1 \gamma^\ast(\zeta)}} 
 \ar[d]_-{\times t^{\int_0^1 \gamma_1^\ast (\zeta)}} & & & \text{Nov}_{\gamma(0)}
 \ar[d]^-{\times t^{\int_0^1 \gamma_0^\ast (\zeta)}} \\
\text{Nov}_{\gamma(1)} \ar[rrr]^-{\times t^{<[\zeta],[\gamma_0\gamma\gamma_1^{-1}]>}} 
& & & \text{Nov}_{\gamma(0)}
}
$$
\proofend

\smallskip\noindent
{\bf Note:} The direction of integration for the 1-form in Definition \ref{zetasystem} is opposite 
that for the 1-form $\eta$ in the local system $e^\eta$ defined in Example \ref{etasystem}. 
The direction of integration for $e^\eta$ was chosen so that $\eta$ is integrated in the direction
of the gradient flow lines in Definition \ref{etatwisted}, whereas the direction of integration 
for the form $\zeta$ in Definition \ref{zetasystem} was chosen to agree with that of
$<[\zeta],[\gamma_0\gamma\gamma_1^{-1}]>$.


\subsubsection{Novikov homology}

Let $\zeta$ be a closed Morse 1-form on a closed smooth Riemannian manifold $(M,\mathsf{g})$,
and let $\Gamma_\xi$ be its group of periods where $[\zeta] = \xi \in H^1_{\text{dR}}(M;\mathbb{R})$.
Since $M$ is compact, the Morse form $\zeta$ has a finite number of zeros of Morse index $k$ for all
$k=0,\ldots, m$. However, the flow associated to $(\zeta,\mathsf{g})$ isn't a gradient flow when 
$\zeta$ is not exact. Thus, the Morse index won't necessarily decrease along the flow, e.g.
there may be flow lines that begin and end at the same zero. 

The lift of the form $\zeta$ and a generic metric $\mathsf{g}$ on $M$ to the covering space 
$\pi:\widetilde{M}_\xi \rightarrow M$ associated to the kernel $\Delta_\xi$ of the period homomorphism
determines a Morse-Smale function $\tilde{f}:\widetilde{M}_\xi \rightarrow \mathbb{R}$ whose critical 
points are the fibers above the zeros of $\zeta$. In fact, the fiber above a zero of $\zeta$ of Morse 
index $k$ consists of critical points of $\tilde{f}$ of Morse index $k$. Moreover, the fiber is in 
bijective correspondence with the group of periods $\Gamma_\xi$ (Claim \ref{periodfiber}). Thus, 
there are an infinite number of critical points of index $k$ above every zero of index $k$ whenever 
$\zeta$ is not exact, i.e. when the rank of $\zeta$ is greater than zero. 

The Novikov complex retains the features of a free finitely generated chain complex of
modules over a ring and the features of a boundary operator defined by a Morse-Smale flow. 
For each $k=0,\ldots, m$, the chain complex $C_k(\zeta;\text{Nov}(\Gamma_\xi))$ is defined to be the free 
finitely generated $\text{Nov}(\Gamma_\xi)$-module on the index $k$ zeros $Z_k(\zeta)$ of $\zeta$, i.e.
$$
C_k(\zeta;\text{Nov}(\Gamma_\xi)) \stackrel{\text{def}}{=} \left.\left\{ 
\sum_{q \in Z_k(\zeta)} gq \right| g \in \text{Nov}(\Gamma_\xi) \right\} \approx 
\bigoplus_{q \in Z_k(\zeta)} \text{Nov}(\Gamma_\xi),
$$
and the boundary operator is defined intuitively using the Morse coefficients of the 
gradient flow of the Morse-Smale function $\tilde{f}:\widetilde{M}_\xi \rightarrow \mathbb{R}$ 
(or of an $\tilde{f}$-gradient  or a ``gradient-like'' flow associated to $\tilde{f})$.  However, 
when the form is of rank greater than one a flow on a related cyclic cover $\bar{M}$ or the universal
cover $\widetilde{M}$ is sometimes used instead.

More specifically, above each zero $q$ of $\zeta$ fix a critical point $\tilde{q}$ in the fiber of 
the cover. The equivariant incidence coefficients are defined to be
$$
<q:p>\ = \sum_{\gamma \in \Gamma} [\tilde{q}:\gamma \tilde{p}] t^\gamma,
$$
where $[\tilde{q}: \gamma \tilde{p}]$ denotes the number of flow lines from $\tilde{q}$ to $\gamma \tilde{p}$ 
counted with sign and we have identified $\Gamma \subseteq \mathbb{R}$ with the group of deck
transformations of whatever cover is being used (see the discussion above Claim \ref{periodfiber}).
After proving that
$$
<q:p>\ \in \text{Nov}(\Gamma_\xi),
$$
$\partial^\zeta_k:C_k(\zeta;\text{Nov}(\Gamma_\xi)) \rightarrow C_{k-1}(\zeta;\text{Nov}(\Gamma_\xi))$
is defined on an elementary chain $gq \in C_k(\zeta;\text{Nov}(\Gamma_\xi))$ by
$$
\partial_k^\zeta(gq) = \sum_{p\in Z_{k-1}(\zeta)} <q:p> g p,
$$
cf. Section 2.2 of \cite{FarTop} or Definition 4.6 of \cite{SchGeo}.

\smallskip
The Novikov Principle relates the homology of the Novikov complex to the homology
of $\widetilde{M}_\xi$ with coefficients in $\text{Nov}(\Gamma_\xi)$, cf. Section 2.2 of \cite{FarTop},
\cite{NovMul}, \cite{NovThe}, or Section 4 of \cite{SchGeo}. Proofs of the Novikov Principle 
(some in special cases) have been given by Farber and Ranicki \cite{FarThe}, Farber \cite{FarMor}, 
Latour \cite{LatExi}, Pazhitnov \cite{PazOnt}, Sch\"utz \cite{SchGeo}, and others. Several of the proofs
are for the special case of a rank one form $\zeta$ or a circle valued Morse function 
$f:M \rightarrow S^1$, where the relevant cover $\widetilde{M}_\xi$ or $\bar{M}$ is a cyclic cover
(Corollary \ref{rank1cyclic}). However, rank one cohomology classes are dense in $H^1(M;\mathbb{R})$
(Remark \ref{rankonedense}), and it is possible to prove the Novikov Principle for a general 
closed Morse 1-form by perturbing the form to a rank one form, cf. Section 4 of \cite{SchGeo}.
 
\begin{theorem}[Novikov Principle]\label{NovPrin}
Given a closed 1-form $\zeta$ with Morse zeros on a finite dimensional closed smooth manifold $M$,
let $\Gamma_\xi \subset \mathbb{R}$ be the group of periods of $\zeta$. There exists a chain complex 
$(C_\ast(\zeta;\text{Nov}(\Gamma_\xi)), \partial^\zeta_\ast)$ of modules over the Novikov ring 
$\text{Nov}(\Gamma_\xi)$ such that
\begin{enumerate}
\item $C_k(\zeta;\text{Nov}(\Gamma_\xi))$ is free and finitely generated with a basis that is in
one-to-one correspondence with the zeros of $\zeta$ of Morse index $k$ for all $k = 0,\ldots ,m$.
\item $(C_\ast(\zeta;\text{Nov}(\Gamma_\xi)), \partial^\zeta_\ast)$ is chain homotopy equivalent to 
the chain complex 
$$
\text{Nov}(\Gamma_\xi) \otimes_{\mathbb{Z}[\Gamma_\xi]} C_\ast(\widetilde{M}_\xi),
$$
where $C_\ast(\widetilde{M}_\xi)$ denotes the chain complex of cellular, simplicial, or singular, 
chains $\widetilde{M}_\xi$.
\end{enumerate}
\end{theorem}

\begin{remark}
The Novikov Principle has been proved for rings more general than $\text{Nov}(\Gamma_\xi)$.
The Novikov Principle has been proved for various completions of the group ring $\mathbb{Z}[\pi_1]$,
where $\pi_1 = \pi_1(M,x_0)$, including the rational part of the Novikov ring, cf. Section 1.3 of 
\cite{FarTop}, the noncommutative {\it Novikov-Sikorav completion} $\widehat{\mathbb{Z}\pi_{1\xi}}$, 
cf. Section 3.1.5 of \cite{FarTop}, and the noncommutative {\it universal Cohn localization} 
$\Sigma^{-1}_\xi(\mathbb{Z}[\pi_1])$, cf. \cite{CohFre} or Section 3.1.7 of \cite{FarTop}. 
The Cohn localization is used for the {\it Universal Novikov Principle}, which can be used to obtain 
the other applications of the Novikov Principle via extension of scalars.
When using more general rings, the universal cover $\widetilde{M}$ is used in place
of $\widetilde{M}_\xi$. That is, the Novikov Principle for a ring $\mathcal{R}$ says that
there is a Novikov complex consisting of free $\mathcal{R}$-modules $C^\zeta_k$, with a basis in 
bijective correspondence with the zeros of index $k$ of a closed Morse 1-form $\zeta$, that is chain
homotopy equivalent to $\mathcal{R} \otimes_{\mathbb{Z}[\pi_1]} C_\ast(\widetilde{M})$. 
\end{remark}

Now suppose that $G \triangleleft \pi_1$ is a normal subgroup, where $\pi_1 = \pi_1(M,x_0)$.
Let $\widetilde{M}_G \approx \widetilde{M}/G$ be the regular covering space of $M$ corresponding 
to $G$, and let $\Gamma \approx \pi_1/G$ be the group of covering transformations of $\widetilde{M}_G$. 

\begin{claim}\label{equivcovers}
For any ring $\mathcal{R}$ and any ring homomorphism $\rho:\mathbb{Z}[\pi_1] \rightarrow \mathcal{R}$
that factors through the group ring $\mathbb{Z}[\Gamma]$ there is a chain equivalence
$$
\mathcal{R} \otimes_{\mathbb{Z}[\pi_1]} C_k(\widetilde{M}) \simeq \mathcal{R} 
\otimes_{\mathbb{Z}[\Gamma]} C_k(\widetilde{M}_G)
$$
for all $k=0, \ldots, m$, where $\pi_1 = \pi_1(M,x_0)$.
\end{claim}

\proofstart
Since $\widetilde{M}_G = \widetilde{M}/G$, an elementary singular chain 
$\sigma:\Delta^k \rightarrow \widetilde{M}$ determines an elementary singular chain
$\pi \circ \sigma: \Delta^k \rightarrow \widetilde{M}_G$, where $\pi:\widetilde{M} \rightarrow
\widetilde{M}_G$ is the projection map. Moreover, every elementary singular chain in 
$\widetilde{M}_G$ is of the form $\pi \circ \sigma$ for some elementary singular chain
$\sigma:\Delta^k \rightarrow \widetilde{M}$, cf. Corollary III.4.2 of \cite{BreTop}.
Thus, the projection map $\pi:\widetilde{M} \rightarrow \widetilde{M}_G$ induces a surjective
chain map $\tilde{\pi}: \mathcal{R} \otimes_\mathbb{Z} C_k(\widetilde{M}) 
\rightarrow \mathcal{R} \otimes_\mathbb{Z} C_k(\widetilde{M}_G)$. 

The assumption that  $\rho$ factors through $\mathbb{Z}[\Gamma]$, where $\Gamma \approx \pi_1/G$,
means that there is a map $\bar{\rho}$ making the following diagram commute.
$$
\xymatrix{
\mathbb{Z}[\pi_1] \ar[rr]^-{\rho} \ar[dr] & & \mathcal{R}\\ 
& \mathbb{Z}[\Gamma] \ar[ur]_-{\bar{\rho}}
}
$$
The chain group $\mathcal{R} \otimes_{\mathbb{Z}[\pi_1]} C_k(\widetilde{M})$ is the quotient of
$\mathcal{R} \otimes_\mathbb{Z} C_k(\widetilde{M})$ by the subgroup generated by elements of the
form
$$
r \rho(\gamma) \otimes \sigma - r \otimes \gamma \cdot \sigma,
$$
where $\gamma \in \pi_1$, $r \in \mathcal{R}$, and $\gamma \cdot \sigma$ denotes the action of
the covering transformation determined by $\gamma$ on the elementary singular chain 
$\sigma:\Delta^k \rightarrow \widetilde{M}$. Similarly, the chain group $\mathcal{R} 
\otimes_{\mathbb{Z}[\Gamma]} C_k(\widetilde{M}_G)$ is the quotient of $\mathcal{R} 
\otimes_\mathbb{Z} C_k(\widetilde{M}_G)$ by the subgroup generated by the elements of the
form
$$
r\bar{\rho}([\gamma]) \otimes \pi\circ \sigma - r \otimes [\gamma] \cdot (\pi \circ \sigma) =
r\rho(\gamma) \otimes \pi\circ \sigma - r \otimes \pi \circ (\gamma \cdot \sigma),
$$
where $[\gamma] \in \pi_1/G \approx \Gamma$, $r \in \mathcal{R}$, and $[\gamma] \cdot (\pi \circ \sigma)$
denotes the action of the covering transformation determined by $[\gamma]$ on the elementary singular chain 
$\pi \circ \sigma:\Delta^k \rightarrow \widetilde{M}_G$. Thus, there is a chain equivalence 
$\bar{\pi}: \mathcal{R} \otimes_{\mathbb{Z}[\pi_1]} C_\ast(\widetilde{M}) \rightarrow
\mathcal{R} \otimes_{\mathbb{Z}[\Gamma]} C_\ast(\widetilde{M}_G)$ making the following diagram commute
for all $k=0,\ldots , m$.

$$
\xymatrix{
\mathcal{R} \otimes_\mathbb{Z} C_k(\widetilde{M}) \ar[r]^-{\tilde{\pi}} \ar[d] &
\mathcal{R} \otimes_\mathbb{Z} C_k(\widetilde{M}_G) \ar[d] \\
\mathcal{R} \otimes_{\mathbb{Z}[\pi_1]} C_k(\widetilde{M}) \ar[r]^-{\bar{\pi}} & 
\mathcal{R} \otimes_{\mathbb{Z}[\Gamma]} C_k(\widetilde{M}_G) 
}
$$

\proofend

The following corollary shows that the Novikov homology of a closed Morse 1-form $\zeta \in
\Omega_{\text{cl}}(M,\mathbb{R})$ is isomorphic to the homology of $M$ with local coefficients 
in the bundle $\mathcal{L}(\Gamma_\xi)$, where $[\zeta] = \xi \in H^1_{\text{dR}}(M;\mathbb{R})$.

\begin{corollary}\label{NovTwist}
If $M$ is a finite dimensional closed smooth manifold and $\xi \in H^1_{\text{dR}}(M;\mathbb{R})$, 
then the following homology groups are isomorphic for all $k=0,\ldots, m$. 
$$
H_k((C_\ast(\zeta;\text{Nov}(\Gamma_\xi)), \partial^\zeta_\ast)) \approx 
H_k((\text{Nov}(\Gamma_\xi) \otimes_{\mathbb{Z}[\pi_1]} C_k(\widetilde{M}),\bar{\partial}_\ast)) \approx
H_k(M;\mathcal{L}(\Gamma_\xi))
$$

\end{corollary}

\proofstart
By the Novikov Principle (Theorem \ref{NovPrin}) and Claim \ref{equivcovers} the
following complexes are chain homotopy equivalent
$$
(C_\ast(\zeta;\text{Nov}(\Gamma_\xi)),\partial^\zeta_\ast) \simeq 
(\text{Nov}(\Gamma_\xi) \otimes_{\mathbb{Z}[\Gamma_\xi]} C_\ast(\widetilde{M}_\xi),\bar{\partial}_\ast)
\simeq (\text{Nov}(\Gamma_\xi) \otimes_{\mathbb{Z}[\pi_1]} C_\ast(\widetilde{M}), \bar{\partial}_\ast),
$$
and by Eilenberg's Theorem (Theorem \ref{Eilenberg}) there is an isomorphism between equivariant 
homology and homology with local coefficients, i.e.
$$
H_k((\text{Nov}(\Gamma_\xi) \otimes_{\mathbb{Z}[\pi_1]} C_\ast(\widetilde{M}), \bar{\partial}_\ast))
\approx H_k(M;\mathcal{L}(\Gamma_\xi))
$$ 
for all $k=0,\ldots, m$.

\proofend

\smallskip
The preceding corollary and the Twisted Morse Homology Theorem (Theorem \ref{twistedsame}) show that
Novikov homology can be viewed as a special case of twisted Morse homology. However, the local
coefficient system $\mathcal{L}(\Gamma_\xi)$ isn't well suited for use with twisted Morse homology, 
because it's only defined up to isomorphism and the exponents of the elements in $\text{Nov}(\Gamma_\xi)$
are limited to the elements in $\Gamma_\xi = \text{Im }\text{Per}_\xi$. For twisted Morse homology it's
more convenient to use the local coefficient system $\mathcal{L}_\zeta$ from Definition \ref{zetasystem},
with fiber $\text{Nov}$ instead of $\text{Nov}(\Gamma_\xi)$, because we can then integrate 
$\zeta$ over the gradient flow lines rather than loops created by arbitrarily picking
homotopy classes of paths from a basepoint to the critical points of the Morse function,
cf. Claim \ref{zetarep}.

\smallskip
The following claim gives the relationship between $H_\ast(M;\mathcal{L}_\xi)$, 
where $\xi = [\zeta]$, and $H_\ast(M;\mathcal{L}(\Gamma_\xi))$.

\begin{claim}
If $M$ is a finite dimensional closed smooth manifold and $\xi \in H^1_{\text{dR}}(M;\mathbb{R})$,
then there is an isomorphism
$$
H_k(M;\mathcal{L}(\Gamma_\xi)) \otimes_{\text{Nov}(\Gamma_\xi)} \text{Nov} 
\approx H_k(M;\mathcal{L}_\xi)
$$
for all $k=0,\ldots, m$.
\end{claim}

\proofstart
Fix a basepoint $x_0 \in M$ and homotopy classes of paths rel endpoints from $x_0$ 
to every point in $M$. Also, pick some $\zeta \in \Omega^1_{\text{cl}}(M,\mathbb{R})$ such
that $[\zeta] = \xi \in H^1_{\text{dR}}(M;\mathbb{R})$. These choices determine a local
coefficient system $\mathcal{L}_\zeta(\Gamma_\xi)$ in the equivalence class $\mathcal{L}(\Gamma_\xi)$
and a local coefficient system $\mathcal{L}_\zeta(\mathbb{R})$ in the equivalence class
$\mathcal{L}_\xi$, 
cf. Theorem VI.1.12 of \cite{WhiEle}. The isomorphism $\gamma_\ast:\text{Nov}(\Gamma_\xi)_{\gamma(1)}
\rightarrow  \text{Nov}(\Gamma_\xi)_{\gamma(0)}$ associated to a path $\gamma:[0,1] \rightarrow M$
by $\mathcal{L}_\zeta(\Gamma_\xi)$ is the restriction of the isomorphism 
$\gamma_\ast:\text{Nov}_{\gamma(1)} \rightarrow \text{Nov}_{\gamma(0)}$ associated to the path by 
$\mathcal{L}_\zeta(\mathbb{R})$. In both cases the isomorphism is given by multiplying by 
$t^{<\xi,[\gamma_0\gamma\gamma_1^{-1}]>}$, where $\gamma_0$ and $\gamma_1$ are paths in the 
chosen homotopy classes of paths rel endpoints from the basepoint $x_0$ to $\gamma(0)$ and 
$\gamma(1)$ respectively.

Note that $\text{Nov}$ is a $\text{Nov}(\Gamma_\xi)$-module because $\text{Nov}(\Gamma_\xi)$
is a subring of $\text{Nov}(\mathbb{R}) = \text{Nov}$. Thus,
$$
\text{Nov}(\Gamma_\xi) \otimes_{\text{Nov}(\Gamma_\xi)} \text{Nov} \approx \text{Nov},
$$
cf. Theorem 5.7 of \cite{HunAlg} or Proposition 2.58 of \cite{RotAnI}. So, picking any Morse-Smale
pair $(f,\mathsf{g})$ on $M$ we have 
$$
C_k(f;\mathcal{L}_\zeta(\Gamma_\xi)) \otimes_{\text{Nov}(\Gamma_\xi)} \text{Nov} 
\approx C_k(f;\mathcal{L}_\zeta(\mathbb{R})),
$$
and hence, 
$$
H_k((C_\ast(f;\mathcal{L}_\zeta(\Gamma_\xi)) \otimes_{\text{Nov}(\Gamma_\xi)} \text{Nov},
\partial_\ast)) \approx H_k((C_\ast(f;\mathcal{L}_\zeta(\mathbb{R})), \partial_\ast))
$$
for all $k=0,\ldots ,m$, since $\partial_\ast^{\mathcal{L}_\zeta(\Gamma_\xi)}$ is simply a 
restriction of $\partial_\ast^{\mathcal{L}_\zeta(\mathbb{R})}$.

Now recall that $\text{Nov}(\Gamma_\xi)$ is a principle ideal domain that is torsion free,
cf. Lemma 1.10 and Lemma 1.12 of \cite{FarTop}. Thus, $(C_\ast(f;\mathcal{L}_\zeta(\Gamma_\xi)),
\partial_\ast)$ is a complex of flat $\text{Nov}(\Gamma_\xi)$-modules whose subcomplex of 
boundaries is also flat, cf. Corollary 3.50 of \cite{RotAnI}. Therefore, the Universal Coefficient
Theorem, cf. Theorem 7.55 of \cite{RotAnI}, gives a short exact sequence
$$
\xymatrix{
0 \ar[r] & H_k((C_\ast(f;\mathcal{L}_\zeta(\Gamma_\xi)),\partial_\ast))
\otimes_{\text{Nov}(\Gamma_\xi)} \text{Nov} \ar@{->}[r]^-{\lambda_k} &  
H_k((C_\ast(f;\mathcal{L}_\zeta(\mathbb{R}),\partial_\ast)) \ar@{-}[r] &   \\ 
\ar[r] & \text{Tor}^{\text{Nov}(\Gamma_\xi)}_1(H_{k-1}((C_\ast(f;\mathcal{L}_\zeta(\Gamma_\xi)),
\partial_\ast)), \text{Nov}) \ar[r] & 0,
}
$$
where $\lambda_k([c] \otimes g) = [c \otimes g]$, and the $\text{Tor}^{\text{Nov}(\Gamma_\xi)}_1$ term is zero
because $\text{Nov}$ is flat as a $\text{Nov}(\Gamma_\xi)$-module, cf. Lemma 1.12 of \cite{FarTop} and 
Theorem 7.2 of \cite{RotAnI}. Therefore,
$$
H_k((C_\ast(f;\mathcal{L}_\zeta(\Gamma_\xi)),\partial_\ast))
\otimes_{\text{Nov}(\Gamma_\xi)} \text{Nov} \approx 
H_k((C_\ast(f;\mathcal{L}_\zeta(\mathbb{R}),\partial_\ast)),
$$
and by the Morse Homology Theorem (Theorem \ref{twistedsame}) we have
$$
H_k(M;\mathcal{L}(\Gamma_\xi)) \otimes_{\text{Nov}(\Gamma_\xi)} \text{Nov} \approx 
H_k(M;\mathcal{L}_\xi)
$$
for all $k=0,\ldots, m$.

\proofend


\subsubsection{Novikov numbers}
The following lemma allows us to define the Novikov numbers
associated to a cohomology class $\xi \in H^1(M;\mathbb{R})$, cf. Section 1.5 of \cite{FarTop}. 

\begin{lemma}
Let $M$ be a connected closed finite dimensional smooth manifold of dimension $m$, and let
$\xi \in H^1(M;\mathbb{R})$. For all $k=0,\ldots, m$, $H_k(M;\mathcal{L}_\xi)$ is a finitely 
generated $\text{Nov}$-module, and hence it is isomorphic to a finitely generated free 
$\text{Nov}$-module and finitely many cyclic torsion modules.
\end{lemma}

\proofstart
Picking any $\zeta\in \Omega^1_{\text{cl}}(M,\mathbb{R})$ such that $[\zeta] = \xi \in H^1(M;\mathbb{R})$
gives a bundle of rank one $\text{Nov}$-modules $\mathcal{L}_\zeta$ in the isomorphism class 
$\mathcal{L}_\xi$. For any Morse-Smale pair $(f,\mathsf{g})$ the chain groups in the twisted 
Morse-Smale-Witten chain complex $(C_\ast(f;\mathcal{L}_\zeta),\partial^{\mathcal{L}_\zeta}_\ast)$ are finitely
generated $\text{Nov}$-modules, and
$$
H_k((C_\ast(f;\mathcal{L}_\zeta),\partial^{\mathcal{L}_\zeta}_\ast)) \approx H_k(M;\mathcal{L}_\xi)
$$
for all $k=0,\ldots, m$ by the Twisted Morse Homology Theorem (Theorem \ref{twistedsame}),
whose proof also holds in the category of modules over a ring. Therefore, $H_k(M;\mathcal{L}_\xi)$
is a finitely generated $\text{Nov}$-module. The second claim follows from the structure theorem for
finitely generated modules over principal ideal domains, cf. Theorem 6.12 of \cite{HunAlg} or 
Corollary 1.5.3 of \cite{SamAlg}.
\proofend

\begin{definition}\label{NovNum}
Let $\xi \in H^1(M;\mathbb{R})$. The {\bf Novikov numbers} $b_k(\xi)$ and $q_k(\xi)$ are 
defined as follows for all $k=0,\ldots, m$.
\begin{eqnarray*}
b_k(\xi) & = & \text{the rank of }H_k(M;\mathcal{L}_\xi) \text{ as a module over Nov}.\\
q_k(\xi) & = & \text{the minimal number of generators of the torsion} \\
         &   & \text{submodule of } H_k(M;\mathcal{L}_\xi).
\end{eqnarray*}
\end{definition}

\noindent
The Novikov numbers generalize the Betti numbers and torsion numbers of a manifold.

\begin{proposition}[Proposition 1.28 \cite{FarTop}]\label{Novexact}
If $\xi = 0\in H^1(M;\mathbb{R})$, then the Novikov number $b_k(\xi)$ coincides with the Betti 
number $b_k(M)$, i.e. the rank of $H_k(M;\mathbb{Z})$, and $q_k(\xi)$ coincides with the minimal
number of generators of the torsion subgroup of $H_k(M;\mathbb{Z})$, for all $k=0,\ldots, m$.
\end{proposition}

\proofstart
The local coefficient system $\mathcal{L}_0$ is constant, since the homomorphisms are given by 
multiplication by $t^0$, and hence $H_k(M;\mathcal{L}_0) = H_k(M;\text{Nov})$ for all $k=0,\ldots, m$. 
By the Universal Coefficient Theorem we have the following split short exact sequence
of abelian groups, cf. Corollary 3A.4 of \cite{HatAlg} or Theorem 7.55 of \cite{RotAnI}
$$
\xymatrix{
0 \ar[r] & H_k(M;\mathbb{Z}) \otimes \text{Nov} \ar@{->}[r]^-{\lambda_k} & H_k(M;\text{Nov}) \ar[r] & 
\text{Tor}(H_{k-1}(M;\mathbb{Z}), \text{Nov}) \ar[r] & 0,
}
$$
where the $\text{Tor}$ term is zero since $\text{Nov}$ is a torsion free abelian group,
cf. Proposition 3A.5 of \cite{HatAlg} or Proposition 3.1.4 of \cite{WeiAnI}. 
Moreover, $\lambda_k$ is defined by $\lambda_k([c]\otimes g) = [c \otimes g]$, 
which is a homomorphism of $\text{Nov}$-modules. Hence $\lambda_k$ gives an isomorphism of 
$\text{Nov}$-modules
$$
H_k(M;\text{Nov}) \approx H_k(M;\mathbb{Z}) \otimes \text{Nov} 
$$
for all $k=0,\ldots, m$.
\nocite{HilAco}
\proofend

\subsubsection{Novikov inequalities}
A closed 1-form $\zeta \in \Omega^1_{\text{cl}}(M,\mathbb{R})$ is locally exact. In particular, 
for every zero $p \in M$ of $\zeta$ there exists a neighborhood $U$ of $p$ such that $\zeta|_U = df$ 
for some function $f:U \rightarrow \mathbb{R}$. A zero $p$ of $\zeta$ is a critical point of $f$,
and hence the notions of nondegenerate and index for critical points of a function carry over to
the zeros of a closed 1-form. That is, a zero $p$ of $\zeta$ is called \textbf{nondegenerate} if and 
only if $p$ is a nondegenerate critical point of $f$, and the \textbf{Morse index} of the zero $p$ is defined to
be the Morse index of $p$ as a critical point of $f$. 

The significance of the Novikov numbers is indicated by the following theorem, which generalizes the
Morse inequalities. The theorem was proved by M. Farber in the case where $\zeta$ has rank one 
in \cite{FarSha}. The theorem can then be proved for a closed 1-form $\zeta$ of higher rank by 
perturbing $\zeta$ to a closed form of rank one having the same zeros as $\zeta$, cf. Lemma 2.5 
of \cite{FarTop}. The theorem also follows from the more general Novikov Principle 
(Theorem \ref{NovPrin}). 

\begin{theorem}[Novikov Inequalities]\label{Novinequalities}
Let $\zeta \in \Omega^1_{\text{cl}}(M,\mathbb{R})$, and assume that all the zeros of $\zeta$
are nondegenerate. If $c_k(\zeta)$ denotes the number of zeros of $\zeta$ with Morse index
$k$, then
$$
c_k(\zeta) \geq b_k(\xi) + q_k(\xi) + q_{k-1}(\xi)
$$
for all $k=0,\ldots, m$, where $\xi = [\zeta] \in H^1_{\text{dR}}(M;\mathbb{R})$.
\end{theorem}


\subsubsection{Novikov numbers and twisted Morse homology}\label{NovNumMor}
The Twisted Morse Homology Theorem (Theorem \ref{twistedsame}) and Claim \ref{zetarep} and imply that 
the Novikov numbers can be computed using a Morse-Smale-Witten chain complex with coefficients 
in the bundle of rank one $\text{Nov}$-modules $\mathcal{L}_\zeta$ defined in 
Definition \ref{zetasystem}.

\begin{proposition}
Let $f:M \rightarrow \mathbb{R}$ be a smooth Morse-Smale function on a closed finite dimensional
smooth Riemannian manifold $(M,\mathsf{g})$. If $\xi \in H^1(M;\mathbb{R})$ and  
$\zeta \in \Omega^1_{\text{cl}}(M,\mathbb{R})$ satisfies $[\zeta] = \xi$, then  
\begin{eqnarray*}
b_k(\xi) & = & \text{the rank of }H_k((C_\ast(f;\mathcal{L}_\zeta),\partial^{\mathcal{L}_\zeta}_\ast)) 
               \text{ as a module over Nov},\\
q_k(\xi) & = & \text{the minimal number of generators of the torsion} \\
         &   & \text{submodule of } H_k((C_\ast(f;\mathcal{L}_\zeta),\partial^{\mathcal{L}_\zeta}_\ast)),
\end{eqnarray*}
for all $k=0,\ldots, m$.
\end{proposition}

\noindent
We now present a few concrete examples where we compute the Novikov numbers using twisted Morse homology.

\begin{example}[The Novikov numbers of a circle]\label{Novcircle}
Consider a closed 1-form $\zeta$ on the unit circle $S^1 \subset \mathbb{R}^2$.
\begin{figure}[h]
\includegraphics{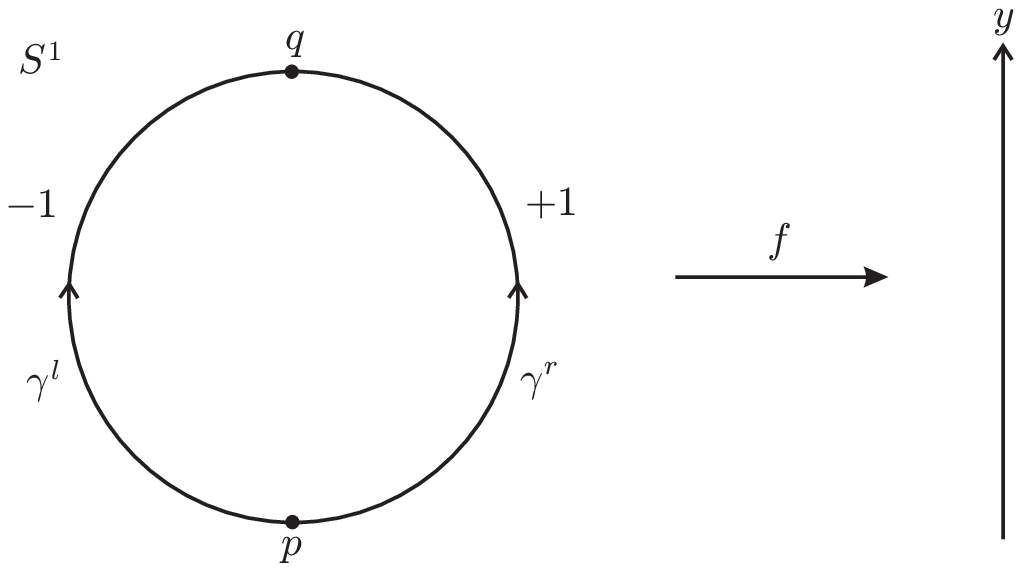}
\end{figure}
Using the Morse-Smale function and orientations from Example \ref{circle} we can compute the Novikov
numbers $b_k(\xi)$ and $q_k(\xi)$ using the following twisted Morse-Smale-Witten complex with coefficients 
in $\mathcal{L}_\zeta$, where $\xi = [\zeta] \in H^1(S^1;\mathbb{R})$.
$$
\xymatrix{
0 \ar[r] & C_1(f;\mathcal{L}_\zeta) \ar[r]^-{\partial_1^{\mathcal{L}_\zeta}}\ar@{<->}[d]^{\approx} & C_0(f;\mathcal{L}_\zeta) \ar@{<->}[d]^{\approx}
\ar[r] & 0\\
0 \ar[r] & \text{Nov}_q \ar[r]^-{\partial_1^{\mathcal{L}_\zeta}} & \text{Nov}_p \ar[r] & 0
}
$$
We have
$$
\partial_1^{\mathcal{L}_\zeta}(g q) = \left(t^{\int_0^1 (\gamma^r)^\ast(\zeta)}g - t^{\int_0^1 
(\gamma^l)^\ast(\zeta)}g \right) p
$$
for all $g \in \text{Nov}$, and hence if $\zeta$ is exact
$$
H_k((C_\ast(f;\mathcal{L}_\zeta),\partial_\ast^{\mathcal{L}_\zeta})) \approx 
\left\{
\begin{array}{ll}
\text{Nov} & \text{if }k = 0,1\\
0          & \text{otherwise}.
\end{array}
\right.
$$
Therefore, $H_k((C_\ast(f;\mathcal{L}_\zeta),\partial_\ast^{\mathcal{L}_\zeta})) \approx H_k(S^1;\mathbb{Z})
\otimes \text{Nov}$ for all $k$ when $[\zeta] = \xi = 0$, since $\mathbb{Z} \otimes \text{Nov}
\approx \text{Nov}$, cf. Theorem 5.7 of \cite{HunAlg}. Thus, $b_k(0)$ and $q_k(0)$ 
agree with the Betti numbers and torsion numbers of $S^1$ for all $k$, cf. Proposition \ref{Novexact}. 

\smallskip
Now consider the case where $\zeta = d\theta$ is the non-exact closed 1-form from Example \ref{circle}.
For any $g \in \text{Nov}$ we have   
$$
\partial_1^{\mathcal{L}_\zeta}(gq) = \left( t^{\int_0^1 (\gamma^r)^\ast(d\theta)}g - t^{\int_0^1 
(\gamma^l)^\ast(d\theta)}g \right) p = \left(\left(t^{\pi} - t^{-\pi}\right) g\right) p,
$$
where $t^\pi - t^{-\pi}$ is invertible since
$$
\frac{1}{t^\pi - t^{-\pi}} = \frac{1}{t^\pi} \frac{1}{1-t^{-2\pi}} = \sum_{i=0}^\infty 
t^{-2\pi i - \pi} \in \text{Nov}.
$$ 
Hence, $\partial_1^{\mathcal{L}_\zeta}: C_1(f;\mathcal{L}_\zeta) \rightarrow C_0(f;\mathcal{L}_\zeta)$
is surjective, and $H_k((C_\ast(f;\mathcal{L}_\zeta),\partial_\ast^{\mathcal{L}_\zeta})) \approx 0$
for all $k$. Thus, $b_k([d\theta]) = q_k([d\theta]) = 0$ for all $k$.
Similarly, one can show that $b_k(\xi) = q_k(\xi) = 0$ for all $k$ whenever $\xi \neq 0$.
\end{example}

\smallskip\noindent
{\bf Note:} If we change the orientation of $W^u(q)$ in the previous example to be counterclockwise
instead of clockwise, then for all $g \in \text{Nov}$
$$
\partial_1^{\mathcal{L}_\zeta}(gq) = \left( t^{\int_0^1 (\gamma^l)^\ast(d\theta)}g - t^{\int_0^1 
(\gamma^r)^\ast(d\theta)}g \right) p = \left(\left(t^{-\pi} - t^{\pi}\right) g\right) p, 
$$
and the element $t^{-\pi} - t^{\pi} \in \text{Nov}$ is still invertible because
$$
\frac{1}{t^{-\pi} - t^{\pi}} = \frac{-1}{t^\pi} \frac{1}{1-t^{-2\pi}} = \sum_{i=0}^\infty 
(-1)t^{-2\pi i - \pi} \in \text{Nov}.
$$
In fact, as noted earlier, an element in $\text{Nov}$ is invertible if and only if its
top coefficient is $\pm 1$, cf. Lemma 1.9 of \cite{FarTop}.

\begin{example}[The Novikov numbers of a torus]\label{Novtorus}
Consider a closed 1-form $\zeta$ on the torus $T^2 = S^1 \times S^1$, viewed as a square
with opposite edges identified. Using the construction detailed in the proof of Theorem
\ref{regularMorse} we can create a Morse-Smale pair $(f,\mathsf{g})$ on $T^2$ with one critical point
$p$ of index $0$, two critical points $q,r$ of index $1$, and one critical point $s$ of index $2$, 
such that there are exactly two gradient flow lines between each pair of critical points of relative
index one. If we orient each $S^1$ clockwise, give $S^1 \times S^1$ the product orientation, and 
follow the orientation conventions in Section \ref{pathcomponents} for the moduli spaces, then the
signs associated to the gradient flow lines are as indicated in the diagram.
\begin{figure}[h]
\includegraphics{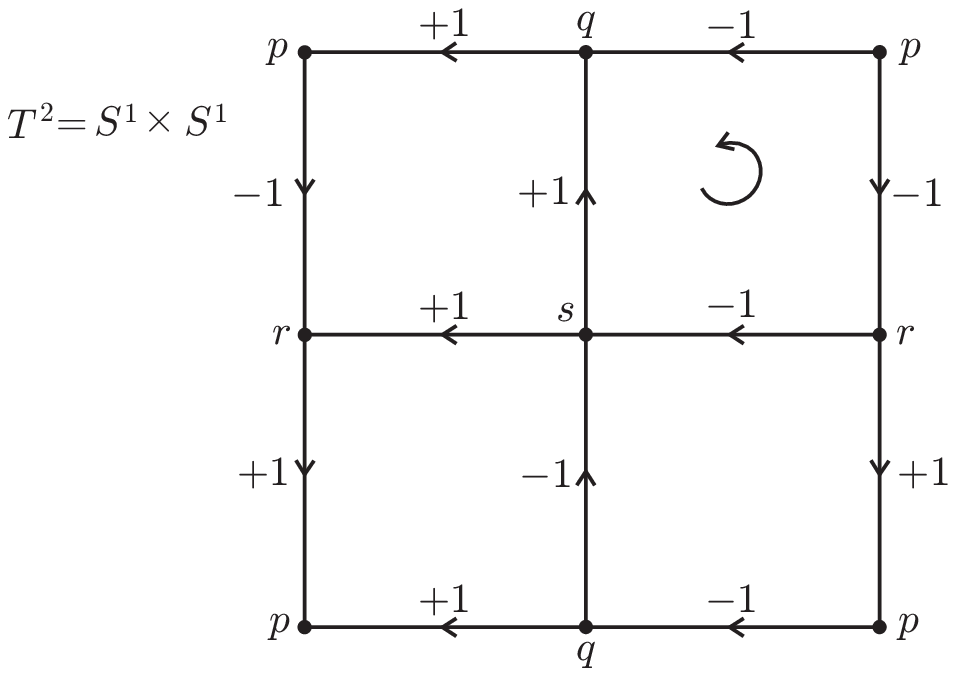}
\end{figure}
For a pair of critical points $(x_1,x_2)$ with $\lambda_{x_2} - \lambda_{x_1} = 1$, let 
$\gamma_{x_1x_2}^+$  and $\gamma_{x_1x_2}^-$ denote parameterizations of the gradient flow lines from 
$x_1$ to $x_2$ with the signs $+1$ and $-1$ respectively. The twisted Morse-Smale-Witten chain complex 
$(C_\ast(f;\mathcal{L}_\zeta), \partial_\ast^{\mathcal{L}_\zeta})$ is  
$$
\xymatrix{
0 \ar[r] & C_2(f;\mathcal{L}_\zeta) \ar[r]^-{\partial_2^{\mathcal{L}_\zeta}}\ar@{<->}[d]^{\approx} & 
 C_1(f;\mathcal{L}_\zeta) \ar[r]^-{\partial_1^{\mathcal{L}_\zeta}}\ar@{<->}[d]^{\approx} & C_0(f;\mathcal{L}_\zeta) \ar@{<->}[d]^{\approx} \ar[r] & 0\\
0 \ar[r] & \text{Nov}_s \ar[r]^-{\partial_2^{\mathcal{L}_\zeta}} & \text{Nov}_r \oplus \text{Nov}_q   \ar[r]^-{\partial_1^{\mathcal{L}_\zeta}} & \text{Nov}_p  \ar[r] & 0
}
$$
where the boundary operator is given by
$$
\partial_k^{\mathcal{L}_\zeta}(g x_2) = \sum_{x_1 \in Cr_{k-1}(f)} \left(t^{\int_0^1
(\gamma_{x_1x_2}^+)^\ast(\zeta)}g - t^{\int_0^1 (\gamma_{x_1x_2}^-)^\ast(\zeta)}g \right)x_1,
$$
for $k = \lambda_{x_2}$ and $g \in \text{Nov}$.
If $\zeta$ is exact, then $\partial^{\mathcal{L}_\zeta}_k = 0$ and
$$
H_k((C_\ast(f;\mathcal{L}_\zeta),\partial_\ast^{\mathcal{L}_\zeta})) \approx H_k(T^2;\mathbb{Z}) \otimes \text{Nov}
$$
for all $k$, since $\mathbb{Z} \otimes \text{Nov} \approx \text{Nov}$, cf. Theorem 5.7 of \cite{HunAlg}.
Hence, $b_k(0)$ and $q_k(0)$ agree with the Betti numbers and torsion numbers 
of $T^2$ for all $k$, cf. Proposition \ref{Novexact}. 

Now consider the case where $\zeta = \pi_1^\ast(d\theta)$ is
the pullback of $d\theta$, the non-exact closed 1-form from Example \ref{circle}, under the projection
onto the first factor of $T^2 = S^1 \times S^1 \subset \mathbb{R}^2 \times \mathbb{R}^2$. 
In this case we have
$$
\partial_1^{\mathcal{L}_\zeta}(g q) = \left(t^{\int_0^1 (\gamma_{pr}^+)^\ast(\zeta)}g - t^{\int_0^1 
(\gamma_{pr}^-)^\ast(\zeta)}g \right)p = ((t^{\pi} - t^{-\pi})g)p,
$$
and
$$
\partial_2^{\mathcal{L}_\zeta}(g s) = \left(t^{\int_0^1 (\gamma_{qs}^+)^\ast(\zeta)}g - t^{\int_0^1 
(\gamma_{qs}^-)^\ast(\zeta)}g \right)q = ((t^\pi - t^{-\pi})g)r,
$$
and $\partial_1^{\mathcal{L}_\zeta}(g r) = 0$ for all $g \in \text{Nov}$. Since $t^{\pi} - t^{-\pi}$ is
invertible in $\text{Nov}$, this implies that
$$
H_k((C_\ast(f;\mathcal{L}_\zeta),\partial_\ast^{\mathcal{L}_\zeta})) \approx 0
$$
for all $k$, and hence, $b_k([\zeta]) = q_k([\zeta]) = 0$ for all $k$ when 
$\zeta = \pi_1^\ast(d\theta)$. 
This is also true when $\zeta = \pi_2^\ast(d\theta)$, and in fact, for any $\xi \in H^1(T^2;\mathbb{R})$
with $\xi \neq 0$. This can be seen directly using a computation similar to the one presented above, 
or by using the fact that $T^2$ is an associative H-space and $t^\gamma - 1$ is invertible in $\text{Nov}$ for
all $\gamma \in \mathbb{R}$ with $\gamma \neq 0$, cf. Example \ref{torusgroup}.
\end{example}

\begin{example}[The Novikov numbers of a Klein bottle]\label{NovKlein}
Consider the Klein bottle $K^2$ as $\mathbb{R}^2$ modulo the properly discontinuous action of the
group $\Delta$ generated by 
\begin{eqnarray*}
\alpha(x,y) & = & (x+1,y)\\
\beta(x,y)  & = & (1-x,y+1).
\end{eqnarray*}
Then $K^2 = \mathbb{R}^2/\Delta$ is a compact smooth manifold with universal cover $\mathbb{R}^2$
and deck transformation group $\Delta$. Moreover, the fundamental group of $K^2$ is
$$
\pi_1(K^2,x_0) \approx \Delta \approx \{ a, b|\ b^{-1} a b = a^{-1} \},
$$
cf. Corollary III.7.3 and Example III.7.5 of \cite{BreTop}. Since $H_1(K^2;\mathbb{Z})$ is 
the abelianization of $\pi_1(K^2,x_0)$, this implies that
$$
H_1(K^2;\mathbb{Z}) \approx \{ a, b|a^2 = 1 \} \approx \mathbb{Z} \oplus \mathbb{Z}_2,
$$
and hence $H_1(K^2;\mathbb{R}) \approx \mathbb{R}$.
If we consider $U = [0,1] \times [0,1] \subset \mathbb{R}^2$ as a fundamental domain for the
action of $\Delta$ with basepoint $x_0 = (0,0)$, then $\alpha \in \pi_1(K^2,x_0)$ corresponds
to a path along the bottom edge of $U$ and $\beta \in \pi_1(K^2,x_0)$ corresponds to a diagonal
in $U$ path from $(0,0)$ to $(1,1)$, cf. Definition III.6.7 of \cite{BreTop}.
The diagonal path from $(0,0)$ to $(1,1)$ is homotopic rel endpoints to a concatenated path along
the boundary of $U$ from $(0,0$) to $(0,1)$ to $(1,1)$, and since homology with real coefficients 
has no torsion, a path $\beta'$ along the left edge of $U$ from $(0,0)$ to $(0,1)$ represents a 
generator of $H_1(K^2;\mathbb{R}) \approx \mathbb{R}$.
\begin{figure}[h]
\includegraphics{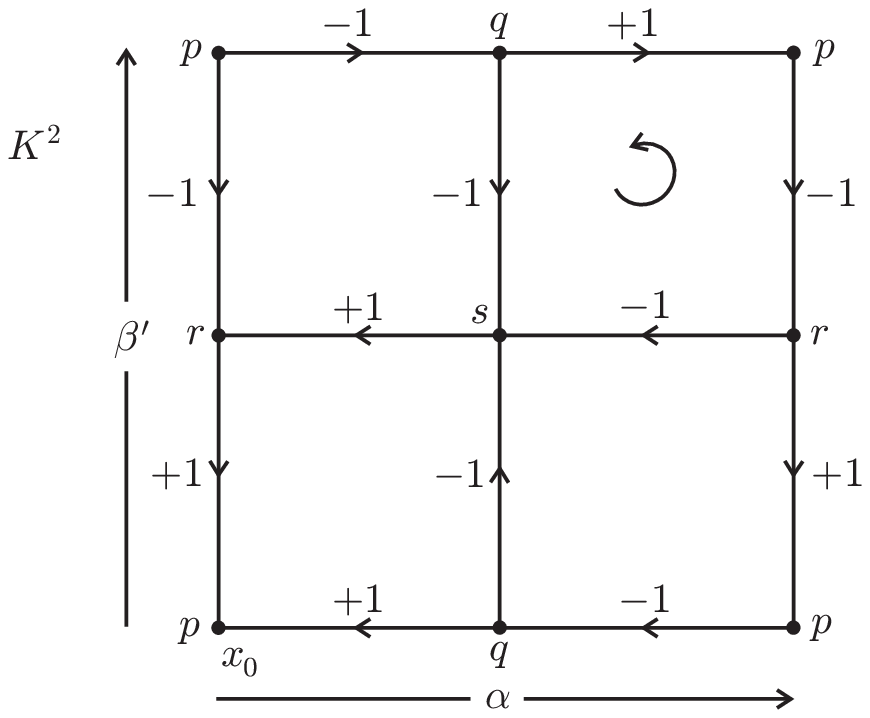}
\end{figure}

Using the construction detailed in the proof of Theorem \ref{regularMorse} we can create
a Morse-Smale pair $(f,\mathsf{g})$ on $K^2$ with one critical point $s$ of index $2$, two critical
points $q,r$ of index $1$, and one critical point $p$ of index $0$, such that there are exactly two
gradient flow lines between each pair of critical points of relative index one. 
If $\zeta \in \Omega^1_{\text{cl}}(K,\mathbb{R})$ is exact, then using the orientations
shown in the diagram we have $\partial_1^{\mathcal{L}_\zeta} = 0$ and for any $g \in \text{Nov}$
$$
\partial_2^{\mathcal{L}_\zeta}(g s) = (t^{\gamma_r}g - t^{\gamma_r}g) r + (-t^{\gamma_q}g - 
t^{\gamma_q}g) q = -2t^{\gamma_q}gq,
$$
for some $\gamma_r,\gamma_q \in \mathbb{R}$, where $-2t^{\gamma_q}$ is not invertible
in $\text{Nov}$. This implies that
$$
H_k((C_\ast(f;\mathcal{L}_\zeta),\partial_\ast^{\mathcal{L}_\zeta})) \approx 
\left\{
\begin{array}{ll}
\text{Nov}                               & \text{if } k=0\\
\text{Nov} \oplus \text{Nov}/2\text{Nov} & \text{if }k = 1\\
0                                        & \text{otherwise}.
\end{array}
\right.
$$
Hence, $b_0(0)= b_1(0) = q_1(0) = 1$, and $b_k(0) = q_k(0) = 0$ for all other $k$, 
cf. Proposition \ref{Novexact}.
 
Now let $\zeta$ be a closed 1-form such that $\int_{\beta'} \zeta = 2\pi$ and
$\int_{\alpha} \zeta = 0$. The form $\zeta$ exists because $\beta'$ represents the generator
of $H_1(K^2;\mathbb{R})$. Explicitly, we can choose $\zeta$ to be the pullback of the closed
1-form $d\theta$ on $S^1$ from Example \ref{circle} along the projection onto the $y$-axis
in the above diagram. With this choice of $\zeta \in \Omega_{\text{cl}}(K^2,\mathbb{R})$ we have
$$
\partial_1^{\mathcal{L}_\zeta}(g r) =  ((t^\pi - t^{-\pi})g)p \quad \text{ and } \quad
\partial_2^{\mathcal{L}_\zeta}(g s) =  ((-t^\pi - t^{-\pi})g)q, 
$$
for any $g \in \text{Nov}$, where both $t^\pi - t^{-\pi}$ and $-t^\pi - t^{-\pi}$ are
invertible in $\text{Nov}$ because their leading coefficients are invertible in $\mathbb{Z}$.
Therefore, $H_k((C_\ast(f;\mathcal{L}_\zeta),\partial_\ast^{\mathcal{L}_\zeta})) \approx 0$ for all $k$,
and hence $b_k([\zeta]) = q_k([\zeta]) = 0$ for all $k$. The same is true when $\zeta$ is any non-exact
closed 1-form. Note that the fundamental difference between the computations when $\zeta$ is exact and
when $\zeta$ is not exact is that $-2t^\gamma$ is not invertible in $\text{Nov}$ for any $\gamma \in
\mathbb{R}$, but $-t^{\gamma} - t^{\gamma- a} $ is invertible in $\text{Nov}$ for all 
$\gamma,a \in \mathbb{R}$ with $a \neq 0$.
\end{example}

\begin{example}[The Novikov numbers of a surface of genus 2]\label{Novsurface2}
Let $S$ be the surface of genus 2 from Example \ref{genus2Lich},
with the same Morse-Smale pair $(f,\mathsf{g})$ and basis of closed 1-forms
$\{\eta_1,\eta_2,\eta_3,\eta_4\}$ for $H^1_{\text{dR}}(S;\mathbb{R})$.
\begin{figure}[h]
\includegraphics{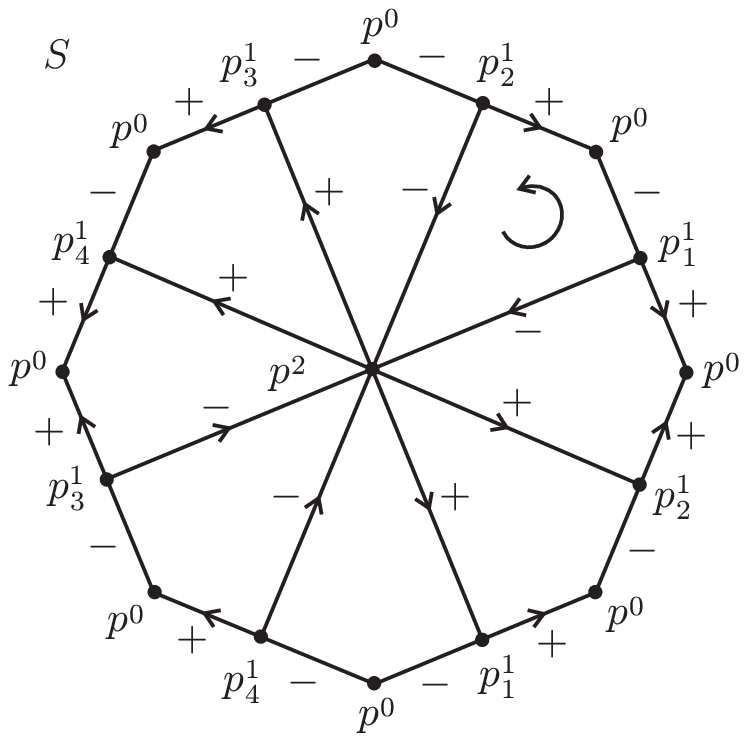}
\end{figure}
For a pair of critical points $(x_1,x_2)$ with $\lambda_{x_2} - \lambda_{x_1} = 1$, let 
$\gamma_{x_1x_2}^+$  and $\gamma_{x_1x_2}^-$ denote parameterizations of the gradient 
flow lines from $x_1$ to $x_2$ with the signs $+1$ and $-1$ respectively.
Let $\zeta$ be a closed 1-form on $S$ and consider the twisted Morse-Smale-Witten chain complex
with coefficients in $\mathcal{L}_\zeta$,
$$
\xymatrix{
0 \ar[r] & C_2(f;\mathcal{L}_\zeta) \ar[r]^-{\partial_2^{\mathcal{L}_\zeta}}\ar@{<->}[d]^{\approx} & 
 C_1(f;\mathcal{L}_\zeta) \ar[r]^-{\partial_1^{\mathcal{L}_\zeta}}\ar@{<->}[d]^{\approx} & 
C_0(f;\mathcal{L}_\zeta) \ar@{<->}[d]^{\approx} \ar[r] & 0\\
0 \ar[r] & \text{Nov}_{p^2} \ar[r]^-{\partial_2^{\mathcal{L}_\zeta}} & \text{Nov}_{p^1_1} \oplus \text{Nov}_{p^1_2} 
\oplus \text{Nov}_{p^1_3} \oplus \text{Nov}_{p^1_4} \ar[r]^-{\partial_1^{\mathcal{L}_\zeta}} & \text{Nov}_{p^0}
\ar[r] & 0
}
$$
where the boundary operator is given by
$$
\partial_k^{\mathcal{L}_\zeta}(g x_2) = \sum_{x_1 \in Cr_{k-1}(f)} \left(t^{\int_0^1 
(\gamma_{x_1x_2}^+)^\ast(\zeta)}g - t^{\int_0^1 (\gamma_{x_1x_2}^-)^\ast(\zeta)}g \right)x_1
$$
for $k = \lambda_{x_2}$ and $g \in \text{Nov}$. If $\zeta$ is exact, then $\partial^{\mathcal{L}_\zeta}_k = 0$
and
$$
H_k((C_\ast(f;\mathcal{L}_\zeta),\partial^{\mathcal{L}_\zeta}_\ast)) \approx H_k(S;\mathbb{Z}) \otimes
\text{Nov}
$$
for all $k$. Thus, $b_0(0) = 1$, $b_1(0) = 4$, $b_2(0) = 1$, $b_k(0) = 0$ for all $k > 2$,
and $q_k(0) = 0$ for all $k$, cf. Proposition \ref{Novexact}.

Now let $\zeta = \eta_1$, the closed 1-form whose integral over $\overline{W^u(p^1_i)}$ is
$1$ when $i=1$ and $0$ when $i \neq 1$. We have
\begin{eqnarray*}
\partial^{\mathcal{L}_{\eta_1}}_1(p^1_i) & = & 0 \text{ for } i\neq 1\\
\partial^{\mathcal{L}_{\eta_1}}_1(p^1_1) & = & (t^{\gamma} - t^{\gamma \pm 1}) p^0 \text{ for some }
\gamma\in \mathbb{R}\\
\partial^{\mathcal{L}_{\eta_1}}_2(p^2) & = & (t^{\gamma'} - t^{\gamma' \pm 1}) p^1_2 \text{ for some }
\gamma'\in \mathbb{R},
\end{eqnarray*}
where $t^{\gamma} - t^{\gamma \pm 1}$ is invertible in $\text{Nov}$ for any $\gamma\in\mathbb{R}$.
Therefore,
$$
H_k((C_\ast(f;\mathcal{L}_{\eta_1}),\partial_\ast^{\mathcal{L}_{\eta_1}})) \approx 
\left\{
\begin{array}{ll}
0                            & \text{if } k=0\\
\text{Nov} \oplus \text{Nov} & \text{if } k = 1\\
0                            & \text{otherwise},
\end{array}
\right.
$$
and we see that $b_1(\eta_1) = 2$, $b_k(\eta_1) = 0$ for all $k \neq 1$, and $q_k(\eta_1) = 0$ 
for all $k$. 

It's clear that same is true for any closed 1-form cohomologous to any of the basis
elements $\{\eta_1,\eta_2,\eta_3,\eta_4\}$ for $H^1_{\text{dR}}(S;\mathbb{R})$.
To compute the Novikov numbers of a general non-exact closed 1-form $\eta \in 
\Omega_{\text{cl}}^1(S,\mathbb{R})$ we will use the Invariance of the Twisted Euler Number
(Corollary \ref{Eulerclassical}), which can be applied to the chain complex 
$(C_\ast(f;\mathcal{L}_\zeta),\partial^{\mathcal{L}_\zeta}_\ast)$ because of 
Theorems \ref{homologyindependence} and \ref{regularMorse} and Lemma \ref{boundarysame}. 

Assume that $\zeta = a_1 \eta_1 + a_2 \eta_2 + a_3 \eta_3 + a_4 \eta_4$ for some 
$a_1,a_2,a_3,a_4 \in \mathbb{R}$ with $a_j \neq 0$ for some $j$, i.e. $\zeta$ is a non-exact closed
1-form, and let $\overline{W(p^2,p^1_l)}$ be the loop on the opposite side of the singular 2-cube
from $\overline{W^u(p^1_j)}$ (see Example \ref{genus2Lich}). For any $g \in G_{p^2}$, 
$$
\partial_2^{\mathcal{L}_\zeta}(g p^2)\ \ = \sum_{p^1_i \in Cr_1(f)}
\left(t^{\int_0^1(\gamma^+_{p^1_ip^2})^\ast(\zeta)}g - t^{\int_0^1(\gamma^-_{p^1_ip^2})^\ast
(\zeta)}g\right)p^1_i,
$$
where the coefficient in front of $p^1_l$ is
\begin{eqnarray*}
\left( t^{\int_0^1 (\gamma^+_{p^1_lp^2})^\ast(a_1 \eta_1 + a_2 \eta_2 + a_3 \eta_3 + a_4 \eta_4)} -
t^{\int_0^1 (\gamma^-_{p^1_lp^2})^\ast(a_1 \eta_1 + a_2 \eta_2 + a_3 \eta_3 + a_4 \eta_4)} \right) g\\
= t^At^Bt^C \left( t^{a_ja}  - t^{a_j(a \pm 1)} \right) g
\end{eqnarray*}
for some $a,A,B,C \in \mathbb{R}$, because
$$
\int_{\gamma^-_{p^1_lp^2}} \eta_i\ = \int_{\gamma^+_{p^1_lp^2}} \eta_i\ \text{ if }i \neq j \quad
\text{and} \quad \int_{\gamma^-_{p^1_lp^2}} \eta_j = a \pm 1\ \text{ if } 
\int_{\gamma^+_{p^1_lp^2}} \eta_j = a.
$$
Thus, $\partial_2^{\mathcal{L}_\zeta}(g p^2) \neq 0$ if $g \neq 0$ and
$H_2((C_\ast(f;\mathcal{L}_\zeta),\partial^{\mathcal{L}_\zeta}_\ast)) \approx 0$.
Now let $g \in G_{p_j^1}$. 
\begin{eqnarray*}
\partial_1^{\mathcal{L}_\zeta}(g p^1_j) & = & \left(t^{\int_0^1(\gamma^+_{p^0p^1_j})^\ast(\zeta)}g - 
t^{\int_0^1(\gamma^-_{p^0p^1_j})^\ast (\zeta)}g\right)p^0\\
& = & \left( t^{\int_0^1 (\gamma^+_{p^0p^1_j})^\ast(a_1 \eta_1 + a_2 \eta_2 + a_3 \eta_3 + a_4 \eta_4)} -
t^{\int_0^1 (\gamma^-_{p^0p^1_j})^\ast(a_1 \eta_1 + a_2 \eta_2 + a_3 \eta_3 + a_4 \eta_4)} \right) gp^0\\
& = & t^At^Bt^C \left( t^{a_ja}  - t^{a_j(a \pm 1)} \right) gp^0
\end{eqnarray*}
for some $a,A,B,C \in \mathbb{R}$, because
$$
\int_{\gamma^-_{p^0p^1_j}} \eta_i\ = \int_{\gamma^+_{p^0p^1_j}} \eta_i\ \text{ if }i \neq j \quad
\text{and} \quad \int_{\gamma^-_{p^0p^1_j}} \eta_j = a \pm 1\ \text{ if } 
\int_{\gamma^+_{p^0p^1_j}} \eta_j = a.
$$
Therefore $\partial_1^{\mathcal{L}_\zeta}$ is surjective when $a_j \neq 0$, and hence 
$H_0((C_\ast(f;\mathcal{L}_\zeta),\partial^{\mathcal{L}_\zeta}_\ast)) \approx 0$.
Thus,
$$
H_k((C_\ast(f;\mathcal{L}_{\zeta}),\partial_\ast^{\mathcal{L}_{\zeta}})) \approx 
\left\{
\begin{array}{ll}
0                            & \text{if } k=0\\
\text{Nov} \oplus \text{Nov} & \text{if } k = 1\\
0                            & \text{otherwise},
\end{array}
\right.
$$
because $\mathcal{X}_{\mathcal{L}_\zeta}(S) = -2$. This shows that for any non-exact closed 1-form $\zeta$
we have $b_1(\zeta) = 2$, $b_k(\zeta) = 0$ for all $k \neq 1$, and $q_k(\zeta) = 0$ for all $k$.

\smallskip
Note that the Novikov inequalities (Theorem \ref{Novinequalities}) are non-trivial and
distinct from the Morse inequalities in this example. That is, the above computation implies that
any non-exact closed Morse 1-form must have at least 2 zeros with Morse index 1. 
\end{example}


\bigskip\noindent
{\bf Acknowledgments}

\smallskip\noindent
We would like to thank the following people for several useful conversations on various
aspects of this project: Pinaki Das, Wojciech Dorabiala, Mark Johnson, and Thomas Krainer.
In particular, the analytic approach to the invariance of the Euler number for Lichnerowicz 
cohomology (Remark \ref{AnalyticEuler}) was brought to our attention by Thomas Krainer, and 
the simple yet elegant argument proving that rank one cohomology classes are dense 
(Remark \ref{rankonedense}) comes from Pinaki Das.


\nocite{KawThe}
\nocite{KroMon}
\nocite{LauAMo}
\nocite{McCAUs}
\nocite{MilMor}
\nocite{OanFib}

\bibliographystyle{amsxport}
\bibliography{references}

\end{document}